# On Deposition of the Product of Demazure Atoms and Demazure Characters

Anna, Ying PUN
Department of Mathematics
University of Pennsylvania, Philadelphia, PA 19104-6395
punying@sas.upenn.edu

June 5, 2016

**Abstract**

This paper studies the properties of Demazure atoms and characters using linear operators and also tableaux-combinatorics. It proves the atom-positivity property of the product of a dominating monomial and an atom, which was an open problem. Furthermore, it provides a combinatorial proof to the key-positivity property of the product of a dominating monomial and a key using skyline fillings, an algebraic proof to the key-positivity property of the product of a Schur function and a key using linear operator and verifies the first open case for the conjecture of key-positivity of the product of two keys using linear operators and polytopes.

## 1 Introduction

Macdonald [Mac96] defined a family of non-symmetric polynomials, called non-symmetric Macdonald polynomials,

$$\{E_\gamma(x_1,\ldots,x_n;q,t)|\gamma \text{ is a weak composition with } n \text{ parts}, n\in\mathbb{N}\}$$

which shares many properties with the family of symmetric Macdonald polynomials [Mac95]

$$\{P_\lambda(x_1,\ldots,x_n;q,t)|\lambda \text{ is a partition with } n \text{ parts}, n\in\mathbb{N}\}.$$

Haglund, Haiman and Loehr [HHL08] obtained a combinatorial formula for $E_\gamma(X;q,t)$ where $X=(x_1,\ldots,x_n)$, using fillings of augmented diagram of shape $\gamma$, called *skyline fillings*, satisfying certain constraints.
Marshall[Mar99] studied the family of non-symmetric Macdonald polynomials using another notation $\hat{E}_\gamma(x_1,\ldots,x_n;q,t) := E_{\overline{\gamma}}(x_n,\ldots,x_1;\frac{1}{q},\frac{1}{t})$. In particular, by setting $q=t=0$ in $\hat{E}_\gamma$, one can obtain Demazure atoms (first studied by Lascoux and Schützenberger[LS90]) $\mathcal{A}_\gamma=\hat{E}_\gamma(x_1,\ldots,x_n;0,0)=E_{\overline{\gamma}}(x_n,\ldots,x_1;\infty,\infty)$. Similarly, one can obtain Demazure characters (key polynomials) by setting $q=t=0$ in $E_\gamma$, *i.e.*, $\kappa_\gamma=E_\gamma(x_1,\ldots,x_n;0,0)=\hat{E}_{\overline{\gamma}}(x_n,\ldots,x_1;\infty,\infty)$. The set of all Demazure atoms forms a basis for the polynomial ring, as does the set of all key polynomials.

Haglund, Luoto, Mason, Remmel and van Willigenburg [HLMvW11], [HMR13] further studied the combinatorial formulas for Demazure atoms and Demazure characters given by the skyline fillings and obtained results which generalized those for Schur functions like the Pieri Rule, the Robinson-Schensted-Knuth (RSK) algorithm, and the Littlewood-Richardson (LR) rule.

It is a classical result in Algebraic Geometry that the product of two Schubert polynomials can be written as a positive sum of Schubert polynomials. A representation theoretic proof is also given recently by using Kráskiewicz-Pragacz modules [Wat14]. However a combinatorial proof of the positivity property of Schubert polynomials has long been open.

Since every Schubert polynomial is a positive sum of key polynomials [LS89], the product of two Schubert polynomials is a positive sum of product of two key polynomials. This suggests one to study the product of two key polynomials. It is known that the product of two key polynomials is not key-positive in general. However, it is still a conjecture that whether the product is atom-positive. This provides a possible approach to a combinatorial proof of the positivity property of Schubert polynomials by trying to recombine the atoms into keys and hence into Schubert polynomials.

Also, since key polynomials are positive sum of atoms [LS90], one can study the atom-positivity properties of the products between atoms and keys or even atoms and atoms to try to prove the conjecture by recombining the atoms back to keys. In this thesis, we prove that the product of a dominating monomial and an atom is always atom positive and that the product of a dominating monomial and a key is always key positive (and hence atom-positive) by using the insertions introduced in [Mas08] and [HLMvW11].

In Section 2, we will give a brief summary on notations and some results in symmetric groups. We will introduce Demazure atoms and Demazure keys in Section 3 by first defining them using linear operators and then define them using semi-standard augmented fillings. We will then study some properties of atoms and characters using both definitions. We also study some properties among the linear operators and obtain certain useful identities for the proofs in later Sections.

In Section 4, we will set up the tools, namely, words and recording tableaux, that we need to prove the main results of this thesis in the first 2 sections in the section and give the proof in Section 4.3. We then give alternative proofs to known results, namely, the key-positivity of the product of a dominating monomial and a key in Section 4.4 using results in Section 6 of [HLMvW11] and the key-positivity of the product of a Schur function and a key in Section 4.5.

We will check the first open case of the conjecture of the key-positivity of the product of two key polynomials in Section 5. We first introduce a geometric interpretation of Demazure atoms and characters in Section 6. We then verify the key-positivity of the product of every pair of keys in this open case in Section 7.

# 2 Symmetric group $S_n$

This section gives a brief summary of the terminologies, notations, lemmas and theorems that will be used in later sections.

Let $[n] = \{1, \dots, n\}$ be the set of all positive integers not greater than $n$. Let $S_n$ be the group of all permutations on $[n]$, i.e. $S_n = \{\sigma : [n] \to [n] | \sigma \text{ is bijective}\}$ with identity $id$ such that $id(j) = j$ for all $j \in [n]$, and the group product is defined as the composition of functions, that is, for all $\sigma_1, \sigma_2 \in S_n$, $\sigma_1\sigma_2(j) = \sigma_1(\sigma_2(j))$ for all $j \in [n]$.

**Definition 2.1.** *Let $n$ be a positive integer and $1 \le k \le n$. A cycle of length $k$, denoted as $(a_1, a_2, \cdots, a_k)$, where $a_1, a_2, \ldots, a_k$ are $k$ distinct integers in $[n]$, is a permutation $\sigma \in S_n$ such that*

$$\begin{cases} \sigma(a_i) = a_{i+1} & \text{for } 1 \le i < k \\ \sigma(a_k) = a_1 & \\ \sigma(j) = j & \text{if } j \neq a_i \text{ for any } 1 \le j \le k \end{cases}$$

.

**Example 1.** *Let $k = 3$ and $n = 5$ , then the 3-cycle $(2, 5, 3)$ represents the permutation*

$$\begin{aligned} \sigma: \quad & 1 \mapsto 1 \\ & 2 \mapsto 5 \\ & 3 \mapsto 2 \\ & 4 \mapsto 4 \\ & 5 \mapsto 3 \end{aligned}$$

.

*Note that $(2, 5, 3)$, $(5, 3, 2)$ and $(3, 2, 5)$ are all treated as the same cycle.*

We say cycles $C_1 = (a_1, \ldots, a_r)$ and $C_2 = (b_1, \ldots, b_k)$ are disjoint if $\{a_1, \ldots, a_r\} \cap \{b_1, \ldots, b_k\} = \emptyset$. For example, $(2, 5, 3)$ and $(1)$ are disjoint cycles while $(2, 5, 3)$ and $(1, 2)$ are not.

**Definition 2.2.** *A cycle of length 2 is called a transposition (or a reflection). In particular, for any positive integer $n$, we call $s_i = (i, i + 1) \in S_n$ a simple transposition (or a simple reflection) for $1 \le i \le n - 1$.*

**Proposition 2.1.** *The simple transpositions in $S_n$ for any integer $n > 1$ satisfy the following relations:*

(i) $s_i^2 = id$ *for* $1 \le i \le n - 2$

(ii) $s_i s_j = s_j s_i$ *for* $|i - j| > 1$

(iii) $s_i s_{i+1} s_i = s_{i+1} s_i s_{i+1} = (i, i + 2)$ *for* $1 \le i \le n - 2$.

**Theorem 2.2.** *Every permutation is a product of disjoint cycles.*

**Theorem 2.3.** *Let $n > 1$ be an integer. The permutation group $S_n$ is generated by simple transpositions, that is ,*

$$S_n = \langle s_1, s_2, \ldots, s_{n-1} \rangle.$$

There are several ways to represent a permutation $\sigma \in S_n$:

1. Two-line notation: $\sigma := \begin{pmatrix} 1 & 2 & \cdots & n \\ \sigma(1) & \sigma(2) & \cdots & \sigma(n) \end{pmatrix}$

2. One-line notation: $\sigma = \sigma(1), \sigma(2), \sigma(3), \cdots, \sigma(n)$

3. Product of disjoint cycles: This follows by Theorem 2.2.

4. Product of simple transpositions: This follows by Theorem 2.3.

**Example 2.** *Consider the permutation $\sigma$ in Example 1, we can write it as:*

1. *Two-line notation:* $\sigma := \begin{pmatrix} 1 & 2 & 3 & 4 & 5 \\ 1 & 5 & 2 & 4 & 3 \end{pmatrix}$

2. *One-line notation:* $\sigma = 1, 5, 2, 4, 3$

3. *Product of disjoint cycles:* $\sigma = (1)(4)(2,5,3)$.

4. *Product of simple transpositions:* $\sigma = (3,4)(4,5)(3,4)(2,3) = s_3 s_4 s_3 s_2$.

From now on, we will use one-line notation to represent a permutation, i.e.

$$\sigma = \sigma(1), \sigma(2), \sigma(3), \cdots, \sigma(n)$$

unless stated otherwise.

Note that applying a transposition $s_i$ on the left of a permutation $\sigma$ means interchanging $i$ and $i+1$ in the one-line notation of $\sigma$ while applying $s_i$ on the right interchanges entries $\sigma(i)$ and $\sigma(i+1)$ in the one line notation of $\sigma$.

By Theorem 2.3, every permutation $\sigma$ can be written as a product of simple transpositions. Hence we can find a decomposition with the shortest length (that is, with the smallest number of transpositions). For $\sigma \neq id$, we call such a decomposition a reduced decomposition of $\sigma$.

**Definition 2.3.** *Let $n \geq 2$ be an integer and $\sigma \in S_n \backslash \{id\}$. Let $\sigma = s_{i_1} s_{i_2} \cdots s_{i_k}$ be a reduced decomposition of $\sigma$. We call $i_1 i_2 \ldots i_k$ a reduced word of $\sigma$.*

**Lemma 2.4.** *Every consecutive substring of a reduced word is also a reduced word.*

**Definition 2.4.** *Let $n$ be any positive integer and a permutation $\sigma \in S_n \backslash \{id\}$. Define the length of $\sigma$, denoted as $l(\sigma)$, as the number of simple transpositions in a reduced decomposition. Define $l(id) = 0$.*

Note that reduced decomposition of a permutation is not unique. For instance, $s_1 s_3 s_2 s_3 = s_1 s_2 s_3 s_2 = s_3 s_1 s_2 s_3$. By Tit's Theorem, any reduced word can be obtained by applying a sequence of *braid relations* (i.e. item (iii) in Proposition 2.1) on any other reduced word representing the same permutation.

**Definition 2.5.** *Let $n$ be a positive integer and $\sigma \in S_n$ be a permutation. The pair $(i,j)$ is called an inversion of $\sigma$ if $i < j$ and $\sigma(i) > \sigma(j)$. Denote $\mathrm{inv}(\sigma)$ as the number of inversions of $\sigma$.*

**Lemma 2.5.** *Let $n \geq 2$ be an integer. For any permutation $\sigma \in S_n$ and a simple transposition $s_i$ ( $1 \leq i \leq n-1$), we have*

$$\mathrm{inv}(s_i\sigma) - \mathrm{inv}(\sigma) = \begin{cases} -1 & \text{if } (\sigma^{-1}(i+1), \sigma^{-1}(i)) \text{ is an inversion pair of } \sigma \\ 1 & \text{else} \end{cases}.$$

**Proposition 2.6.** *Let* $\sigma = s_{i_1}s_{i_2}\dots s_{i_k}$ *(not necessarily reduced). Then* $k \equiv \text{inv}(\sigma) \pmod 2$*.*

**Lemma 2.7.** *Let* $n > 1$ *be an integer and* $\sigma \in S_n$ *be a permutation. Let* $s_i$ *be a transposition in* $S_n$*, where* $1 \le i \le n-1$*. Then* $|l(s_i\sigma) - l(\sigma)| = 1$*.*

**Lemma 2.8.** *Let* $n > 1$ *be an integer and* $\sigma \in S_n$*.* $l(s_i\sigma) = l(\sigma) - 1$ *if and only if there exists a reduced decomposition* $s_{r_1}s_{r_2}\dots s_{r_{l(\sigma)}}$ *such that* $r_1 = i$*.*

*Proof.* If $\sigma$ has a reduced decomposition $s_{r_1}s_{r_2}\dots s_{r_{l(\sigma)}}$ such that $r_1 = i$, then

$$s_i\sigma = s_is_{r_1}s_{r_2}\dots s_{r_{l(\sigma)}} = s_i^2s_{r_2}\dots s_{r_{l(\sigma)}} = s_{r_2}\dots s_{r_{l(\sigma)}}$$

by item (i) in Proposition 2.1. By Lemma 2.4, we know that $s_{r_2}\dots s_{r_{l(\sigma)}}$ is reduced and hence $l(s_i\sigma) = l(\sigma) - 1$.

If $l(s_i\sigma) = l(\sigma)-1$, then consider a reduced decomposition of $s_i\sigma$, say $s_i\sigma = s_{i_1}\cdots s_{i_{l(\sigma)-1}}$, by item (i) in Proposition 2.1, applying $s_i$ on both sides gives $\sigma = s_is_{i_1}\cdots s_{i_{l(\sigma)-1}}$ with exactly $l(\sigma)$ transpositions, which implies $s_is_{i_1}\cdots s_{i_{l(\sigma)-1}}$ is a reduced decomposition of $\sigma$. Hence $\sigma$ has a decomposition with $s_i$ as the leftmost simple transposition. □

**Proposition 2.9.** $l(\sigma) = \text{inv}(\sigma)$ *for any permutation* $\sigma$*.*

*Proof.* We first consider $\sigma^{-1}(1)$. If $\sigma^{-1}(1) \neq 1$, then all the integers before 1 in $\sigma$, i.e. $\sigma(1), \dots, \sigma(\sigma^{-1}(1)-1)$, are all larger than 1, and hence $(r, \sigma^{-1}(1))$ are inversions of $\sigma$ for all $1 \le r < \sigma^{-1}(1)$. So by interchanging 1 with $\sigma(\sigma^{-1}(1)-1)$, and then with $\sigma(\sigma^{-1}(1)-2)$ until with $\sigma(1)$, we can put 1 to the leftmost of the new $\sigma$ (the sigma after interchanging 1 with the $\sigma^{-1}(1)-1$ integers). Indeed, by a previous note, the procedure described above is exactly applying transpositions on the right of $\sigma$, resulting in a new permutation, call it $\sigma^{(1)} = \sigma s_{\sigma^{-1}(1)-1}s_{\sigma^{-1}(1)-2}\cdots s_1$.

Note that each of the above procedure of moving 1 to the front decreases the the number of inversions by exactly 1.

We then use the same procedure by moving 2 to the second leftmost position of $\sigma^{(1)}$ by applying ${\sigma^{(1)}}^{-1}(2) - 1$ simple transpositions on the right of $\sigma^{(1)}$ and get $\sigma^{(2)}$.

Continue this process until we get the $\sigma^{(n-1)}$ which has no inversion, i.e. $\sigma^{(n-1)} = id$. Since each time we apply the interchanging procedure, we are actually applying a simple transposition on the right and also decrease the number of inversion by exactly 1, we have performed exactly $\text{inv}(\sigma)$ interchanging procedures from $\sigma$ to $id$. As a result, we get $\sigma s_{i_1}s_{i_2}\cdots s_{i_{\text{inv}(\sigma)}} = id$ and hence $\sigma = s_{i_{\text{inv}(\sigma)}}\cdots s_{i_1}$. (This also proves Theorem 2.3) which implies $l(\sigma) \le \text{inv}(\sigma)$.

Let $s_{r_1}\cdots s_{r_{l(\sigma)}}$ be a reduced decomposition of $\sigma$. By Lemma 2.5, we know $\text{inv}(s_{r_1}\cdots s_{r_{l(\sigma)}}) \le \text{inv}(s_{r_2}\cdots s_{r_{l(\sigma)}}) + 1 \le \cdots \le \text{inv}(s_{r_{l(\sigma)}}) + l(\sigma) - 1 = l(\sigma)$ and hence we get $\text{inv}(\sigma) \le l(\sigma)$.

As a result, $l(\sigma) = \text{inv}(\sigma)$. □

Note that $n, n-1, \cdots, 1$ has the longest length in $S_n$ as it has the maximum number (namely, $\binom{n}{2}$) of inversions.

**Corollary 2.10.** *Let* $\sigma = n, n-1, \cdots, 1$*. Then for any* $i \in [n-1]$*, there is a reduced word of* $\sigma$ *starting with* $i$*.*

*Proof.* Let $i \in [n-1]$. Since $|l(s_i\sigma) - l(\sigma)| = 1$ and $l(\sigma) > l(s_i\sigma)$ as $\sigma$ has the longest length among all permutations in $S_n$, we have $l(s_i\sigma) = l(\sigma) - 1$. Hence result follows by Lemma 2.8. □

There are several equivalent definitions of Bruhat order on $S_n$ and we will use the reduced word definition. See [Deo77] for further discussion.

**Definition 2.6.** *Let $n$ be a positive integer. Define a partial ordering $\leq$ on $S_n$ such that $\sigma \leq \gamma$ if and only if there exists a reduced word of $\sigma$ which is a substring (not necessarily consecutive) of some reduced word of $\gamma$.*

**Lemma 2.11.** *Let $k \geq 2$ be a positive integer and $i_1 i_2 \dots i_k$ be a reduced word. Let $\sigma' = s_{i_2} \cdots s_{i_k}$ and $\sigma = s_{i_1}\sigma'$. Then $\{\tau | \tau \leq \sigma\} = \{s_{i_1}\gamma, \gamma | \gamma \leq \sigma', l(s_{i_1}\gamma) = l(\gamma) + 1\}$.*

*Proof.* First note that as $i_1 i_2 \dots i_k$ is reduced, by Lemma 2.4, $i_2 \dots i_k$ is a reduced word of $\sigma'$. Hence $\sigma' \leq \sigma$.

Consider $\gamma$ such that $\gamma \leq \sigma'$ such that $l(s_{i_1}\gamma) = l(\gamma) + 1$.

Since $\gamma \leq \sigma'$ and $\sigma' \leq \sigma$, we have $\gamma \leq \sigma$. Also, $\gamma \leq \sigma'$ implies $\gamma$ has a reduced word which is a substring of $i_2 \dots i_k$, say $i_{r_1} i_{r_2} \dots i_{r_{l(\gamma)}}$ where $2 \leq r_1 < r_2 < \cdots < r_{l(\gamma)} \leq k$. Then $s_{i_1}\gamma = s_{i_1} s_{i_{r_1}} s_{i_{r_2}} \cdots s_{i_{r_{l(\gamma)}}}$ which is reduced as $l(s_{i_1}\gamma) = l(\gamma) + 1$. Hence $s_{i_1}\gamma \leq \sigma$ (as $i_1 i_{r_1} \dots i_{r_{l(\gamma)}}$ is a substring of $i_1 i_2 \dots i_k$).

We thus have $\{\tau | \tau \leq \sigma\} \supseteq \{s_{i_1}\gamma, \gamma | \gamma \leq \sigma', l(s_{i_1}\gamma) = l(\gamma) + 1\}$.

Now consider $\tau$ such that $\tau \leq \sigma$. Let $i_{j_1} i_{j_2} \dots i_{j_{l(\tau)}}$ be a reduced word of $\tau$, where $1 \leq j_1 < j_2 < \cdots < j_{l(\gamma)} \leq k$. If $j_1 \neq 1$ for any reduced word of $\tau$, then by Lemma 2.8, $l(s_{i_1}\tau) \neq l(\tau) - 1$. By Lemma 2.7, $l(s_{i_1}\tau) = l(\tau) + 1$. Also, $2 \leq j_1 < j_2 < \cdots < j_{l(\gamma)} \leq k$ implies $\tau \leq \sigma'$. As a result, $\tau \leq \sigma'$ and $l(s_{i_1}\tau) = l(\tau) + 1$ if $j_1 \neq 1$ for any reduced word $i_{j_1} i_{j_2} \dots i_{j_{l(\tau)}}$ of $\tau$.

If $j_1 = 1$, then we can write $\tau = s_{i_1}\tau'$ where $\tau' = s_{i_{j_2}} s_{i_{j_3}} \dots s_{i_{j_{l(\tau)}}}$ which is also a reduced decomposition by Lemma 2.4. Since $1 = j_1 < j_2$, we have $\tau' \leq \sigma'$. Also $l(\tau') = l(\sigma) - 1$ which implies $l(s_{i_1}\tau') = l(\sigma) = l(\tau') + 1$. Hence $\tau = s_{i_1}\tau'$ where $\tau' \leq \sigma', l(s_{i_1}\tau') = l(\tau') + 1$ if $j_1 = 1$ for some reduced word $i_{j_1} i_{j_2} \dots i_{j_{l(\tau)}}$ of $\tau$.

Therefore $\{\tau | \tau \leq \sigma\} \subseteq \{s_{i_1}\gamma, \gamma | \gamma \leq \sigma', l(s_{i_1}\gamma) = l(\gamma) + 1\}$ and result follows.

□

# 3 Demazure atoms and characters

## 3.1 Linear operators

Let $P$ be the polynomial ring $\mathbb{Z}[x_1, x_2, \dots]$ and $S_\infty$ be the permutation group of the positive integers, acting on $P$ by permuting the indices of the variables. For any positive integer $i$, define linear operators

$$\begin{aligned} \partial_i &:= \frac{1 - s_i}{x_i - x_{i+1}} \\ \pi_i &:= \partial_i x_i \\ \theta_i &:= x_{i+1}\partial_i \end{aligned}$$

where $s_i$ is the elementary transposition $(i, i+1)$ and 1 is the identity element in $S_\infty$. Note that for $f \in P$, $(x_i - x_{i+1})$ is a factor of $(1 - s_i)f$ and hence $\partial_i f \in P$. Therefore, $\pi_i f, \theta_i f \in P$.

**Example 3.** *Let $i = 2$ and consider the monomials $x_1^5 x_2^4 x_3$, $x_1^3 x_3^2$ and $x_1 x_2^2 x_3^2$, we have*

1 a) $\partial_2(x_1^5 x_2^4 x_3) = \dfrac{x_1^5 x_2^4 x_3 - s_2(x_1^5 x_2^4 x_3)}{x_2 - x_3} = \dfrac{x_1^5 x_2^4 x_3 - x_1^5 x_2 x_3^4}{x_2 - x_3} = x_1^5(x_2^3 x_3 + x_2^2 x_3^2 + x_2 x_3^3)$

b) $\pi_2(x_1^5 x_2^4 x_3) = \partial_2 x_2(x_1^5 x_2^4 x_3) = \partial_2(x_1^5 x_2^5 x_3) = x_1^5(x_2^4 x_3 + x_2^3 x_3^2 + x_2^2 x_3^3 + x_2 x_3^4)$

c) $\theta_2(x_1^5 x_2^4 x_3) = x_3 \partial_2(x_1^5 x_2^4 x_3) = x_1^5(x_2^3 x_3^2 + x_2^2 x_3^3 + x_2 x_3^4)$

2 a) $\partial_2(x_1^3 x_3^2) = \dfrac{x_1^3 x_3^2 - s_2(x_1^3 x_3^2)}{x_2 - x_3} = \dfrac{x_1^3 x_3^2 - x_1^3 x_2^2}{x_2 - x_3} = -x_1^3(x_2 + x_3)$

b) $\pi_2(x_1^3 x_3^2) = \partial_2 x_2(x_1^3 x_3^2) = \partial_2(x_1^3 x_2 x_3^2) = -x_1^3 x_2 x_3$

c) $\theta_2(x_1^3 x_3^2) = x_3 \partial_2(x_1^3 x_3^2) = -x_1^3(x_2 x_3 + x_3^2)$

3 a) $\partial_2(x_1 x_2^2 x_3^2) = \dfrac{x_1 x_2^2 x_3^2 - s_2(x_1 x_2^2 x_3^2)}{x_2 - x_3} = \dfrac{x_1 x_2^2 x_3^2 - x_1 x_2^2 x_3^2}{x_2 - x_3} = 0$

b) $\pi_2(x_1 x_2^2 x_3^2) = \partial_2 x_2(x_1 x_2^2 x_3^2) = \partial_2(x_1 x_2^3 x_3^2) = x_1 x_2^2 x_3^2$

c) $\theta_2(x_1 x_2^2 x_3^2) = x_3 \partial_2(x_1 x_2^2 x_3^2) = 0$

**Proposition 3.1.** *For any positive integer $i$, we have*

1. $\pi_i = \theta_i + 1$;

2. $\pi_i \theta_i = \theta_i \pi_i$;

3. $s_i \partial_i = -\partial_i s_i = \partial_i$, $s_i \pi_i = \pi_i$, $\pi_i s_i = -\partial_i x_{i+1}$, $s_i \theta_i = x_i \partial_i$, $\theta_i s_i = -\theta_i$;

4. $\partial_i \partial_i = 0$, $\partial_i \pi_i = \theta_i \partial_i = 0$, $\pi_i \partial_i = -\partial_i \theta_i = \partial_i$;

5. $\pi_i \pi_i = \pi_i$, $\theta_i \theta_i = -\theta_i$, $\pi_i \theta_i = \theta_i \pi_i = 0$.

*Proof.* Let $f \in P$. Then

$$\begin{aligned}
\pi_i f \quad = \partial_i x_i f \quad &= \frac{x_i f - s_i(x_i f)}{x_i - x_{i+1}} \\
&= \frac{x_i f - x_{i+1} s_i f}{x_i - x_{i+1}} \\
&= \frac{x_{i+1} f - x_{i+1} s_i f + (x_i - x_{i+1}) f}{x_i - x_{i+1}} \\
&= x_{i+1} \frac{f - s_i(f)}{x_i - x_{i+1}} + f \\
&= x_{i+1} \partial_i f + f \\
&= (\theta_i + 1) f
\end{aligned}$$

and hence $\pi_i = \theta_i + 1$.

By item 1, we have $\pi_i\theta_i = (\theta_i+1)\theta_i = \theta_i\theta_i + \theta_i = \theta_i(\theta_i+1) = \theta_i\pi_i$, proving item 2.

$s_i\partial_i f = s_i(\dfrac{f - s_if}{x_i - x_{i+1}}) = \dfrac{s_if - s_is_if}{x_{i+1} - x_i} = \dfrac{s_if - f}{x_{i+1} - x_i} = \dfrac{f - s_if}{x_i - x_{i+1}} = \partial_i f.$

$\partial_i s_i f = \dfrac{s_if - s_i(s_if)}{x_i - x_{i+1}} = -\dfrac{f - s_if}{x_i - x_{i+1}} = -\partial_i f.$

$s_i\pi_i f = s_i(\partial_i x_i f) = s_i\partial_i(x_if) = \partial_i(x_if) = \pi_i f.$

$\pi_i s_i f = \partial_i(x_is_if) = \partial_i s_i(x_{i+1}f) = -\partial_i x_{i+1} f.$

$s_i\theta_i f = s_i(x_{i+1}\partial_i f) = x_is_i\partial_i f = x_i\partial_i f.$

$\theta_i s_i f = x_{i+1}\partial_i s_i f = -x_{i+1}\partial_i f = -\theta_i f.$

Hence, $s_i\partial_i = \partial_i, \partial_i s_i = -\partial_i, s_i\pi_i = \pi_i, \pi_i s_i = -\partial_i x_{i+1}, s_i\theta_i = x_i\partial_i$, $\theta_i s_i = -\theta_i$ and item 3 follows.

By item 3, $\partial_i\partial_i f = \dfrac{\partial_i f - s_i\partial_i f}{x_i - x_{i+1}} = \dfrac{\partial_i f - \partial_i f}{x_i - x_{i+1}} = 0.$

Hence $\partial_i\pi_i f = \partial_i\partial_i x_i f = 0$ and $\theta_i\partial_i f = x_{i+1}\partial_i\partial_i f = 0$.

Then $\pi_i\partial_i f = (1+\theta_i)\partial_i f = \partial_i f + \theta_i\partial_i f = \partial_i f$.

Also, $\partial_i\theta_i f = \partial_i(\pi_i - 1)f = \partial_i\pi_i f - \partial_i f = -\partial_i f$.

As a result, $\partial_i\partial_i = \partial_i\pi_i = \theta_i\partial_i = 0$, $\pi_i\partial_i = \partial_i$, $\partial_i\theta_i = -\partial_i$, proving item 4.

Now by item 4, we have

$\pi_i\pi_i f = (\pi_i\partial_i)(x_if) = \partial_i(x_if) = \pi_i f,$

$\theta_i\theta_i f = (x_{i+1}\partial_i)\theta_i f = x_{i+1}(\partial_i\theta_i)f = -x_{i+1}\partial_i f = -\theta_i f,$

Therefore $\theta_i\pi_i f = \theta_i(\theta_i+1)f = \theta_i\theta_i f + \theta_i f = -\theta_i f + \theta_i f = 0$.

Furthermore, $\pi_i\theta_i f = (\theta_i+1)\theta_i f = \theta_i\theta_i f + \theta_i f = 0$ and item 5 follows.

□

**Proposition 3.2.** *$\partial_i\partial_j = \partial_j\partial_i$, $\pi_i\pi_j = \pi_j\pi_i$ and $\theta_i\theta_j = \theta_j\theta_i$ for $|i-j| \geq 2$.*

*Proof.* Let $f \in P$ and $i, j \in \mathbb{N}$ such that $|i-j| \geq 2$. Then

$$\partial_i\partial_j f = \frac{\partial_j f - s_i\partial_j f}{x_i - x_{i+1}}.$$

Since $s_is_j = s_js_i$, we have

$$s_i\partial_j f = s_i\frac{f - s_jf}{x_j - x_{j+1}} = \frac{s_if - s_is_jf}{x_j - x_{j+1}} = \frac{s_if - s_js_if}{x_j - x_{j+1}} = \partial_j s_i f.$$

Hence

$$\begin{aligned}\partial_i\partial_j f &= \frac{\partial_j f - \partial_j(s_if)}{x_i - x_{i+1}} \\ &= \frac{\partial_j(f - s_if)}{x_i - x_{i+1}} \\ &= \frac{(f - s_if) - s_j(f - s_if)}{(x_j - x_{j+1})(x_i - x_{i+1})} \\ &= \frac{f - s_if - s_jf - s_is_jf}{(x_i - x_{i+1})(x_j - x_{j+1})}.\end{aligned}$$

Similarly $\partial_j\partial_i f = \dfrac{f - s_jf - s_if - s_js_if}{(x_j - x_{j+1})(x_i - x_{i+1})} = \dfrac{f - s_if - s_jf - s_is_jf}{(x_i - x_{i+1})(x_j - x_{j+1})}.$

Thus $\partial_i \partial_j = \partial_j \partial_i$.

Note that $x_i \partial_j f = x_i \dfrac{f - s_j f}{x_j - x_{j+1}} = \dfrac{x_i f - x_i s_j f}{x_j - x_{j+1}} = \dfrac{x_i f - s_j x_i f}{x_j - x_{j+1}} = \partial_j x_i f$, that is, $x_i \partial_j = \partial_j x_i$. Similarly, $x_j \partial_i = \partial_i x_j$.

Hence we have

$$\pi_i \pi_j f = \partial_i x_i \partial_j x_j f = \partial_i \partial_j x_i x_j f = \partial_j \partial_i x_i x_j f = \partial_j \partial_i x_j x_i f = \partial_j x_j \partial_i x_i f = \pi_j \pi_i f.$$

$$\theta_i \theta_j = (\pi_i - 1)(\pi_j - 1) = \pi_i \pi_j - \pi_i - \pi_j + 1 = \pi_j \pi_i - \pi_i - \pi_j + 1 = (\pi_j - 1)(\pi_i - 1) = \theta_j \theta_i.$$

□

**Proposition 3.3.** *For any positive integer $i$, linear operators $\partial_i, \pi_i$ and $\theta_i$ satisfy the braid relation. That is, $\partial_i \partial_{i+1} \partial_i = \partial_{i+1} \partial_i \partial_{i+1}$, $\pi_i \pi_{i+1} \pi_i = \pi_{i+1} \pi_i \pi_{i+1}$ and $\theta_i \theta_{i+1} \theta_i = \theta_{i+1} \theta_i \theta_{i+1}$.*

*Proof.* By direct computation. □

For any permutation $\sigma \neq id$, define $\partial_\sigma = \partial_{i_1} \dots \partial_{i_j}$ where $s_{i_1} \dots s_{i_j}$ is a reduced decomposition of $\sigma$. By the same argument in Proposition 3.2 and Proposition 3.3, we also have $\theta_\sigma = \theta_{i_1} \dots \theta_{i_j}$ and $\pi_\sigma = \pi_{i_1} \dots \pi_{i_j}$. We define $\partial_{id} = \theta_{id} = \pi_{id} = id$.

**Lemma 3.4.** *Let $n > 1$ be an integer and consider a permutation $\gamma \in S_n$. For $1 \leq i \leq n-1$,*
$$\theta_i \theta_\gamma = \begin{cases} -\theta_\gamma & \text{if } l(s_i \gamma) = l(\gamma) - 1 \\ \theta_{s_i \gamma} & \text{if } l(s_i \gamma) = l(\gamma) + 1 \end{cases}$$

*Proof.* By Lemma 2.8, if $l(s_i \gamma) = l(\gamma) - 1$, then there exists a reduced decomposition of $\gamma = s_i s_{r_2} \cdots s_{r_{l(\gamma)}}$ and hence $\theta_\gamma = \theta_i \theta_{r_2} \cdots \theta_{r_{l(\gamma)}}$. By item 5. in Proposition 3.1, we have $\theta_i \theta_\gamma = (\theta_i \theta_i) \theta_{r_2} \cdots \theta_{r_{l(\gamma)}} = (-\theta_i) \theta_{r_2} \cdots \theta_{r_{l(\gamma)}} = -\theta_i \theta_{r_2} \cdots \theta_{r_{l(\gamma)}} = -\theta_\gamma$. Otherwise if $l(s_i \gamma) = l(\gamma) + 1$, $s_i s_{i_1} s_{i_2} \cdots s_{i_{l(\gamma)}}$ is a reduced decomposition of $s_i \gamma$ for any reduced decomposition $s_{i_1} s_{i_2} \cdots s_{i_{l(\gamma)}}$ of $\gamma$. Thus $\theta_i \theta_\gamma = \theta_i \theta_{i_1} \theta_{i_2} \cdots \theta_{i_{l(\gamma)}} = \theta_{s_i \gamma}$. □

**Lemma 3.5.** *For any permutation $\sigma$, $\pi_\sigma = \sum_{\gamma \leq \sigma} \theta_\gamma$.*

*Proof.* Let $k = l(\sigma)$. We prove the statement by induction on $k$.
For $\sigma = id$ (i.e. $k = 0$), $\pi_{id} = \theta_{id} = id$.
When $k = 1$, then $\sigma = s_i$ for some positive integer $i$. By item 1. in Proposition 3.1, we have $\pi_\sigma = \pi_{s_i} = \pi_i = 1 + \theta_i = \theta_{id} + \theta_{s_i}$. Hence the statement is true for $l(\sigma) = 1$.

Assume the statement is true for all non-negative integers $k \leq m$ for some $m \geq 1$.

Let $l(\sigma) = m + 1$. Let $s_{i_1} s_{i_2} \cdots s_{i_{m+1}}$ be a reduced decomposition of $\sigma$. Let $\sigma' = s_{i_2} \cdots s_{i_{m+1}}$ (which implies $l(\sigma') \leq m$ by definition) and hence $\sigma = s_{i_1} \sigma'$. Note that by Lemma 2.4, $s_{i_2} \cdots s_{i_{m+1}}$ is a reduced decomposition of $\sigma'$.

By induction assumption, $\pi'_\sigma = \sum_{\gamma \le \sigma'} \theta_\gamma$. Now we have

$$\begin{aligned}
& \pi_\sigma \\
= \;& \pi_{i_1}\pi_{i_2}\cdots\pi_{i_{m+1}} \\
= \;& \pi_{i_1}\pi_{\sigma'} \\
= \;& (1+\theta_{i_1})\sum_{\gamma\le\sigma'}\theta_\gamma \\
= \;& (1+\theta_{i_1})\Bigg(\sum_{\substack{\gamma\le\sigma' \\ l(s_{i_1}\gamma)=l(\gamma)+1}}\theta_\gamma + \sum_{\substack{\gamma\le\sigma' \\ l(s_{i_1}\gamma)=l(\gamma)-1}}\theta_\gamma\Bigg) \qquad \text{(by Lemma 2.7)} \\
= \;& \sum_{\substack{\gamma\le\sigma' \\ l(s_{i_1}\gamma)=l(\gamma)+1}}\theta_\gamma + \sum_{\substack{\gamma\le\sigma' \\ l(s_{i_1}\gamma)=l(\gamma)-1}}\theta_\gamma + \sum_{\substack{\gamma\le\sigma' \\ l(s_{i_1}\gamma)=l(\gamma)+1}}\theta_{i_1}\theta_\gamma + \sum_{\substack{\gamma\le\sigma' \\ l(s_{i_1}\gamma)=l(\gamma)-1}}\theta_{i_1}\theta_\gamma \\
= \;& \sum_{\substack{\gamma\le\sigma' \\ l(s_{i_1}\gamma)=l(\gamma)+1}}\theta_\gamma + \sum_{\substack{\gamma\le\sigma' \\ l(s_{i_1}\gamma)=l(\gamma)-1}}\theta_\gamma + \sum_{\substack{\gamma\le\sigma' \\ l(s_{i_1}\gamma)=l(\gamma)+1}}\theta_{s_{i_1}\gamma} + \sum_{\substack{\gamma\le\sigma' \\ l(s_{i_1}\gamma)=l(\gamma)-1}}-\theta_\gamma \\
& \quad \text{(by Lemma 3.4)} \\
= \;& \sum_{\substack{\gamma\le\sigma' \\ l(s_{i_1}\gamma)=l(\gamma)+1}}(\theta_\gamma+\theta_{s_{i_1}\gamma}) \\
= \;& \sum_{\tau\le\sigma}\theta_\tau \qquad \text{(by Lemma 2.11)}
\end{aligned}$$

and result follows by induction. □

**Lemma 3.6.** *For any $f, g \in P$ and $i \in \mathbb{N}$, we have*

1. $\partial_i(fg) = (\partial_i f)g + (s_i f)(\partial_i g)$;
2. $\theta_i(fg) = (\theta_i f)g + (s_i f)(\theta_i g)$;
3. $\pi_i(fg) = (\pi_i f)g + (s_i f)(\theta_i g)$.

*Proof.*

$$\begin{aligned}
& \partial_i(fg) \\
= \;& \frac{fg - s_i(fg)}{x_i - x_{i+1}} \\
= \;& \frac{fg - (s_i f)g + (s_i f)g - s_i f s_i g}{x_i - x_{i+1}} \\
= \;& g\frac{f - s_i f}{x_i - x_{i+1}} + (s_i f)\frac{g - s_i g}{x_i - x_{i+1}} \\
= \;& g(\partial_i f) + (s_i f)(\partial_i g) \\
= \;& (\partial_i f)g + (s_i f)(\partial_i g).
\end{aligned}$$

Therefore

$$\begin{aligned} & \theta_i(fg) \\ = \ & x_{i+1}\partial_i(fg) \\ = \ & x_{i+1}\big((\partial_i f)g + (s_i f)(\partial_i g)\big) \\ = \ & (x_{i+1}\partial_i f)g + (s_i f)(x_{i+1}\partial_i g) \\ = \ & (\theta_i f)g + (s_i f)(\theta_i g). \end{aligned}$$

By item 1. of Proposition 3.1,

$$\begin{aligned} & \pi_i(fg) \\ = \ & (1+\theta_i)(fg) \\ = \ & fg + \theta_i(fg) \\ = \ & fg + (\theta_i f)g + (s_i f)(\theta_i g) \\ = \ & \big((1+\theta_i)f\big)g + (s_i f)(\theta_i g) \\ = \ & (\pi_i f)g + (s_i f)(\theta_i g). \end{aligned}$$

□

## 3.2 Semi-Standard Augmented filling

Let $\mathbb{N}$ (or $\mathbb{Z}^+$) be the set of all positive integers and $\mathbb{Z}_{\geq 0}$ be the set of non-negative integers. Also we denote $\epsilon_k = 12\cdots k = id$ as the identity element (we write permutations in one line notation). For $n \in \mathbb{Z}_{\geq 0}$ and $k \in \mathbb{N}$, we say $\alpha = (\alpha_1, \alpha_2, \ldots, \alpha_k) \in (\mathbb{Z}_{\geq 0})^k$ is a weak composition $n$ (denoted as $\alpha \vDash n$) with $k$ parts if $\sum_{i=1}^{k} \alpha_i = n$ and write $l(\alpha) = k$ to denote the length (the number of parts) of $\alpha$. Furthermore, if $\alpha_1 \geq \alpha_2 \geq \cdots \geq \alpha_k \geq 0$, we call $\alpha$ a partition of $n$ with $k$ parts and write $\alpha \vdash n$ (usually we denote $l(\alpha) = \max\{i : \alpha_i > 0\}$ for $\alpha$ being a partition). We use $\mathrm{Par}(n)$ to denote the set of all partitions of a nonnegative integer $n$ and Par to denote the set of all partitions. For a weak composition $\alpha$ with $k$ parts, define $x^\alpha := x_1^{\alpha_1} x_2^{\alpha_2} \cdots x_k^{\alpha_k}$. We call $x^\alpha$ a dominating monomial if $\alpha$ is a partition.

We denote $\overline{\alpha}$ as the *reverse* of $\alpha$, that is , $\overline{\alpha} = (\alpha_k, \ldots, \alpha_1)$. (Note that in [HMR13], they use $\alpha^*$ instead of $\overline{\alpha}$.) Similarly, we write $\overline{X}$ as the *reverse* of $X$ for any finite string of alphabets $X$. For example, $\overline{cacdba} = abdcac$ and $\overline{14D9c7} = 7c9D41$.

Define $\omega_\alpha$ as the permutation of minimal length such that

$$\omega_\alpha(\alpha) := (\alpha_{\omega_\alpha(1)}, \alpha_{\omega_\alpha(2)}, \ldots, \alpha_{\omega_\alpha(k)})$$

is a partition.

Given two weak compositions $\alpha$ and $\beta$, we write $\beta \geq \alpha$ if and only if $\omega_\beta \leq \omega_\alpha$ in the strong Bruhat order.

Let $\alpha$ be a weak composition. The augmented diagram of shape $\alpha$ is the figure with $|\alpha| + l(\alpha)$ cells (or boxes) where column $i$ has $\alpha_i + 1$ cells. The bottom row is called the *basement* of the augmented diagram.

For example, if $\alpha = (1, 0, 1, 0, 0, 4, 0, 6, 5)$, then the augmented diagram of $\alpha$ is

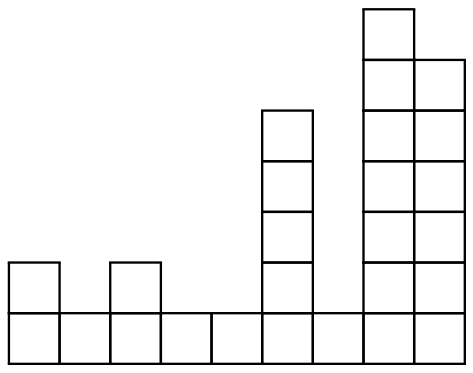

Also we impose an order, called the *reading order*, on the cells of the diagram which starts from left to right, top to bottom. So the order of the above diagram is:

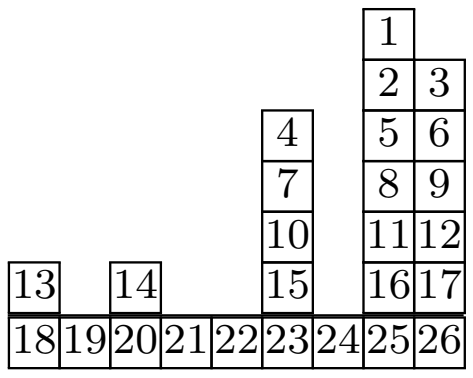

where the number in each cell represents the order of that cell in reading order.

A filling of an augmented diagram is an assignment of a positive integer to each cell in the diagram.

From now on we only consider fillings whose entries in each column are weakly decreasing from the bottom to the top.

For any two columns (including the basement cells) $i$ and $j$ with $i < j$, we pick three cells $X, Y$ and $Z$, where cell $X$ is immediately above cell $Y$ in the 'taller' column $k$ , where

$$k = \begin{cases} i & \text{if} \quad \alpha_i \geq \alpha_j \\ j & \text{if} \quad \alpha_i < \alpha_j \end{cases},$$

and cell $Z$ from the 'shorter' column to form a *triple* $(X, Y, Z)$ in the following way:

$$\begin{cases} \text{Type A triple: cell } Z \text{ is in the same row as cell } X & \text{if} \quad \alpha_i \geq \alpha_j \\ \text{Type B triple: cell } Z \text{ is in the same row as cell } Y & \text{if} \quad \alpha_i < \alpha_j \end{cases}.$$

Here are some examples of triples:

(The first two are type A triples and the last two are type B triples.)

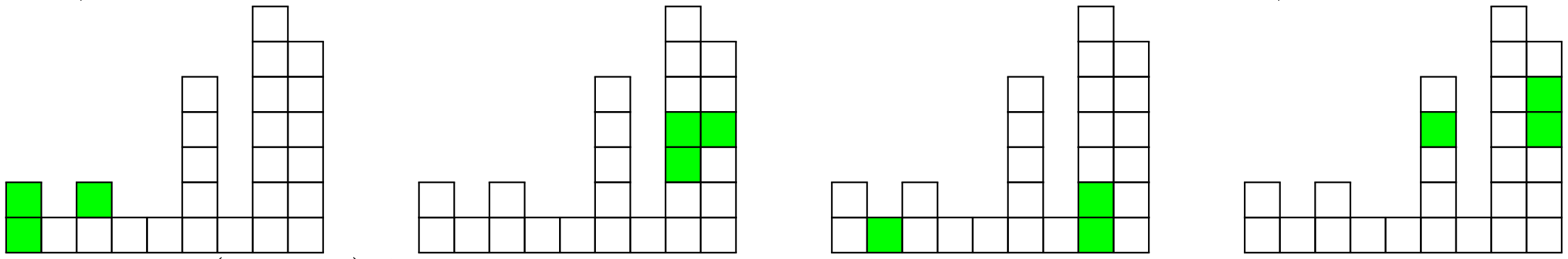

We say $(X, Y, Z)$ forms a *coinversion triple* if the filling $F$ of the diagram assigns each cell in the triple a positive integer, say $F(X), F(Y), F(Z)$ respectively, in such a way that $F(X) \leq F(Z) \leq F(Y)$. Otherwise we call $(X, Y, Z)$ an *inversion triple.*

For instance in the following examples, the second and the third ones are coinversion triples while the first and the last one are inversion triples.

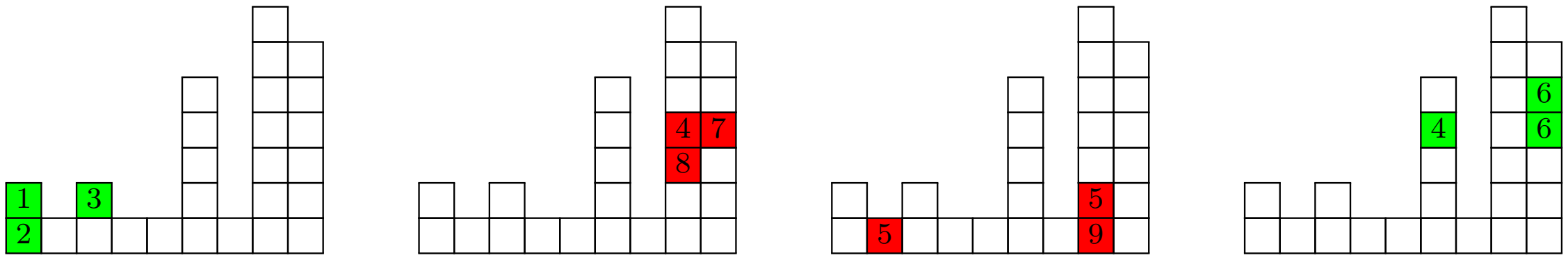

**Definition 3.1.** *A semi-standard augmented filling (SSAF) of an augmented diagram with shape being a weak composition* $\alpha = (\alpha_1, \alpha_2, \dots, \alpha_k)$ *is a filling satisfying:*

1. *the basement entries form a permutation* $\sigma$ *(in one line notation) of* $\{1, ..., k\}$*, i.e.* $\sigma \in S_k$*;*

2. *every (Type A or B) triple is an inversion triple.*

We denote $SSAF(\sigma, \alpha)$ the set of all SSAF of an augmented diagram of shape $\alpha = (\alpha_1, \dots, \alpha_k)$ with basement entries (from left to right) being $\sigma \in S_k$ (i.e. basement of column $i$ has entry $\sigma(i)$).

**Example 4.** *The following SSAFs are all the elements in the set* $SSAF(4132, 1032)$*:*

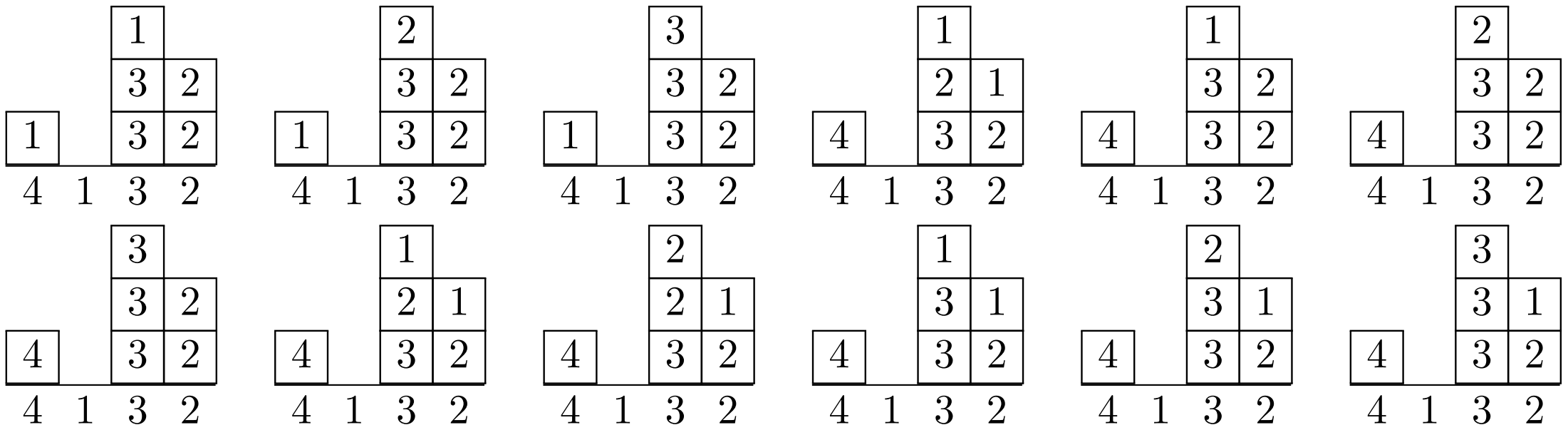

Given an SSAF $F$ with basement of length $n$ for some $n \in \mathbb{N}$, define the weight of $F$ as $x^F := \prod_{i=1}^{n} x_i^{m_i(F)-1}$, where $m_i(F)$ is the number of $i$ appearing in $F$ for $1 \leq i \leq n$.

**Example 5.** $x^F = x_1^2 x_2^2 x_3^2$ *for* $F =$

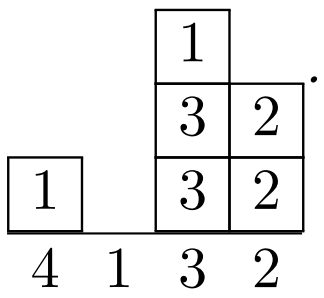

.

### 3.3 Demazure atoms and Demazure characters

**Definition 3.2.** *A Demazure atom of shape $\alpha$, where $\alpha$ is a weak composition, is defined as*

$$\mathcal{A}_\alpha := \sum_{F \in SSAF(\epsilon_k, \alpha)} x^F.$$

**Definition 3.3.** *A Demazure character (key polynomial, or key) of shape $\alpha$, where $\alpha$ is a weak composition, is defined as*

$$\kappa_{\overline{\alpha}} := \sum_{F \in SSAF(\overline{\epsilon_k}, \alpha)} x^F.$$

Remark: See [Mas09] for further discussion on key polynomials.

The following theorem gives different equivalent definitions of Demazure atoms ($\mathcal{A}_\alpha$) and Demazure characters ($\kappa_{\overline{\alpha}}$) of shape $\alpha$ and is proved in [HHL08, LS90, Mas09].

**Theorem 3.7.** *For $k \in \mathbb{N}$ and composition $\alpha$ with $l(\alpha) = k$,*

$$E_{\overline{\alpha}}(x_k, \ldots, x_1; \infty, \infty) = \sum_{F \in SSAF(\epsilon_k, \alpha)} x^F = \theta_{\omega_\alpha^{-1}} x^{\omega_\alpha(\alpha)}$$

$$E_{\overline{\alpha}}(x_1, \ldots, x_k; 0, 0) = \sum_{F \in SSAF(\overline{\epsilon_k}, \alpha)} x^F = \pi_{\overline{\epsilon_k}\omega_\alpha^{-1}} x^{\omega_\alpha(\alpha)}$$

*where $E_{\overline{\alpha}}(x_k, \ldots, x_1; \infty, \infty)$ and $E_{\overline{\alpha}}(x_1, \ldots, x_k; 0, 0)$ are the nonsymmetric Macdonald polynomials of shape $\overline{\alpha}$ with $q = t = \infty$ and $X = (x_k, \ldots, x_1)$ and with $q = t = 0$ and $X = (x_1, \ldots, x_k)$ respectively.*

**Example 6.** *Let* $\alpha = (1,0,3)$. *Then* $\omega_\alpha = 231 = (12)(23) = s_1s_2$.

*Hence* $\mathcal{A}_{(1,0,3)} = \theta_2\theta_1 x_1^3x_2 = \theta_2(x_1^2x_2^2 + x_1x_2^3) = x_1^2x_2x_3 + x_1^2x_3^2 + x_1x_2^2x_3 + x_1x_2x_3^2 + x_1x_3^3$.

*The following are all the SSAFs of* $SSAF(123,(1,0,3))$ :

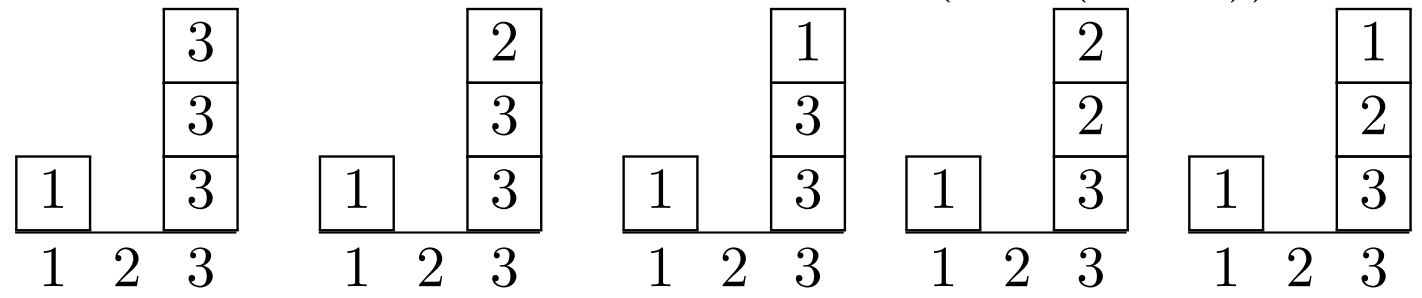

*weights :* $x_1x_3^3$ $\quad x_1x_2x_3^2$ $\quad x_1^2x_3^2$ $\quad x_1x_2^2x_3$ $\quad x_1^2x_2x_3$

*As a result, the sum of all weights of the SSAFs gives* $\mathcal{A}_{(1,0,3)}$.

*Since* $\overline{\epsilon_3}\omega_\alpha^{-1} = (s_2s_1s_2)(s_2s_1) = s_2$, *we have* $\kappa_{(3,0,1)} = \pi_{\overline{\epsilon_3}\omega_\alpha^{-1}}(x_1^3x_2) = \pi_2(x_1^3x_2) = x_1^3x_2 + x_1^3x_3$.

*The following are all the SSAFs of* $SSAF(321,(1,0,3))$ :

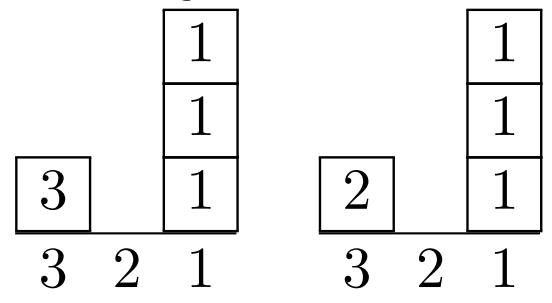

*weights :* $x_1^3x_3$ $\quad x_1^3x_2$

*As a result, the sum of all weights of the SSAFs gives* $\kappa_{(3,0,1)}$.

The following theorem gives some properties of Demazure atoms and characters.

**Theorem 3.8.**

1. *A key polynomial is a positive sum of Demazure atoms. In fact,*

$$\kappa_{\overline{\alpha}} = \sum_{\beta \geq \alpha} \mathcal{A}_\beta.$$

2. *A key polynomial with a partition shape* $\lambda$, *with* $l(\lambda) = k$, *is the Schur polynomial* $s_\lambda$, *i.e.,* $\kappa_{\overline{\lambda}} = s_\lambda(x_1, \ldots, x_k)$.

3. *The set of all Demazure atoms* $\{\mathcal{A}_\gamma : \gamma \vDash n, n \in \mathbb{Z}_{\geq 0}\}$ *forms a basis for the polynomial ring, and so does the set of all key polynomials* $\{\kappa_\gamma : \gamma \vDash n, n \in \mathbb{Z}_{\geq 0}\}$.

*Proof.* Item 1 follows directly from Theorem 3.5. We can describe combinatorially how to get the atoms from the key (a particular case of Proposition 6.1 in [Mas09]):

Consider a filling $F \in SSAF(\overline{\epsilon_k}, \alpha)$ and an empty filling $G_0$ with basement $\epsilon_k$. Consider the entries of the first row (from the bottom above the basement) of $F$, namely $a_{11} < a_{12} < \cdots < a_{1r_1}$ where $r_1$ is the length of the first row. Create the first row of $G_0$ by placing $a_{1i}$ in the cell right above $a_i$ in the basement of $G_0$ (that is, $a_{1i}$ is placed in the first row above the basement and also in the $a_{1i}^{th}$ column of $G$) for $1 \leq i \leq r_1$. We call the new filling $G_1$

Now consider the entries $a_{21} < a_{22} < \cdots < a_{2r_2}$ of the second row of $F$ where $r_2$ is the length of the second row of $F$. Search in the top row of $G_1$ for the leftmost number not less that $a_{2r_2}$ and place $a_{2r_2}$ in the cell right above it. Then search for the leftmost available number (i.e. not chosen yet) in the top row of $G_1$ not less than $a_{2,r_2-1}$ and place $a_{2,r_2-1}$ in the cell right above it, and so on until $a_{21}$ is placed. We now get a new filling with 2 rows above the basement and call it $G_2$.

By repeating the same process until all entries of $F$ are placed and we get a filling $G_r$ with basement $\epsilon_k$ with shape less than or equal to $\alpha$, where $r$ is the number of rows in $F$.

Item 2 follows from Theorem 3.7 as $s_\lambda = E_\lambda(x_1, \ldots, x_k; 0, 0)$ ([HHL08]). It is also proved in [LS90]. A combinatorial proof can be found in Theorem 4.1 in [HLMvW11] which uses the insertion algorithm discussed in [Mas08, Mas09].

Item 3 also follows from Theorem 3.7 as $\{E_\alpha(x; q, t) : \alpha \vDash n, n \in \mathbb{Z}_{\geq 0}\}$ forms a basis for the polynomial ring over $\mathbb{Q}(q, t)$. Again, it is also proved in [LS90].

□

# 4 Decomposition of products of Demazure atoms and characters

In this section, we study the decomposition of the products of Demazure atoms and characters with respect to the atom-basis $\{\mathcal{A}_\gamma : \gamma \vDash n, n \in \mathbb{Z}_{\geq 0}\}$ and key-basis $\{\kappa_\gamma : \gamma \vDash n, n \in \mathbb{Z}_{\geq 0}\}$.

Let $\lambda, \mu$ be partitions and $\alpha, \beta$ be weak compositions. Let $+_\mathcal{A}$ and $+_\kappa$ denote the property of being able to be decomposed into a positive sum of atoms and keys respectively. Note that by item 1 in Theorem 3.8, $+_\kappa$ implies $+_\mathcal{A}$. Otherwise, we put an $\times$ in the cell. For example, a partition ($\mu$)-shaped atom times a key of any shape ($\alpha$) is key positive and hence we put $+_\kappa$ in the corresponding box.

| | | Atoms | | Keys | |
|---|---|---|---|---|---|
| | shape | $\lambda$ | $\alpha$ | $\lambda$ | $\alpha$ |
| Atoms | $\mu$ | $+_\mathcal{A}$ | $+_\mathcal{A}$ ① | $+_\mathcal{A}$ | $+_\kappa$ ② |
| | $\beta$ | | $\times$ | $+_\mathcal{A}$ | $\times$ |
| Keys | $\mu$ | | | $+_\kappa$ | $+_\kappa$ |
| | $\beta$ | | | | open ③ |

Table 1: Decomposition of products of atoms and keys into atoms

The positive results in the table can be found in [HLMvW11], except for the cells marked ①, ② and ③.

We will prove ① (which was previously open) in this section using words and insertion algorithm introduced in [HLMvW11, Mas09]:

**Theorem 4.1.** *The product $\mathcal{A}_\mu \cdot \mathcal{A}_\alpha$ is atom-positive for any partition $\mu$ and weak composition $\alpha$.*

The coefficients in the decomposition into atoms are actually counting the number of ways to insert words arising from an SSAF of shape $\alpha$ into an SSAF of shape $\mu$ and we will discuss properties of words and how to record different ways of insertion in Section 4.1 and Section 4.2. Also note that the product in the theorem is not key positive. A simple counter example would be just putting $\mu$ as the empty partition, that is, with all entries 0 and $\alpha = (0, 1)$ and hence $\mathcal{A}_\mu \cdot \mathcal{A}_\alpha = \mathcal{A}_\alpha = \theta_1(x_1) = (\pi_1 - id)(x_1) = \kappa_{(0,1)} - \kappa_{(1,0)}$.

② is proved in [Jos03] (the proof involves crystals but does not involve SSAF). Both results ① and ② imply $+_{\mathcal{A}}$ for the $\mathcal{A}_\mu \cdot \kappa_{\overline{\alpha}}$ cell. We will apply the bijection in the proof of Theorem 6.1 in [HLMvW11] to Theorem 4.1 to give a tableau-combinatorial proof of ② in Section 4.4.

As for the product of two keys of arbitrary shapes, that is, the cell marked with ③, there are examples showing that such a product is not a positive sum of keys. For example, $\kappa_{(0,2)} \cdot \kappa_{(1,0,2)} = \kappa_{(1,2,2)} + \kappa_{(1,3,1)} + \kappa_{(1,4,0)} + \kappa_{(2,3,0)} + \kappa_{(3,0,2)} - \kappa_{(3,2,0)} + \kappa_{(4,0,1)} - \kappa_{(4,1,0)}$. Thus it remains to check whether it is a positive sum of atoms, which is still open. Hence ③ gives the following conjecture (first appearing in an unpublished work of Victor Reiner and Mark Shimozono).

**Conjecture 1.** *Let $\alpha, \beta$ be weak compositions. Then the product of the key polynomials of shape $\overline{\alpha}$ and $\overline{\beta}$ can be written as a positive sum of atoms, i.e.,*

$$\kappa_\alpha \cdot \kappa_\beta = \sum_{\gamma \vDash |\alpha|+|\beta|} c_{\alpha\beta}^{\gamma} \mathcal{A}_\gamma$$

*for some nonnegative integers $c_{\alpha\beta}^{\gamma}$.*

We will verify Conjecture 1 for $l(\alpha), l(\beta) \leq 3$ in Section 5.

### 4.1 Convert a column word to a row word

**Definition 4.1.** *A **word** is a sequence of positive integers.*

**Definition 4.2.** *Let $a, b, c \in \mathbb{N}$ and $u, v$ be some fixed (can be empty) words. Define twisted Knuth relation $\backsim^*$ by:*

1. *$ubacv \backsim^* ubcav$ if $c \leq b < a$*
2. *$uacbv \backsim^* ucabv$ if $c < b \leq a$.*

*Then we say two words $w$ and $w'$ are twisted Knuth equivalent if $w$ can be transformed to $w'$ by repeated use of 1. and 2. and we write $w \backsim^* w'$.*

**Definition 4.3.** *A word $w$ is a **column word** if it can be broken down into $k$ weakly decreasing subsequences of weakly decreasing lengths*

$$w = a_{11} \dots a_{1c_1} | a_{21} \dots a_{2c_2} | \cdots | a_{k1} \dots a_{kc_k}$$

*where $c_1 \geq c_2 \geq \cdots \geq c_k > 0$, $c_1, \dots, c_k \in \mathbb{N}$*
*such that* $\begin{cases} a_{ij} \geq a_{i,j+1} & \forall 1 \leq j < c_i, 1 \leq i \leq k \\ a_{i+1,c_{i+1}-j} > a_{i,c_i-j} & \forall 0 \leq j < c_i, 1 \leq i < k \end{cases}.$

**Definition 4.4.** *A word $w$ is a **row word** if it can be broken down into $k$ strictly increasing subsequences of weakly decreasing lengths*

$$w = a_{11} \dots a_{1r_1} | a_{21} \dots a_{2r_2} | \cdots | a_{k1} \dots a_{kr_k}$$

*where $r_1 \geq r_2 \geq \cdots \geq r_k > 0$, $r_1, \dots, r_k \in \mathbb{N}$*
*such that* $\begin{cases} a_{ij} < a_{i,j+1} & \forall 1 \leq j < r_i, 1 \leq i \leq k \\ a_{i+1,r_{i+1}-j} \leq a_{i,c_i-j} & \forall 0 \leq j < r_i, 1 \leq i < k \end{cases}.$

Given an SSAF, one can get its column word by using the algorithm described in [HLMvW11], while one can get its row word (which is the reverse reading word defined in [HMR13] by reading the entries of each row in ascending order, starting from the bottom row to the top row. We call a word a column (resp. row) word because when we insert each subsequence of the word using the insertion in [Mas08], a new column (resp. row) will be created.

**Example 7.** 886531|97643|9764|5|6 *is a column word whose corresponding SSAF is:*

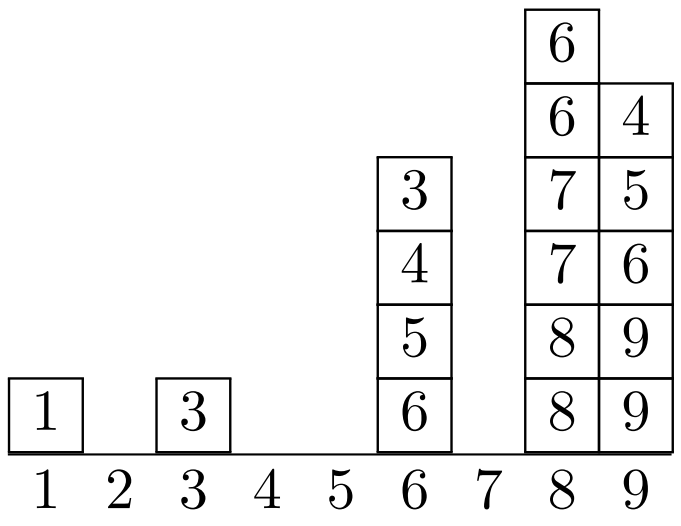

*and the corresponding row word is* 13689|589|467|357|46|6*.*

We describe an insertion algorithm for inserting an integer $c \leq k$ into an SSAF with basement $\epsilon_k$. A detailed description can be found in [Mas08].

Given an SSAF $F$ with basement $\epsilon_k$ and an integer $c$ to be inserted, we write $c \to F$ or $F \leftarrow c$ to represent the insertion of $c$ into $F$ (and similarly we denote $c_2c_1 \to F$ and $F \leftarrow c_1c_2$ as first inserting $c_1$ to $F$ and then insert $c_2$, and so on).

To insert $c$, we first find the cell $A$ in $F$ with the smallest order, say $m$, such that $F(A') < c \leq F(A)$, where $A'$ is the cell immediately above $A$ if it exists and assign $F(A') = 0$ if $A$ is the top cell of a column and just treat $A'$ as an empty cell to be filled in. If cell $A$ is the top cell of a column, then we create a new cell immediately above $A$, i.e. $A'$, and assign $c$ to the new cell and we are done. Otherwise, we replace the entry $y = F(A')$ by $c$ and now insert $y$ as we do for $c$, but now finding a cell $B$ of the smallest order **larger than** $m$ such that $F(B') < y \leq F(B)$ as treating $B$ as $A$ in the previous step. Repeat the process until a new cell is finally created.

**Example 8.** *Let $F$ be the SSAF in Example 7 where $k = 9$. We find the new SSAF created by inserting* 7*, i.e.* $7 \to F$*, as follows:*

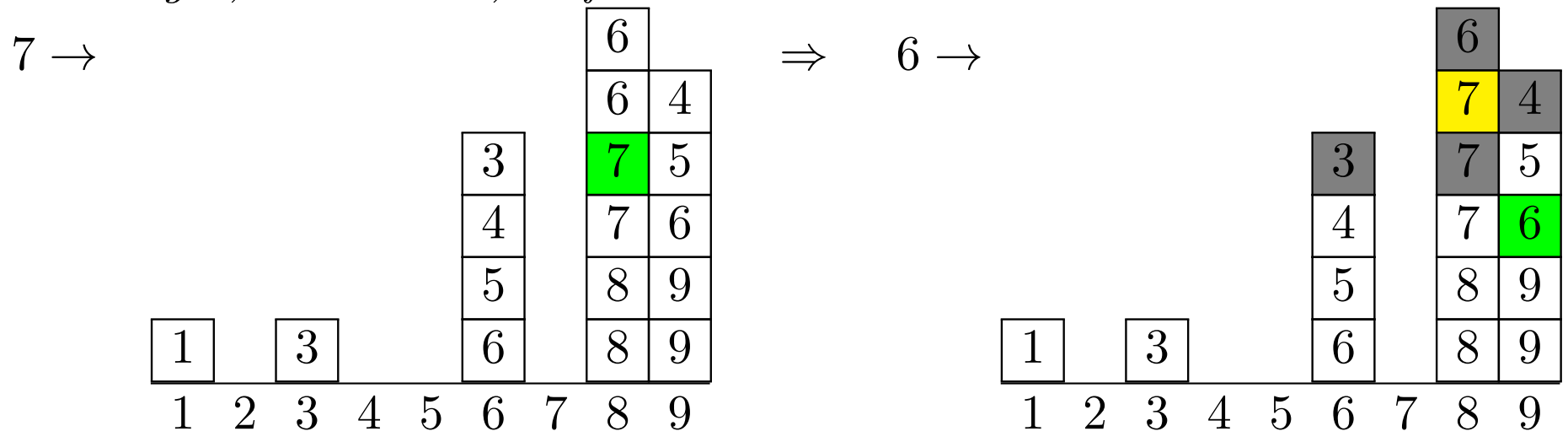

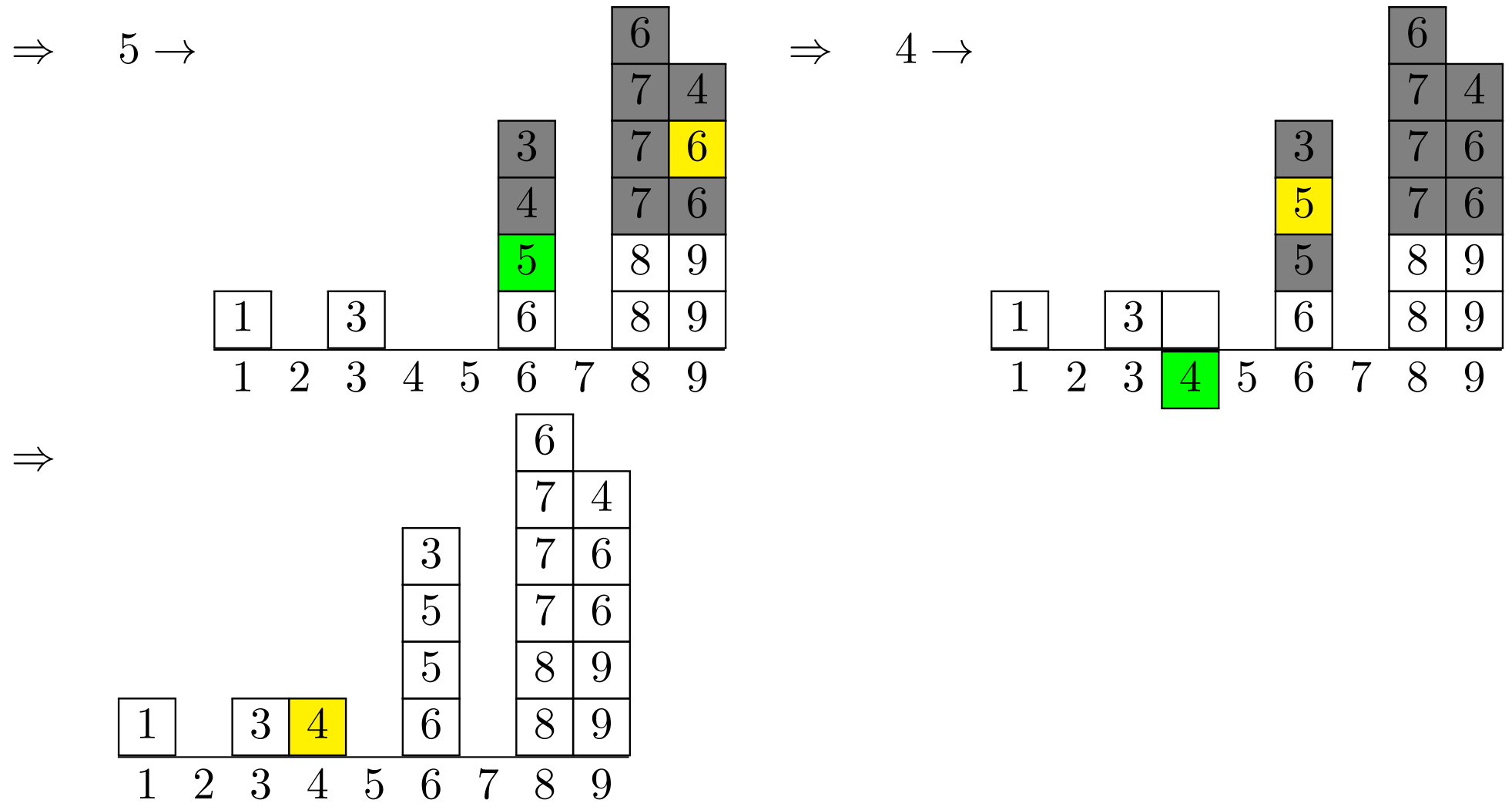

*The green cell represents the cell $A$ in each step. The yellow cell represents the cell $A'$ whose value is changed after the insertion step (i.e. the original entry is bumped out by the number being inserted). The white cells are the cells under consideration in each step (i.e. candidates for the position of $B$, that is, the new $A$).*

**Lemma 4.2.** *Let $u := a_1a_2\dots a_n|b$ and $v := a_1a_2\dots a_n|b|c$ be two column words, where $c > b > a_n$, $n \in \mathbb{N}$. Then*

1. *$u \backsim^* b'a_1'...a_n'$ where $b' < a_1'$, $a_1' \geq a_2' \geq \cdots \geq a_n'$, and $b' := a_t$ where $t = \min\{j : b > a_j\}$*
2. *$v \backsim^* b'c'a_1''a_2''\dots a_n''$ where $b'$ is defined as in 1., $b' < c' < a_1''$ and $a_1'' \geq a_2'' \geq \cdots \geq a_n''$.*

*Proof.* One can check that for $t = \min\{j : b > a_j\}$,

$$u \backsim^* \begin{cases} a_1ba_2\dots a_n & \text{if} \quad t = 1 \\ a_ta_1\dots a_{t-1}ba_{t+1}\dots a_n & \text{if} \quad 1 < t < n \\ a_na_1\dots a_{n-1}b & \text{if} \quad t = n \end{cases}$$

and *1.* follows. Also note that $b' = a_t$.

Indeed, we can visualize by applying the insertion in [Mas08]. When we insert the first $n$ integers of $u$, we get a column above the basement entry $a_1$:

If $b > a_1$, we have

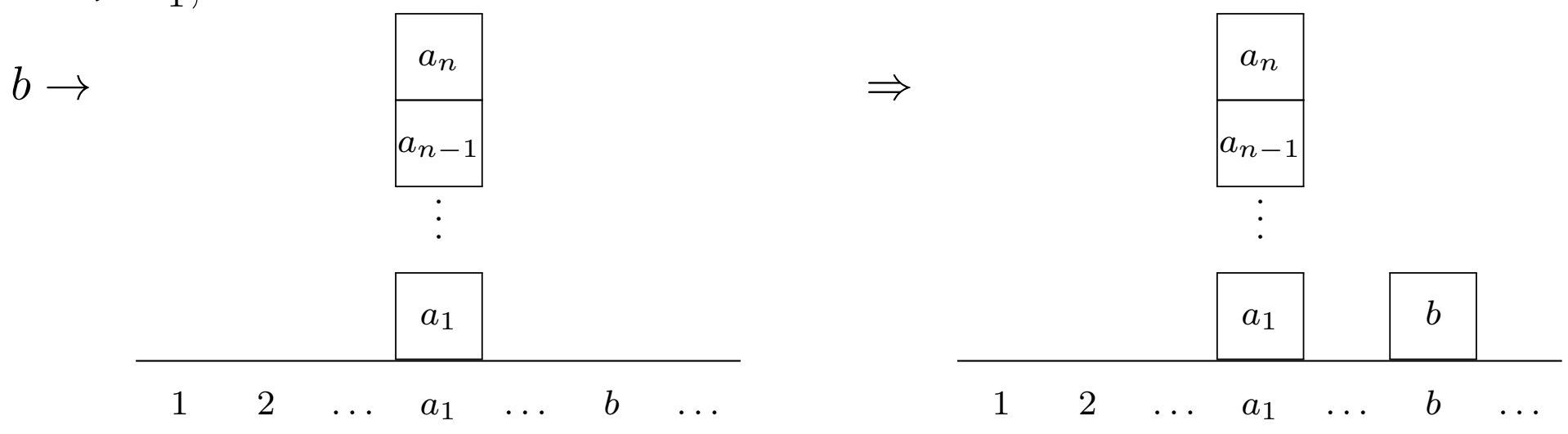

giving the row word $a_1ba_2\dots a_n$.

If $b \leq a_1$, then for $t = \min\{j : b > a_j\}$, we have:
For $t < n$

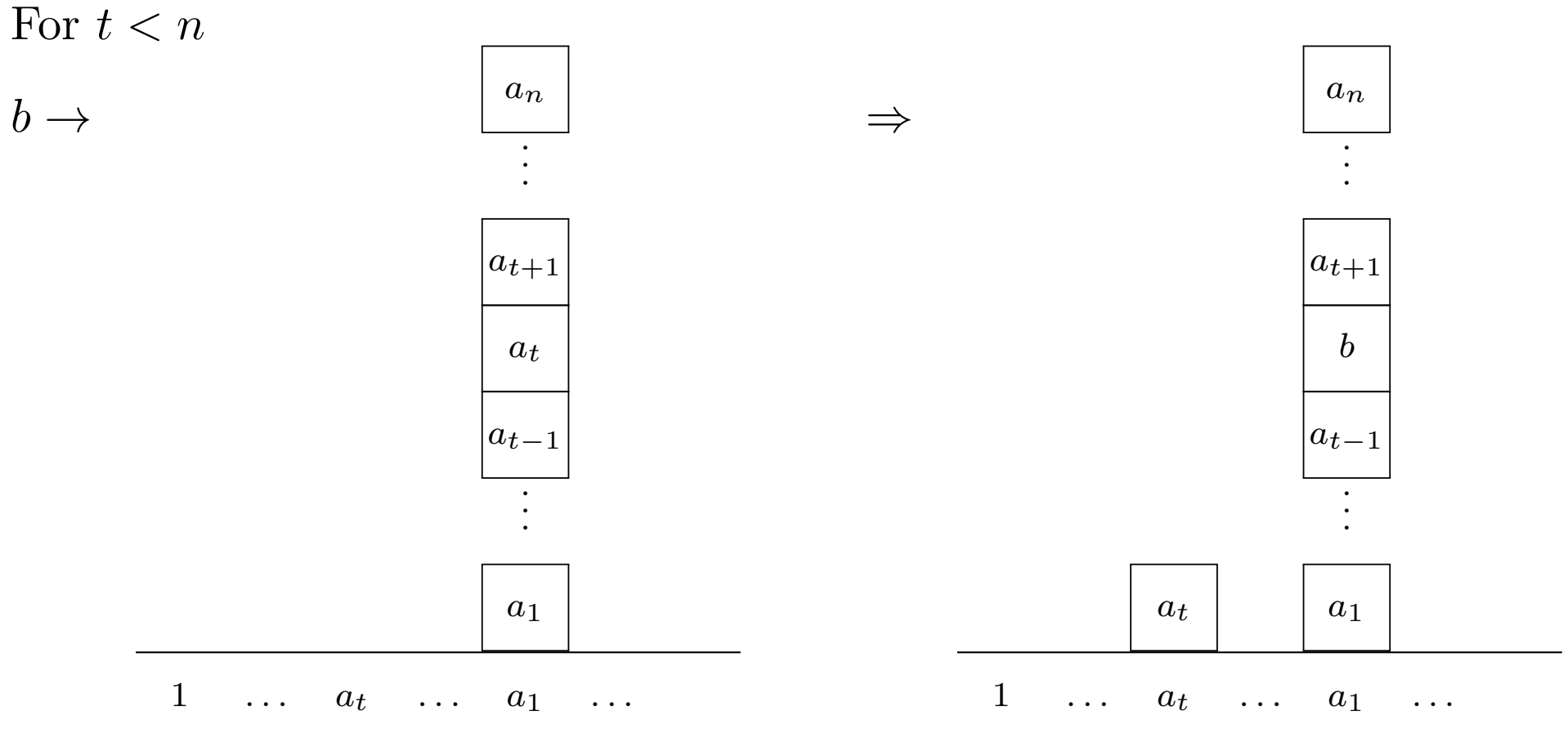

giving the row word $a_t a_1 \ldots a_{t-1} b a_{t+1} \ldots a_n$.

For $t = n$

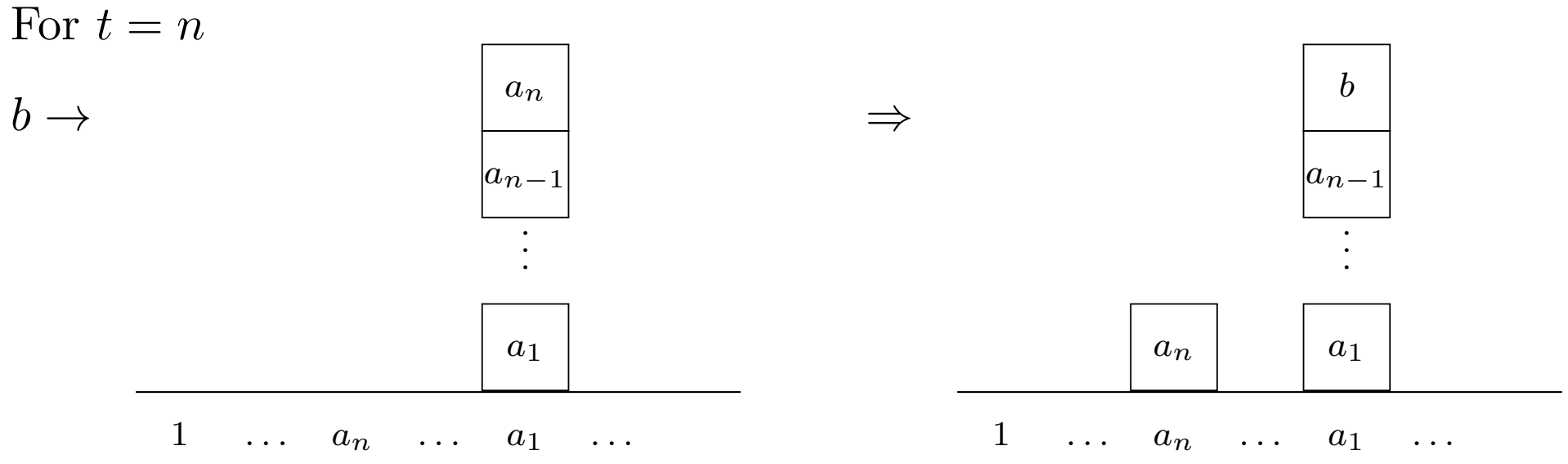

giving the row word $a_n a_1 \ldots a_{n-1} b$.

Since inserting $v$ to an empty atom is the same as inserting $c$ to the atom created by inserting $u$ to an empty atom, we can inset $c$ to the tableaux above and get the following:

If $c > b > a_1$, we have

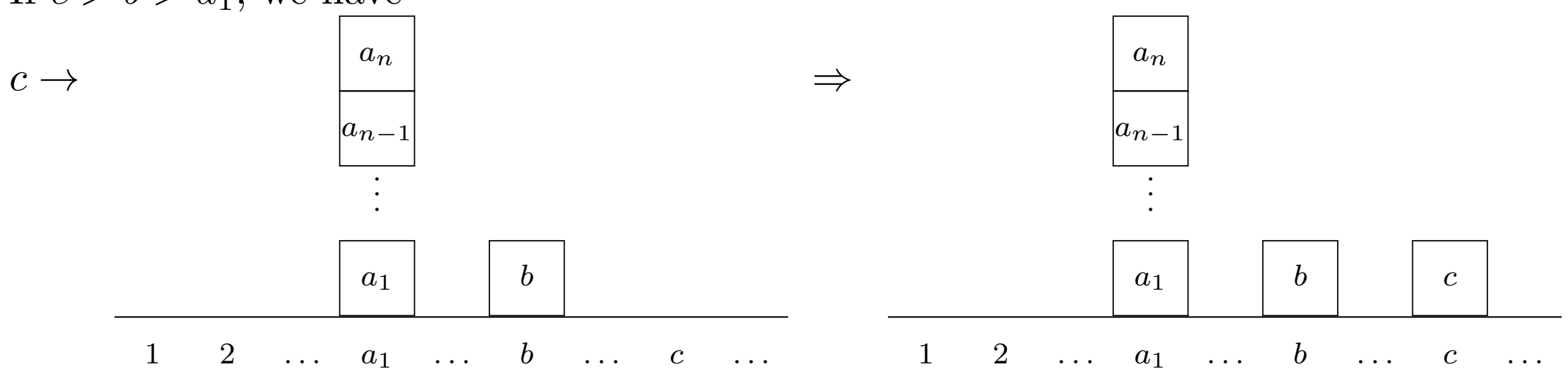

giving the row word $a_1 b c a_2 \ldots a_n$.

If $b \leq a_1$ and $a_{t-1} \geq c > b > a_t$, we have:
For $t < n$

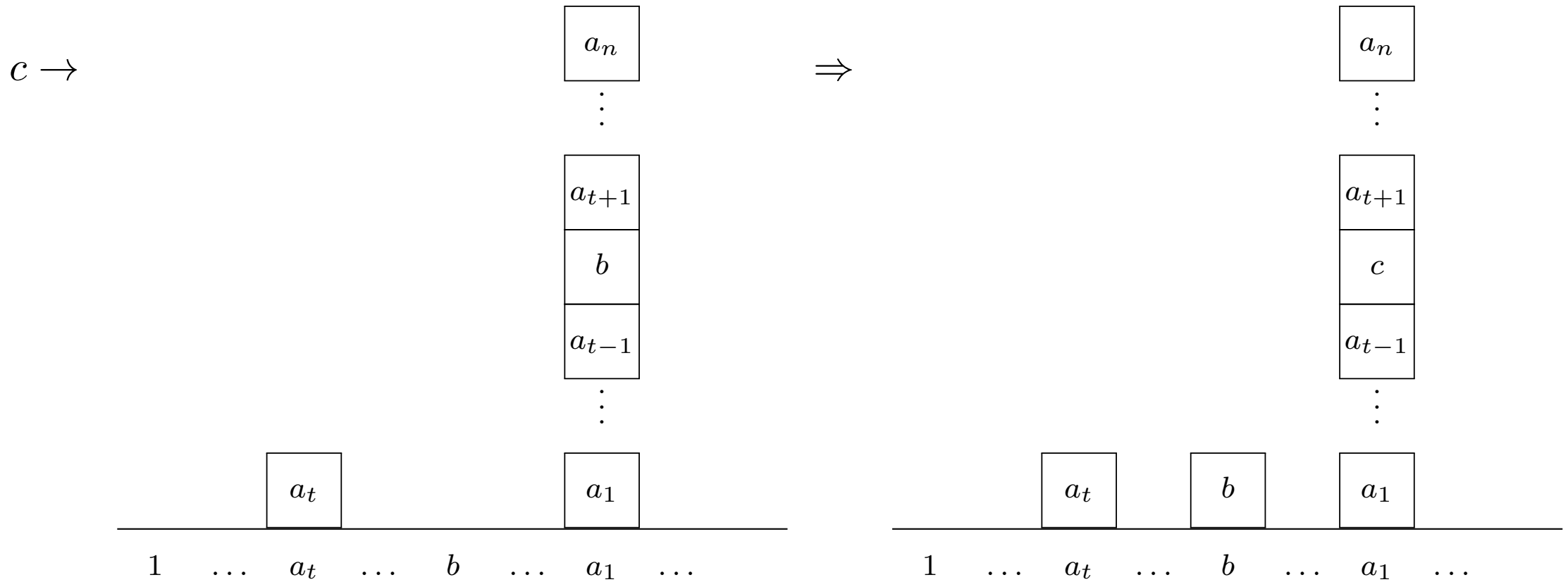

giving the row word $a_t b a_1 \dots a_{t-1} c a_{t+1} \dots a_n$.

For $t = n$

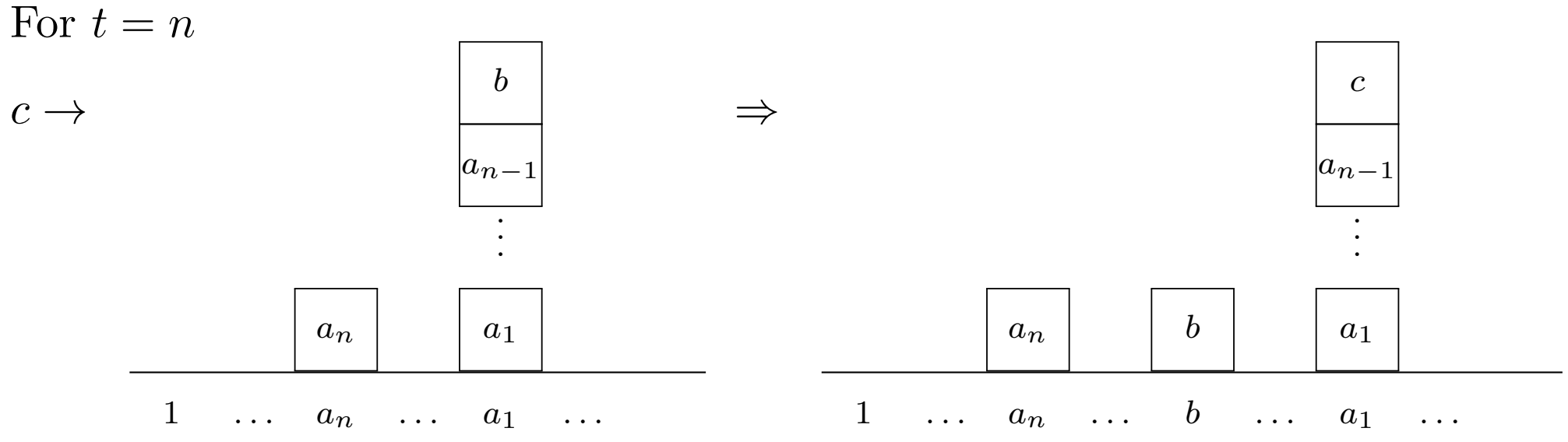

giving the row word $a_n b a_1 \dots a_{n-1} c$.
If $b \leq a_1$ and $c > a_{t-1}$, we have:
For $t < n$, let $t' = \min\{j : c > a_j, 1 \leq j \leq t-1\}$, we have:
For $t' = 1$,

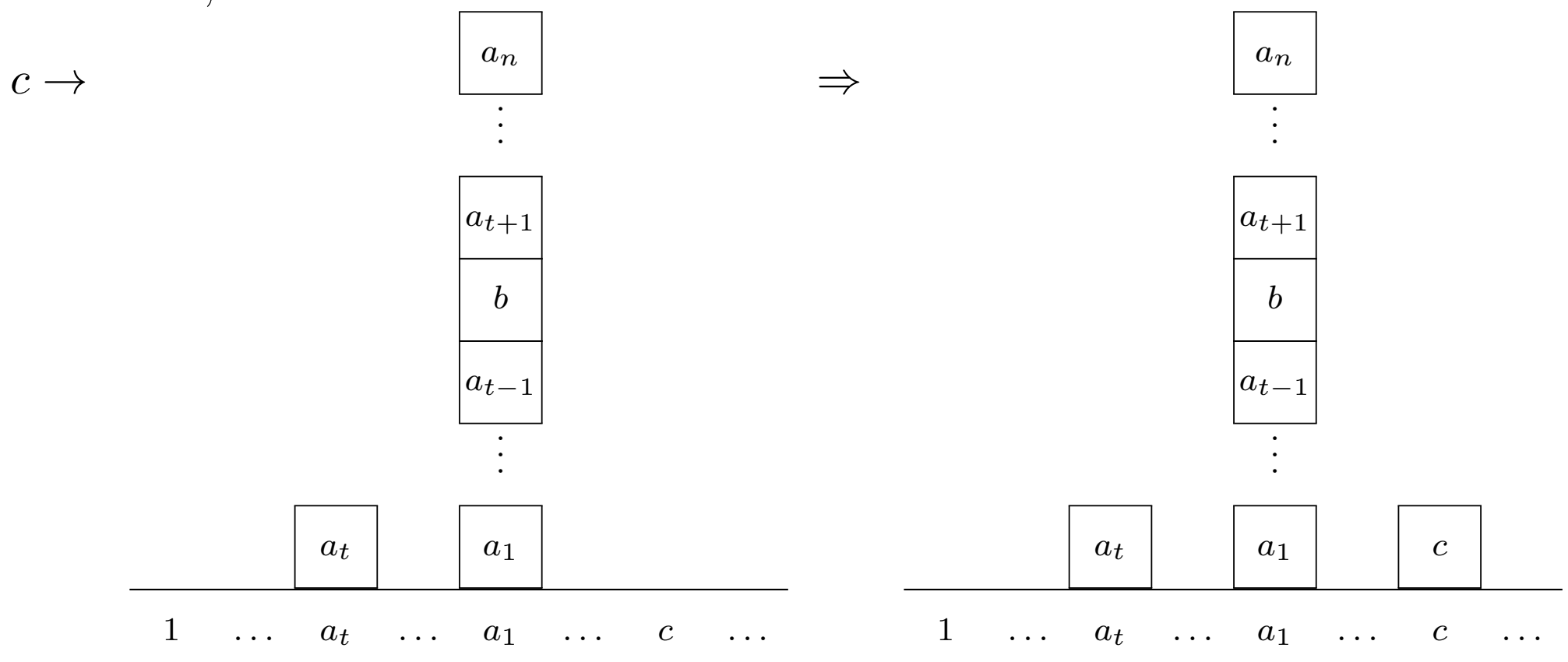

giving the row word $a_t a_1 c a_2 \dots a_{t-1} b a_{t+1} \dots a_n$.

For $1 < t' < t-1$,

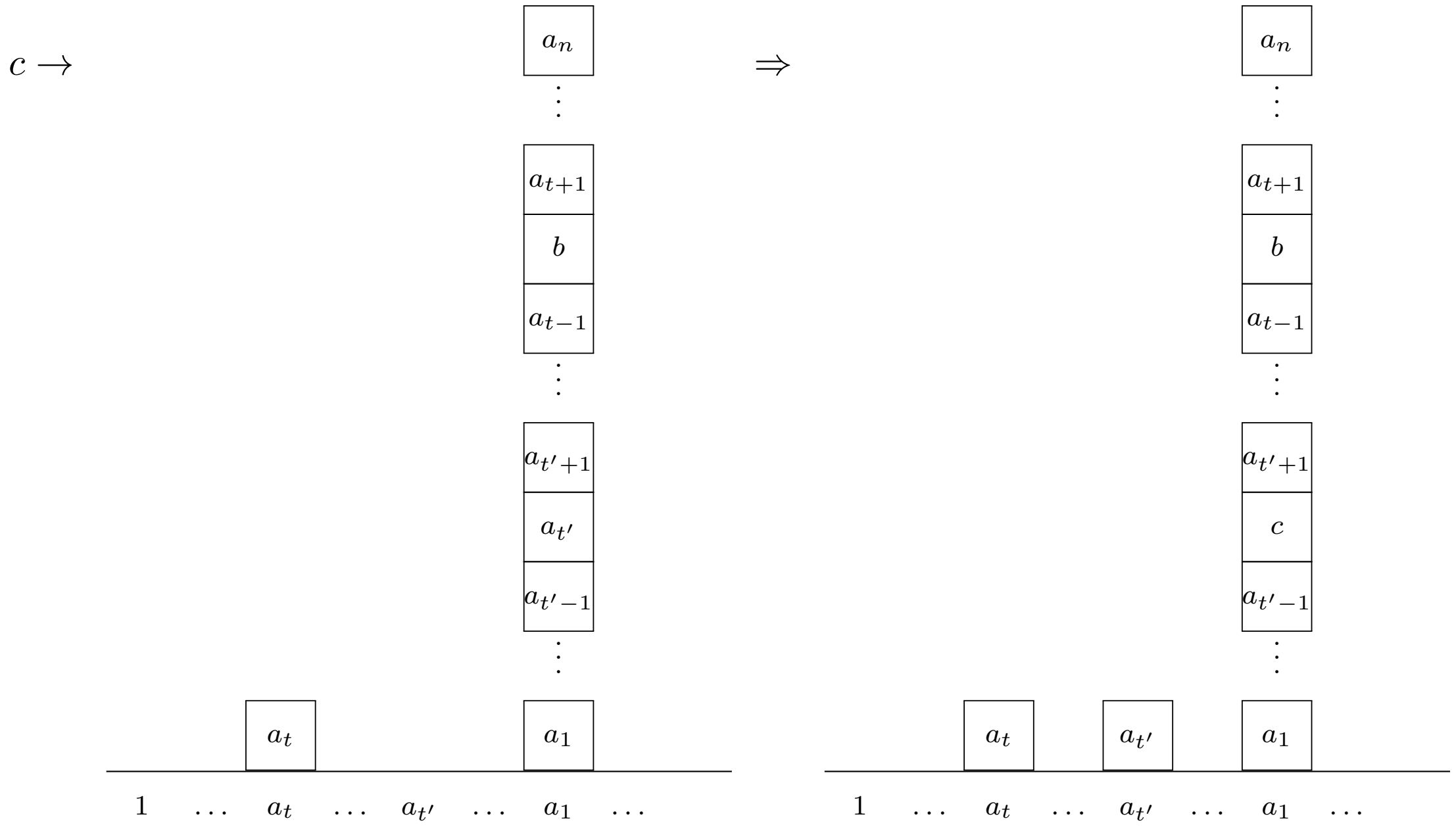

giving the row word $a_t a_{t'} a_1 \dots a_{t'-1} c a_{t'+1} \dots a_{t-1} b a_{t+1} \dots a_n$.

For $t' = t-1$,

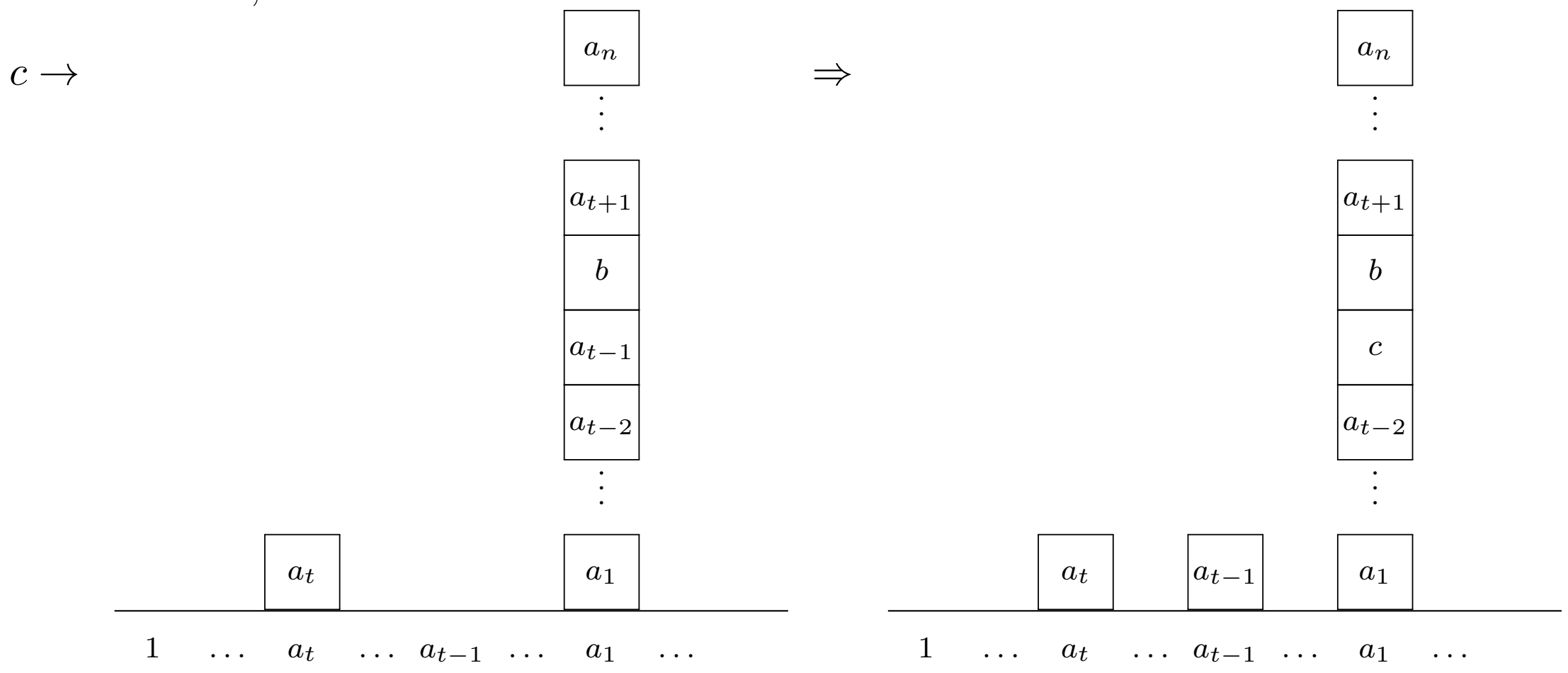

giving the row word $a_t a_{t-1} a_1 \dots a_{t-2} c b a_{t+1} \dots a_n$.

For $t = n$, let $t' = \min\{j : c > a_j\}$, we have:
For $t' = 1$,

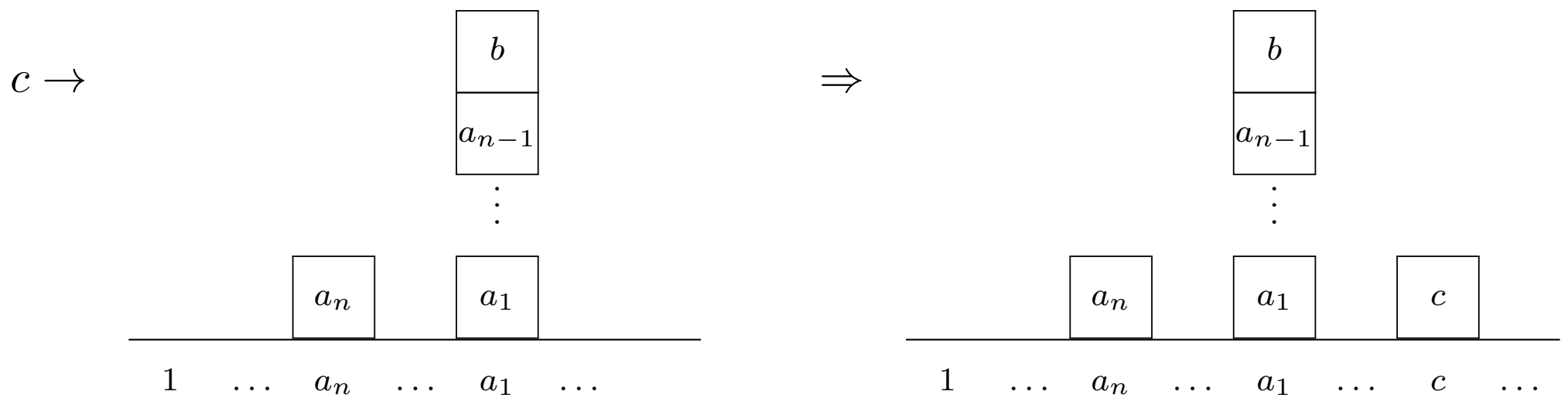

giving the row word $a_n a_1 c a_2 \dots a_{n-1} b$.
For $1 < t' < n-1$,

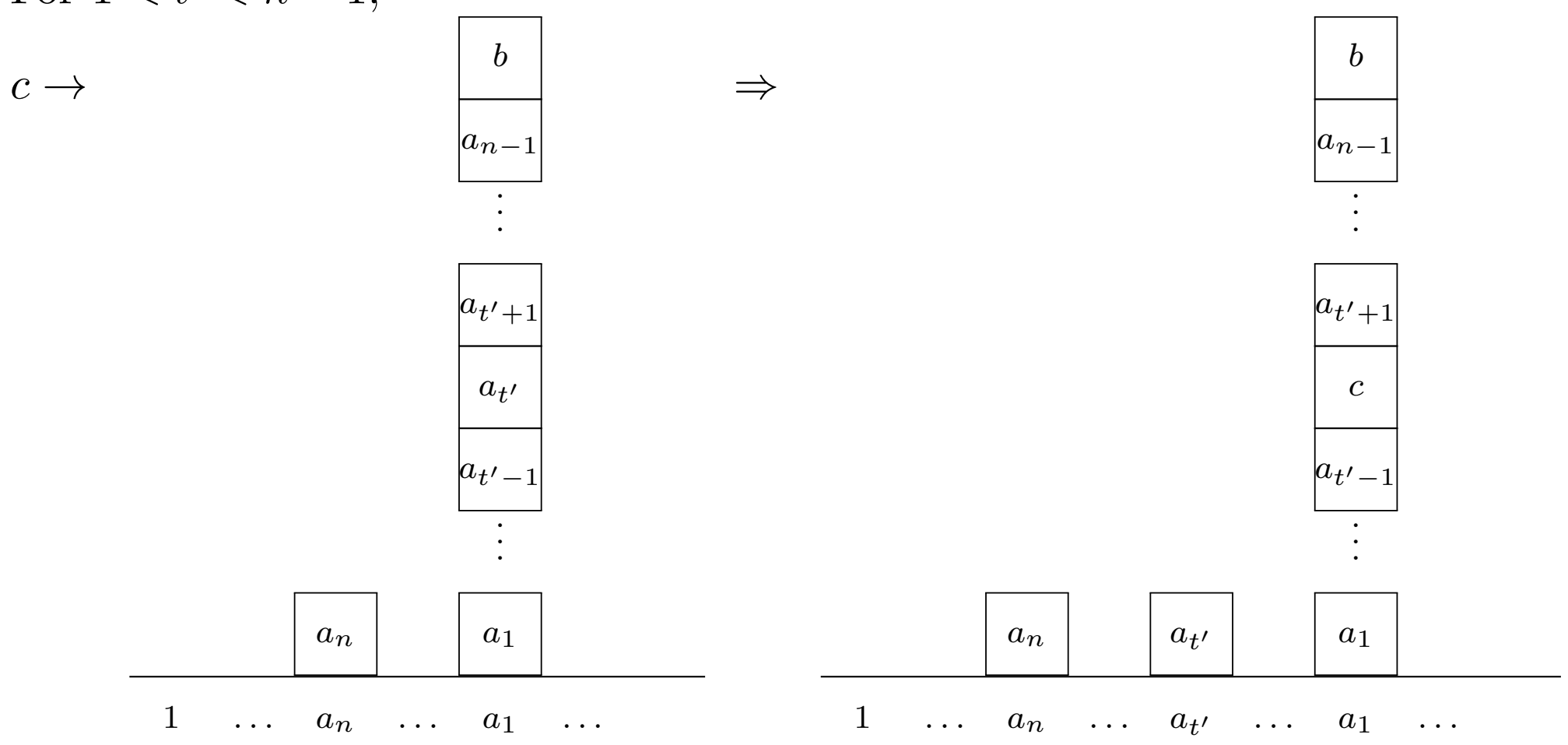

giving the row word $a_n a_{t'} a_1 a_2 \dots a_{t'-1} c a_{t'+1} \dots a_{n-1} b$.
For $t' = n-1$,

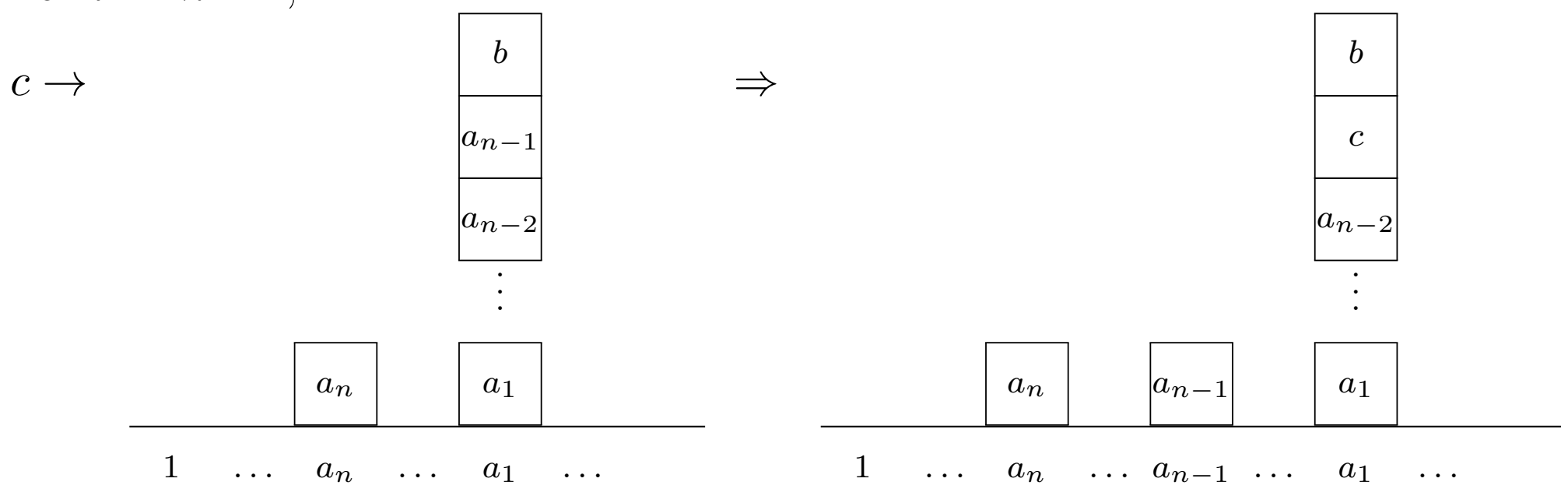

giving the row word $a_n a_{n-1} a_1 a_2 \dots a_{n-2} c b$.

Notice that for all the cases of inserting $c$, the first row has the least entry unaffected, meaning that the row word has the same first entry after inserting $c$. The result follows by setting the row word of the tableau of inserting $v$ into an empty atom as $b'c'a''_1 \dots a''_n$. $\square$

**Example 9.** *One can also prove Lemma 4.2 by repeated use of 1. and 2. in Definition 4.2, for example:*

$9764|5 \backsim^* 97465 \backsim^* 94765 \backsim^* 4|9765$ *by using 2. repeatedly, and similarly, we have* $9764|5|6 \backsim^* 4|9765|6 \backsim^* 45|9766$.

*Using SSAF, it means*

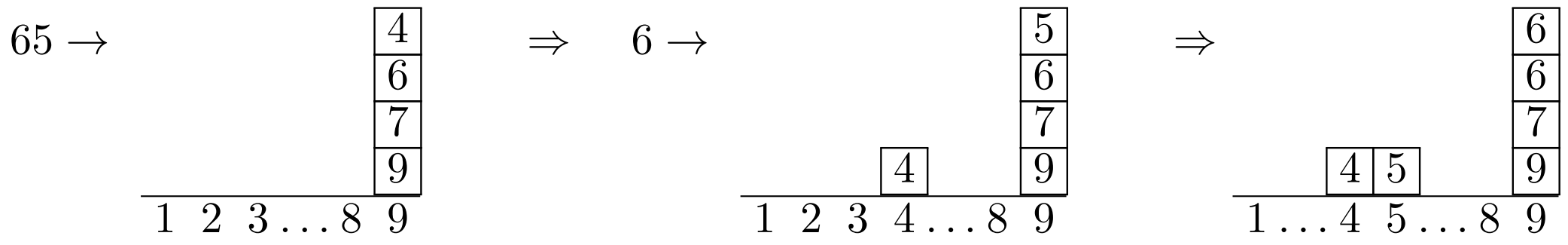

*where* $65 \to$ *means* 5 *is inserted before* 6.

**Lemma 4.3.** *Let* $k \geq 2$ *be an integer and* $i_1, \dots, i_k$ *be nonnegative integers. Let*

$$w_0 = a_{11}^{(i_1)} \dots a_{1c_1}^{(i_1)} | a_{21}^{(i_2)} \dots a_{2c_2}^{(i_2)} | \cdots | a_{k1}^{(i_k)} \dots a_{kc_k}^{(i_k)}$$

*be a column word. Then there exist* $a_{rs}^{(i_r+1)}$ *where* $1 \leq s \leq c_r$, $1 \leq r < k$ *and* $a_{k1}^{(i_k+m)}$ *where* $1 \leq m < k$ *satisfying*

$$\begin{cases} a_{k-j-1,c_{k-j-1}}^{(i_{k-j-1})} < a_{k1}^{(i_k+j)} < a_{k-j,1}^{(i_{k-j}+1)} & \forall 1 \leq j < k-1 \quad (k > 3) \\ a_{k1}^{(i_k+k-1)} < a_{11}^{(i_1+1)} & \\ a_{1,1}^{(i_1+1)} \dots a_{1,c_1}^{(i_1+1)} | \cdots | a_{k-1,1}^{(i_{k-1}+1)} \dots a_{k-1,c_{k-1}}^{(i_{k-1}+1)} | a_{k2}^{(i_k)} \dots a_{kc_k}^{(i_k)} & \textit{is a column word} \end{cases}$$

*such that* $w_0 \backsim^* w_1 \backsim^* \cdots \backsim^* w_{k-1}$, *where*

$$\begin{cases} w_j := a_{11}^{(i_1)} \dots a_{1c_1}^{(i_1)} \cdots a_{k-j-1,1}^{(i_{k-j-1})} \dots a_{k-j-1,c_{k-j-1}}^{(i_{k-j-1})} a_{k1}^{(i_k+j)} a_{k-j,1}^{(i_{k-j}+1)} \dots a_{k-j,c_{k-j}}^{(i_{k-j}+1)} \cdots \\ \qquad \cdots a_{k-1,1}^{(i_{k-1}+1)} \dots a_{k-1,c_{k-1}}^{(i_{k-1}+1)} a_{k2}^{(i_k)} \dots a_{kc_k}^{(i_k)} \qquad \forall 1 \leq j < k-1, \\ \\ w_{k-1} := a_{k1}^{i_k+k-1} a_{11}^{(i_1+1)} \dots a_{1c_1}^{(i_1+1)} \cdots a_{k-1,1}^{(i_{k-1}+1)} \dots a_{k-1,c_{k-1}}^{(i_{k-1}+1)} a_{k2}^{(i_k)} \dots a_{kc_k}^{(i_k)}. \end{cases}$$

*Proof.* We prove this by induction on $k$.

When $k = 2$, we have $w_0 := a_{11}^{(i_1)} a_{12}^{(i_1)} \dots a_{1c_1}^{(i_1)} a_{21}^{(i_2)} a_{22}^{(i_2)} \dots a_{2c_2}^{(i_2)}$.

Let $t_{12} = \min_{1 \leq j \leq c_1} \{ j : a_{21}^{(i_2)} > a_{ij}^{(i_1)} \}$.

By applying the proof of *1.* in Lemma 4.2 on $a_{11}^{(i_1)} a_{12}^{(i_1)} \dots a_{1c_1}^{(i_1)} a_{21}^{(i_2)}$, we have

$$\begin{aligned} w_0 \backsim^* w_1 :=& a_{21}^{(i_2+1)} a_{11}^{(i_1+1)} a_{12}^{(i_1+1)} \dots a_{1c_1}^{(i_1+1)} a_{22}^{(i_2)} \dots a_{2c_2}^{(i_2)} \\ :=& \begin{cases} a_{11}^{(i_1)} a_{21}^{(i_2)} a_{12}^{(i_1)} \dots a_{1c_1}^{(i_1)} a_{22}^{(i_2)} \dots a_{2c_2}^{(i_2)} & \text{if} \quad t_{12} = 1 \\ a_{t_{12}}^{(i_1)} a_{11}^{(i_1)} \dots a_{1,t_{12}-1}^{(i_1)} a_{21}^{(i_2)} a_{1,t_{12}+1}^{(i_1)} \dots a_{1c_1}^{(i_1)} a_{22}^{(i_2)} \dots a_{2c_2}^{(i_2)} & \text{if} \quad t_{12} > 1 \end{cases}. \end{aligned}$$

where $a_{21}^{(i_2+1)} < a_{11}^{(i_1+1)}$ and $a_{11}^{(i_1+1)} \geq \cdots \geq a_{1c_1}^{(i_1+1)}$. Also note that $a_{21}^{(i_2+1)} = a_{1t_{12}}^{(i_1)}$.

For $t_{12} = 1$, we know that $a_{12}^{(i_1)} \dots a_{1c_1}^{(i_1)} a_{22}^{(i_2)} \dots a_{2c_2}^{(i_2)}$ is a column word and hence $a_{12}^{(i_1+1)} \dots a_{1c_1}^{(i_1+1)} a_{22}^{(i_2)} \dots a_{2c_2}^{(i_2)}$.

For $t_{12} > 1$, if $c_2 - (c_1 - t_{12} + 1) > 0$, then

$$a_{21}^{(i_2)} < a_{1,t_{12}-1}^{(i_1)} = a_{1,c_1-(c_1-t_{12}+1)}^{(i_1)} < a_{2,c_2-(c_1-t_{12}+1)}^{(i_2)} \text{ (by Definition 4.3)}$$

which leads to a contradiction as $a_{21}^{(i_2)} \geq a_{2j}^{(i_2)}$ for $1 \leq j \leq c_2$.

As a result, we have $c_2 - (c_1 - t_{12} + 1) \leq 0$ which implies $c_2 - 1 \leq c_1 - t_{12}$.

Thus $a_{1,t_{12}+1}^{(i_1)} \dots a_{1c_1}^{(i_1)} a_{22}^{(i_2)} \dots a_{2c_2}^{(i_2)}$ is a column word, and hence $a_{12}^{(i_1)} \dots a_{1,t_{12}-1}^{(i_1)} a_{21}^{(i_2)} a_{1,t_{12}+1}^{(i_1)} \dots a_{1c_1}^{(i_1)} a_{22}^{(i_2)} \dots a_{2c_2}^{(i_2)}$ is also a column word.

Therefore $a_{12}^{(i_1+1)} \dots a_{1c_1}^{(i_1+1)} a_{22}^{(i_2)} \dots a_{2c_2}^{(i_2)}$ is a column word.

Note that $a_{11}^{(i_1+1)} \geq a_{12}^{(i_1+1)}$, we can conclude $a_{11}^{(i_1+1)} \dots a_{1c_1}^{(i_1+1)} a_{22}^{(i_2)} \dots a_{2c_2}^{(i_2)}$ is a column word.

Hence the statement is true for $k = 2$.

For $k = 3$, we have $w_0 := a_{11}^{(i_1)} a_{12}^{(i_1)} \dots a_{1c_1}^{(i_1)} a_{21}^{(i_2)} a_{22}^{(i_2)} \dots a_{2c_2}^{(i_2)} a_{31}^{(i_3)} a_{32}^{(i_3)} \dots a_{3c_3}^{(i_3)}$.

Consider $a_{21}^{(i_2)} a_{22}^{(i_2)} \dots a_{2c_2}^{(i_2)} a_{31}^{(i_3)} a_{32}^{(i_3)} \dots a_{3c_3}^{(i_3)}$, by $k = 2$ case, we know

$$a_{21}^{(i_2)} a_{22}^{(i_2)} \dots a_{2c_2}^{(i_2)} a_{31}^{(i_3)} a_{32}^{(i_3)} \dots a_{3c_3}^{(i_3)} \backsim^* a_{31}^{(i_3+1)} a_{21}^{(i_2+1)} a_{22}^{(i_2+1)} \dots a_{2c_2}^{(i_2+1)} a_{32}^{(i_3)} \dots a_{3c_3}^{(i_3)}$$

where $a_{31}^{(i_3+1)} < a_{21}^{(i_2+1)}$ and $a_{21}^{(i_2+1)} \dots a_{2c_2}^{(i_2+1)} a_{32}^{(i_3)} \dots a_{3c_3}^{(i_3)}$ is a column word. Also, $a_{31}^{(i_3+1)} = a_{2t_{23}}^{(i_2)} \geq a_{2c_2}^{(i_2)} > a_{1c_1}^{(i_1)}$, $t_{23} := \min\limits_{1 \leq j \leq c_2} \{j : a_{31}^{(i_3)} > a_{2j}^{(i_2)}\}$.

Hence, $a_{1c_1}^{(i_1)} < a_{31}^{(i_3+1)} < a_{21}^{(i_2+1)}$.

Furthermore, $a_{11}^{(i_1)} a_{12}^{(i_1)} \dots a_{1c_1}^{(i_1)} a_{31}^{(i_3+1)}$ is also a column word and by *1.* in Lemma 4.2,

$$a_{11}^{(i_1)} a_{12}^{(i_1)} \dots a_{1c_1}^{(i_1)} a_{31}^{(i_3+1)} \backsim^* a_{31}^{(i_3+2)} a_{11}^{(i_1+1)} a_{12}^{(i_1+1)} \dots a_{1c_1}^{(i_1+1)}$$

where $a_{31}^{(i_3+2)} < a_{11}^{(i_1+1)}$ and $a_{11}^{(i_1+1)} a_{12}^{(i_1+1)} \dots a_{1c_1}^{(i_1+1)}$ is a column word for $a_{11}^{(i_1+1)} \geq a_{12}^{(i_1+1)} \geq \dots \geq a_{1c_1}^{(i_1+1)}$. Also $a_{31}^{(i_3+2)} = a_{1t_{13}}^{(i_1)}$ where $t_{13} := \min\limits_{1 \leq j \leq c_1} \{j : a_{31}^{(i_3+1)} > a_{1j}^{(i_1)}\}$.

It remains to check that $a_{11}^{(i_1+1)} a_{12}^{(i_1+1)} \dots a_{1c_1}^{(i_1+1)} a_{21}^{(i_2+1)} a_{22}^{(i_2+1)} \dots a_{2c_2}^{(i_2+1)}$ is a column word.

First note that by *1.* in Lemma 4.2,

$$a_{21}^{(i_2+1)} a_{22}^{(i_2+1)} \dots a_{2c_2}^{(i_2+1)} = \begin{cases} a_{31}^{(i_3)} a_{22}^{(i_2)} \dots a_{2c_2}^{(i_2)} & \text{if} \quad t_{23} = 1 \\ a_{21}^{(i_2)} \dots a_{2,t_{23}-1}^{(i_2)} a_{31}^{(i_3)} a_{2,t_{23}+1}^{(i_2)} \dots a_{2c_2}^{(i_2)} & \text{if} \quad t_{23} > 1 \end{cases},$$

$$\begin{aligned} a_{11}^{(i_1+1)} a_{12}^{(i_1+1)} \dots a_{1c_1}^{(i_1+1)} &= \begin{cases} a_{31}^{(i_3+1)} a_{12}^{(i_1)} \dots a_{1c_1}^{(i_1)} & \text{if} \quad t_{13} = 1 \\ a_{11}^{(i_1)} \dots a_{1,t_{13}-1}^{(i_1)} a_{31}^{(i_3+1)} a_{1,t_{13}+1}^{(i_1)} \dots a_{1c_1}^{(i_1)} & \text{if} \quad t_{13} > 1 \end{cases} \\ &= \begin{cases} a_{2t_{23}}^{(i_2)} a_{12}^{(i_1)} \dots a_{1c_1}^{(i_1)} & \text{if} \quad t_{13} = 1 \\ a_{11}^{(i_1)} \dots a_{1,t_{13}-1}^{(i_1)} a_{2t_{23}}^{(i_2)} a_{1,t_{13}+1}^{(i_1)} \dots a_{1c_1}^{(i_1)} & \text{if} \quad t_{13} > 1 \end{cases}. \end{aligned}$$

By definition of a column word, we have $a_{2t_{23}}^{(i_2)} > a_{1,c_1-(c_2-t_{23})}^{(i_1)}$. By definition of $t_{13}$, we have $t_{13} \leq c_1 - (c_2 - t_{23})$ which implies $c_2 - t_{23} \leq c_1 - t_{13}$.

Combining, we have the following four cases:

(i) $\boldsymbol{t_{13} = t_{23} = 1}$ **:**

$a_{11}^{(i_1+1)} a_{12}^{(i_1+1)} \ldots a_{1c_1}^{(i_1+1)} a_{21}^{(i_2+1)} a_{22}^{(i_2+1)} \ldots a_{2c_2}^{(i_2+1)} = a_{21}^{(i_2)} a_{12}^{(i_1)} \ldots a_{1c_1}^{(i_1)} a_{31}^{(i_3)} a_{22}^{(i_2)} \ldots a_{2c_2}^{(i_2)}$.

Since $a_{31}^{(i_3)} > a_{21}^{(i_2)} = a_{11}^{(i_1+1)} \geq a_{1j}^{(i_1+1)} = a_{1j}^{(i_1)}$ for $1 < j \leq c_1$, and $a_{12}^{(i_1)} \ldots a_{1c_1}^{(i_1)} a_{22}^{(i_2)} \ldots a_{2c_2}^{(i_2)}$ is a column word, we can conclude that $a_{21}^{(i_2)} a_{12}^{(i_1)} \ldots a_{1c_1}^{(i_1)} a_{31}^{(i_3)} a_{22}^{(i_2)} \ldots a_{2c_2}^{(i_2)}$ is also a column word.

(ii) $\boldsymbol{t_{13} = 1, t_{23} > 1}$ **:**

$a_{11}^{(i_1+1)} \ldots a_{1c_1}^{(i_1+1)} a_{21}^{(i_2+1)} \ldots a_{2c_2}^{(i_2+1)}$
$= a_{2t_{23}}^{(i_2)} a_{12}^{(i_1)} \ldots a_{1c_1}^{(i_1)} a_{21}^{(i_2)} \ldots a_{2,t_{23}-1}^{(i_2)} a_{31}^{(i_3)} a_{2,t_{23}+1}^{(i_2)} \ldots a_{2c_2}^{(i_2)}$.

Since $a_{31}^{(i_3)} > a_{2t_{23}}^{(i_2)}$ and $a_{12}^{(i_1)} \ldots a_{1c_1}^{(i_1)} a_{22}^{(i_2)} \ldots a_{2,t_{23}-1}^{(i_2)} a_{2t_{23}}^{(i_2)} a_{2,t_{23}+1}^{(i_2)} \ldots a_{2c_2}^{(i_2)}$ is a column word, we know $a_{12}^{(i_1)} \ldots a_{1c_1}^{(i_1)} a_{22}^{(i_2)} \ldots a_{2,t_{23}-1}^{(i_2)} a_{31}^{(i_3)} a_{2,t_{23}+1}^{(i_2)} \ldots a_{2c_2}^{(i_2)}$ is a column word.

Also $a_{21}^{(i_2)} \geq a_{2,t_{23}-1}^{(i_2)} \geq a_{31}^{(i_3)} > a_{2t_{23}}^{(i_3)} = a_{11}^{(i_1+1)} \geq a_{1j}^{(i_1+1)} = a_{1j}^{(i_1)}$ for $1 < j \leq c_1$,

thus $a_{2t_{23}}^{(i_2)} a_{12}^{(i_1)} \ldots a_{1c_1}^{(i_1)} a_{21}^{(i_2)} \ldots a_{2,t_{23}-1}^{(i_2)} a_{31}^{(i_3)} a_{2,t_{23}+1}^{(i_2)} \ldots a_{2c_2}^{(i_2)}$ is a column word.

(iii) $\boldsymbol{t_{13} > 1, t_{23} = 1}$ **:**

$a_{11}^{(i_1+1)} a_{12}^{(i_1+1)} \ldots a_{1c_1}^{(i_1+1)} a_{21}^{(i_2+1)} a_{22}^{(i_2+1)} \ldots a_{2c_2}^{(i_2+1)}$
$= a_{11}^{(i_1)} \ldots a_{1,t_{13}-1}^{(i_1)} a_{21}^{(i_2)} a_{1,t_{13}+1}^{(i_1)} \ldots a_{1c_1}^{(i_1)} a_{31}^{(i_3)} a_{22}^{(i_2)} \ldots a_{2c_2}^{(i_2)}$.

Since $c_2 - 1 \leq c_1 - t_{13}$ which implies $c_2 \leq c_1 - (t_{13} - 1)$, we just need to consider

$$a_{21}^{(i_2)} a_{1,t_{13}+1}^{(i_1)} \ldots a_{1c_1}^{(i_1)} a_{31}^{(i_3)} a_{22}^{(i_2)} \ldots a_{2c_2}^{(i_2)}$$

to conclude that

$$a_{11}^{(i_1)} \ldots a_{1,t_{13}-1}^{(i_1)} a_{21}^{(i_2)} a_{1,t_{13}+1}^{(i_1)} \ldots a_{1c_1}^{(i_1)} a_{31}^{(i_3)} a_{22}^{(i_2)} \ldots a_{2c_2}^{(i_2)}$$

is a column word.

As $c_2 - 1 \leq c_1 - t_{13}$, $a_{1,t_{13}+1}^{(i_1)} \ldots a_{1c_1}^{(i_1)} a_{22}^{(i_2)} \ldots a_{2c_2}^{(i_2)}$ is a column word. Also $a_{31}^{(i_3)} > a_{21}^{(i_2)}$ for $t_{23} = 1$, and $a_{21}^{(i_2)} = a_{1t_{13}}^{(i_1+1)} \geq a_{1j}^{(i_1+1)} = a_{1j}^{(i_1)}$ for $t_{13} < j \leq c_1$ and hence $a_{21}^{(i_2)} a_{1,t_{13}+1}^{(i_1)} \ldots a_{1c_1}^{(i_1)} a_{31}^{(i_3)} a_{22}^{(i_2)} \ldots a_{2c_2}^{(i_2)}$ is a column word and result follows.

(iv) $\boldsymbol{t_{13} > 1, t_{23} > 1}$ **:**

$a_{11}^{(i_1+1)} a_{12}^{(i_1+1)} \ldots a_{1c_1}^{(i_1+1)} a_{21}^{(i_2+1)} a_{22}^{(i_2+1)} \ldots a_{2c_2}^{(i_2+1)}$
$= a_{11}^{(i_1)} \ldots a_{1,t_{13}-1}^{(i_1)} a_{2t_{23}}^{(i_2)} a_{1,t_{13}+1}^{(i_1)} \ldots a_{1c_1}^{(i_1)} a_{21}^{(i_2)} \ldots a_{2,t_{23}-1}^{(i_2)} a_{31}^{(i_3)} a_{2,t_{23}+1}^{(i_2)} \ldots a_{2c_2}^{(i_2)}$.

Since $a_{11}^{(i_1)} \ldots a_{1,t_{13}-1}^{(i_1)} a_{1t_{13}}^{(i_1)} a_{1,t_{13}+1}^{(i_1)} \ldots a_{1c_1}^{(i_1)} a_{21}^{(i_2)} \ldots a_{2,t_{23}-1}^{(i_2)} a_{2t_{23}}^{(i_2)} a_{2,t_{23}+1}^{(i_2)} \ldots a_{2c_2}^{(i_2)}$ is a column word, it remains to check if $a_{1,c_1-(c_2-t_{23})}^{(i_1+1)} < a_{2t_{23}}^{(i_2+1)}$ and $a_{1t_{13}}^{(i_1+1)} < a_{2,c_2-(c_1-t_{13})}^{(i_2+1)}$.
(We need to check the latter condition only when $c_2 - (c_1 - t_{13}) > 0$.)

If $c_2 - t_{23} = c_1 - t_{13}$, then

$$a_{1,c_1-(c_2-t_{23})}^{(i_1+1)} = a_{1t_{13}}^{(i_1+1)} = a_{2t_{23}}^{(i_2)} < a_{31}^{(i_3)} = a_{2t_{23}}^{(i_2+1)} = a_{2,c_2-(c_1-t_{13})}^{(i_2+1)}.$$

If $c_2 - t_{23} < c_1 - t_{13}$, then $c_2 - (c_1 - t_{13}) < t_{23}$.
Hence $a^{(i_1+1)}_{1t_{13}} = a^{(i_2)}_{2t_{23}} < a^{(i_2)}_{2,c_2-(c_1-t_{13})} = a^{(i_2+1)}_{2,c_2-(c_1-t_{13})}$.

Similarly, $t_{13} < c_1 - (c_2 - t_{23})$ implies

$$a^{(i_1+1)}_{1,c_1-(c_2-t_{23})} = a^{(i_1)}_{1,c_1-(c_2-t_{23})} < a^{(i_2)}_{2,c_2-(c_2-t_{23})} = a^{(i_2)}_{2t_{23}} < a^{(i_3)}_{31} = a^{(i_2+1)}_{2t_{23}}.$$

Set $\begin{cases} w_1 := a^{(i_1)}_{11} a^{(i_1)}_{12} \dots a^{(i_1)}_{1c_1} a^{(i_3+1)}_{31} a^{(i_2+1)}_{21} a^{(i_2+1)}_{22} \dots a^{(i_2+1)}_{2c_2} a^{(i_3)}_{32} \dots a^{(i_3)}_{3c_3} \\ w_2 := a^{(i_3+2)}_{31} a^{(i_1+1)}_{11} a^{(i_1+1)}_{12} \dots a^{(i_1+1)}_{1c_1} a^{(i_2+1)}_{21} a^{(i_2+1)}_{22} \dots a^{(i_2+1)}_{2c_2} a^{(i_3)}_{32} \dots a^{(i_3)}_{3c_3} \end{cases}$,
and we thus have $w_0 \backsim^* w_1 \backsim^* w_2$ and result follows.

Therefore the statement is true for $k = 3$.

Assume the statement is true for all $k = 2, \dots, m, m+1$ for some $m \geq 2$.

When $k = m + 2$,

By considering $a^{(i_m)}_{m,1} \dots a^{(i_m)}_{m,c_m} a^{(i_{m+1})}_{m+1,1} \dots a^{(i_{m+1})}_{m+1,c_{m+1}} a^{(i_{m+2})}_{m+2,1} \dots a^{(i_{m+2})}_{m+2,c_{m+2}}$ and apply the result in $k = 3$, we get

$$\begin{aligned} & a^{(i_1)}_{11} \dots a^{(i_1)}_{1c_1} a^{(i_2)}_{21} \dots a^{(i_2)}_{2c_2} \cdots a^{(i_{m+1})}_{m+1,1} \dots a^{(i_{m+1})}_{m+1,c_{m+1}} a^{(i_{m+2})}_{m+2,1} \dots a^{(i_{m+2})}_{m+2,c_{m+2}} \\ \backsim^* & a^{(i_1)}_{11} \dots a^{(i_1)}_{1c_1} a^{(i_2)}_{21} \dots a^{(i_2)}_{2c_2} \cdots a^{(i_{m+2}+1)}_{m+2,1} a^{(i_{m+1}+1)}_{m+1,1} \dots a^{(i_{m+1}+1)}_{m+1,c_{m+1}} a^{(i_{m+2})}_{m+2,2} \dots a^{(i_{m+2})}_{m+2,c_{m+2}} \\ \backsim^* & a^{(i_1)}_{11} \dots a^{(i_1)}_{1c_1} a^{(i_2)}_{21} \dots a^{(i_2)}_{2c_2} \cdots a^{(i_{m+2}+2)}_{m+2,1} a^{(i_m+1)}_{m,1} \dots \\ & \cdots a^{(i_m+1)}_{m,c_m} a^{(i_{m+1}+1)}_{m+1,1} \dots a^{(i_{m+1}+1)}_{m+1,c_{m+1}} a^{(i_{m+2})}_{m+2,2} \dots a^{(i_{m+2})}_{m+2,c_{m+2}} \end{aligned}$$

where $a^{(i_m+1)}_{m,1} \dots a^{(i_m+1)}_{m,c_m} a^{(i_{m+1}+1)}_{m+1,1} \dots a^{(i_{m+1}+1)}_{m+1,c_{m+1}} a^{(i_{m+2})}_{m+2,2} \dots a^{(i_{m+2})}_{m+2,c_{m+2}}$ is a column word, and $a^{(i_m)}_{m,c_m} < a^{(i_{m+2}+1)}_{m+2,1} < a^{(i_{m+1}+1)}_{m+1,1}$.

As $a^{(i_m)}_{m,c_m} < a^{(i_{m+2}+1)}_{m+2,1}$, we can apply the induction assumption on

$$a^{(i_1)}_{11} \dots a^{(i_1)}_{1c_1} \cdots a^{(i_m)}_{m1} \dots a^{(i_m)}_{mc_2} a^{(i_{m+2}+1)}_{m+2,1}$$

and hence we have

$$\begin{aligned} & a^{(i_1)}_{11} \dots a^{(i_1)}_{1c_1} a^{(i_2)}_{21} \dots a^{(i_2)}_{2c_2} \cdots a^{(i_m)}_{m1} \dots a^{(i_m)}_{mc_m} a^{(i_{m+2}+1)}_{m+2,1} \\ \backsim^* & a^{(i_1)}_{11} \dots a^{(i_1)}_{1c_1} a^{(i_2)}_{21} \dots a^{(i_2)}_{2c_2} \cdots a^{(i_{m-1})}_{m-1,1} \dots a^{(i_{m-1})}_{m-1,c_{m-1}} a^{(i_{m+2}+2)}_{m+2,1} a^{(i_m+1)}_{m1} \dots a^{(i_m+1)}_{mc_m} \\ \backsim^* & a^{(i_1)}_{11} \dots a^{(i_1)}_{1c_1} a^{(i_2)}_{21} \dots a^{(i_2)}_{2c_2} \cdots a^{(i_{m-2})}_{m-2,1} \dots a^{(i_{m-2})}_{m-2,c_{m-2}} a^{(i_{m+2}+3)}_{m+2,1} a^{(i_{m-1}+1)}_{m-1,1} \cdots \\ & \cdots a^{(i_{m-1}+1)}_{m-1,c_{m-1}} a^{(i_m+1)}_{m1} \dots a^{(i_m+1)}_{mc_m} \\ & \vdots \\ \backsim^* & a^{(i_{m+2}+m+1)}_{m+2,1} a^{(i_1+1)}_{11} \dots a^{(i_1+1)}_{1c_1} a^{(i_2+1)}_{21} \dots a^{(i_2+1)}_{2c_2} \cdots a^{(i_m+1)}_{m1} \dots a^{(i_m+1)}_{mc_m} \end{aligned}$$

where $\begin{cases} a^{(i_{m-j})}_{m-j,c_{m-j}} < a^{(i_{m+2}+1+j)}_{m+2,1} < a^{(i_{m+1-j}+1)}_{m+1-j,1} & \forall 1 \leq j < m \\ a^{(i_{m+2}+m+1)}_{m+2,1} < a^{(i_1+1)}_{11} & \\ a^{(i_1+1)}_{11} \dots a^{(i_1+1)}_{1c_1} a^{(i_2+1)}_{21} \dots a^{(i_2+1)}_{2c_2} \cdots a^{(i_m+1)}_{m1} \dots a^{(i_m+1)}_{mc_m} & \text{is a column word} \end{cases}$.

Combining the above results, we have

$$\begin{aligned}
& a^{(i_1)}_{11} \dots a^{(i_1)}_{1c_1} \cdots a^{(i_{m+1})}_{m+1,1} \dots a^{(i_{m+1})}_{m+1,c_{m+1}} a^{(i_{m+2})}_{m+2,1} \dots a^{(i_{m+2})}_{m+2,c_{m+2}} \\
\backsim^* & a^{(i_1)}_{11} \dots a^{(i_1)}_{1c_1} \cdots a^{(i_m)}_{m1} \dots a^{(i_m)}_{mc_m} a^{(i_{m+2}+1)}_{m+2,1} a^{(i_{m+1}+1)}_{m+1,1} \dots a^{(i_{m+1}+1)}_{m+1,c_{m+1}} a^{(i_{m+2})}_{m+2,2} \dots a^{(i_{m+2})}_{m+2,c_{m+2}} \\
\backsim^* & a^{(i_1)}_{11} \dots a^{(i_1)}_{1c_1} \cdots a^{(i_{m-1})}_{m-1,1} \dots a^{(i_{m-1})}_{m-1,c_{m-1}} a^{(i_{m+2}+2)}_{m+2,1} a^{(i_m+1)}_{m,1} \dots a^{(i_m+1)}_{m,c_m} a^{(i_{m+1}+1)}_{m+1,1} \dots \\
& \dots a^{(i_{m+1}+1)}_{m+1,c_{m+1}} a^{(i_{m+2})}_{m+2,2} \dots a^{(i_{m+2})}_{m+2,c_{m+2}} \\
\backsim^* & a^{(i_1)}_{11} \cdots a^{(i_{m-2})}_{m-2,c_{m-2}} a^{(i_{m+2}+3)}_{m+2,1} a^{(i_{m-1}+1)}_{m-1,1} \dots \\
& \dots a^{(i_{m-1}+1)}_{m-1,c_{m-1}} a^{(i_m+1)}_{m1} \dots a^{(i_m+1)}_{mc_m} a^{(i_{m+1}+1)}_{m+1,1} \dots a^{(i_{m+1}+1)}_{m+1,c_{m+1}} a^{(i_{m+2})}_{m+2,2} \dots a^{(i_{m+2})}_{m+2,c_{m+2}} \\
& \qquad \vdots \\
\backsim^* & a^{(i_{m+2}+m+1)}_{m+2,1} a^{(i_1+1)}_{11} \dots a^{(i_1+1)}_{1c_1} \cdots a^{(i_m+1)}_{m1} \dots a^{(i_m+1)}_{mc_m} a^{(i_{m+1}+1)}_{m+1,1} \dots \\
& \dots a^{(i_{m+1}+1)}_{m+1,c_{m+1}} a^{(i_{m+2})}_{m+2,2} \dots a^{(i_{m+2})}_{m+2,c_{m+2}}
\end{aligned}$$

Both $a^{(i_1+1)}_{11} \dots a^{(i_1+1)}_{1c_1} a^{(i_2+1)}_{21} \dots a^{(i_2+1)}_{2c_2} \cdots a^{(i_m+1)}_{m1} \dots a^{(i_m+1)}_{mc_m}$ and $a^{(i_m+1)}_{m,1} \dots a^{(i_m+1)}_{m,c_m} a^{(i_{m+1}+1)}_{m+1,1} \dots a^{(i_{m+1}+1)}_{m+1,c_{m+1}} a^{(i_{m+2})}_{m+2,2} \dots a^{(i_{m+2})}_{m+2,c_{m+2}}$ are column words, so $a^{(i_1+1)}_{11} \dots a^{(i_1+1)}_{1c_1} \cdots a^{(i_m+1)}_{m1} \dots a^{(i_m+1)}_{mc_m} a^{(i_{m+1}+1)}_{m+1,1} \dots a^{(i_{m+1}+1)}_{m+1,c_{m+1}} a^{(i_{m+2})}_{m+2,2} \dots a^{(i_{m+2})}_{m+2,c_{m+2}}$ is a column word.

Also,
$\begin{cases} a^{(i_{m-j})}_{m-j,c_{m-j}} < a^{(i_{m+2}+1+j)}_{m+2,1} < a^{(i_{m+1-j}+1)}_{m+1-j,1} & \forall 1 \leq j < m \\ a^{(i_m)}_{m,c_m} < a^{(i_{m+2}+1)}_{m+2,1} < a^{(i_{m+1}+1)}_{m+1,1} & \end{cases}$
implies $a^{(i_{m+1-j})}_{m+1-j,c_{m+1-j}} < a^{(i_{m+2}+j)}_{m+2,1} < a^{(i_{m+2-j}+1)}_{m+2-j,1} \quad \forall 1 \leq j < m+1$.
Therefore the statement is true for $k = m+2$.
By Mathematical Induction, the statement is true for all integers $k \geq 2$.

$\square$

Lemma 4.3 shows how to create the first entry of the corresponding row word, leaving the remaining part as a column word . We illustrate the Lemma by an example with $a^{(i_k+j)}_{k1}$ circled for all $0 \leq j \leq k-1$.

**Example 10.** *Let $w_0$ be the word in Example 7, hence $k = 5$, $c_1 = 6, c_2 = 5, c_3 = 4, c_4 = c_5 = 1$. Set $i_j = 0$ for all $j$. Then we have*

$$\begin{cases} w_0 = a_{11}^{(0)} \dots a_{16}^{(0)} | a_{21}^{(0)} \dots a_{25}^{(0)} | a_{31}^{(0)} \dots a_{34}^{(0)} | a_{41}^{(0)} | \textcircled{a_{51}^{(0)}} = 886531|97643|9764|5|\textcircled{6} \\ w_1 = a_{11}^{(0)} \dots a_{16}^{(0)} | a_{21}^{(0)} \dots a_{25}^{(0)} | a_{31}^{(0)} \dots a_{34}^{(0)} | \textcircled{a_{51}^{(1)}} | a_{41}^{(1)} = 886531|97643|9764|\textcircled{5}|6 \\ w_2 = a_{11}^{(0)} \dots a_{16}^{(0)} | a_{21}^{(0)} \dots a_{25}^{(0)} | \textcircled{a_{51}^{(2)}} | a_{31}^{(1)} \dots a_{34}^{(1)} | a_{41}^{(1)} = 886531|97643|\textcircled{4}|9765|6 \\ w_3 = a_{11}^{(0)} \dots a_{16}^{(0)} | \textcircled{a_{51}^{(3)}} | a_{21}^{(1)} \dots a_{25}^{(1)} | a_{31}^{(1)} \dots a_{34}^{(1)} | a_{41}^{(1)} = 886531|\textcircled{3}|97644|9765|6 \\ w_4 = \textcircled{a_{51}^{(4)}} | a_{11}^{(1)} \dots a_{16}^{(1)} | a_{21}^{(1)} \dots a_{25}^{(1)} | a_{31}^{(1)} \dots a_{34}^{(1)} | a_{41}^{(1)} = \textcircled{1}|886533|97644|9765|6 \end{cases} .$$

**Lemma 4.4.** *Let $w = a_{11} \dots a_{1c_1} \cdots a_{k1} \dots a_{kc_k}$ be a column word. Then there exists a sequence $\{b_{ij}\}_{1 \le j \le c_i, 1 \le i \le k}$ such that*

$$w \backsim^* b_{k1} b_{k-1,1} \dots b_{11} b_{12} \dots b_{1c_1} b_{22} \dots b_{2c_2} \cdots b_{k2} \dots b_{kc_k},$$

*where $b_{11} > b_{21} > \cdots > b_{k1}$ and $b_{11} \dots b_{1c_1} | b_{22} \dots b_{2c_2} | \cdots | b_{k2} \dots b_{kc_k}$ is a column word (with lengths $c_1 > c_2 - 1 \ge \cdots \ge c_k - 1$).*

*Proof.* Apply Lemma 4.3 on $a_{11} \dots a_{1c_2} \cdots a_{k1} \dots a_{kc_k}$, we have

$a_{11} \dots a_{1c_1} a_{21} \dots a_{2c_2} \cdots a_{k1} \dots a_{kc_k}$
$\backsim^* a_{k1}^{(k-1)} a_{11}^{(1)} \dots a_{1c_1}^{(1)} \cdots a_{k-1,1}^{(1)} \dots a_{k-1,c_{k-1}}^{(1)} a_{k2} \dots a_{kc_k}$
where $a_{11}^{(1)} \dots a_{1c_1}^{(1)} \cdots a_{k-1,1}^{(1)} \dots a_{k-1,c_{k-1}}^{(1)} a_{k2} \dots a_{kc_k}$ is a column word.

Apply Lemma 4.3 on $a_{11}^{(1)} \dots a_{1c_1}^{(1)} \cdots a_{k-1,1}^{(1)} \dots a_{k-1,c_{k-1}}^{(1)}$, we have

$a_{11}^{(1)} \dots a_{1c_1}^{(1)} \cdots a_{k-1,1}^{(1)} \dots a_{k-1,c_{k-1}}^{(1)}$
$\backsim^* a_{k-1,1}^{(k-1)} a_{11}^{(2)} \dots a_{1c_1}^{(2)} \cdots a_{k-2,1}^{(2)} \dots a_{k-2,c_{k-2}}^{(2)} a_{k-1,2}^{(1)} \dots a_{k-1,c_{k-1}}^{(1)}$,
where $a_{11}^{(2)} \dots a_{1c_1}^{(2)} \cdots a_{k-2,1}^{(2)} \dots a_{k-2,c_{k-2}}^{(2)} a_{k-1,2}^{(1)} \dots a_{k-1,c_{k-1}}^{(1)}$ is a column word.

Since by the $k = 2$ case in the proof of Lemma 4.3, we have $a_{k1}^{(1)} < a_{k-1,1}^{(1)}$, then by applying *2.* of Lemma 4.2 inductively, we have $a_{k1}^{(k-1)} < a_{k-1,1}^{(k-1)}$.

Combining the above, we have

$$\begin{aligned} & a_{11} \dots a_{1c_1} a_{21} \dots a_{2c_2} \cdots a_{k1} \dots a_{kc_k} \\ \backsim^* & a_{k1}^{(k-1)} a_{11}^{(1)} \dots a_{1c_1}^{(1)} \cdots a_{k-1,1}^{(1)} \dots a_{k-1,c_{k-1}}^{(1)} a_{k2} \dots a_{kc_k} \\ \backsim^* & a_{k1}^{(k-1)} a_{k-1,1}^{(k-1)} a_{11}^{(2)} \dots a_{1c_1}^{(2)} \cdots a_{k-2,1}^{(2)} \dots a_{k-2,c_{k-2}}^{(2)} a_{k-1,2}^{(1)} \cdots a_{k-1,c_{k-1}}^{(1)} a_{k2} \dots a_{kc_k} \end{aligned}$$

with $a_{k1}^{(k-1)} < a_{k-1,1}^{(k-1)}$ and $a_{11}^{(2)} \dots a_{1c_1}^{(2)} \cdots a_{k-2,1}^{(2)} \dots a_{k-2,c_{k-2}}^{(2)} a_{k-1,2}^{(1)} \cdots a_{k-1,c_{k-1}}^{(1)} a_{k2} \dots a_{kc_k}$ is a column word.

Hence by induction, we have

$$
\begin{aligned}
& a_{11}\dots a_{1c_1}a_{21}\dots a_{2c_2}\cdots a_{k1}\dots a_{kc_k} \\
\backsim^* & a^{(k-1)}_{k1}a^{(1)}_{11}\dots a^{(1)}_{1c_1}\cdots a^{(1)}_{k-1,1}\dots a^{(1)}_{k-1,c_{k-1}}a_{k2}\dots a_{kc_k} \\
\backsim^* & a^{(k-1)}_{k1}a^{(k-1)}_{k-1,1}a^{(2)}_{11}\dots a^{(2)}_{1c_1}\cdots a^{(2)}_{k-2,1}\dots a^{(2)}_{k-2,c_{k-2}}a^{(1)}_{k-1,2}\dots a^{(1)}_{k-1,c_{k-1}}a_{k2}\dots a_{kc_k} \\
\backsim^* & a^{(k-1)}_{k1}a^{(k-1)}_{k-1,1}a^{(k-1)}_{k-2,1}a^{(3)}_{11}\dots a^{(3)}_{1c_1}\cdots a^{(3)}_{k-3,1}\dots a^{(3)}_{k-3,c_{k-3}}a^{(2)}_{k-2,1}\dots a^{(2)}_{k-2,c_{k-2}}a^{(1)}_{k-1,2}\cdots \\
& \cdots a^{(1)}_{k-1,c_{k-1}}a_{k2}\dots a_{kc_k} \\
& \qquad\vdots \\
\backsim^* & a^{(k-1)}_{k1}a^{(k-1)}_{k-1,1}\dots a^{(k-1)}_{11}a^{(k-1)}_{12}\dots a^{(k-1)}_{1c_1}a^{(k-2)}_{22}\dots a^{(k-2)}_{2c_2}\cdots \\
& \cdots a^{(1)}_{k-1,2}\dots a^{(1)}_{k-1,c_{k-1}}a_{k2}\dots a_{kc_k}
\end{aligned}
$$

where $a^{(k-1)}_{k1} < a^{(k-1)}_{k-1,1} < \cdots < a^{(k-1)}_{11}$ and $a^{(k-1)}_{11}\dots a^{(k-1)}_{1c_1}a^{(k-2)}_{22}\dots a^{(k-2)}_{2c_2}\cdots a^{(1)}_{k-1,2}\dots a^{(1)}_{k-1,c_{k-1}}a_{k2}\dots a_{kc_k}$ is a column word.

Set $b_{i1} := a^{(k-1)}_{i1}$ and $b_{ij} := a^{(k-i)}_{ij}$ for $1 < j \le c_i$, $1 \le i \le k$ and result follows.

□

**Example 11.** *We use the same word in Example 10 to illustrate Lemma 4.4 by repeated use of Lemma 4.3. We first have*

$w_0 = a^{(0)}_{11}\dots a^{(0)}_{16}|a^{(0)}_{21}\dots a^{(0)}_{25}|a^{(0)}_{31}\dots a^{(0)}_{34}|a^{(0)}_{41}|a^{(0)}_{51} = 886531|97643|9764|5|6$

*and by Example 10, we have*

$a^{(4)}_{51}\|\boxed{a^{(1)}_{11}\dots a^{(1)}_{16}|a^{(1)}_{21}\dots a^{(1)}_{25}|a^{(1)}_{31}\dots a^{(1)}_{34}|a^{(1)}_{41}} = 1\|\boxed{886533|97644|9765|6}$

*and we apply Example 10 again on the latter part (marked by the rectangle), we have*

$$
\begin{cases}
a^{(4)}_{51}|a^{(1)}_{11}\dots a^{(1)}_{16}|a^{(1)}_{21}\dots a^{(1)}_{25}|\enclose{circle}{a^{(2)}_{41}}|a^{(2)}_{31}\dots a^{(2)}_{34} = 1|886533|97644|\enclose{circle}{5}|9766 \\
a^{(4)}_{51}|a^{(1)}_{11}\dots a^{(1)}_{16}|\enclose{circle}{a^{(3)}_{41}}|a^{(2)}_{21}\dots a^{(2)}_{25}|a^{(2)}_{31}\dots a^{(2)}_{34} = 1|886533|\enclose{circle}{4}|97654|9766 \\
a^{(4)}_{51}\enclose{circle}{a^{(4)}_{41}}\|\boxed{a^{(2)}_{11}\dots a^{(2)}_{16}|a^{(2)}_{21}\dots a^{(2)}_{25}|a^{(2)}_{31}\dots a^{(2)}_{34}} = 1\enclose{circle}{3}\|\boxed{886543|97654|9766}
\end{cases}
$$

*and again we apply Example 10 again on the latter part (marked by the rectangle), we have*

$$
\begin{cases}
a^{(4)}_{51}a^{(4)}_{41}|a^{(2)}_{11}\dots a^{(2)}_{16}|a^{(2)}_{21}\dots a^{(2)}_{25}|\enclose{circle}{a^{(2)}_{31}}a^{(2)}_{32}a^{(2)}_{33}a^{(2)}_{34} = 13|886543|97654|\enclose{circle}{9}766 \\
a^{(4)}_{51}a^{(4)}_{41}|a^{(2)}_{11}\dots a^{(2)}_{16}|\enclose{circle}{a^{(3)}_{31}}|a^{(3)}_{21}\dots a^{(3)}_{25}|a^{(2)}_{32}a^{(2)}_{33}a^{(2)}_{34} = 13|886543|\enclose{circle}{7}|99654|766 \\
a^{(4)}_{51}a^{(4)}_{41}\enclose{circle}{a^{(4)}_{31}}\|\boxed{a^{(3)}_{11}\dots a^{(3)}_{16}|a^{(3)}_{21}\dots a^{(3)}_{25}|a^{(2)}_{32}a^{(2)}_{33}a^{(2)}_{34}} = 13\enclose{circle}{6}\|\boxed{887543|99654|766}
\end{cases}
$$

*and by repeating the same process again, we get*

$$
\begin{cases}
a^{(4)}_{51}a^{(4)}_{41}a^{(4)}_{31}|a^{(3)}_{11}\dots a^{(3)}_{16}|\enclose{circle}{a^{(3)}_{21}}a^{(3)}_{22}\dots a^{(3)}_{25}|a^{(2)}_{32}a^{(2)}_{33}a^{(2)}_{34} = 136|887543|\enclose{circle}{9}9654|766 \\
a^{(4)}_{51}a^{(4)}_{41}a^{(4)}_{31}\enclose{circle}{a^{(4)}_{21}}\|\boxed{a^{(4)}_{11}\dots a^{(4)}_{16}|a^{(3)}_{22}\dots a^{(3)}_{25}|a^{(2)}_{32}a^{(2)}_{33}a^{(2)}_{34}} = 136\enclose{circle}{8}\|\boxed{987543|9654|766}
\end{cases}
$$

*and finally we get*

$a^{(4)}_{51}a^{(4)}_{41}a^{(4)}_{31}a^{(4)}_{21}|a^{(4)}_{11}a^{(4)}_{12}\dots a^{(4)}_{16}|a^{(3)}_{22}\dots a^{(3)}_{25}|a^{(2)}_{32}a^{(2)}_{33}a^{(2)}_{34} = 1368|987543|9654|766.$

*So we have:*
$b_{51}b_{41}b_{31}b_{21}|b_{11}b_{12}\dots b_{16}|b_{22}\dots b_{25}|b_{31}\dots b_{34} = 1368|987543|9654|766$*, where*
$b_{11}b_{12}\dots b_{16}|b_{22}\dots b_{25}|b_{32}b_{33}b_{34} = 987543|9654|766$ *is a column word with lengths* $6, 4, 3$.

Using the notation in Lemma 4.4, since

$$w \backsim^* b_{k1}b_{k-1,1}\dots b_{11}b_{12}\dots b_{1c_1}b_{22}\dots b_{2c_2}\cdots b_{k2}\dots b_{kc_k},$$

where $b_{11} > b_{21} > \cdots > b_{k1}$, the SSAF with basement being $\epsilon_n$ representing the word $w$ is the same as that representing $b_{k1}b_{k-1,1}\dots b_{11}b_{12}\dots b_{1c_1}b_{22}\dots b_{2c_2}\cdots b_{k2}\dots b_{kc_k}$. Denote $F(w)$ as the SSAF created.

Since $b_{k1}b_{k-1,1}\dots b_{11}$ is strictly increasing, by Lemma 15 in [HMR13], they create new cells in ascending reading order (i.e. one after another) and hence is exactly the first entire row of $F(w)$ as the entries are fixed when inserting $b_{k1}b_{k-1,1}\dots b_{11}$ into an empty atom, and there are $k$ columns (as $w$ has $k$ subsequences) and so the first row has length $k$ and hence the row reading word has exactly $b_{k1}, b_{k-1,1}, \dots, b_{11}$ as the first subsequence. That means we can apply Lemma 4.4 to find the first subsequence of the row reading word of $F(w)$.

Since $b_{12}\dots b_{1c_1}b_{22}\dots b_{2c_2}\cdots b_{k2}\dots b_{kc_k}$ is a column word, we can apply Lemma 4.4 again and get the second subsequence of the row reading word of $F(w)$, and we can apply Lemma 4.4 repeatedly on the remaining $b_{ij}$'s until we get all the subsequences of the row reading word of $F(w)$. As a result, we can convert $w$ into the row reading word of $F(w)$ by applying Lemma 4.4 repeatedly as described.

We illustrate this by using the example in Example 11:

**Example 12.** *From Example 7 and Example 11, we notice that* $b_{51}b_{41}b_{31}b_{21}b_{11} = 13689$ *which is exactly the first subsequence of the row word of the SSAF representing the column word* $w = 886531|97643|9764|5|6$. *Also one can check that* $87543|9654|766$ *is a column word (by Definition 4.3).*

*We can apply Lemma 4.4 on* $87543|9654|766$ *as we did in Example 11 on* $886531|97643|9764|5|6$ *and get* $589|7643|754|66$ *which gives the second subsequence of the row word*: $589$*, and by repeating the same process, we get* $467|753|64|6$*,* $357|64|6$*,* $46|6$ *and lastly* $6$*.*

*As a result, we get* $467$*,* $357$*,* $46$ *and* $6$ *as the third to the last subsequence of the row word* (*see Example 7*).

## 4.2 Convert a Column Recording Tableau to a Row Recording Tableau

This section gives an interpretation of the twisted Knuth equivalence using recording tableaux. We use the insertion in [Mas08] and the generalized Littlewood-Richardson rule in [HLMvW11]. We also use the notation $(U \leftarrow W)$ for an SSAF $U$ and a biword $W = \begin{pmatrix} x_1 & x_2 & \dots x_n \\ y_1 & y_2 & \dots y_n \end{pmatrix}$ for some positive integer $n$ to denote the pair $(U', L)$ where $U'$ is the SSAF obtained by $(U \leftarrow y_1y_2\dots y_n)$ while $L$ is the recording tableau, i.e. by putting $x_i$ into the cell created when $y_i$ is being inserted. In particular, if $y_1y_2\dots y_n$ is a column (resp. row) word, then we call $L$ as a column (resp. row) recording tableau. By abuse of notation, we sometimes refer $(U \leftarrow W)$ to either $V$ or $L$ (depending on the context).

(Note that if we change the basement $\epsilon_n$ into the large basement in [HLMvW11], a column recording tableau is the same as an LRS defined in Section 4 of [HLMvW11].)

**Lemma 4.5.** *Let $U$ be an SSAF with basement $\epsilon_n$ and shape $\alpha$ for some positive integer $n$ and $l(\alpha) \leq n$. Consider the biword $W = \begin{pmatrix} 2 & 2 & 1 \\ a & b & c \end{pmatrix}$, $a \geq b, c > b$ (i.e. $ab|c$ is a column word.). Let $L = U \leftarrow W$ be the recording tableau. Let $V$ be the SSAF representing the word $abc$ and $a'b'c'$ be the row reading word of $V$ (so $a' < b'$ and $c' \leq b'$ and $abc \backsim^* a'b'c'$). Consider the biword $\widetilde{W} = \begin{pmatrix} 1 & 1 & 2 \\ a' & b' & c' \end{pmatrix}$ and let $\widetilde{L} = U \leftarrow \widetilde{W}$ be the recording tableau. Then $L$ determines $\widetilde{L}$.*

*Proof.* There are two cases to consider: $a \geq c$ or $a < c$.

**Case(I): $\boldsymbol{a \geq c > b}$**

Hence we have $a'b'c' = bac$.

We first consider $L$. Since $a > b$, when we insert $ab$ into $U$, i.e. $\left(U \leftarrow \begin{pmatrix} 2 & 2 \\ a & b \end{pmatrix}\right)$, the cell created by $b$ is strictly above that created by inserting $a$ (by Lemma 15 in [HMR13]):

(i) The cell appears when inserting $b$ is immediately above that of inserting $a$:

$$\begin{array}{c} \boxed{2} \\ \boxed{2} \end{array}$$

(ii) The two cells are the top cells of two distinct columns (by Theorem 16 in [HMR13]) and the cells appear in ascending reading order, (i.e. one after another):

Since $c > b$, by Lemma 15 in [HMR13], the cell created when inserting $c$ after inserting $ab$ into $U$ is after the first cell created by inserting $b$ in reading order.

For (i), we have $L =$

(a) $\begin{array}{ccc} \boxed{2} & & \boxed{1} \\ \boxed{2} & & \end{array}$ or (b) $\begin{array}{ccc} & & \boxed{2} \\ \boxed{1} & & \boxed{2} \end{array}$ or (c) $\begin{array}{ccc} \boxed{2} & & \\ \boxed{2} & & \boxed{1} \end{array}$ or (d) $\begin{array}{ccc} \boxed{2} & & \\ \boxed{2} & & \\ & & \boxed{1} \end{array}$

(a) is impossible as the cell created by inserting $c$, i.e. $\boxed{1}$, is not a removable cell (defined in [HMR13]).

For (b), $\widetilde{L} = \begin{array}{ccc} & & \boxed{2} \\ \boxed{1} & & \boxed{1} \end{array}$ as the cell created by $c'$ must be the top cell of some column and is strictly above that created by $b'$ but there is only one such cell.

(c) and (d) can be considered as the same case by viewing the cell created by inserting $c$ is after that created by inserting $a$ (the second $\boxed{2}$ in reading order). By the same argument as ($b$), we have $\widetilde{L} = \begin{array}{cc} \boxed{2} & \\ \boxed{1} & \\ & \boxed{1} \end{array}$

For (ii), since the reading word of $L$ is contre-lattice by [HLMvW11] which means it should be 221 or 212, the cell created by inserting $c$ is either the middle in reading order or the last one in reading order.

Suppose the cell created by inserting $c$ is the middle cell among the three in reading order, i.e. the reading word of $L$ is 212.

Since $c' = c \leq a = b' > a' = b$, so when performing $(U \leftarrow a'b'c' = bac)$ to get $\widetilde{L}$, the cell created by inserting $a$ is after both of that inserted by $b$ and by $c$, by Lemma 15 of [HMR13], the last cell must be created when $a$ is inserted and hence marked as 1.

If the middle cell (in reading order) is immediately above the last cell, then it cannot be created by $(U \leftarrow b)$ as the last cell is not created if $b$ is inserted before $a$, and so it must be marked as 2, i.e. when $c' = c$ is inserted.

$$\begin{array}{cc} \boxed{\phantom{1}} & \\ & \boxed{2} \\ & \boxed{1} \end{array} \quad \text{or} \quad \begin{array}{cc} & \boxed{\phantom{1}} \\ \boxed{2} & \\ \boxed{1} & \end{array}$$

Suppose the middle cell (in reading order) is not immediately above the last cell, then all three cells are the top cell of three distinct columns. Suppose the middle cell is marked as 1 in $\widetilde{L}$, meaning that when we insert $b$ before inserting $a$, the cell created is after the cell created by $b$ when we first insert $a$ into $U$.

$$\begin{array}{ccc} \boxed{\phantom{1}} & & \\ & & \boxed{1} \\ & \boxed{1} & \end{array} \quad \text{or} \quad \begin{array}{ccc} & \boxed{\phantom{1}} & \\ \boxed{1} & & \\ & & \boxed{1} \end{array} \quad \text{or} \quad \begin{array}{ccc} & & \boxed{\phantom{1}} \\ \boxed{1} & & \\ & \boxed{1} & \end{array} \quad \text{or} \quad \begin{array}{ccc} & \boxed{\phantom{1}} & \\ & & \boxed{1} \\ \boxed{1} & & \end{array}$$

That means the bumping sequence for $(U \leftarrow a)$ involves some cell in the bumping sequence in $(U \leftarrow b)$ and this implies the cell created by $(U \leftarrow a)$ is the same as $(U \leftarrow b)$ (because the bumping sequence has the same ending subsequence starting from the common cell and hence creating the same final cell), and this leads to a contradiction.

Therefore we have the reading order (from first to last in reading order) of $\widetilde{L} = 121$.

Now consider the case when $L$ has the reading word 221.

By the same argument, we know that the last cell in $\widetilde{L}$ is created when $b' = a$ is inserted. If $a' = b$ creates the first cell among the three in reading order as $b$ does in $L$ after $a$ is inserted, then the insertion of $(U \leftarrow a)$ has no common cell in the bumping sequence of $(U \leftarrow b)$. That means when we insert $a$ after inserting $b$, the bumping sequence is the same as inserting $a$ before inserting $b$ (this is true because the bumping sequence of $b$ is a decreasing sequence, so inserting $b$ into $U$ would not affect the first (in reading order) cell in $U$ containing an entry larger than $a$ as $a > b$.). That means $(U \leftarrow ab)$ is the same as $(U \leftarrow ba)$. This leads to a contradiction as we assumed $a$ create different cells in the two cases (the middle and the last cell respectively).

So we know the middle cell is marked as 1 in $\widetilde{L}$. So the reading word of $\widetilde{L}$ is 211.

**Case(II):** $\boldsymbol{c > a \geq b}$

Then $a'b'c' = acb$. Hence the first cells created in $L$ and $\widetilde{L}$ are always the same. Note that $c > a$ meaning the cell created by inserting $a = a'$ is not the last (in reading order) cell, and $a \leq b$ meaning that the cell created by inserting $a$ is not the first (in reading order) cell either. As a result, we know that the middle cell must be created when $a$ is inserted. Since the first (in reading order) is created in $L$ when inserting $b$ while the last (in reading order) cell is created in $\widetilde{L}$ when inserting $b' = c$ is inserted, we know the first cell of $L$ is marked as 2 while the last cell in $\widetilde{L}$ must be marked as 1 and that means the reading word of $L$ is 221 and that of $\widetilde{L}$ is 211.

This shows that if the reading word of $L$ is 221 or when $L$ is of the form as **Case(I)**(i) then the reading word of $\widetilde{L}$ is 211, otherwise the reading word of $\widetilde{L}$ is 121.

Hence $L$ determines $\widetilde{L}$.

□

**Lemma 4.6.** *Let $U$ be an SSAF with basement $\epsilon_n$ and shape $\alpha$ for some positive integer $n$ and $l(\alpha) \leq n$. Consider the biword $W = \begin{pmatrix} 2 & 2 & \cdots & 2 & 1 \\ a_1 & a_2 & \cdots & a_k & b \end{pmatrix}$, $a_1 \geq a_2 \geq \cdots \geq a_k, b > a_k$ (i.e. $a_1a_2\ldots a_k|b$ is a column word ). Let $L = U \leftarrow W$ be the recording tableau. If $b > a_1$, then $L$ has reading word $22\ldots21$, i.e. $b$ creates the last cell in reading order and the insertion of $a_2 \ldots a_k$ is independent of the insertion of $b$.*

*Proof.* Since $a_1a_2\ldots a_kb \backsim^* a_1ba_2\ldots a_k$, we have

$$(U \leftarrow a_1a_2\ldots a_kb) = (U \leftarrow a_1ba_2\ldots a_k).$$

When $k = 2$, the result follows by **Case(II)** in the proof of Lemma 4.5.

For $k > 2$, since $a_1a_2\ldots a_kb \backsim^* a_1ba_2\ldots a_k$, we have $(U \leftarrow a_1a_2\ldots a_kb) = (U \leftarrow a_1ba_2\ldots a_k) = ((U \leftarrow a_1ba_2) \leftarrow a_3\ldots a_k)$. By the case when $k = 2$, we know $b$ creates a cell after both $a_1$ and $a_2$ (and $a_3\ldots a_k$ create cells with decreasing reading order each of which has an order smaller than that of the cell created by $a_2$) and result follows. □

**Example 13.** *Let $W = \begin{pmatrix} 2 & 2 & 2 & 2 & 2 & 2 & 1 \\ 8 & 8 & 6 & 5 & 3 & 1 & 9 \end{pmatrix}$ where $a_1 = 8 < 9 = b$. Let $n = 9$ and $U$ be an SSAF with basement $\epsilon_9 = 123456789$ and shape $\alpha = (1, 1, 0, 0, 3, 0, 2, 0, 6)$ :*

$U =$

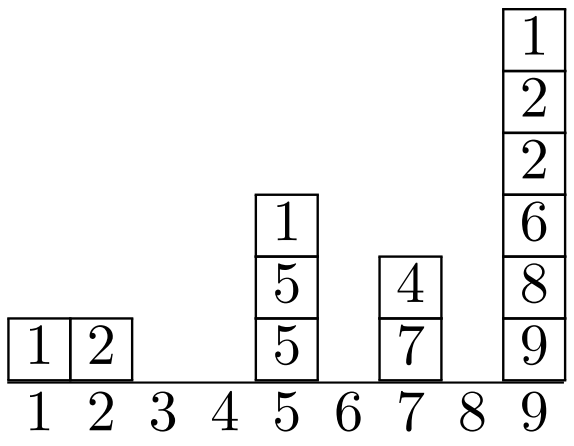

*Therefore* , $(U \leftarrow W)$ *means:*

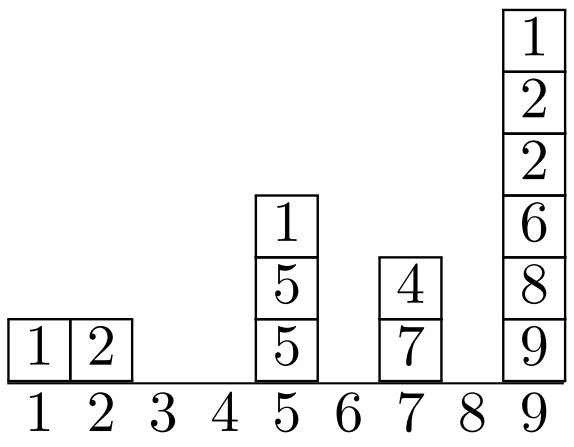

$$\leftarrow \begin{pmatrix} 2 & 2 & 2 & 2 & 2 & 2 & 1 \\ 8 & 8 & 6 & 5 & 3 & 1 & 9 \end{pmatrix}$$

*which results in*

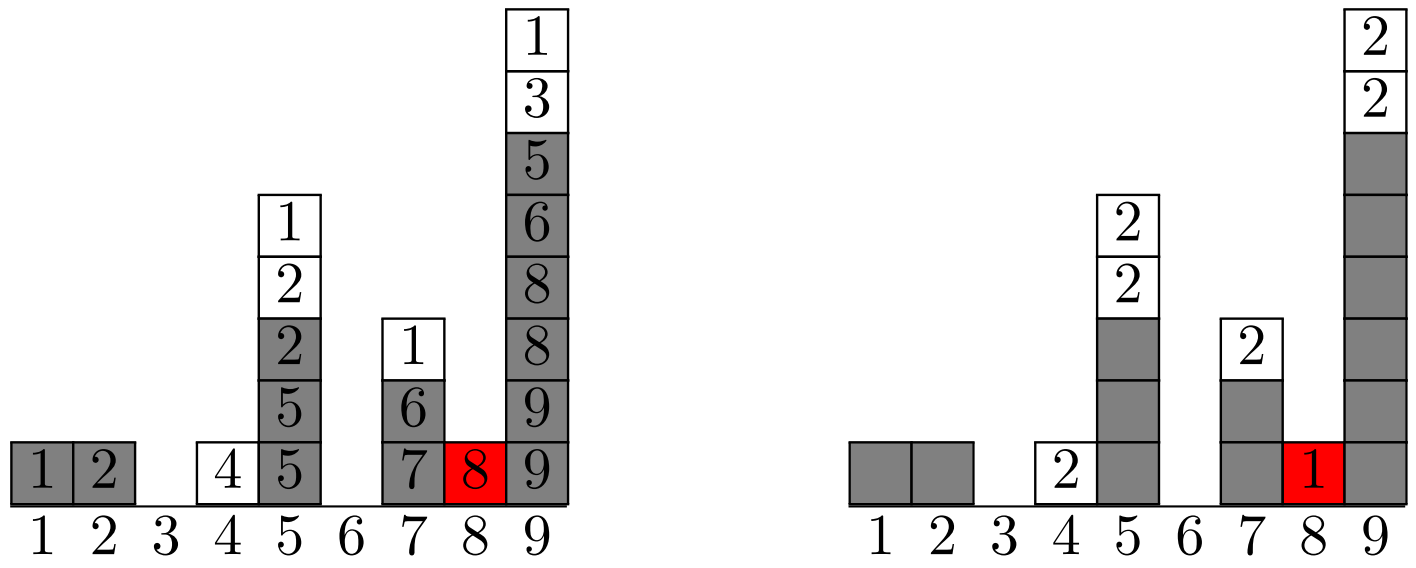

*When* 9 *is inserted, it creates the last cell (we marked the cell red) in reading order among all cells created, and one can see this by reading the recording tableau where the cell is marked as* 1*, which is the last cell in reading order among all cells.*

*Also, note that the bumping route of* 9 *starts from the second row of* $U$ *while those of* $a_2, \dots, a_6 = 8, 6, 5, 3, 1$ *starts from the third row. Hence the insertion of* 9 *does not affect the bumping routes of inserting* $a_2, \dots, a_6$.

**Lemma 4.7.** *Let* $U$ *be an SSAF with basement* $\epsilon_n$ *and shape* $\alpha$ *for some positive integer* $n$ *and* $l(\alpha) \leq n$. *Consider the biword* $W = \begin{pmatrix} 2 & 2 & \cdots & 2 & 1 \\ a_1 & a_2 & \cdots & a_k & b \end{pmatrix}$, $a_1 \geq a_2 \geq \cdots \geq a_k, b > a_k$ *(i.e.* $a_1 a_2 \dots a_k | b$ *is a column word ). Let* $L = (U \leftarrow W)$ *be the recording tableau. If* $i = \min_{1 \leq j \leq k} \{j : b > a_j\}$, *then the cell created by inserting* $a_m$ *for* $m > i$ *is not affected by the insertion of* $b$.

*Proof.* Since $(U \leftarrow a_1 a_2 \dots a_k b) = ((U \leftarrow a_1 a_2 \dots a_{i-1}) \leftarrow a_i \dots a_k b)$. By applying Lemma 4.6 with the $U$ in the lemma being $(U \leftarrow a_1 \dots a_{i-1})$ and result follows.

□

**Lemma 4.8.** *Let* $U$ *be an SSAF with basement* $\epsilon_n$ *and shape* $\alpha$ *for some positive integer* $n$ *and* $l(\alpha) \leq n$. *Consider the biword* $W = \begin{pmatrix} 2 & 2 & \cdots & 2 & 1 \\ a_1 & a_2 & \cdots & a_k & b \end{pmatrix}$, $a_1 \geq a_2 \geq \cdots \geq a_k, b > a_k$ (*i.e.* $a_1a_2\dots a_k|b$ *is a column word* ). *Let* $L = (U \leftarrow W)$ *be the recording tableau. Suppose* $b \le a_{k-1}$. *Let* $\widetilde{L} = (U \leftarrow \widetilde{W})$ *where* $\widetilde{W} = \begin{pmatrix} 1 & 1 & 2 & \cdots & 2 & 2 \\ a_k & a_1 & a_2 & \cdots & a_{k-1} & b \end{pmatrix}$, *then* $L$ *determines* $\widetilde{L}$.

*Proof.* When $k = 2$, this is proved by **Case(I)** in Lemma 4.5.

For $k > 2$, since $a_1 \dots a_k b \backsim^* a_k a_1 \dots a_{k-1} b$, $L$ and $\widetilde{L}$ has the same shape.

Also, $a_1 \dots a_k b \backsim^* a_1 \dots a_{k-1} a_k b$, we have

$$(U \leftarrow a_1 \dots a_{k-1} a_k b) = ((U \leftarrow a_1 \dots a_{k-2}) \leftarrow a_{k-1} a_k b).$$

By applying the case when $k = 2$ on $(U \leftarrow a_1 \dots a_{k-2})$ and $a_{k-1}a_kb$ being the length-3-word inserted, we know the order of the cell being created which represents $((U \leftarrow a_1 \dots a_{k-2}) \leftarrow a_{k-1}a_kb) = ((U \leftarrow a_1 \dots a_{k-2}) \leftarrow a_k a_{k-1} b)$ and so we know which one among the three is the last cell being created (by Lemma 4.5, the last cell being created may either be the first or second cell in reading order among the three cells created). By marking that cell as 2 in $\widetilde{L}$ and then consider the first two cells being created we know how the cells are being created by the insertion on $((U \leftarrow a_1 \dots a_{k-2}) \leftarrow a_k a_{k-1}) = (U \leftarrow a_1 \dots a_{k-2} a_k a_{k-1}) = (U \leftarrow a_k a_1 \dots a_{k-2} a_{k-1})$ by induction (as $a_{k-2} \le a_{k-1} < a_k$ for $a_k < b \le a_{k-1}$) and together with the last cell marked by 2 as mentioned, we know how to label the entries of the recording tableau $\widetilde{L}$ for $(U \leftarrow a_k a_1 \dots a_{k-2} a_{k-1} b)$.

□

**Example 14.** *Let* $U$ *as in Example 13 and let* $W = \begin{pmatrix} 2 & 2 & 2 & 2 & 1 \\ 9 & 7 & 6 & 4 & 5 \end{pmatrix}$.

*Then* $\widetilde{W} = \begin{pmatrix} 2 & 2 & 2 & 2 & 1 \\ 4 & 9 & 7 & 6 & 5 \end{pmatrix}$.

$(U \leftarrow W)$ *means:*

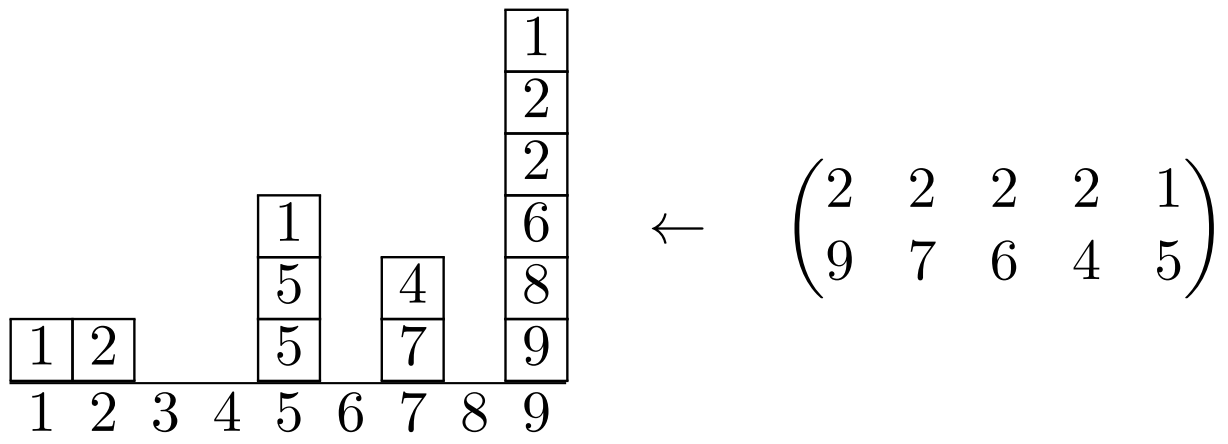

*which results in*

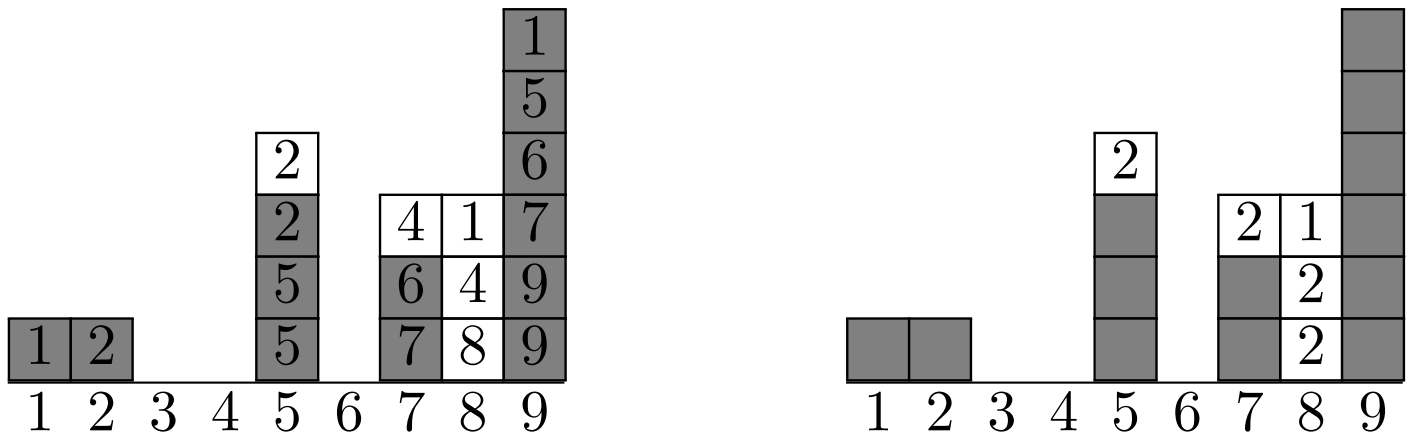

$=: L$

*We illustrate the proof of Lemma 4.8 to get* $\widetilde{L}$ *from* $L$.

*Consider the last three cells created which are marked green in* $L$:

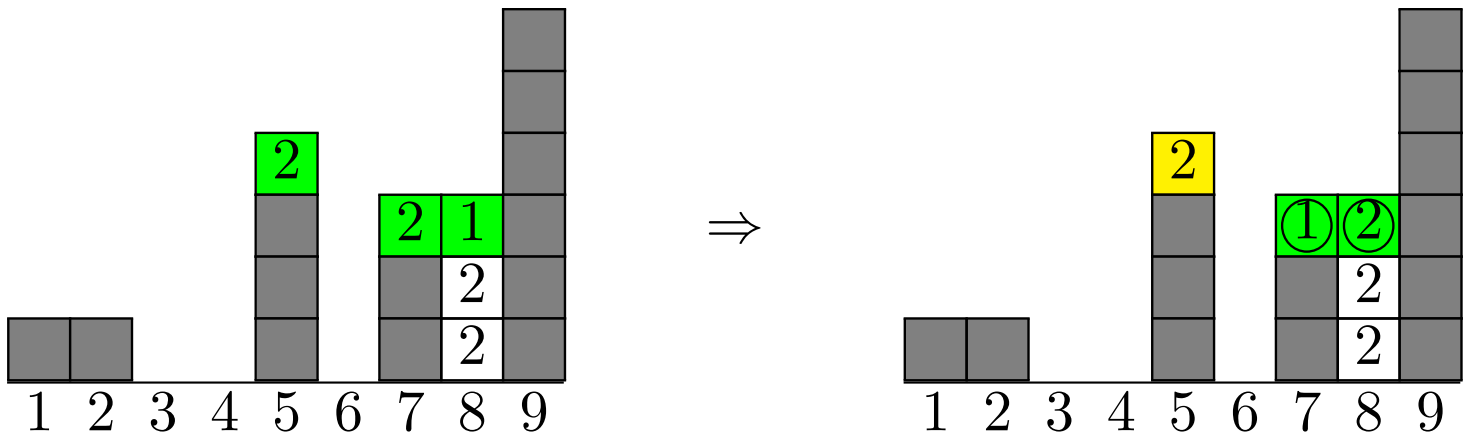

*We get the first cell (in reading order) in* $\widetilde{L}$ *(the yellow cell) marked as* 2 *by applying Case* (ii) *in the proof of Lemma 4.5 as we know the green cells would have reading* 211 *in* $\widetilde{L}$ *which implies the first cell in reading order among the three green cells in* $L$ *must be marked with* 2 *in* $\widetilde{L}$ *while the two cells appear in the order as circled (i.e. the one marked as* ① *appear before the one marked with* ② *when we insert* $a_4a_3b = 465$ *to get* $\widetilde{L}$*.)*

*Hence we can treat the cell marked with* ② *as the last cell created among the remaining 4 cells in* $\widetilde{L}$ *to be filled, and so we mark that cell as* 1 *( treating it as the new* $b$ *being inserted as* $4 < 6 < 7$ *). Now we have the three new lastly-created cells marked green:*

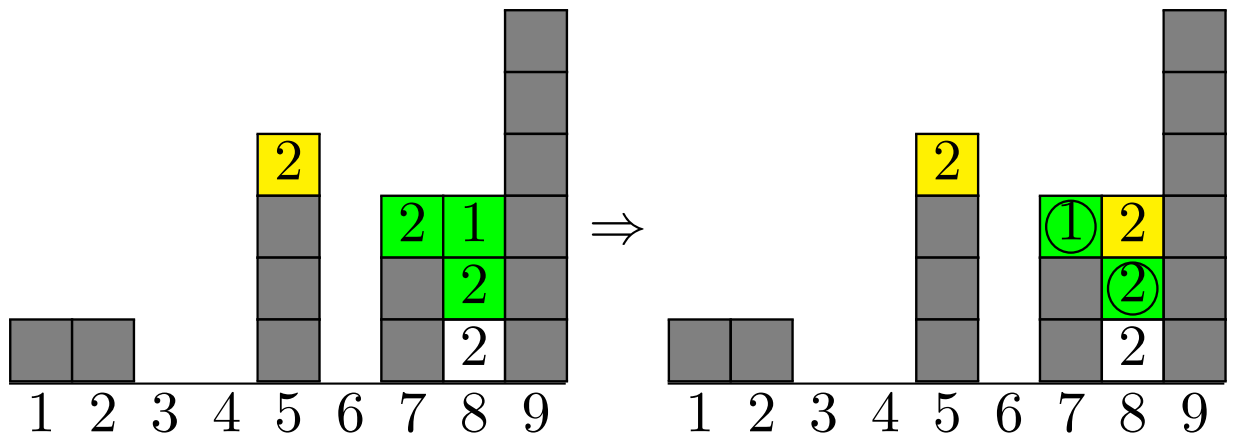

*Again, we get the second (in reading order) cell among the three green cells (the new yellow cell ) marked as* 2 *by applying Case* (ii) *in the proof of Lemma 4.5 as we know the green cells would have reading* 121 *in* $\widetilde{L}$ *which implies the second cell in reading order among the three green cells in* $L$ *must be marked* 2 *in* $\widetilde{L}$ *while the two cells appear in the order as circled (i.e. the one marked with* ① *appears before the one marked with* ② *when we insert* $a_4a_2a_3 = 476$ *to get* $\widetilde{L}$*.)*

*Hence we can, again, treat the cell marked with* ② *as the last cell created among the remaining three cells in* $\widetilde{L}$ *to be filled, and so mark it as* 1 *and the other two cells as* 2 *by a similar argument as the previous step. Now we have the three new lastly-created cells marked green:*

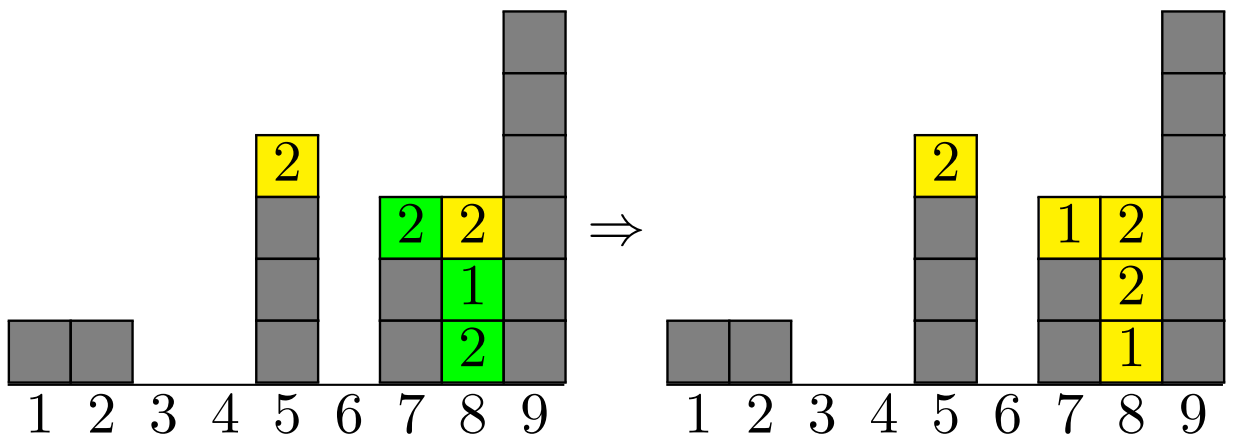

$= \widetilde{L}$

*We get the remaining entries of* $\widetilde{L}$ *by applying Lemma 4.5 again as we know the reading word among these three cells would change from* 212 *to* 121*.*

*One can verify that* $\widetilde{L}$ *is indeed the recording tableau of* $(U \leftarrow \widetilde{W})$*.*

**Lemma 4.9.** *Let $U$ be an SSAF with basement $\epsilon_n$ and shape $\alpha$ for some positive integer $n$ and $l(\alpha) \leq n$. Consider the biword $W = \begin{pmatrix} 2 & 2 & \cdots & 2 & 1 \\ a_1 & a_2 & \cdots & a_k & b \end{pmatrix}$, $a_1 \geq a_2 \geq \cdots \geq a_k, b > a_k$ (i.e. $a_1a_2\ldots a_k|b$ is a column word ). Let $L = (U \leftarrow W)$ be the recording tableau. Let $\widetilde{L} = (U \leftarrow \widetilde{W})$ where $\widetilde{W} = \begin{pmatrix} 1 & 1 & 2 & \cdots & 2 \\ b' & a'_1 & a'_2 & \cdots & a'_k \end{pmatrix}$, such that $b' = a_i$, where $i = \min_{1\leq j\leq k}\{j : b > a_j\}$ and $b' < a'_1$, $a'_1 \geq a'_2 \geq \cdots \geq a'_k$ (one can verify that $b'a'_1\ldots a'_k$ is indeed the row reading word of the SSAF representing the word $a_1a_2...a_kb$ by 1. in Lemma 4.2), then $L$ determines $\widetilde{L}$.*

*Proof.* If $i = k$, then we are done by Lemma 4.8.

If $i = 1$, then by *1.* in Lemma 4.2, we have $b'a'_1\ldots a'_k = a_1ba_2\ldots a_k$ and hence by Lemma 4.6, $b$ would create the last cell in reading order in $\widetilde{L}$ when being inserted and hence we know that which two cells are created by $a_1$ and $b$ in $\widetilde{L}$ and hence we know how to label the entries by marking those two cells as 1 and the rest as 2.

Suppose $i < k$ then by *1.* in Lemma 4.2, we have

$$b'a'_1\ldots a'_k = a_ia_1\ldots a_{i-1}ba_{i+1}\ldots a_k.$$

Now by Lemma 4.7, we know all $a_j, j > i$ creates create the same cells in $L$ and $\widetilde{L}$ and since $a_1\ldots a_{i-1}a_ib \backsim^* a_ia_1\ldots a_{i-1}b$, we just need to first apply Lemma 4.8 on $(U \leftarrow a_1\ldots a_ib)$ as the insertion recorded by $L$ in the lemma in order to find the first $i+1$ cells created in $\widetilde{L}$ by the insertion $(U \leftarrow a_ia_1\ldots a_{i-1}b)$ and then label the rest of cells created by $((U \leftarrow a_ia_1\ldots a_{i-1}b) \leftarrow a_{i+1}\ldots a_k)$ as 2 and get the entries of $\widetilde{L}$ and result follows. □

Lemma 4.9 gives a recording tableau interpretation of *1.* in Lemma 4.2 using $L$ and $\widetilde{L}$.

**Lemma 4.10.** *Let $U$ be an SSAF with basement $\epsilon_n$ and shape $\alpha$ for some positive integer $n$ and $l(\alpha) \leq n$. Let $w_0 = a_{11}\ldots a_{1c_1}|a_{21}\ldots a_{2c_2}|\cdots|a_{k1}\ldots a_{kc_k}$ with $k \leq n$ be a column word and using the notation in Lemma 4.3 while we assume $(i_j) = 0$ for $1 \leq j \leq c_r$, $1 \leq r \leq k$, the recording tableau $L$ of $(U \leftarrow W_0)$ determines the recording tableau $\widetilde{L}_m$ to $(U \leftarrow W_m)$ for $0 \leq m \leq k-1$, where $L := \widetilde{L}_0$ and $W_m$ is a biword with the lower word being $w_m$ and the upper word has entry $k+1-t$ if the lower word entry just below it is $a^{(s)}_{tj}$ for $s = i_t$, $i_t+1$, and $1_{(1)}$ if $m > 0$ and the entry is $a^{i_k+m}_{k1}$.*

*Proof.* Suppose $k = 2$ then

$$(U \leftarrow w_0) = ((U \leftarrow a_{11}\ldots a_{1c_1}a_{21}) \leftarrow a_{22}\ldots a_{2c_2}) = \left(((U \leftarrow a^{(1)}_{21}) \leftarrow a^{(1)}_{11}\ldots a^{(1)}_{1c_1}) \leftarrow a_{22}\ldots a_{2c_2}\right)$$

and by Lemma 4.9 we know how the cells are created in $((U \leftarrow a^{(1)}_{21}) \leftarrow a^{(1)}_{11}\ldots a^{(1)}_{1c_1})$. So giving the recording tableau $L$ of $(U \leftarrow a_{11}\ldots a_{1c_1}a_{21})$, we know how the cells are created by $(U \leftarrow a^{(1)}_{21}a^{(1)}_{11}\ldots a^{(1)}_{1c_1})$. Now inserting the rest of the sequence $a_{22}\ldots a_{2c_2}$ and which creates in $\widetilde{L}_1$ the same last $(c_2-1)$ cells as in $L$, we know how to label the recording tableau $\widetilde{L}_1$.

Suppose $k > 2$, then

$(U \leftarrow w_0) = (U \leftarrow a_{11}\dots a_{1c_1}a_{21}\dots a_{2c_2}\cdots a_{k1}\dots a_{kc_k})$
$= ((U \leftarrow a_{11}\dots a_{1c_1}\cdots a_{k-2,1}\dots a_{k-2,c_{k-2}}) \leftarrow a_{k-1,1}\dots a_{k-1,c_{k-1}}a_{k1}\dots a_{kc_k})$ and by the argument with

$$(U \leftarrow a_{11}\dots a_{1c_1}\cdots a_{k-2,1}\dots a_{k-2,c_{k-2}})$$

being the $U$ to be considered for $k=2$, we know how the cells are created in $\widetilde{L}_1$ for

$$((U \leftarrow a_{11}\dots a_{1c_1}\cdots a_{k-2,1}\dots a_{k-2,c_{k-2}}) \leftarrow a^{(1)}_{k1}a^{(1)}_{k-1,1}\dots a^{(1)}_{k-1,c_{k-1}}a_{k2}\dots a_{kc_k})$$

and by the same argument, we know how the cells are created in $\widetilde{L}_2$ for
$((U \leftarrow a_{11}\dots a_{1c_1}\cdots a_{k-3,1}\dots a_{k-3,c_{k-3}}) \leftarrow a^{(2)}_{k1}a^{(1)}_{k-2,1}\dots a^{(1)}_{k-2,c_{k-2}}a^{(1)}_{k-1,1}\cdots a^{(1)}_{k-1,c_{k-1}}a_{k2}\dots a_{kc_k})$
and result follows by induction.

□

**Example 15.** *We use the words in Example 10 and the same $U$ as in Example 13 to illustrate Lemma 4.10.*

*Recall that in Example 10, we have*

$$\begin{cases} w_0 = a^{(0)}_{11}\dots a^{(0)}_{16}|a^{(0)}_{21}\dots a^{(0)}_{25}|a^{(0)}_{31}\dots a^{(0)}_{34}|a^{(0)}_{41}|\textcircled{$a^{(0)}_{51}$} = 886531|97643|9764|5|\textcircled{6} \\ w_1 = a^{(0)}_{11}\dots a^{(0)}_{16}|a^{(0)}_{21}\dots a^{(0)}_{25}|a^{(0)}_{31}\dots a^{(0)}_{34}|\textcircled{$a^{(1)}_{51}$}|a^{(1)}_{41} = 886531|97643|9764|\textcircled{5}|6 \\ w_2 = a^{(0)}_{11}\dots a^{(0)}_{16}|a^{(0)}_{21}\dots a^{(0)}_{25}|\textcircled{$a^{(2)}_{51}$}|a^{(1)}_{31}\dots a^{(1)}_{34}|a^{(1)}_{41} = 886531|97643|\textcircled{4}|9765|6 \\ w_3 = a^{(0)}_{11}\dots a^{(0)}_{16}|\textcircled{$a^{(3)}_{51}$}|a^{(1)}_{21}\dots a^{(1)}_{25}|a^{(1)}_{31}\dots a^{(1)}_{34}|a^{(1)}_{41} = 886531|\textcircled{3}|97644|9765|6 \\ w_4 = \textcircled{$a^{(4)}_{51}$}|a^{(1)}_{11}\dots a^{(1)}_{16}|a^{(1)}_{21}\dots a^{(1)}_{25}|a^{(1)}_{31}\dots a^{(1)}_{34}|a^{(1)}_{41} = \textcircled{1}|886533|97644|9765|6 \end{cases}.$$

*Therefore, we now have*

$$\begin{cases} W_0 = \left(\begin{array}{cccccc|ccccc|cccc|c|c} 5 & 5 & 5 & 5 & 5 & 5 & 4 & 4 & 4 & 4 & 4 & 3 & 3 & 3 & 3 & 2 & 1 \\ 8 & 8 & 6 & 5 & 3 & 1 & 9 & 7 & 6 & 4 & 3 & 9 & 7 & 6 & 4 & 5 & \textcircled{6} \end{array}\right) \\ W_1 = \left(\begin{array}{cccccc|ccccc|cccc|c|c} 5 & 5 & 5 & 5 & 5 & 5 & 4 & 4 & 4 & 4 & 4 & 3 & 3 & 3 & 3 & 1_{(1)} & 2 \\ 8 & 8 & 6 & 5 & 3 & 1 & 9 & 7 & 6 & 4 & 3 & 9 & 7 & 6 & 4 & \textcircled{5} & 6 \end{array}\right) \\ W_2 = \left(\begin{array}{cccccc|ccccc|c|cccc|c} 5 & 5 & 5 & 5 & 5 & 5 & 4 & 4 & 4 & 4 & 4 & 1_{(1)} & 3 & 3 & 3 & 3 & 2 \\ 8 & 8 & 6 & 5 & 3 & 1 & 9 & 7 & 6 & 4 & 3 & \textcircled{4} & 9 & 7 & 6 & 5 & 6 \end{array}\right) \\ W_3 = \left(\begin{array}{cccccc|c|ccccc|cccc|c} 5 & 5 & 5 & 5 & 5 & 5 & 1_{(1)} & 4 & 4 & 4 & 4 & 4 & 3 & 3 & 3 & 3 & 2 \\ 8 & 8 & 6 & 5 & 3 & 1 & \textcircled{3} & 9 & 7 & 6 & 4 & 4 & 9 & 7 & 6 & 5 & 6 \end{array}\right) \\ W_4 = \left(\begin{array}{c|cccccc|ccccc|cccc|c} 1_{(1)} & 5 & 5 & 5 & 5 & 5 & 5 & 4 & 4 & 4 & 4 & 4 & 3 & 3 & 3 & 3 & 2 \\ \textcircled{1} & 8 & 8 & 6 & 5 & 3 & 3 & 9 & 7 & 6 & 4 & 4 & 9 & 7 & 6 & 5 & 6 \end{array}\right) \end{cases}.$$

*We now show how to get $\widetilde{L}_1$ to $\widetilde{L}_4$ from $L = \widetilde{L}_0$.*

*For $U \leftarrow W_0$:*

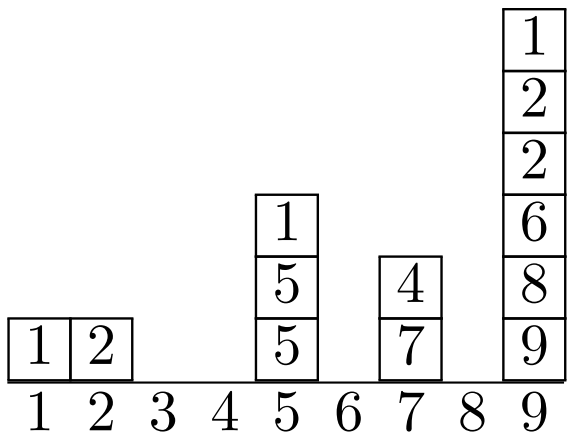

$$\leftarrow \left(\begin{array}{cccccc|ccccc|cccc|c|c} 5 & 5 & 5 & 5 & 5 & 5 & 4 & 4 & 4 & 4 & 4 & 3 & 3 & 3 & 3 & 2 & 1 \\ 8 & 8 & 6 & 5 & 3 & 1 & 9 & 7 & 6 & 4 & 3 & 9 & 7 & 6 & 4 & 5 & 6 \end{array}\right)$$

*which results in*

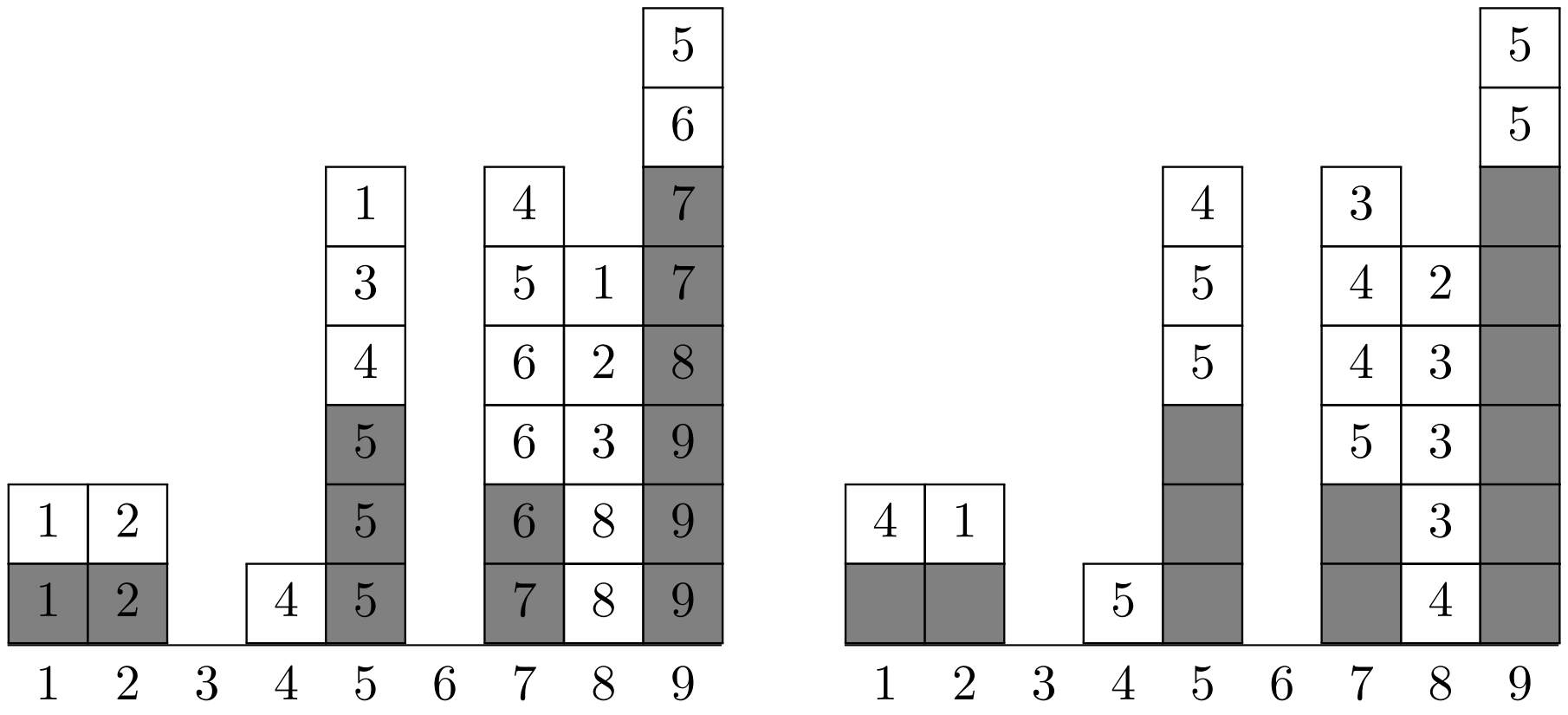

$= L = \widetilde{L}_0.$

*To get $\widetilde{L}_0$, we just consider the last (in reading order) cell containing* 1 *and all cells containing* 2 *and* 1 *and apply Lemma 4.10. Note that in this case since there are only two cells involved, we can simply move* 1 *to the cell with a smaller reading order cell and mark as* $1_{(1)}$ :

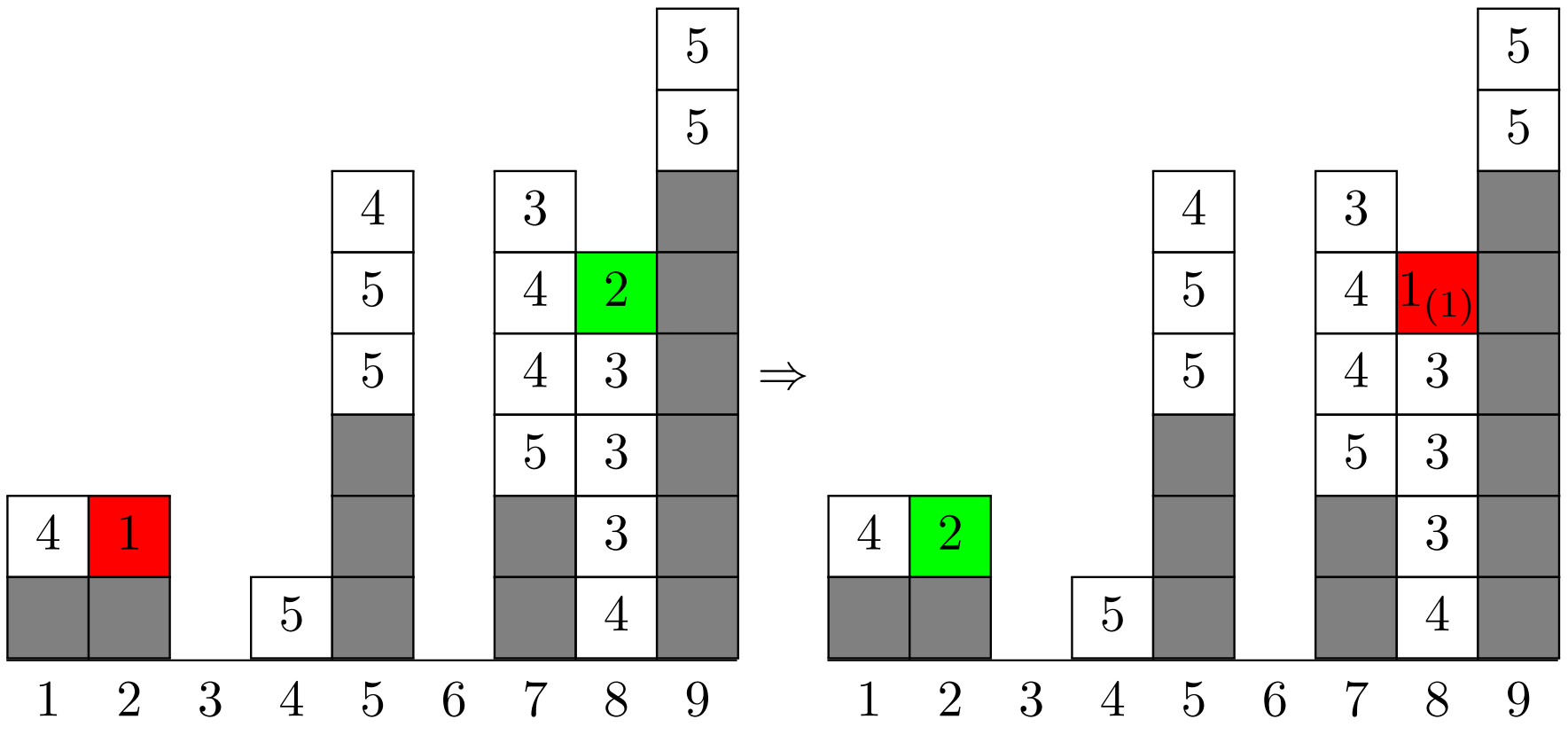

$= \widetilde{L}_1.$

*From $\widetilde{L}_1$ we get $\widetilde{L}_2$ by considering the cell containing $1_{(1)}$ and all the cells containing* 3 *and apply Lemma 4.10 to find the new position of $1_{(1)}$:*

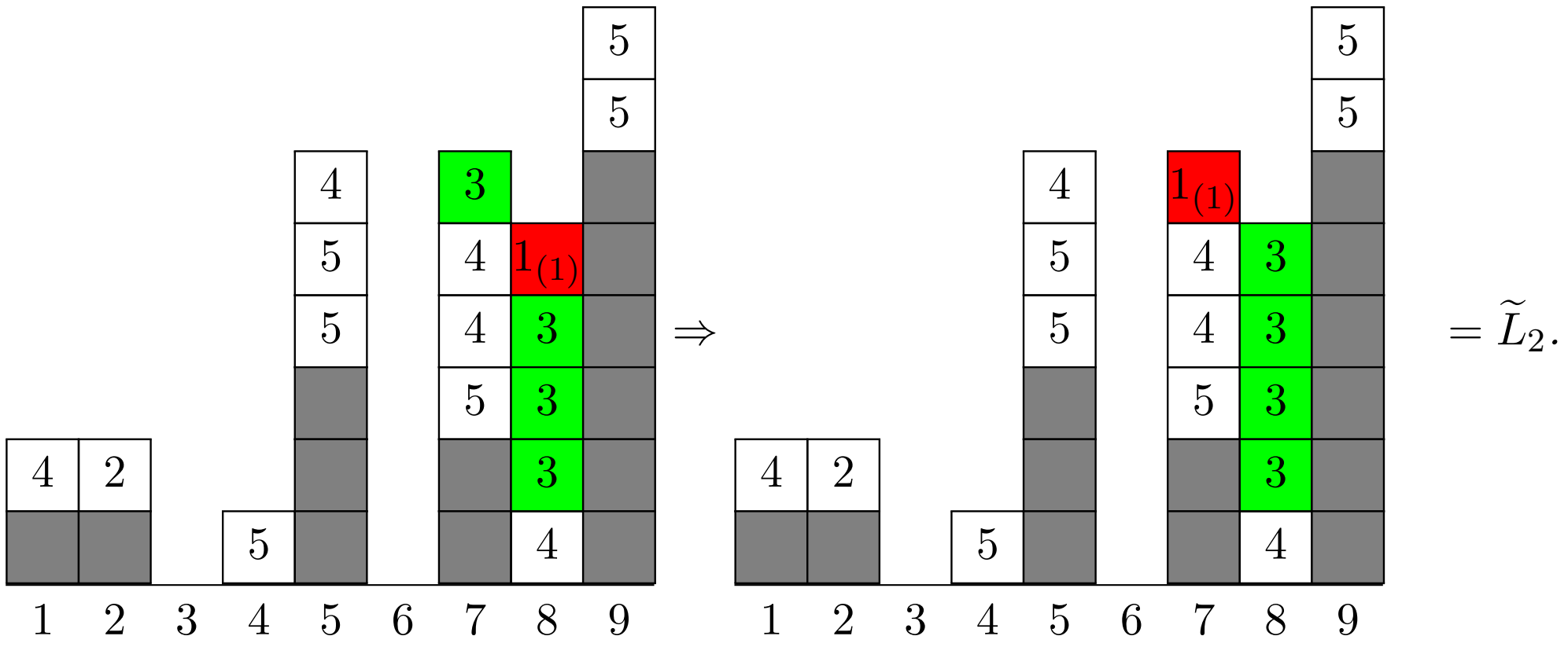

*Then consider the cells with* 4 *and also the cell with* $1_{(1)}$ *in* $\widetilde{L}_2$*, we get* $\widetilde{L}_3$ *by applying Lemma 4.10 to find the new position of* $1_{(1)}$*:*

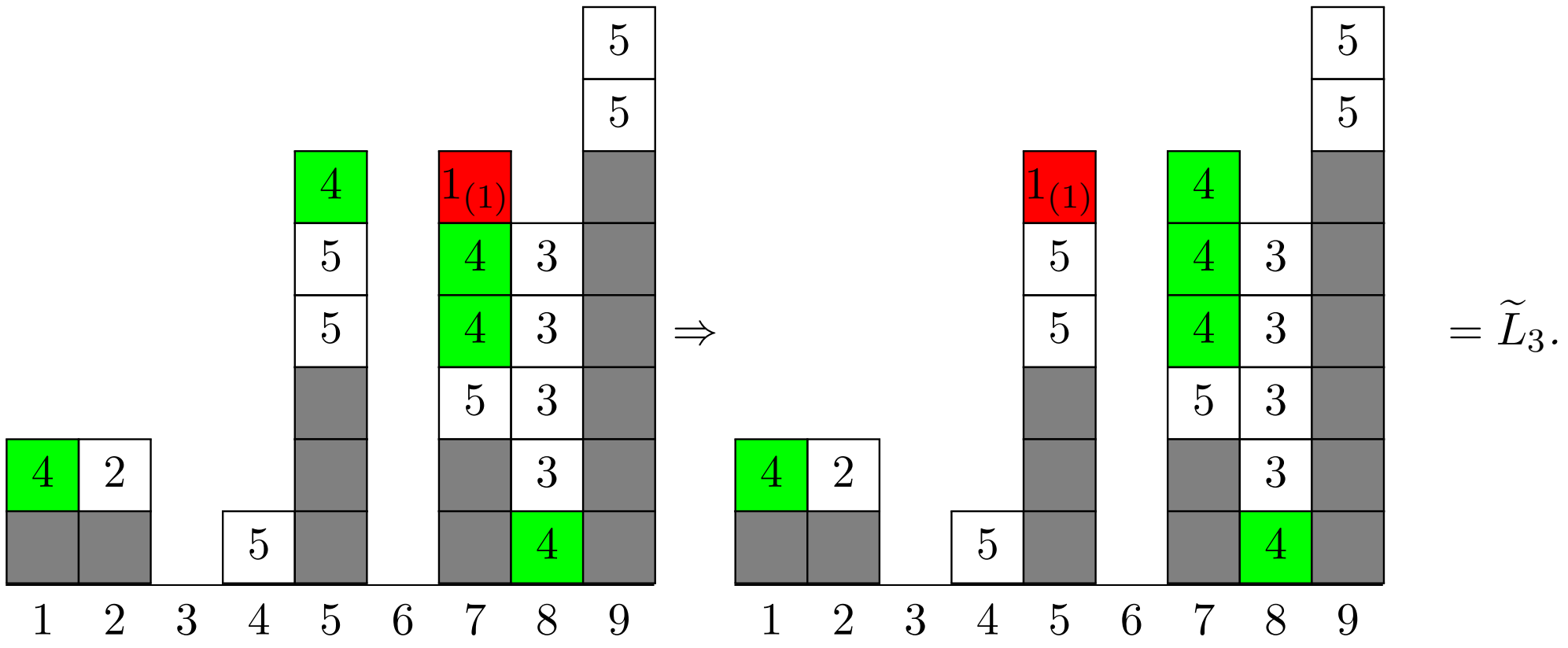

*By the same argument and consider the cells with* 5 *and also the cell with* $1_{(1)}$*, we get* $\widetilde{L}_4$*:*

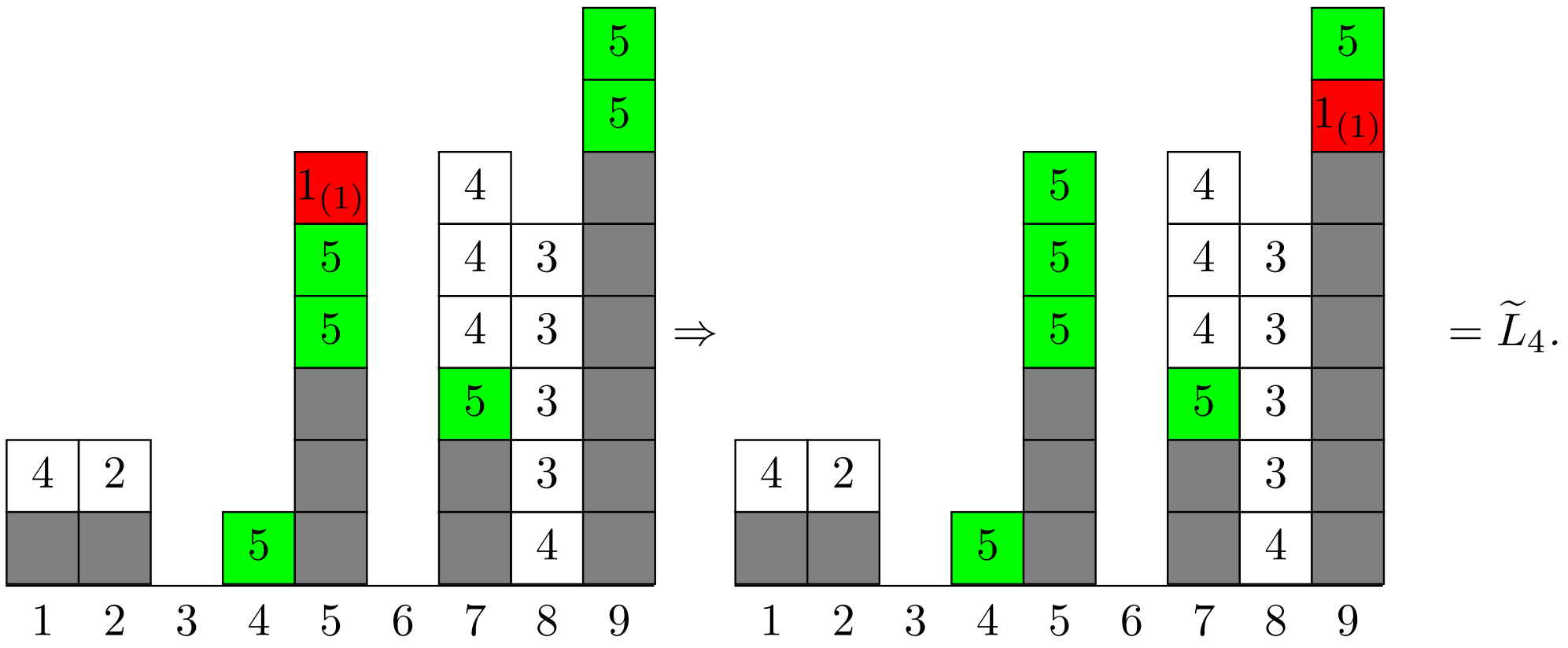

*We describe a more direct way to get* $L = \widetilde{L}_0$ *from* $\widetilde{L}_4$*. Note that the reading word of* $L$ *is a contre-lattice word, we can always find a* 2 *before (in reading order)* 1*. By the last* *part in the proof of Lemma 4.5, we know that, except for the case when $L$ is of the form of* **Case(I)**(i) *, the first cell created in $\widetilde{L}$ (which is the first cell in reading order containing a $1$) is the cell with largest reading order containing $2$ in $L$ before the cell containing $1$ in $L$, and so we put $1_{(1)}$ into that cell in $\widetilde{L}$.*

*In short, the red cell with $1_{(1)}$ is always interchanged with the last green cell before it in reading order, except for the case when that green cell is immediately above the next green cell (which must be after the red cell if it exists) in which the red cell stays the same and compare to the next set of green cells (if available).*

**Lemma 4.11.** *Let $U$ be an SSAF with basement $\epsilon_n$ and shape $\alpha$ for some positive integer $n$ and $l(\alpha) \leq n$. Let $w = a_{11} \dots a_{1c_1} | a_{21} \dots a_{2c_2} | \cdots | a_{k1} \dots a_{kc_k}$ with $k \leq n$ be a column word and using the notation in Lemma 4.4, the recording tableau $L$ of $(U \leftarrow W)$ determines the recording tableau $\widetilde{L}$ of $(U \leftarrow \widetilde{W})$, where*

$$W = \begin{pmatrix} k & k & \dots & k & \cdots & 1 & 1 & \dots & 1 \\ a_{11} & a_{12} & \dots & a_{1c_1} & \cdots & a_{k1} & a_{k2} & \dots & a_{kc_k} \end{pmatrix}$$

*and*

$$\widetilde{W} = \begin{pmatrix} 1_{(1)} & 2_{(1)} & \dots & k_{(1)} & k & k & \dots & k & \cdots & 1 & 1 & \dots & 1 \\ b_{k1} & b_{k-1,1} & \dots & b_{11} & b_{12} & b_{13} & \dots & b_{1c_1} & \cdots & b_{k2} & b_{k3} & \dots & b_{kc_k} \end{pmatrix}.$$

*Proof.* As in the proof of Lemma 4.4 which depends mostly on Lemma 4.3, we apply Lemma 4.10 (which is like the tableaux version of Lemma 4.3) to get the result.

Given $L$, by Lemma 4.10, we can get $\widetilde{L}_{k-1}$ and if we ignore all the 1 entries (including $1_{(1)}$) in $\widetilde{L}_{k-1}$, we get a new column recording tableau $L'$ with entries $2, \dots k$, and we apply Lemma 4.10 again (treating $r$ in $L'$ as $r-1$ in $L$ when we apply Lemma 4.10) and get a $\widetilde{L'}_{k-2}$ with an entry $2_{(1)}$. By changing the corresponding entries of $\widetilde{L}$ with those in $\widetilde{L'}_{k-2}$, we have $2_{(2)}$ in a cell after the cell containing $1_{(1)}$ in reading order (because $a_{k1}^{(k-1)} < a_{k-1,1}^{(k-1)}$ and apply Lemma 15 of [HMR13] ). Result follows by repeating this process until $(k-1)_{(1)}$ is formed and then convert the first $k$ as $k_{(1)}$. □

**Example 16.** *We use the same word in Example 7 and the same $U$ in Example 13, which is also the $U$ and the word we used as $w_0$ in Example 15, to illustrate Lemma 4.11.*

*By Example 11,*

$$\widetilde{W} = \begin{pmatrix} 1_{(1)} & 2_{(1)} & 3_{(1)} & 4_{(1)} & 5_{(1)} & 5 & 5 & 5 & 5 & 5 & 4 & 4 & 4 & 4 & 3 & 3 & 3 \\ 1 & 3 & 6 & 8 & 9 & 8 & 7 & 5 & 4 & 3 & 9 & 6 & 5 & 4 & 7 & 6 & 6 \end{pmatrix}.$$

*We continue Example 15 to illustrate Lemma 4.11 to get $\widetilde{L}$.*

*We mark the cell which is being moved and added subscript (1) as* ***green*** *and the the cells under consideration to get the new position of the green cell as* ***yellow*** *and mark the final cell of the subscripted green cell as* ***red****.*

*By Example 15, we get the position $1_{(1)}$ and we continue using the same process:*

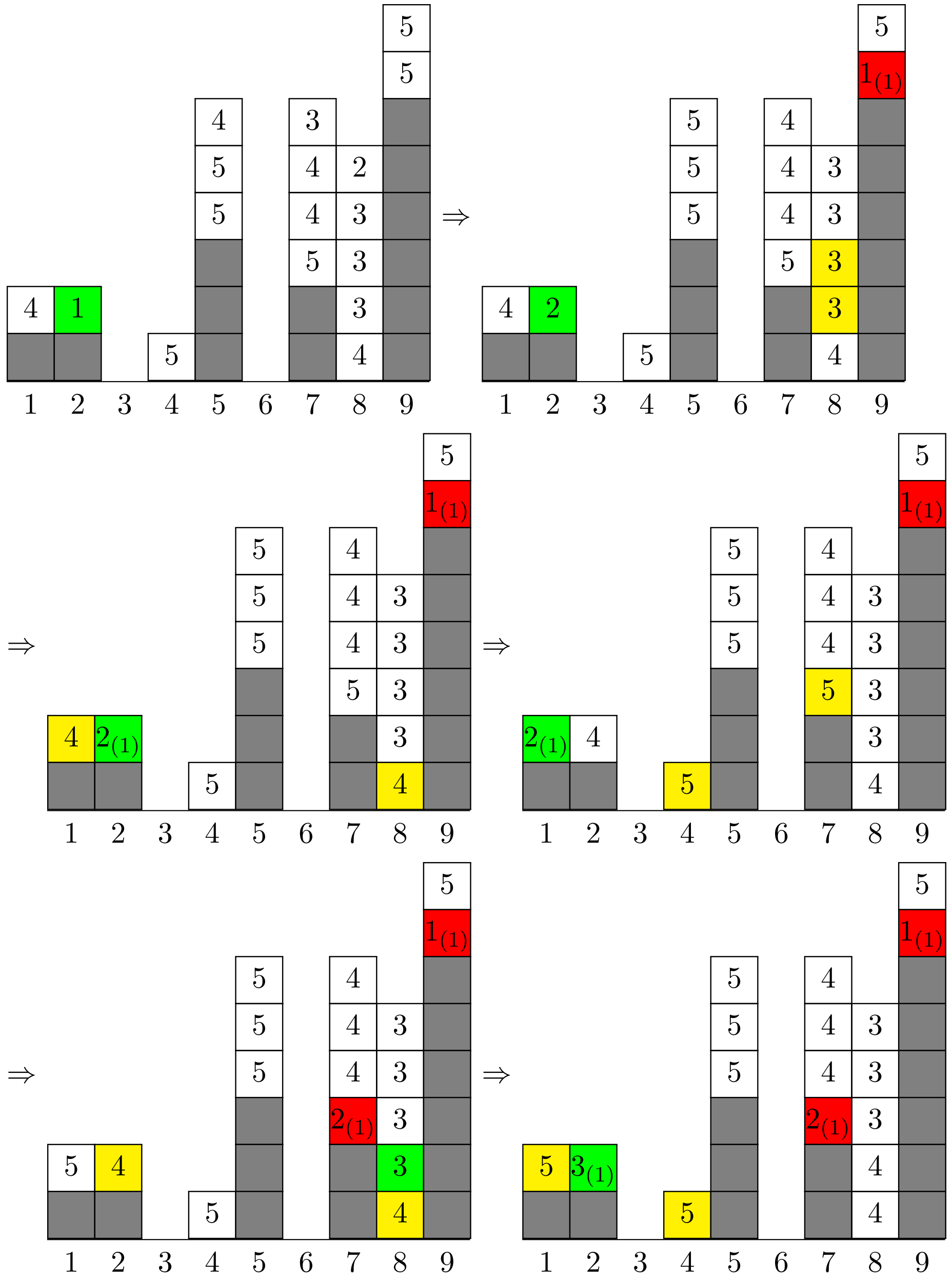

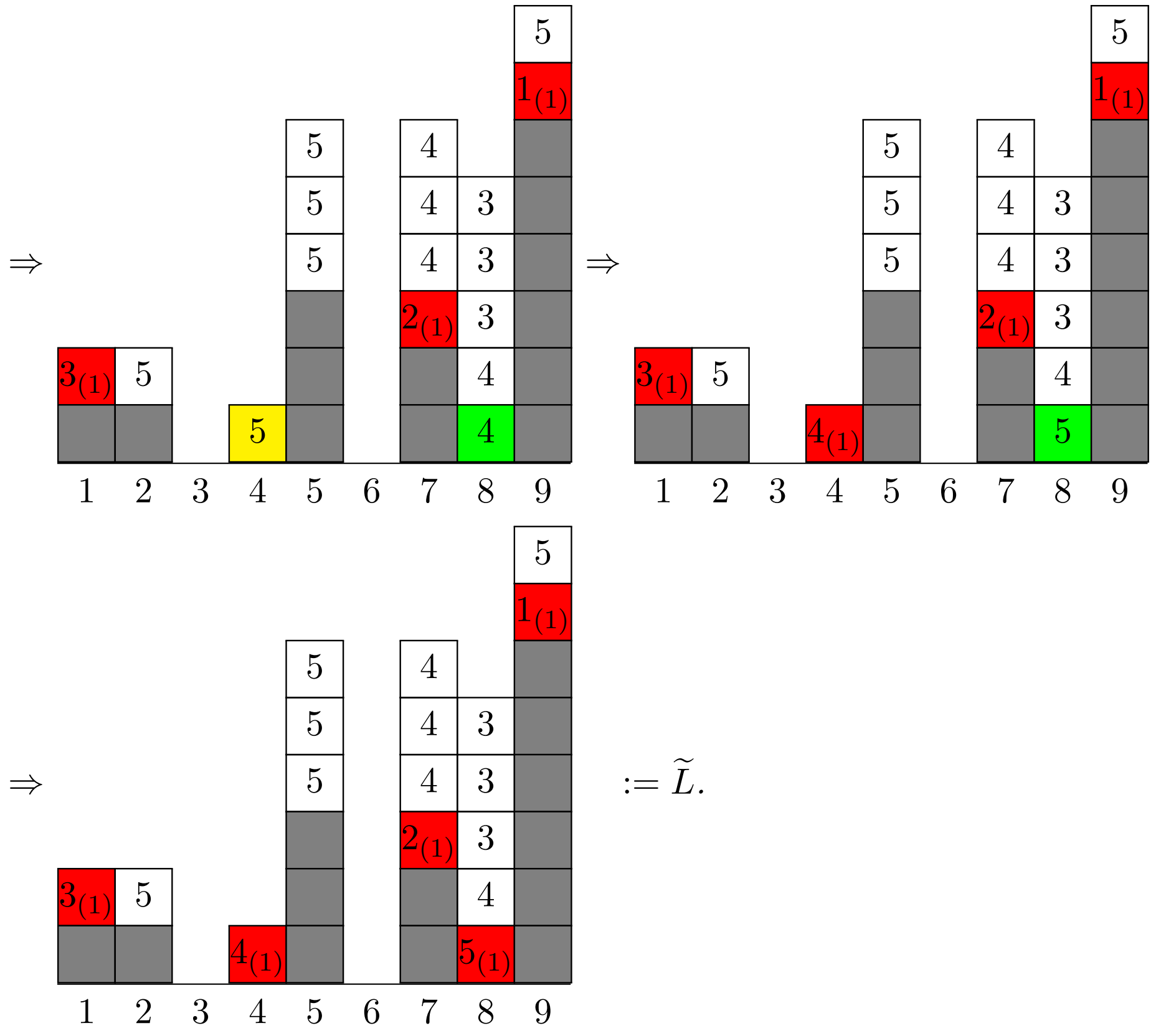

$:= \widetilde{L}.$

*One can verify that $\widetilde{L}$ is actually the recording tableau of $(U \leftarrow \widetilde{W})$.*

**Lemma 4.12.** *Let $U$ be an SSAF with basement $\epsilon_n$ and shape $\alpha$ for some positive integer $n$ and $l(\alpha) \leq n$. Let $w = a_{11} \dots a_{1c_1} | a_{21} \dots a_{2c_2} | \cdots | a_{k1} \dots a_{kc_k}$ with $k \leq n$ be a column word. Let $V$ be the SSAF representing $w$ (i.e. inserting $w$ into an empty atom). Let $\tilde{w} = x_{11} x_{12} \dots x_{r1} \cdots x_{c_1,1} \dots x_{c_1,r_{c_1}}$ be the row reading word of $V$ (as $V$ has $c_1$ rows and $k$ columns and so $r_1$ is $k$), where $k = r_1 \geq r_2 \geq \cdots \geq r_{c_1} > 0$ are the row lengths of $V$ from bottom to top. Then the recording tableau $L$ of $(U \leftarrow W)$ determines the recording tableau $R$ of $(U \leftarrow \widetilde{W})$, where*

$$\begin{cases} W = \begin{pmatrix} k & k & \dots & k & \cdots & 1 & 1 & \dots & 1 \\ a_{11} & a_{12} & \dots & a_{1c_1} & \cdots & a_{k1} & a_{k2} & \dots & a_{kc_k} \end{pmatrix} \\ \\ \widetilde{W} = \begin{pmatrix} 1 & 1 & \dots & 1 & \cdots & c_1 & c_1 & \dots & c_1 \\ x_{11} & x_{12} & \dots & x_{1r_1} & \cdots & x_{c_1,1} & x_{c_1,2} & \dots & x_{c_1,r_{c_1}} \end{pmatrix} \end{cases}.$$

*Proof.* By Lemma 4.11 and the corresponding word in Lemma 4.4, given $L$, we know how to enter all the 1's in $R$, which are those cells marked $1_{(1)}, 2_{(1)}, \dots, k_{(1)}$ after applying Lemma 4.11 as the $b_{k1} b_{k-1,1} \dots b_{11} = x_{11} x_{12} \dots x_{1r_1}$ (by the argument after Lemma 4.4).

We remove the cells from $\widetilde{L}$ in Lemma 4.11 to create a new $L$ to apply Lemma 4.11 on, we can get the second row entries (as described in the paragraphs after Lemma 4.4),and hence we know how to put all the 2's into $R$. By the same argument, we can fill in all entries in $R$ and hence $L$ determines $R$.

□

**Example 17.** *We continue with the* $\widetilde{L}$ *in Example 16 to illustrate Lemma 4.12.*

*We mark the final position of the subscripted cells a different color for a different subscript. Starting with what we get in Example 16:*

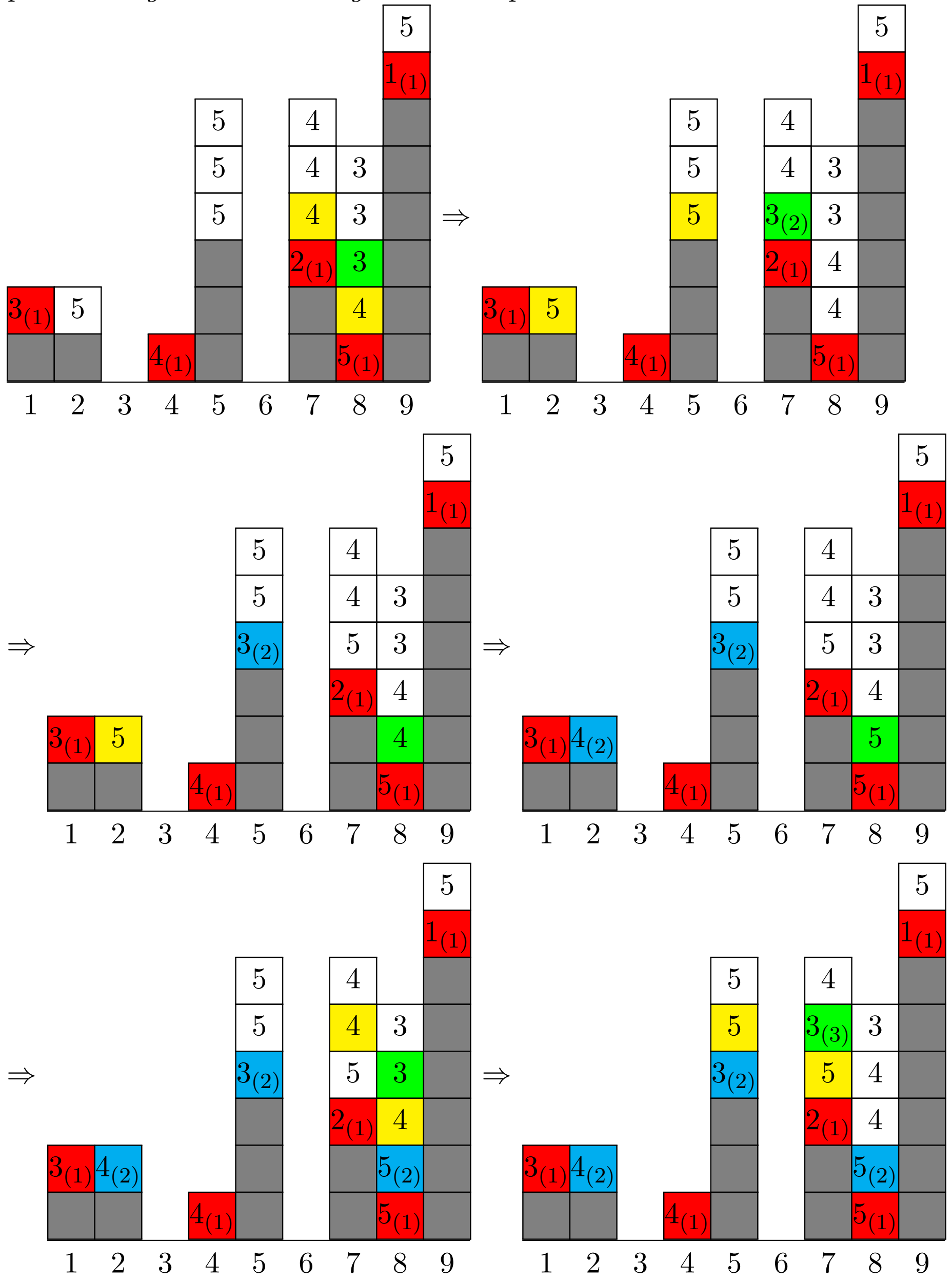

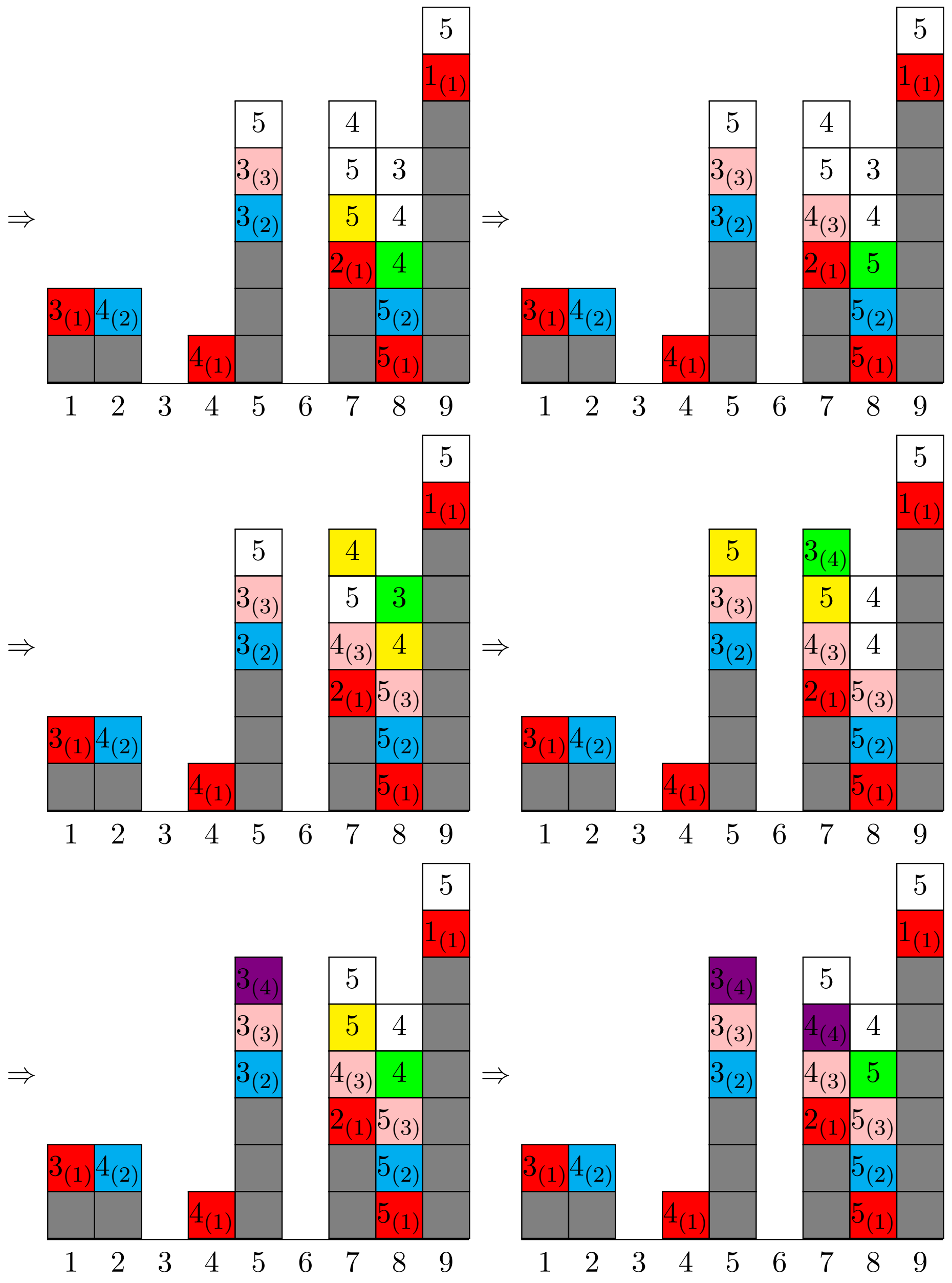

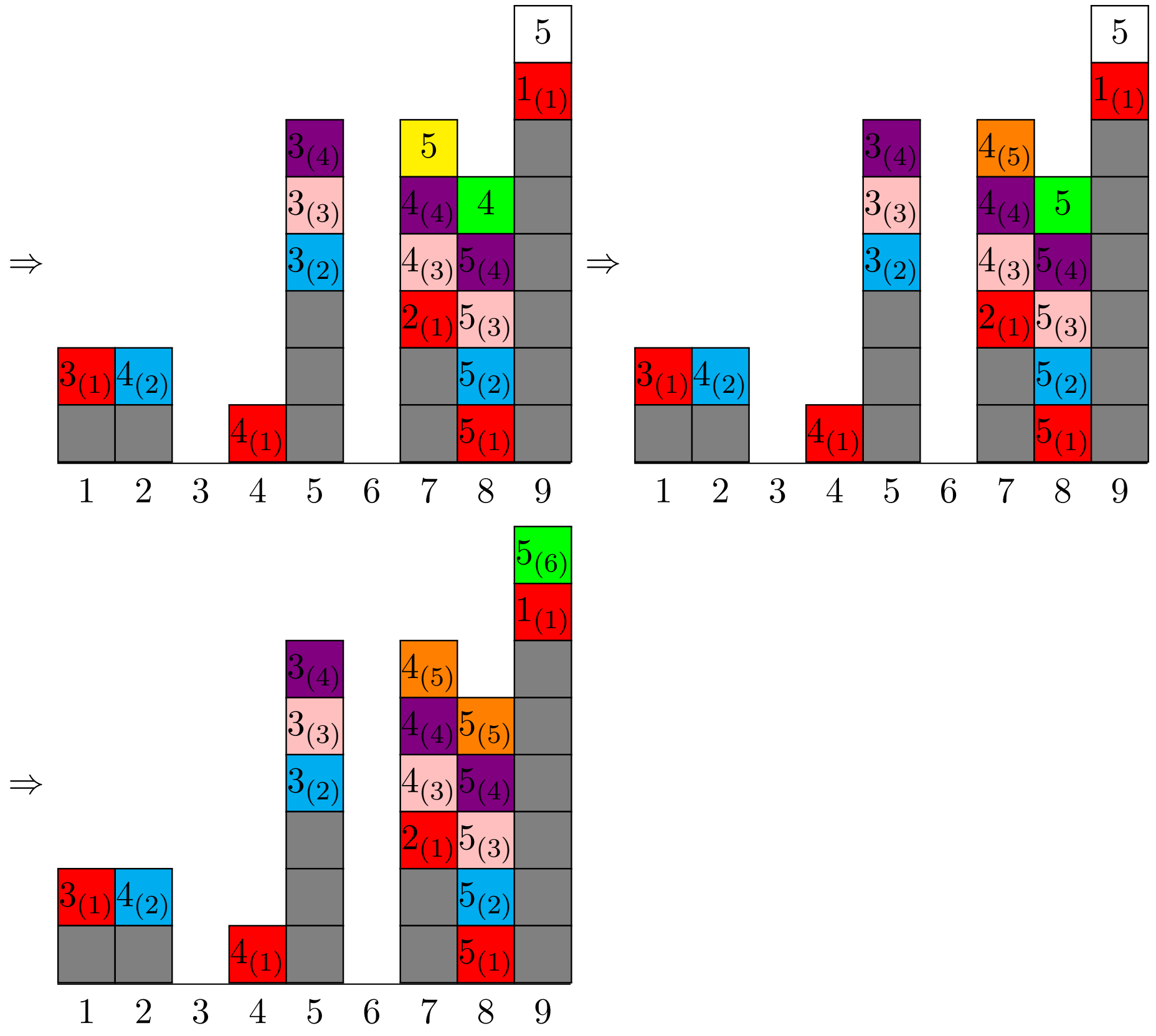

*Now by replacing the entry of each cell with the subscript number, we get*

$R =$ 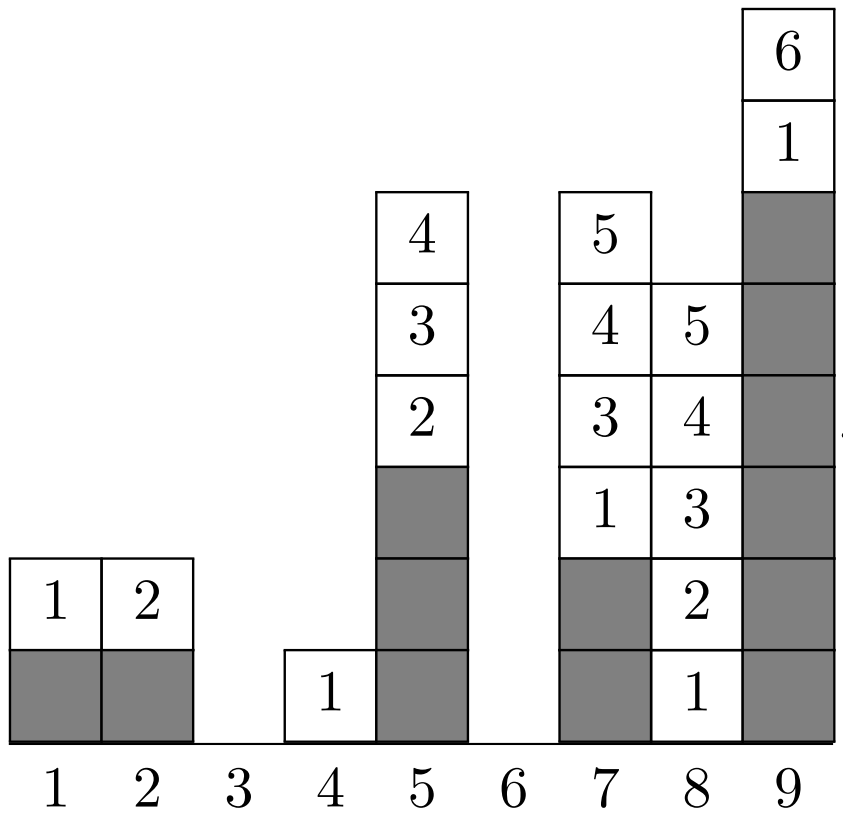 .

*One can verify that $R$ can be obtained by using the row word mentioned in Example7.*

### 4.3 Decomposition of the product of a dominating monomial and an atom into a positive sum of atoms

We now prove Theorem 4.1 mentioned in the beginning of this Section. We rephrase the Theorem in a more precise form as follows.

**Theorem 4.13.** *Let $\lambda$ be a partition and $\alpha$ be a weak composition. Let $\mathcal{A}_\lambda$ and $\mathcal{A}_\alpha$ be atoms of shape $\lambda$ and $\alpha$ respectively. Then*

$$\mathcal{A}_\lambda \cdot \mathcal{A}_\alpha = \sum_{\beta \vDash |\lambda|+|\alpha|, \lambda \subseteq \beta} c^{\beta}_{\lambda\alpha} \mathcal{A}_\beta$$

*where $c^{\beta}_{\lambda\alpha}$ is the number of distinct LRS of shape $\beta/\lambda$ created by column words whose corresponding SSAF has shape $\alpha$.*

*Proof.* Since $\mathcal{A}_\lambda = x^\lambda$ (there is exactly one SSAF with shape $\lambda$, denoted by $U_\lambda$) and $\mathcal{A}_\alpha = \sum_{F \in SSAF(\alpha)} x^F$, we have

$$\mathcal{A}_\lambda \cdot \mathcal{A}_\alpha = \sum_{F \in SSAF(\alpha)} x^\lambda x^F.$$

To prove the theorem, we only need to check that given a LRS $L$ of shape $\beta/\lambda$ created by column words whose corresponding SSAF has shape $\alpha$, if there is some column word

$$w = a_{11} \dots a_{1c_1} | a_{21} \dots a_{2c_2} | \cdots | a_{k1} \dots a_{kc_k}$$

and biword

$$W = \begin{pmatrix} k & k & \dots & k & \cdots & 1 & 1 & \dots & 1 \\ a_{11} & a_{12} & \dots & a_{1c_1} & \cdots & a_{k1} & a_{k2} & \dots & a_{kc_k} \end{pmatrix}$$

such that $(U_\lambda \leftarrow W)$ creates the same $L$, then the SSAF corresponding to $w$ also has shape $\alpha$.

First consider the last cell in reading order among all those containing entry $k$ in $L$, that means it is the very first entry inserted. Since $a_{11}$ must be inserted in a cell immediately above the cell (including those in basement) containing $a_{11}$, the position of that cell fixes the value of $a_{11}$.

Now consider all the cells with $k$ and also the last cell in reading order among all those containing the entry $k-1$, then these cells are created by $U_\lambda \leftarrow a_{11} \dots a_{1c_1} a_{21}$. By Lemma 4.12, we know the corresponding row recording tableau and hence we know which two cells are the first two entries being inserted using the corresponding row word. Note that the first row consists of distinct entries and is inserted in ascending order using the row word, then by Lemma 15 in [HMR13], we know the cells are created in ascending reading order (one after another in reading order) and so they must be the cells immediately above those in $U_\lambda$, and hence the value inside each of those cells in the SSAF created by inserting the row word into $U_\lambda$ is exactly the value inside the cell just below it. Hence we know what the first two row entries of the corresponding SSAF of $a_{11}a_{12} \dots a_{1c_1} a_{21}$ are. Since we already know the first row entry, which is the lowest entry of the column corresponding to $a_{11} \dots a_{1c_1}$ is, we now know what the lowest entry of the second column (corresponding to $a_{11} \dots a_{1c_1} a_{21}$ and hence the same for $a_{11} \dots a_{1c_1} a_{21} \dots a_{2c_2}$) is.

We can repeat the same process until we get all the last entries of the $k$ columns and hence fix the shape of the SSAF corresponding to $w$. Since we read those entries just by considering $L$, this shows that $L$ fixes the shape of the corresponding SSAF of $w$ and result follows.

□

**Example 18.** *Pick* $\lambda = 4332221$ *and* $\alpha = (1, 0, 1, 0, 0, 4, 0, 6, 5)$.
*Let* $\beta = (5, 3, 5, 2, 4, 6, 1, 6, 2)$.

*Let* $U =$ 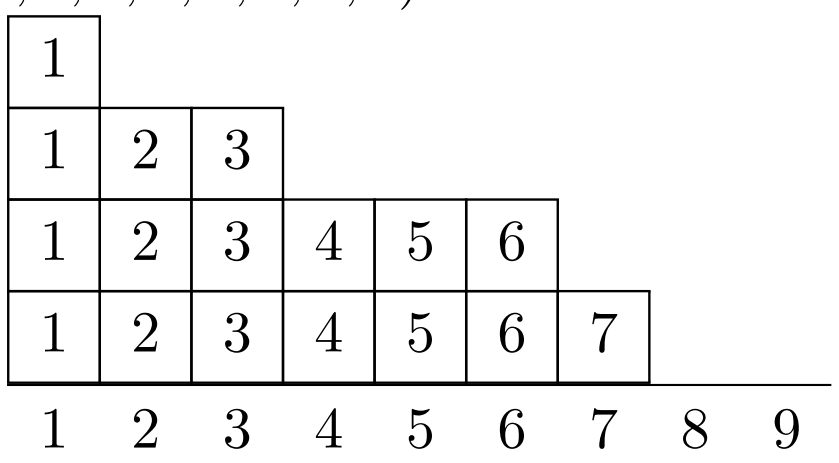

*One can first check that a recording tableaux* $L$ *(which is also an LRS if we change the basement) of shape* $\beta/\lambda$ *is created when the column word in Example 7 whose corresponding SSAF has shape* $\alpha$ *is inserted into* $U$*. Indeed, we have*

$L =$ 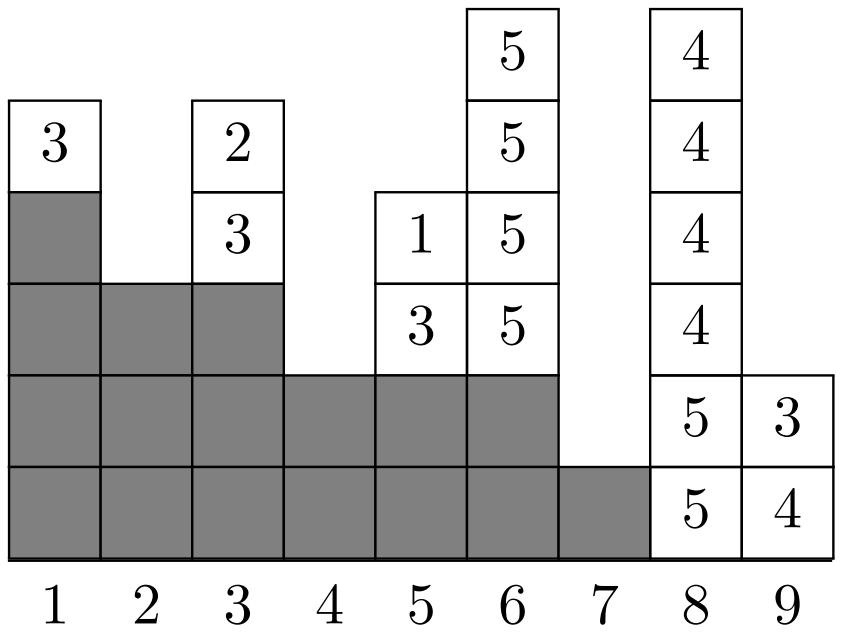

*If a column word* $w$ *would create the same* $L$ *when inserted in* $U$*, we know that it has column lengths* $6, 5, 4, 1, 1$*, and hence we can break it into 5 subsequences:*

$$w = a_{11} \dots a_{16} | a_{21} \dots a_{25} | a_{31} \dots a_{34} | a_{41} | a_{51}$$

*satisfying the conditions of being a column word in Definition 4.3. Let* $F(w)$ *be the SSAF corresponding to* $w$*, i.e. the SSAF created when inserting* $w$ *into an empty SSAF with basement being* $1 \quad 2 \quad 3 \quad 4 \quad \dots$ .

*Consider the cells with entry* 5 *in* $L$*:*

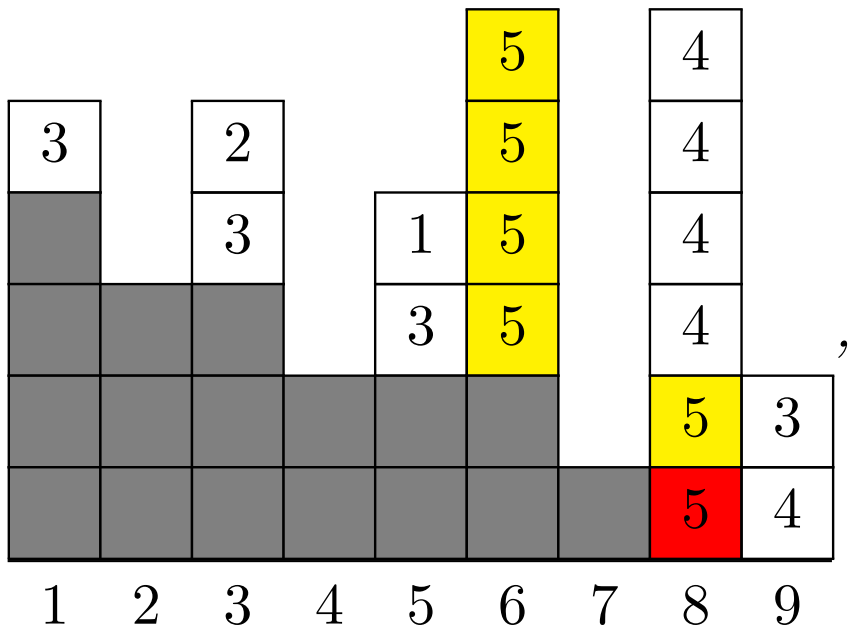 ,

*the largest in reading order (marked as red) is created when $a_{11}$ is inserted in $U$. Note that as $U$ is a partition shaped SSAF, $a_{11}$ must be placed immediately above the cell (including basement) containing the entry $a_{11}$. Hence $a_{11}$ must be 8. That means the column in $F(w)$ corresponding to $a_{11} \dots a_{16}$ is above the basement entry 8.*

*Next consider the cells with entry 5 and the last cell in reading order containing 4:*

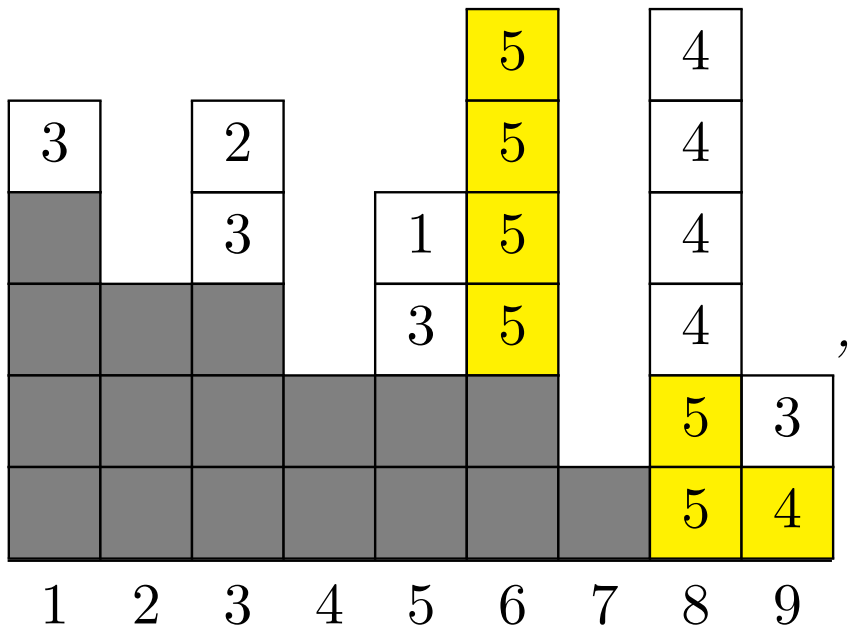

*by applying Lemma 4.11, we know the first two cells created when inserting the row word corresponding to the column word $a_{11} \dots a_{16}a_{21}$ must be the ones marked red as below:*

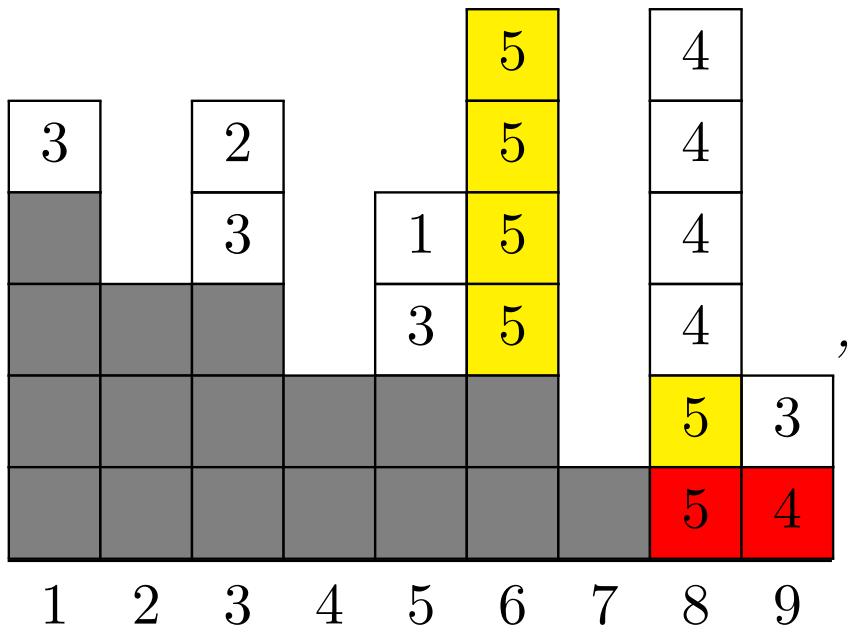

*and so we know the two numbers in the row word must be 8 and 9, and since we already know 8 is the basement entry which the new cell created is above, when inserting $a_{11} \dots a_{16}$ and hence we know the second column of $F(w)$ when inserting $a_{21} \dots a_{25}$ is above the basement entry 9.*

*Apply Lemma 4.11 on the cells with entries 4, 5 and also the last cell in reading order with entry 3, we know the first three cells created when inserting the corresponding row word of the column word $a_{11} \dots a_{16}a_{21} \dots a_{25}a_{31}$ into $U$, and we mark them red as shown:*

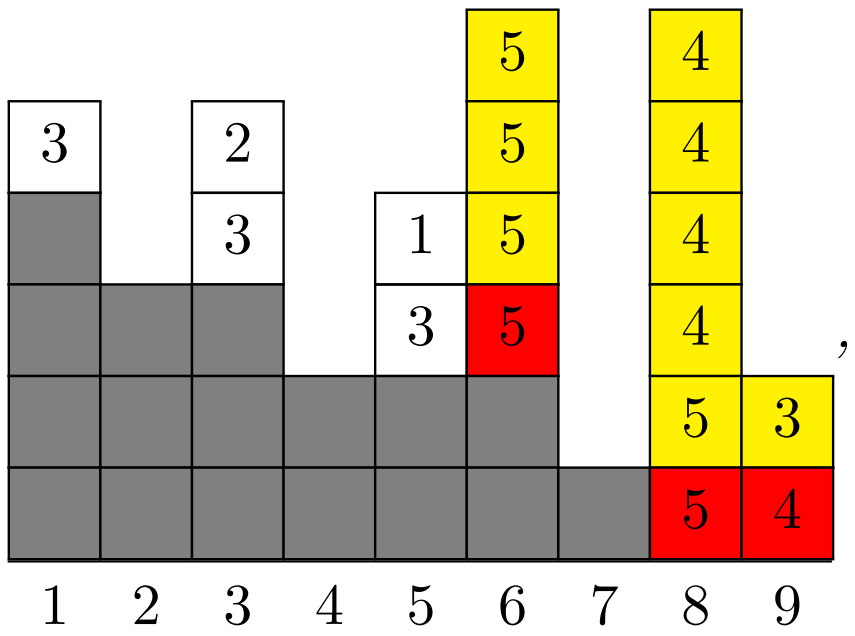

*That means the first three numbers in the corresponding row word is* $6, 8, 9$ *and since we already know the first two columns are above basement entries* $8$ *and* $9$*, we can conclude that the third column created in* $F(w)$ *when inserting* $a_{31} \dots a_{34}$ *is above the basement entry* $6$*.*

*Continue with the same process and we have:*

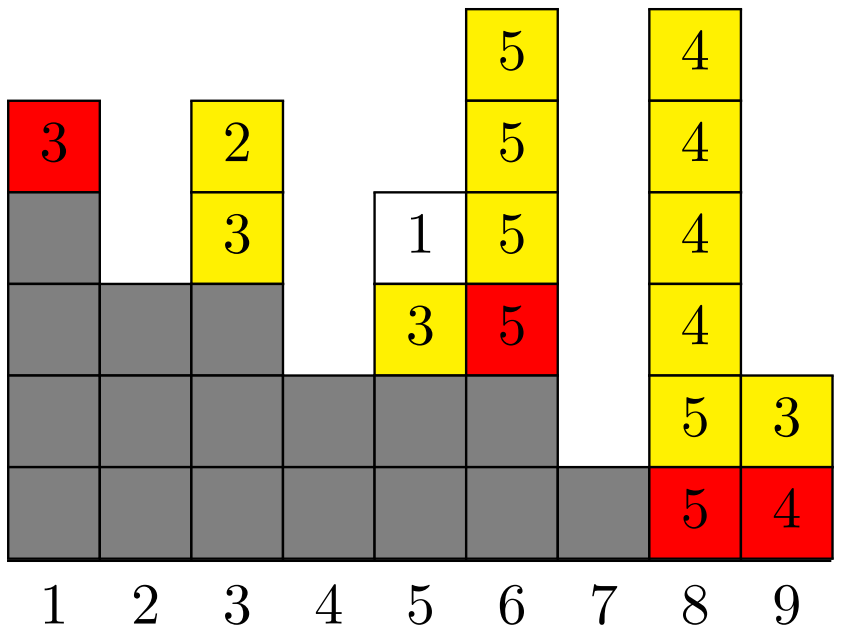

*and*

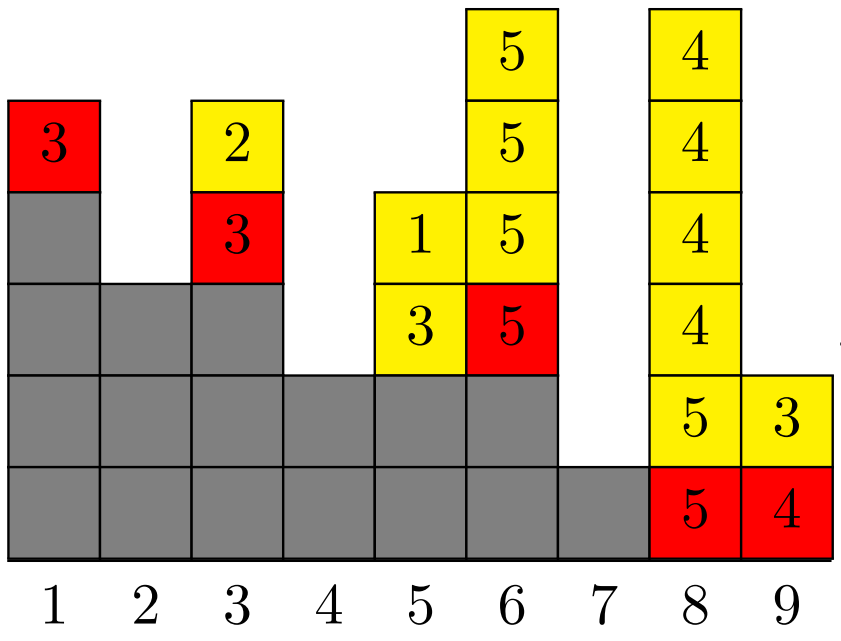

.

*Therefore we know the shape of* $F(w)$ *is* $(1, 0, 1, 0, 0, 4, 0, 6, 5)$ *which is exactly* $\alpha$*.*

## 4.4 Decomposition of the product of a dominating monomial and a key into a positive sum of keys

We would adapt the notations and apply the results in [HLMvW11] to prove the key-positivity property of the product of a dominating monomial and a key. However, we will still use $\overline{\alpha}$ instead of $\alpha^*$ which [HLMvW11] uses to denote the *reverse* of $\alpha$.

**Lemma 4.14.** *For any given partition* $\lambda$ *and weak compositions* $\beta, \gamma$ *such that* $\lambda \subseteq \beta$*, Let*

$$S_1 := \left\{ K \;\middle|\; \begin{array}{c} K \text{ is a LRK of shape } \overline{\delta}/\lambda \text{ with content created } \overline{\omega_\gamma(\gamma)} \text{ and } \phi(K) \text{ is an LRS} \\ \text{created by a column word whose corresponding SSAF has a shape} \geq \overline{\gamma} \text{ and} \\ \text{basement shape } \lambda,\ \overline{\delta} \leq \beta \end{array} \right\},$$

$$S_2 := \left\{ L \;\middle|\; \begin{array}{c} L \text{ is an LRS of shape } \beta/\lambda \text{ created by a column word whose corresponding} \\ \text{SSAF has shape } \alpha, \alpha \geq \overline{\gamma} \end{array} \right\}$$

*(where* $\phi$ *is defined in the proof in Theorem 6.1 in [HLMvW11]. Then* $\phi|_{S_1} : S_1 \to S_2$ *is a bijection.)*

*Proof.* Let $K \in S_1$. Then by fixing $\beta$, we know $\phi(K)$ has an overall shape $\beta$. Hence $\phi(K)$ has shape $\beta/\lambda$ and thus $\phi(K) \in S_2$.

Therefore we have $\phi(S_1) \subseteq S_2$.

Since $\phi$ and hence $\phi|_{S_1}$ is injective, it remains to check $\phi|_{S_1} : S_1 \to S_2$ is surjective.

Let $L \in S_2$. Since $L$ is created by a column word whose corresponding SSAF has shape $\alpha$ and $\alpha \geq \overline{\gamma}$, we know $L$ has content $\overline{\omega_\gamma(\gamma)}$. Then by the proof in Theorem 6.1 in [HLMvW11], $\phi^{-1}(L)$ is an LRK of shape $\overline{\delta}/\lambda$ for some $\overline{\delta} \leq \beta$ and has the same column sets as $L$. Hence $\phi^{-1}(L) \in S_1$. Therefore for any $L \in S_2$, we can find an LRK, namely, $\phi^{-1}(L) \in S_1$ such that $\phi|_{S_1}(\phi^{-1}(L)) = L$.

We thus have $\phi|_{S_1}$ is surjective and result follows. $\square$

**Theorem 4.15.** *Let $\lambda$ be a partition and $\gamma$ be a weak composition. Let $\mathcal{A}_\lambda$ and $\kappa_\gamma$ be an atom of shape $\lambda$ and a key of shape $\gamma$ respectively. Then*

$$\mathcal{A}_\lambda \cdot \kappa_\gamma = \sum_{\alpha} b^{\alpha}_{\lambda\gamma}\kappa_\alpha$$

*where $b^{\alpha}_{\lambda\gamma}$ is the number of distinct LRK of shape $\overline{\alpha}/\lambda$ with content $\overline{\omega_\gamma(\gamma)}$ and the image under $\phi$ is an LRS with basement shape $\lambda$ created by a column word whose corresponding SSAF has shape $\rho$ for some $\rho \geq \overline{\gamma}$.*

*Proof.* By Theorem 1 and Theorem 4.13, we have

$$\mathcal{A}_\lambda \cdot \kappa_\gamma = \mathcal{A}_\lambda \cdot \sum_{\alpha \geq \overline{\gamma}} \mathcal{A}_\alpha = \sum_{\alpha \geq \overline{\gamma}} (\mathcal{A}_\lambda \cdot \mathcal{A}_\alpha) = \sum_{\alpha \geq \overline{\gamma}} \sum_{\substack{\beta \vDash |\lambda|+|\alpha| \\ \beta \supseteq \lambda}} c^{\beta}_{\lambda\alpha}\mathcal{A}_\beta = \sum_{\beta \supseteq \lambda} \sum_{\substack{\alpha \geq \overline{\gamma} \\ \alpha \vDash |\beta|-|\lambda|}} c^{\beta}_{\lambda\alpha}\mathcal{A}_\beta.$$

Also, by Theorem 1, we have

$$\sum_{\delta} b^{\delta}_{\lambda\gamma}\kappa_\delta = \sum_{\delta}\sum_{\beta \geq \overline{\delta}} b^{\delta}_{\lambda\gamma}\mathcal{A}_\beta = \sum_{\beta}(\sum_{\overline{\delta} \leq \beta} b^{\delta}_{\lambda\gamma})\mathcal{A}_\beta.$$

Hence to prove the theorem, we only need to prove

$$\sum_{\substack{\alpha \geq \overline{\gamma} \\ \alpha \vDash |\beta|-|\lambda|}} c^{\beta}_{\lambda\alpha} = \sum_{\beta \geq \overline{\delta}} b^{\delta}_{\lambda\gamma}$$

which follows from Lemma 4.14 as $|S_2| = \sum_{\substack{\alpha \geq \overline{\gamma} \\ \alpha \vDash |\beta|-|\lambda|}} c^{\beta}_{\lambda\alpha}$ and $|S_1| = \sum_{\beta \geq \overline{\delta}} b^{\delta}_{\lambda\gamma}$

□

## 4.5 Decomposition of the product of a Schur function and a Demazure character

It is proved in [HLMvW11] that the product of a Schur function and a Demazure character is key-positive using tableaux-combinatorics. In this section, we give another proof using linear operators.

**Lemma 4.16.** *If $f$ is atom-positive, then $\pi_i f$ is also atom positive for all positive integers $i$.*

*Proof.* Let $f = A_0(x) + \sum_{\substack{I, \\ I \text{is a reduced word}}} \theta_I A_I(x)$ where $A_0(x) = \sum_{\lambda \in Par} a_\lambda x^\lambda$ with $a_\lambda \in \mathbb{Z}_{\geq 0}$ and $A_I(x) = \sum_{\lambda \in Par} a^I_\lambda x^\lambda$ with $a^I_\lambda \in \mathbb{Z}_{\geq 0}$ for any reduced word $I$.

Let $\sigma_I$ be the permutation corresponding to a reduced word $I$. Then by Lemma 2.7, $|l(s_i\sigma_I) - l(\sigma)| = 1$. If $l(s_i\sigma_I) = l(\sigma_I) - 1$, then by Lemma 2.8, there exists a reduced word of $\sigma_I$ starting with $i$. Hence $\pi_i\theta_I = 0$ by item 5. of Proposition 3.1. If $l(s_i\sigma_I) = l(\sigma_I) + 1$, then $iI$ is also a reduced word and hence $\pi_i\theta_I = (1 + \theta_i)\theta_I = \theta_I + \theta_{iI}$.

As a result,

$$\pi_i f = \pi_i(A_0(x) + \sum_{\substack{I \\ I\text{is a reduced word}}} A_I(x)) = A_0(x) + \theta_i(A_0(x)) + \sum_{\substack{I, \\ iI\text{is a reduced word}}} (\theta_I A_I(x) + \theta_{iI} A_I(x)).$$

Since $A_0(x)$ and $A_I(x)$ are sums of dominating monomials with nonnegative integer coefficients, $\pi_i f$ is also atom positive.

□

**Lemma 4.17.** *Let $\sigma = n, n-1, \ldots, 1 \in S_n$. Then $\pi_i \pi_\sigma = \pi_\sigma$ and hence $\theta_i \sigma = 0$.*

*Proof.* By Corollary 2.10, there exists a reduced word of $\sigma$ starting with $i$ and hence by item 5. of Proposition 3.1, we have $\pi_i \pi_i = \pi_i$ which implies $\pi_i \pi_\sigma = \pi_\sigma$.

Hence $\theta_i \pi_\sigma = (\pi_i - 1)\pi_\sigma = \pi_i \pi_\sigma - \pi_\sigma = \pi_\sigma - \pi_\sigma = 0$. □

**Theorem 4.18.** *The product of a key and a Schur function is key positive.*

*Proof.* Let $\sigma = n, n-1, \ldots, 1 \in S_n$ and $\lambda, \mu$ be partitions with length at most $n$. By item 2 of Theorem 3.8, we have $s_\lambda = \pi_\sigma(x^\lambda)$.

Let $I = i_1 i_2 \ldots i_k$ be a reduced word of some permutation in $S_n$. We prove that $\pi_I(x^\mu) \times \pi_\sigma(x^\lambda)$ is key positive.

By Lemma 3.6,

$$\begin{aligned} & \pi_{i_k}(x^\mu \times \pi_\sigma(x^\lambda)) \\ = \; & \pi_{i_k}(x^\mu) \times \pi_\sigma(x^\lambda) + s_{i_k}(x^\mu) \times (\theta_{i_k} \pi_\sigma(x^\lambda)) \\ = \; & \pi_{i_k}(x^\mu) \times \pi_\sigma(x^\lambda) \end{aligned}$$

.

$$\begin{aligned} & \pi_{i_{k-1}} \pi_{i_k}(x^\mu \times \pi_\sigma(x^\lambda)) \\ = \; & \pi_{i_{k-1}}(\pi_{i_k}(x^\mu) \times \pi_\sigma(x^\lambda)) \\ = \; & (\pi_{i_{k-1}} \pi_{i_k}(x^\mu)) \times \pi_\sigma(x^\lambda) + s_{i_{k-1}} \pi_{i_k}(x^\mu) \times \theta_{i_{k-1}} \pi_\sigma(x^\lambda) \\ = \; & (\pi_{i_{k-1}} \pi_{i_k}(x^\mu)) \times \pi_\sigma(x^\lambda) \end{aligned}$$

Inductively, we get

$$\pi_I(x^\mu) \times \pi_\sigma(x^\lambda) = \pi_{i_1} \pi_{i_2} \cdots \pi_{i_k}(x^\mu) \times \pi_\sigma(x^\lambda) = \pi_{i_1} \pi_{i_2} \cdots \pi_{i_k}(x^\mu \times \pi_\sigma(x^\lambda)).$$

By Lemma 3.5, $\pi_\sigma(x^\lambda) = \sum_{\gamma \leq \sigma} \theta_\gamma(x^\lambda) = \sum_{\gamma \in S_n} \theta_\gamma(x^\lambda)$ which implies $x^\mu \times \pi_\sigma(x^\lambda) = \sum_{\gamma \in S_n} x^\mu \times \theta_\gamma(x^\lambda)$. Therefore by Theorem 4.13, $x^\mu \times \pi_\sigma(x^\lambda)$ is atom positive. By Lemma 4.16, we know $\pi_{i_k}(x^\mu \times \pi_\sigma(x^\lambda))$ is atom-positive. Inductively applying Lemma 4.16, $\pi_{i_1} \pi_{i_2} \cdots \pi_{i_k}(x^\mu \times \pi_\sigma(x^\lambda))$ is also atom-positive. □

# 5 Atom positivity of the product of two key polynomials whose basements have length at most 3

In this section, we prove Conjectue 1 for $l(\alpha), l(\beta) \leq 3$. Note that if $l(\alpha) < 3$, we can always add zero parts to it to increase the length to 3. Hence, we can assume $l(\alpha) = 3$. Similarly, we can assume $l(\beta) = 3$.

Let $\lambda_\alpha = \omega_\alpha(\alpha)$ and $\lambda_\beta = \omega_\beta(\beta)$. We claim that we can consider $l(\lambda_\alpha), l(\lambda_\beta) \leq 2$, i.e. both $\alpha$ and $\beta$ have at least one zero part.

First note that for integers $r \geq 0$ and $a \geq b \geq c \geq 0$, $(x_1x_2x_3)^r\theta_\tau(x_1^a x_2^b x_3^c) = \theta_\tau(x_1^{a+r}x_2^{b+r}x_3^{c+r})$ and $\pi_\tau(x_1^a x_2^b x_3^c) = (x_1x_2x_3)^c\pi_\tau(x_1^{a-c}x_2^{b-c})$ for any $\tau \in S_3$. That means the monomial $(x_1x_2x_3)^r$ times any atom is still an atom, and same for the case for key. We can also interpret this by considering fillings, as multiplying $(x_1x_2x_3)^r$ to an atom or a key is just adding $r$ bottom rows to the diagram and there is only one way to fill in these cells in an atoms or a key.

Suppose Conjectue 1 is true for $l(\lambda_\alpha), l(\lambda_\beta) \leq 2$. Then for any $\omega, \tau \in S_3$ and integers $a \geq b \geq c \geq 0$, $s \geq t \geq u \geq 0$,

$$\begin{aligned}
\pi_\omega(x_1^a x_2^b x_3^c) \cdot \pi_\tau(x_1^s x_2^t x_3^u) &= (x_1x_2x_3)^{c+u}\big(\pi_\omega(x_1^{a-c}x_2^{b-c}) \cdot \pi_\tau(x_1^{s-u}x_2^{t-u})\big) \\
&= (x_1x_2x_3)^{c+u}\big(\sum_\gamma c_\gamma \mathcal{A}_\gamma\big) \\
&= \sum_\gamma c_\gamma((x_1x_2x_3)^{c+u}\mathcal{A}_\gamma) \\
&= \sum_\gamma c_\gamma \mathcal{A}_{\gamma'}
\end{aligned}$$

where $c_\gamma$ are all nonnegative integers and $(x_1x_2x_3)^{c+u}\mathcal{A}_\gamma = \mathcal{A}_{\gamma'}$ with $\gamma' = \gamma + (c+u, c+u, c+u)$ (i.e. $\gamma'$ can be obtained by adding $c+u$ to each part of $\gamma$).

As a result, we now only consider compositions of length 3 and with at most two nonzero parts.

# 6 Polytopes

In this section, we introduce another way to view Demazure atoms and characters.

For each weak composition with length $k$, we can view it as a lattice point in $\mathbb{Z}^k_{\geq 0}$. We will focus on the case $k = 3$. (Everything in this section applies for any positive integer $k$.) Hence we have the bijection:

$$\alpha \leftrightarrow (\alpha_1, \alpha_2, \alpha_3) \in \mathbb{Z}^3_{\geq 0} \leftrightarrow x^\alpha = x_1^{\alpha_1}x_2^{\alpha_2}x_3^{\alpha_3}.$$

Consider the Coxeter arrangement (of type $A_2$) in $\mathbb{R}^3$:

$$\mathbf{Cox}(3) = \{a_i - a_j : 1 \leq i < j \leq 3\}$$

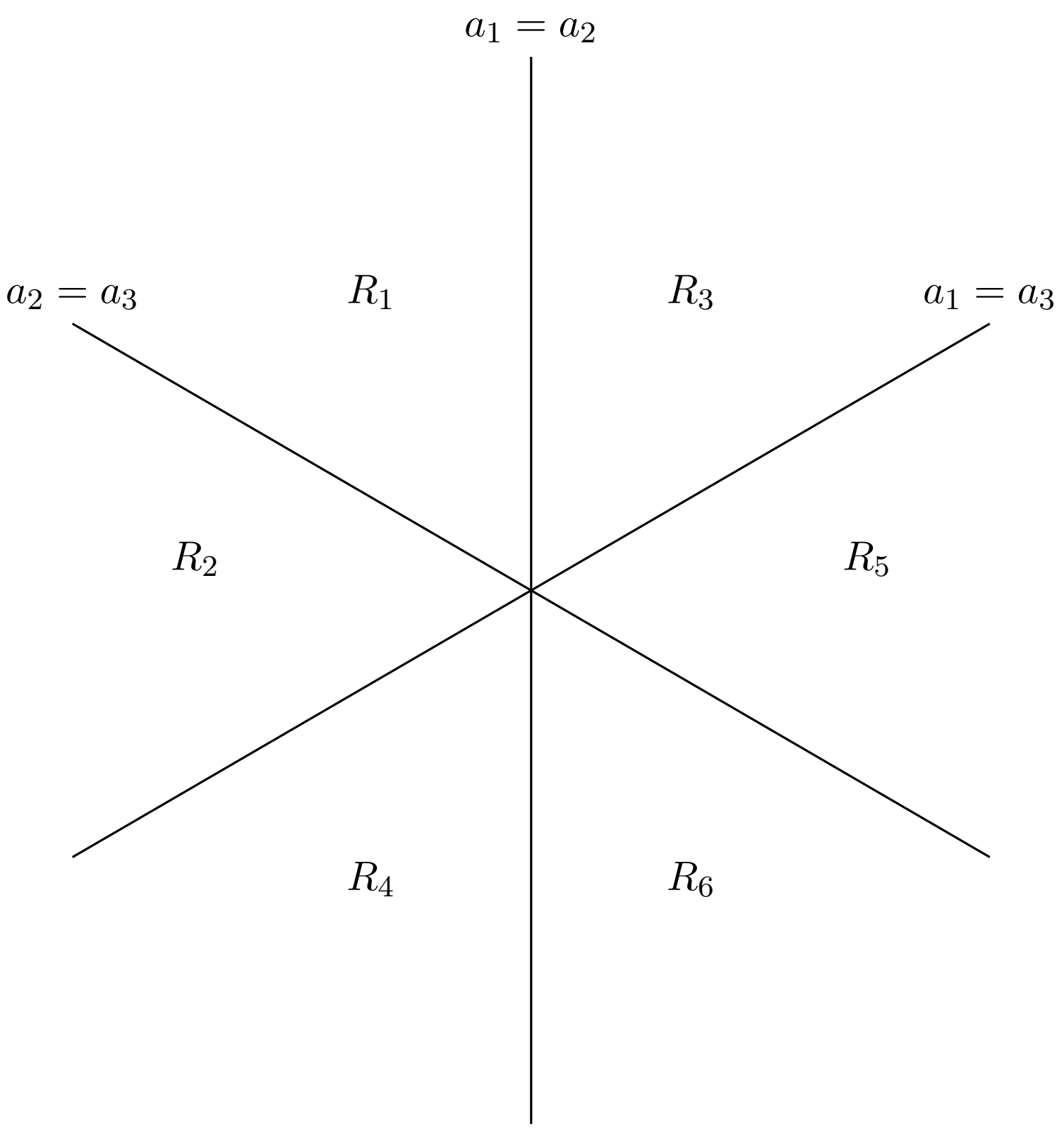

$$R_1:\ a_1 \geq a_2 \geq a_3;$$

$$R_2:\ a_1 \geq a_3 \geq a_2;$$

$$R_3:\ a_2 \geq a_1 \geq a_3;$$

$$R_4:\ a_3 \geq a_1 \geq a_2;$$

$$R_5:\ a_2 \geq a_3 \geq a_1;$$

$$R_6:\ a_3 \geq a_2 \geq a_1;$$

## 6.1 Demazure characters and polytopes

Let $\alpha$ be a weak composition with $\lambda_\alpha = (m, n, 0)$, for some integers $m \geq n \geq 0$. Then there are exactly 6 key polynomials obtained from $\lambda_\alpha$, namely:

$$x^{\lambda_\alpha}, \pi_1 x^{\lambda_\alpha}, \pi_2 x^{\lambda_\alpha}, \pi_{21} x^{\lambda_\alpha}, \pi_{12} x^{\lambda_\alpha}, \pi_{121} x^{\lambda_\alpha}.$$

We now plot each of these 6 key polynomials in the Coxeter arrangement:

Case 1. $x^{\lambda_\alpha} = x_1^m x_2^n$, so it corresponds to the point $(m, n, 0)$ in $R_1$:

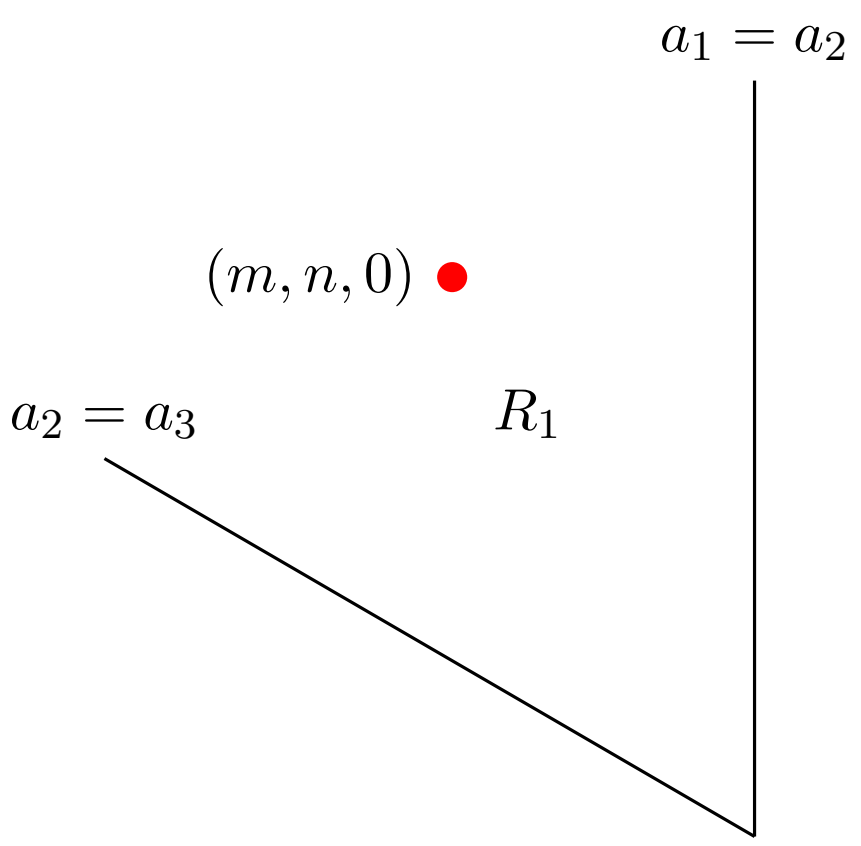

Case 2. $\pi_1 x^{\lambda_\alpha} = \pi_1(x_1^m x_2^n) = x_1^m x_2^n + x_1^{m-1} x_2^{n+1} + \cdots + x_1^n x_2^m$

Therefore it corresponds to the line joining $(m, n, 0)$ and $(n, m, 0)$. That is, each monomial corresponds to a lattice point (and vice versa) on the line obtained by joining $(m, n, 0)$ and its reflection along $a_1 = a_2$.

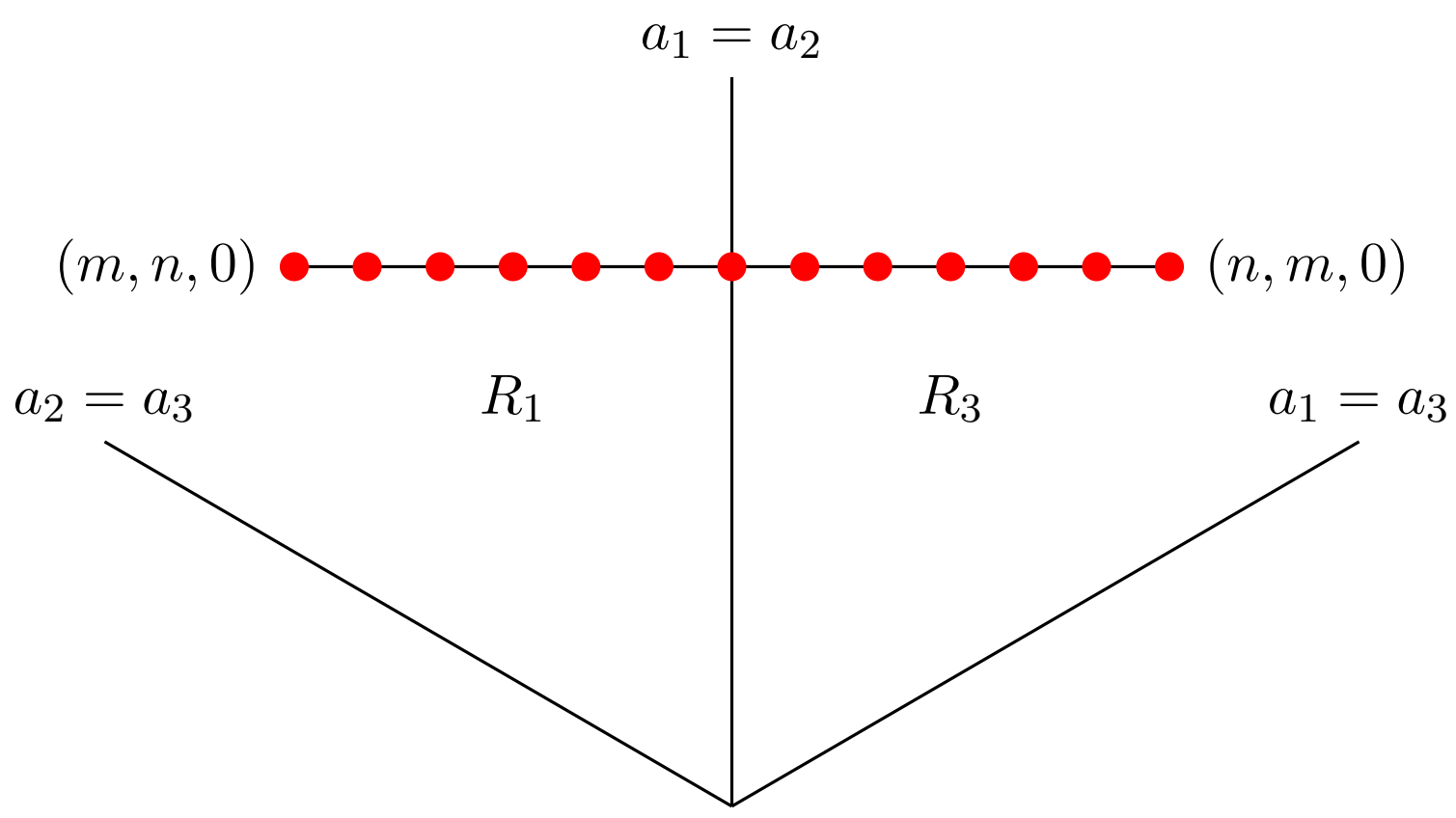

Case 3. $\pi_2 x^{\lambda_\alpha} = \pi_2(x_1^m x_2^n) = x_1^m x_2^n + x_1^m x_2^{n-1} x_3 + \cdots + x_1^m x_3^n$

Therefore it corresponds to the line joining $(m, n, 0)$ and $(m, 0, n)$. That is, each monomial corresponds to a lattice point (and vice versa) on the line obtained by joining $(m, n, 0)$ and its reflection along $a_2 = a_3$.

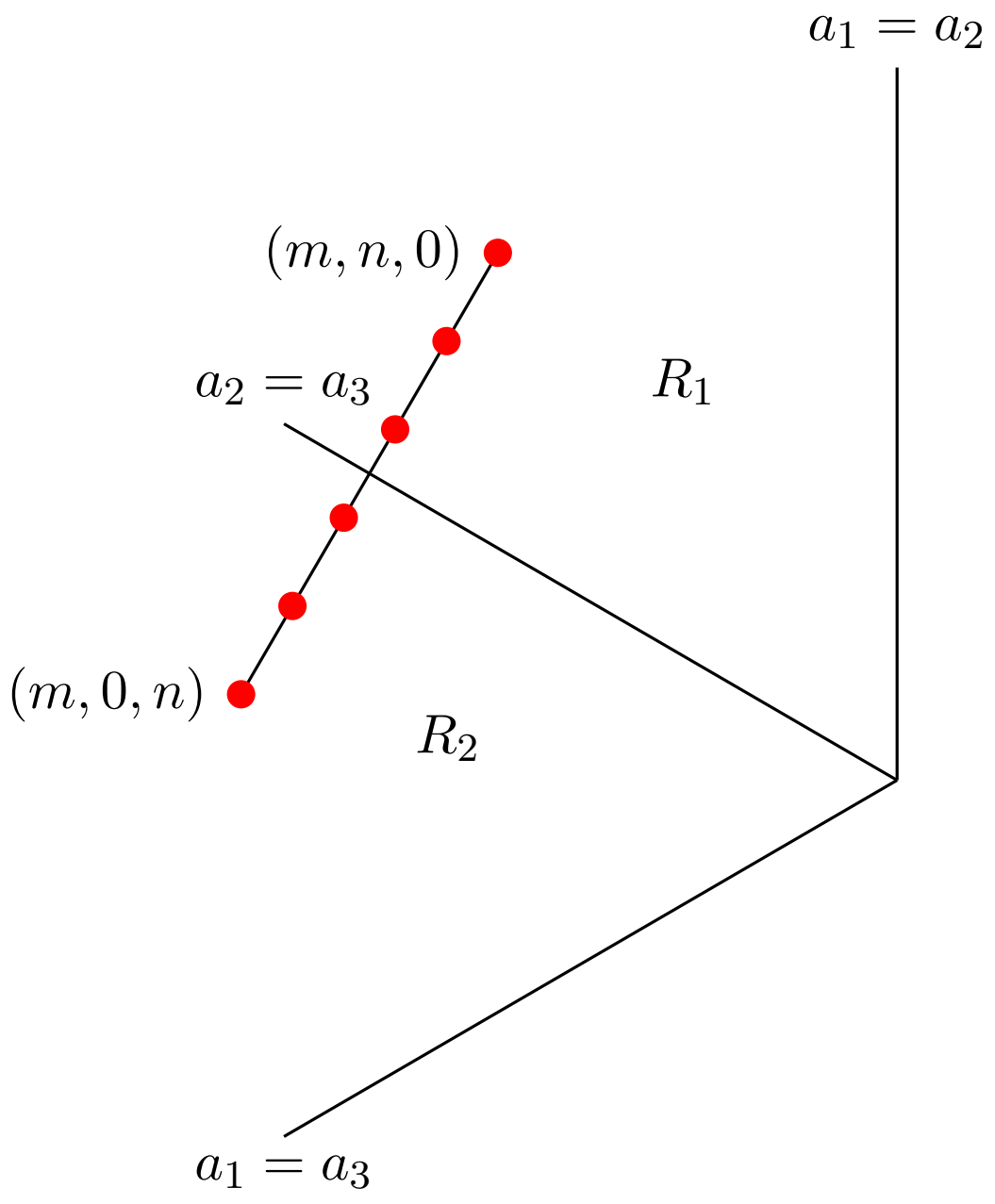

Case 4. $\pi_{21} x^{\lambda_\alpha} = \pi_2 \pi_1 (x_1^m x_2^n) = \pi_2 (x_1^m x_2^n + x_1^{m-1} x_2^{n+1} + \cdots + x_1^n x_2^m)$

Since $\pi_2(x_1^m x_2^n + x_1^{m-1} x_2^{n+1} + \cdots + x_1^n x_2^m) = \pi_2(x_1^m x_2^n) + \pi_2(x_1^{m-1} x_2^{n+1}) + \cdots + \pi_2(x_1^n x_2^m)$, we can apply the same correspondence in Case 3. for each key in the summand. Therefore $\pi_{21} x^{\lambda_\alpha}$ corresponds to the $m - n + 1$ lines obtained by reflecting each lattice point on the line joining $(m, n, 0)$ and $(m, 0, n)$ along $a_2 = a_3$. That is, each monomial corresponds to a lattice point (and vice versa) in the trapezoid obtained by first reflecting $(m, n, 0)$ along the line $a_1 = a_2$ followed by reflecting the resulting line along $a_2 = a_3$ as shown below:

First reflect along $a_1 = a_2$ and get a line joining $(m, n, 0)$ and $(m, 0, n)$:

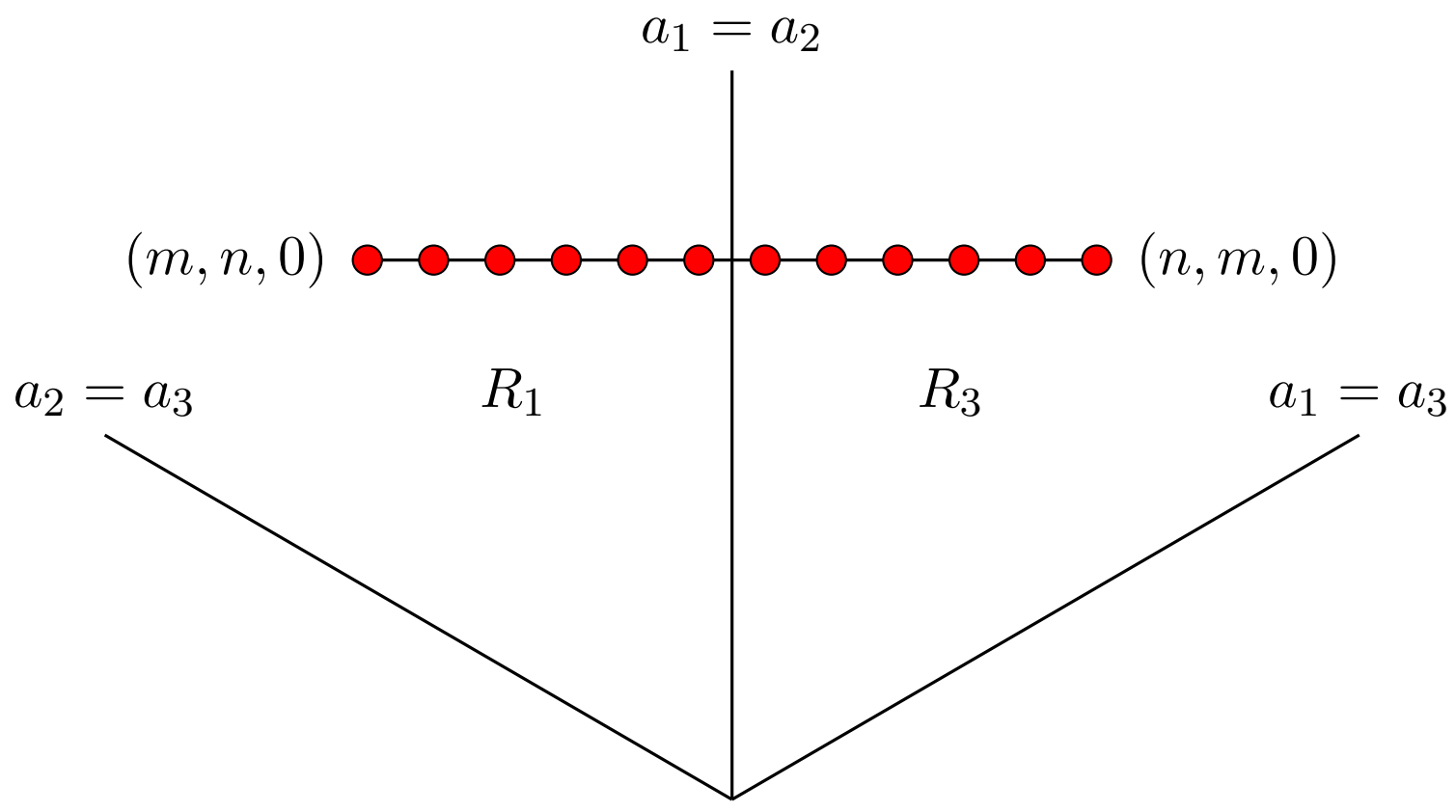

Then reflect the line along $a_2 = a_3$:

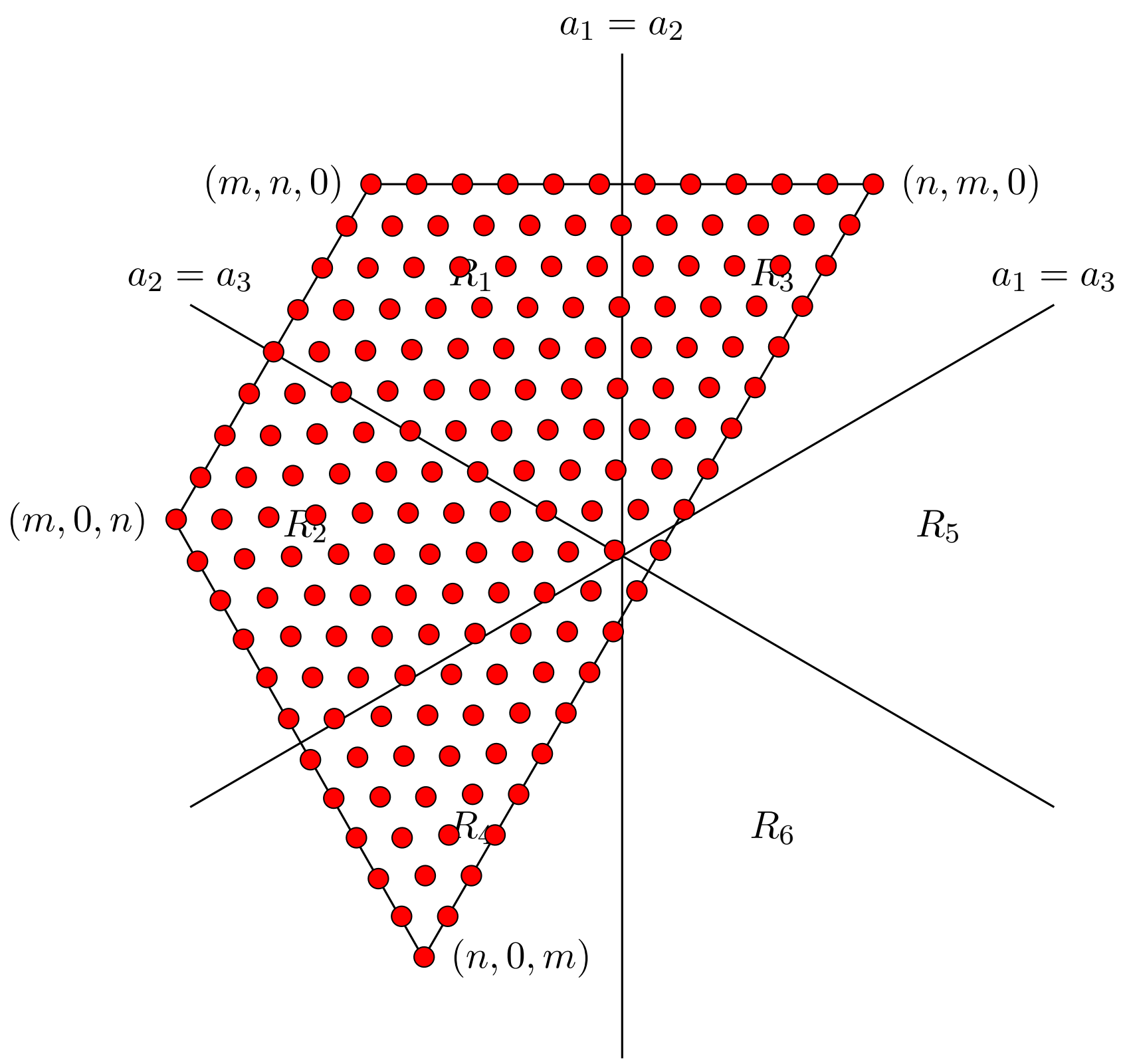

Note that there is another case where the trapezoid does not have any point in $R_5$ or $R_6$:

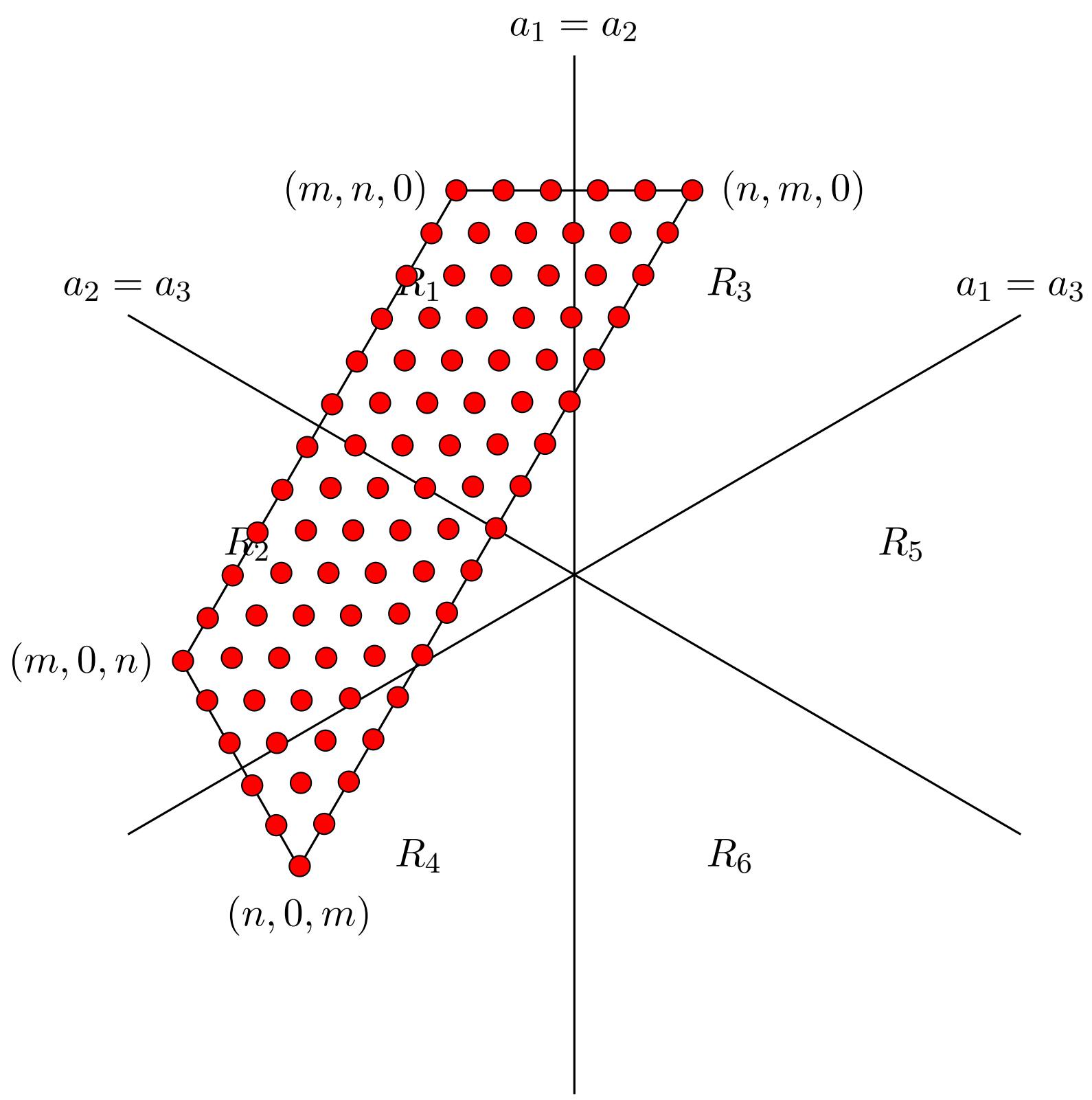

Indeed, the trapezoid has at least a point in $R_5$ or $R_6$ if and only if $m \geq 2n$ (one can prove that by locating the mid-point of the line joining $(n, m, 0)$ and $(n, 0, m)$.)

Case 5. $\pi_{12}x^{\lambda_\alpha} = \pi_1\pi_2(x_1^m x_2^n) = \pi_1(x_1^m x_2^n + x_1^m x_2^{n-1} x_3 + \cdots + x_1^m x_3^n)$

Similar to Case 4., it corresponds to the lattice points in the trapezoid formed by first reflecting $(m, n, 0)$ along $a_2 = a_3$ followed by reflecting the resulting line along $a_1 = a_2$. The trapezoid has at least a point in $R_4$ or $R_6$ if and only if $2n \geq m$. The trapezoid is as follows (we just shade the region using dotted pattern for convenience, but the actual correspondence should be lattice points in the shaded region including the boundary):

For $2n < m$:

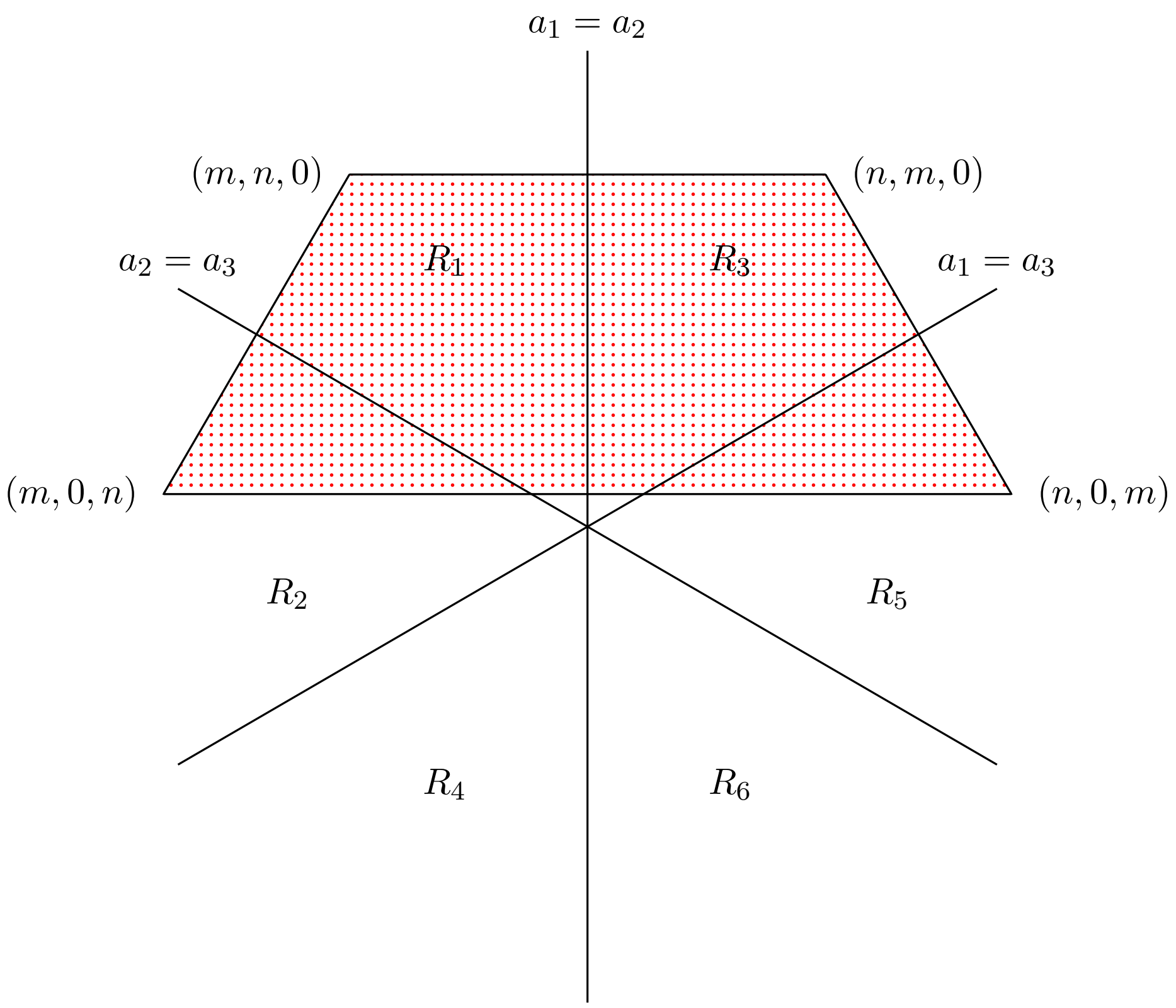

For $2n \geq m$:

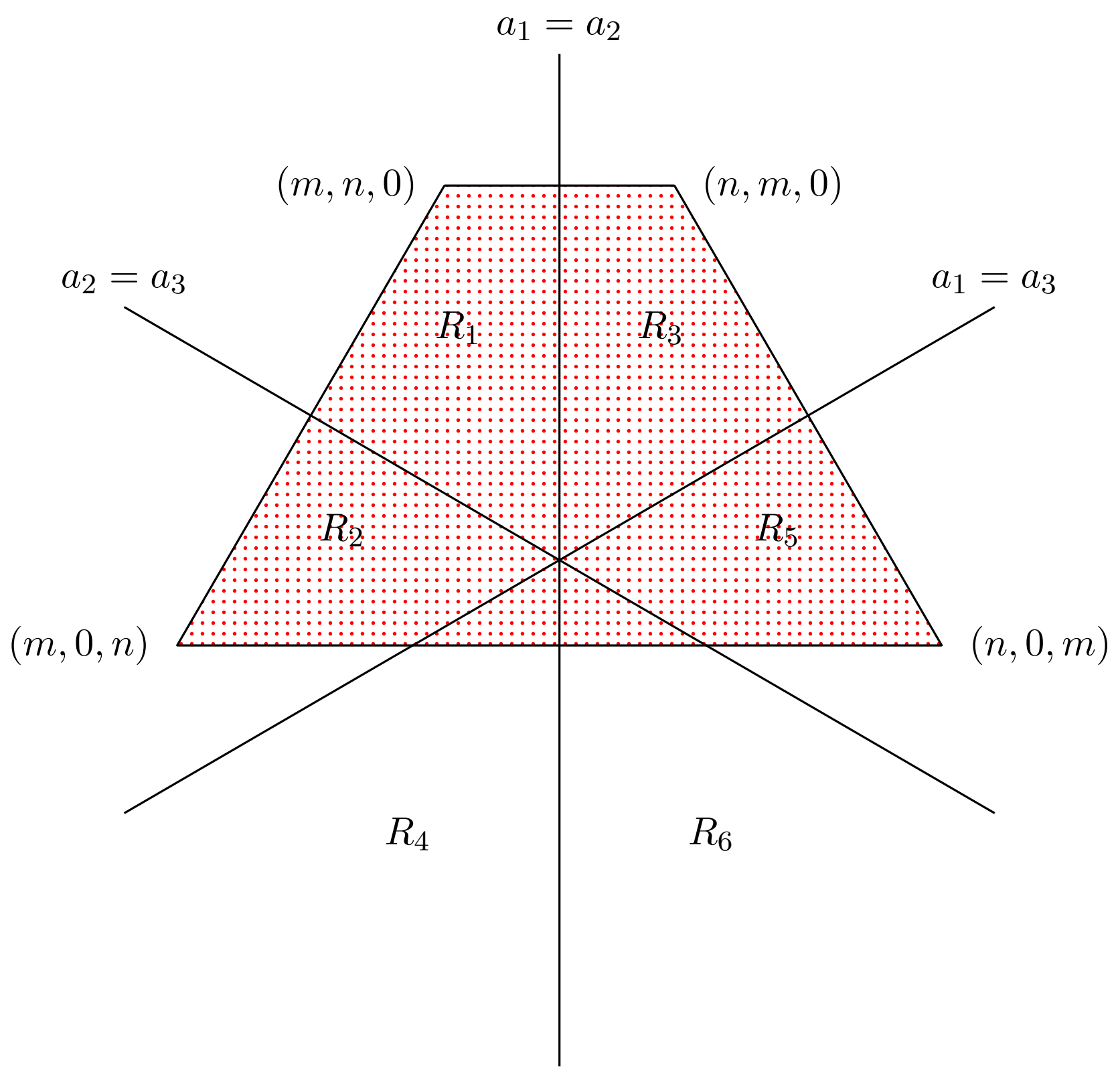

Case 6. $\pi_{121}x^{\lambda_\alpha} = \pi_1\pi_2\pi_1(x_1^m x_2^n) = \pi_1(\pi_2\pi_1(x_1^m x_2^n))$

First recall that by Proposition 3.1 $\pi_i s_i = (1+\theta_i)s_i = s_i + \theta_i s_i = s_i - \theta_i$, hence $\pi_i s_i + \pi_i = s_i - \theta_i + (1+\theta_i) = s_i + 1$.

For example, when $i = 1$, $\pi_1(x_1^a x_2^b) + \pi_1(x_1^b x_2^a) = x_1^a x_2^b + x_1^b x_2^a$ for any and any integers $a \geq b \geq 0$. Hence $\pi_1(x_1^b x_2^a) = x_1^a x_2^b + x_1^b x_2^a - \pi_1(x_1^a x_2^b)$.

We can plot $\pi_1(x_1^b x_2^a)$ as:

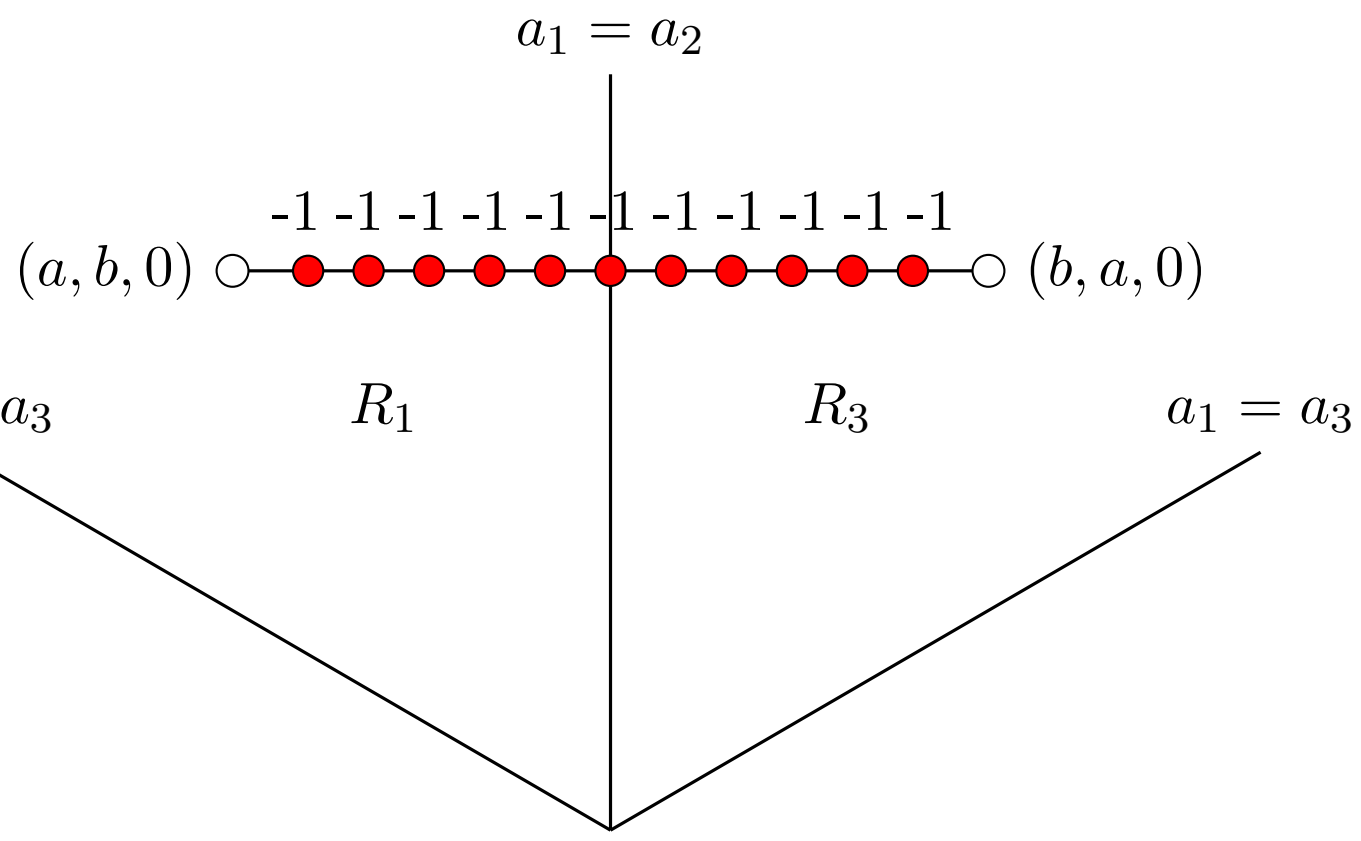

and $\pi_1(x_1^a x_2^b) + \pi_1(x_1^b x_2^a)$ as

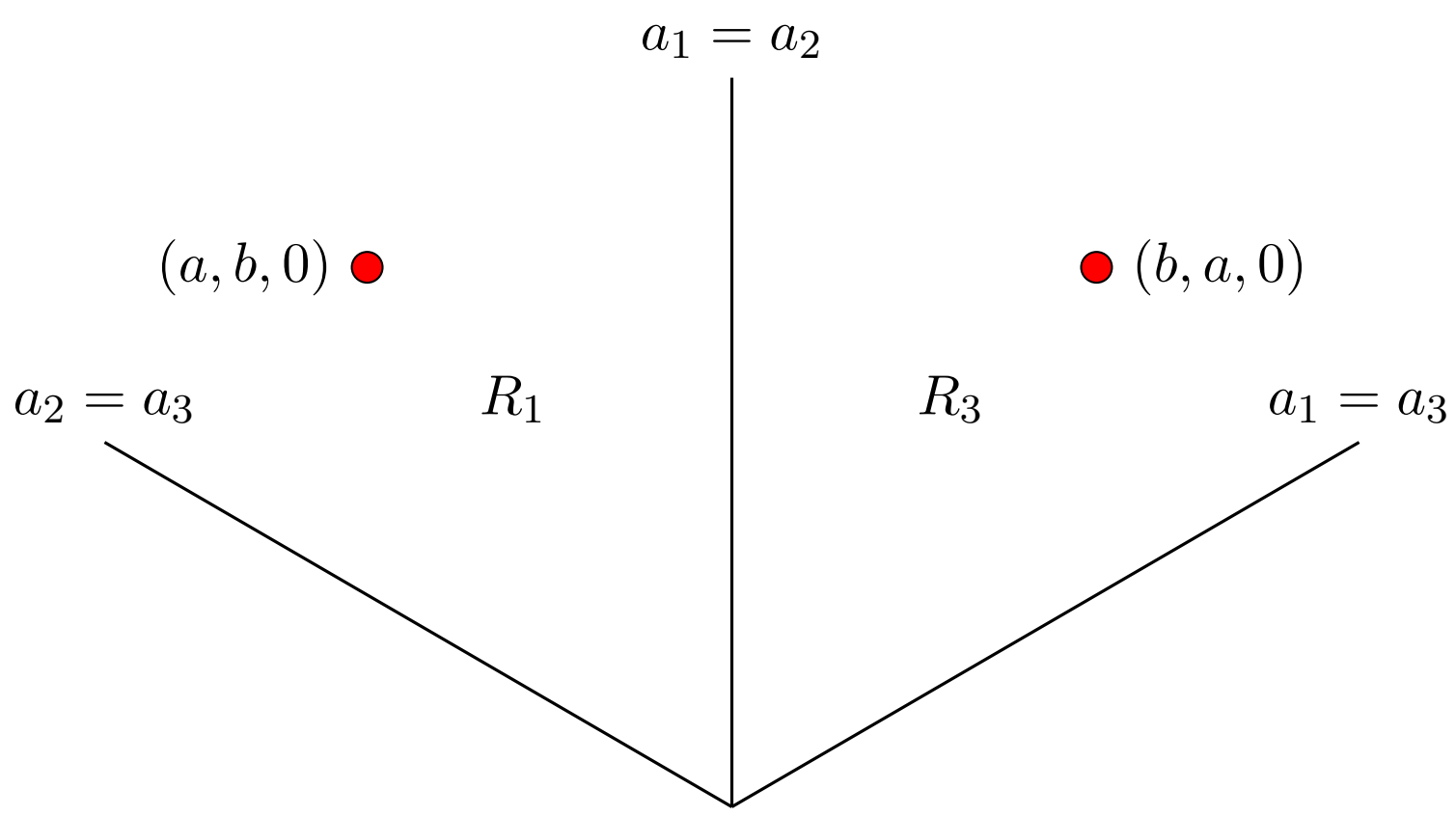

Now consider $\pi_{121}(x_1^m x_2^n) = \pi_1(\pi_{21}(x_1^m x_2^n))$.

By Case 4., we have $\pi_{21}(x_1^m x_2^n)$ as a trapezoid as follows:

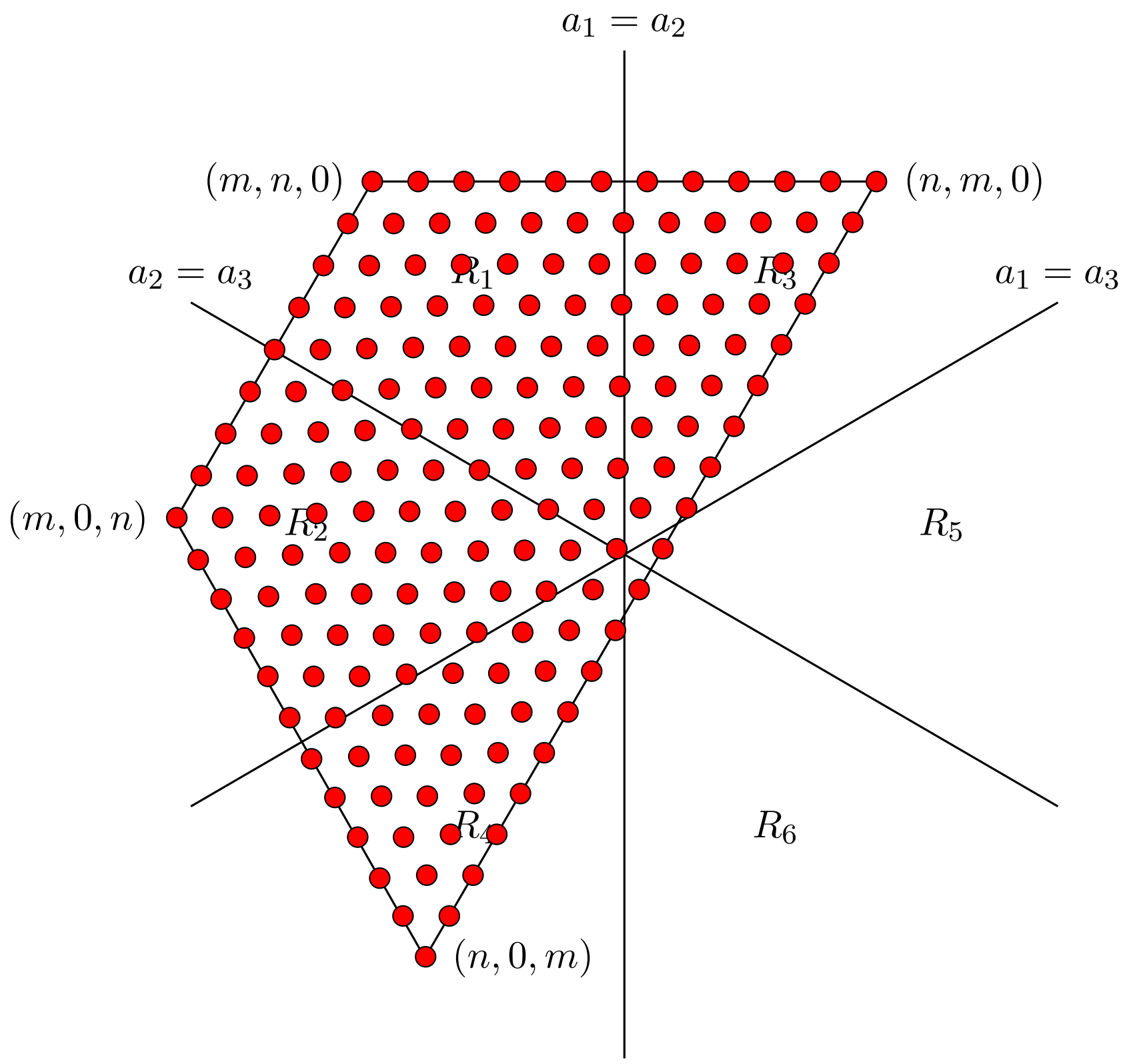

Apply $\pi_1$ to each lattice points in the trapezoid is equivalent to reflecting the trapezoid along $a_1 = a_2$ and get a hexagon with multiplicities on the points. Consider the multiplicity of a lattice point $(a_1, a_2, a_3)$ in the hexagon with $a_1 \geq a_2$. The multiplicity of its '$a_1 = a_2$'-reflection, namely, $(a_2, a_1, a_3)$ is the same. Recall that for each '$a_1 = a_2'$ - reflection pair in the trapezoid region (i.e. $(a.b, c)$ and $(b, a, c)$), applying $\pi_1$ to both of them results in the two points themselves. So the multiplicity of the lattice points increases by 1 'horizontally' from the boundary, starting from 1, until it first hits the line joining $(m,n,0)$ and $(0,n,m)$ or the line joining $(n,m,0)$ and $(n,0,m)$ (i.e. hit either of the lines in region $a_1 > a_2$), then it becomes stable. Here is an example for $m \geq 2n$:

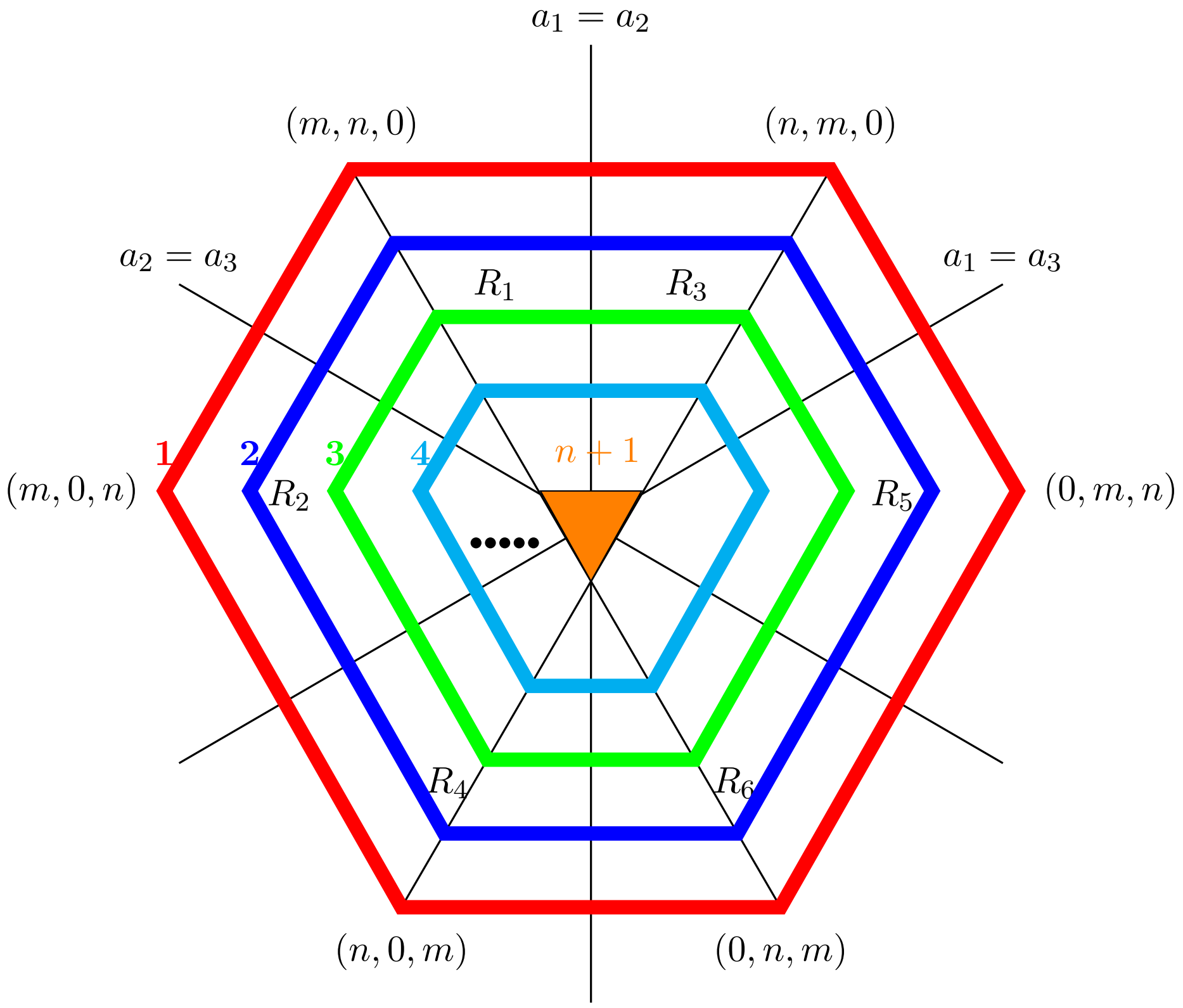

The monomials corresponding to the lattice points on the red boundary have coefficient 1 in the key polynomial, and those monomial corresponding to the lattice points on the blue boundary have coefficient 2 and so on, while those corresponding to the lattice points in the inner most triangle (the orange region) have coefficient $n+1$ for $m \geq 2n$ (and $m-n+1$ if $m \leq 2n$).

If we also plot the multiplicity (with $xy$ -plane being the Coxeter arrangement and $z$-axis represents the multiplicity), we get a polytope as follows:

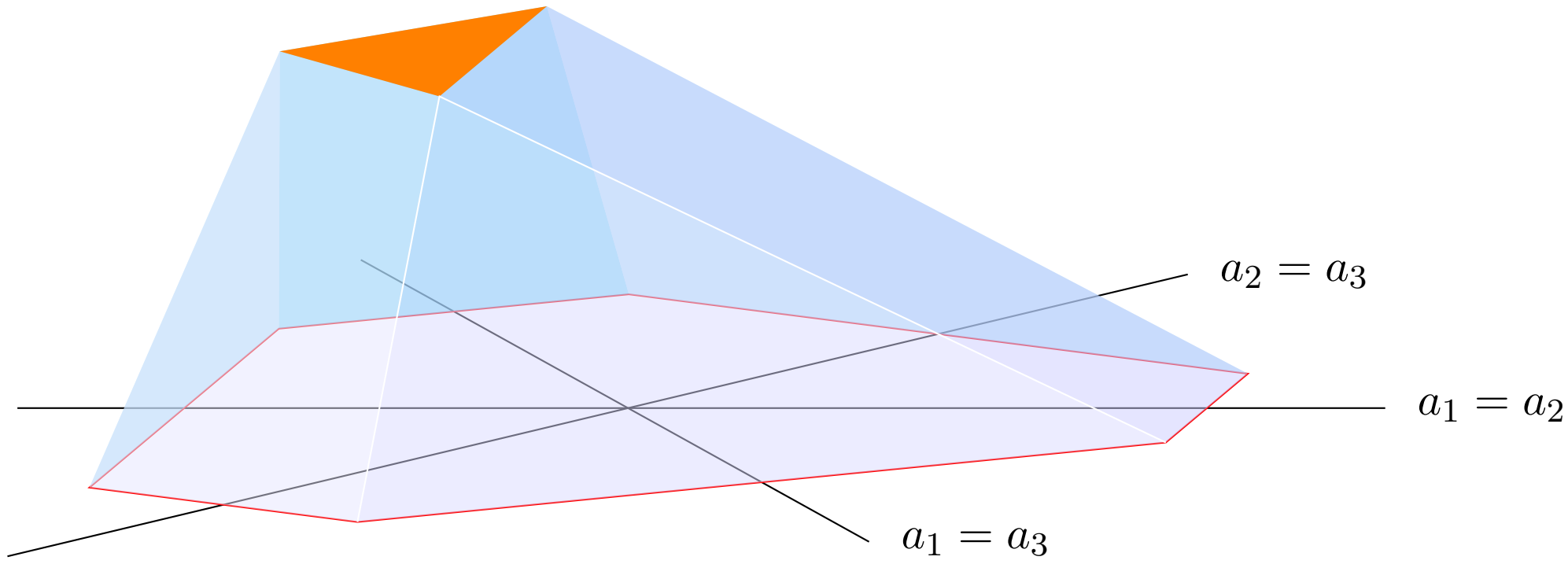

Note: One can also start with $\pi_{12}$ and apply $\pi_2$ on each lattice point in the trapezoid corresponding to $\pi_{12}$ and form $\pi_{212}$ to check that $\pi_{121} = \pi_{212}$.

### 6.2 Demazure atoms and polytopes

In this section, we will discuss how one can obtain a polytope from a Demazure atom. Again, we focus on the case where the shape of the atom is a weak composition of length 3 with at least one zero part.

Let $\alpha$ be a weak composition with $\lambda_\alpha = (m, n, 0)$, for some integers $m \geq n \geq 0$. Then there are exactly 6 key polynomials obtained from $\lambda_\alpha$, namely:

$$x^{\lambda_\alpha}, \theta_1 x^{\lambda_\alpha}, \theta_2 x^{\lambda_\alpha}, \theta_{21} x^{\lambda_\alpha}, \theta_{12} x^{\lambda_\alpha}, \theta_{121} x^{\lambda_\alpha}.$$

We now plot each of these 6 Demazure atoms in the Coxeter arrangement:

Case 1. $x^{\lambda_\alpha} = x_1^m x_2^n$, so it corresponds to the point $(m, n, 0)$ in $R_1$:

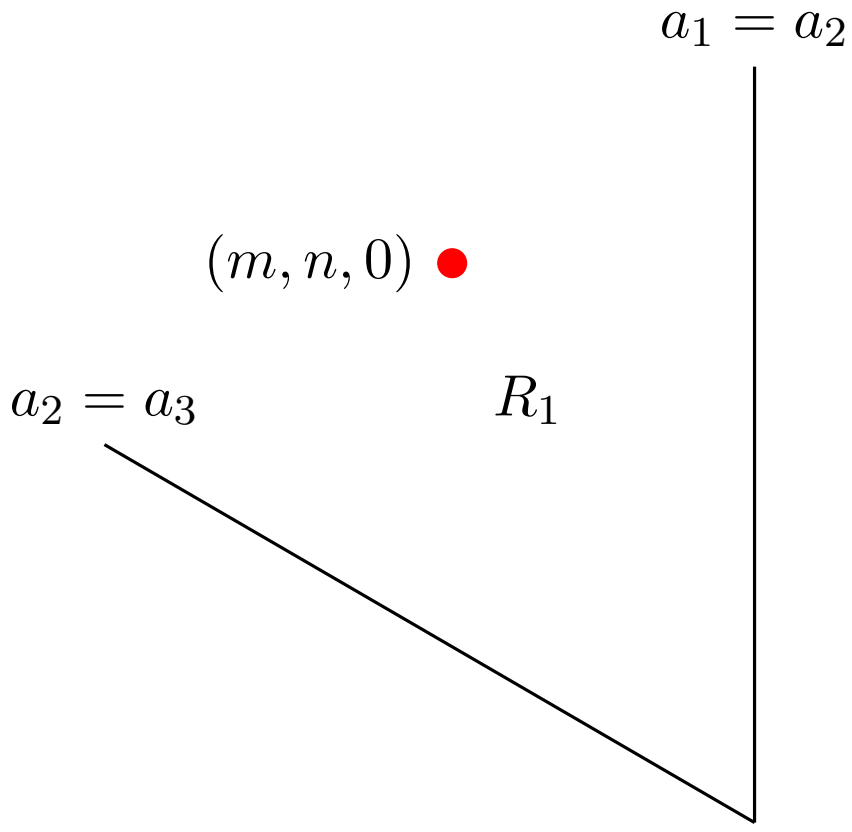

Case 2. $\theta_1 x^{\lambda_\alpha} = \theta_1(x_1^m x_2^n) = x_1^{m-1} x_2^{n+1} + x_1^{m-2} x_2^{n+2} + \cdots + x_1^n x_2^m$

Hence each monomial corresponds to a lattice point except $(m, n, 0)$ (and vice versa) on the line obtained by joining $(m, n, 0)$ and its reflection along $a_1 = a_2$.

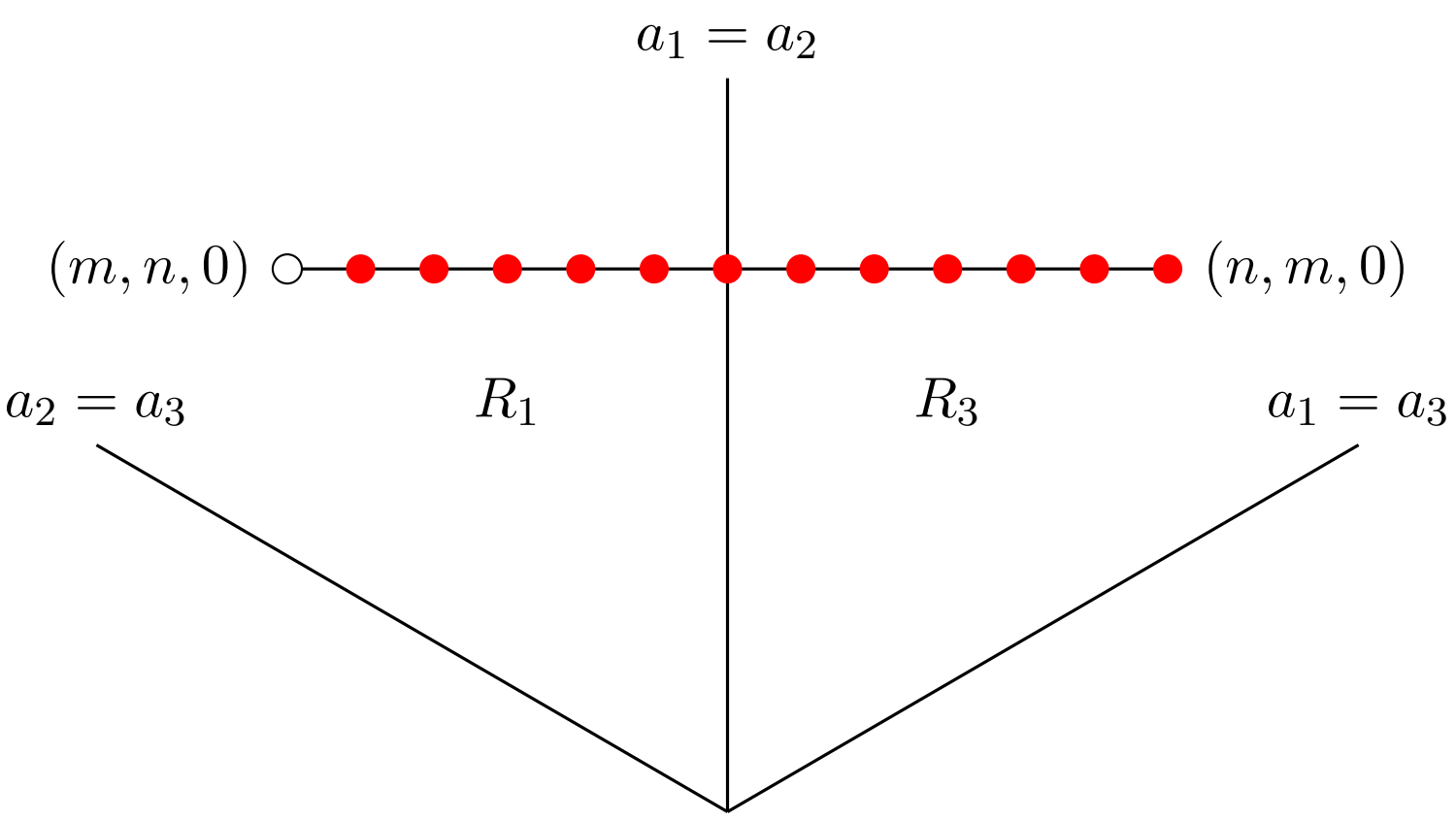

This also shows that $\pi_1 = \theta_1 + 1$

Case 3. $\theta_2 x^{\lambda_\alpha} = \theta_2(x_1^m x_2^n) = x_1^m x_2^{n-1} x_3 + x_1^m x_2^{n-2} x_3^2 + \cdots + x_1^m x_3^n$

Therefore it corresponds to the line joining $(m, n, 0)$ and $(m, 0, n)$ excluding $(m, n, 0)$. That is, each monomial corresponds to a lattice point except $(m, n, 0)$ (and vice versa) on the line obtained by joining $(m, n, 0)$ and its reflection along $a_2 = a_3$.

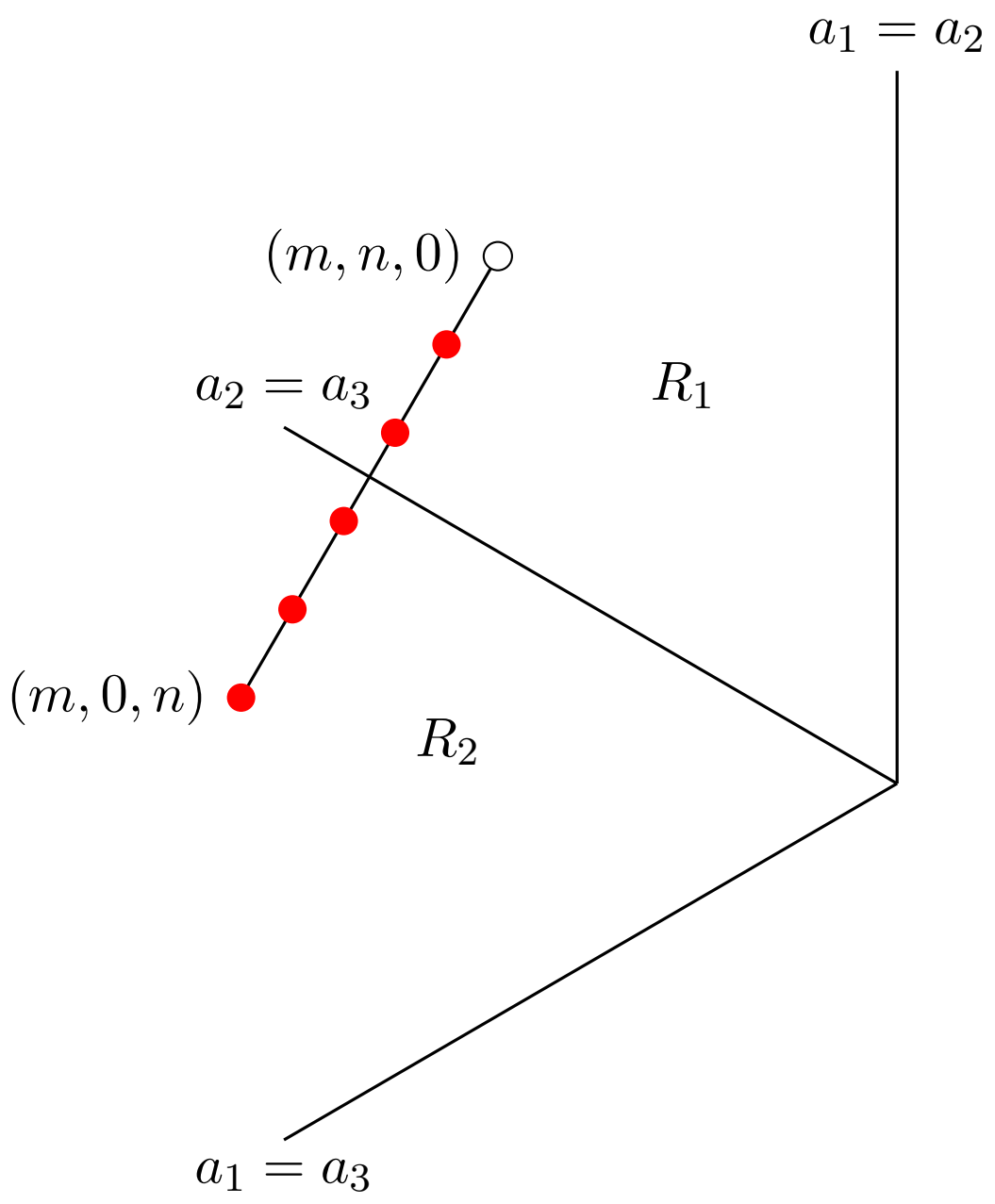

This also shows $\pi_2 = \theta_2 + 1$

Case 4. $\theta_{21} x^{\lambda_\alpha} = \theta_2 \theta_1 (x_1^m x_2^n) = \theta_2(x_1^{m-1} x_2^{n+1} + x_1^{m-2} x_2^{n+2} + \cdots + x_1^n x_2^m)$

Since $\theta_2(x_1^{m-1} x_2^{n+1} + x_1^{m-2} x_2^{n+2} + \cdots + x_1^n x_2^m) = \theta_2(x_1^{m-1} x_2^{n+1}) + \theta_2(x_1^{m-2} x_2^{n+2}) + \cdots + \theta_2(x_1^n x_2^m)$, we can apply the same correspondence in Case 3. for each atom in the summand. Therefore $\theta_{21} x^{\lambda_\alpha}$ corresponds to the $m - n + 1$ lines obtained by reflecting each lattice point except $(m, n, 0)$ on the line joining $(m, n, 0)$ and $(m, 0, n)$ along $a_2 = a_3$. That is, each monomial corresponds to a lattice point (and vice versa) in the 'semi-open' trapezoid obtained by first reflecting $(m, n, 0)$ along the line $a_1 = a_2$ followed by reflecting the resulting line along $a_2 = a_3$ as shown below:

First reflect along $a_1 = a_2$ and get a line joining $(m, n, 0)$ and $(m, 0, n)$:

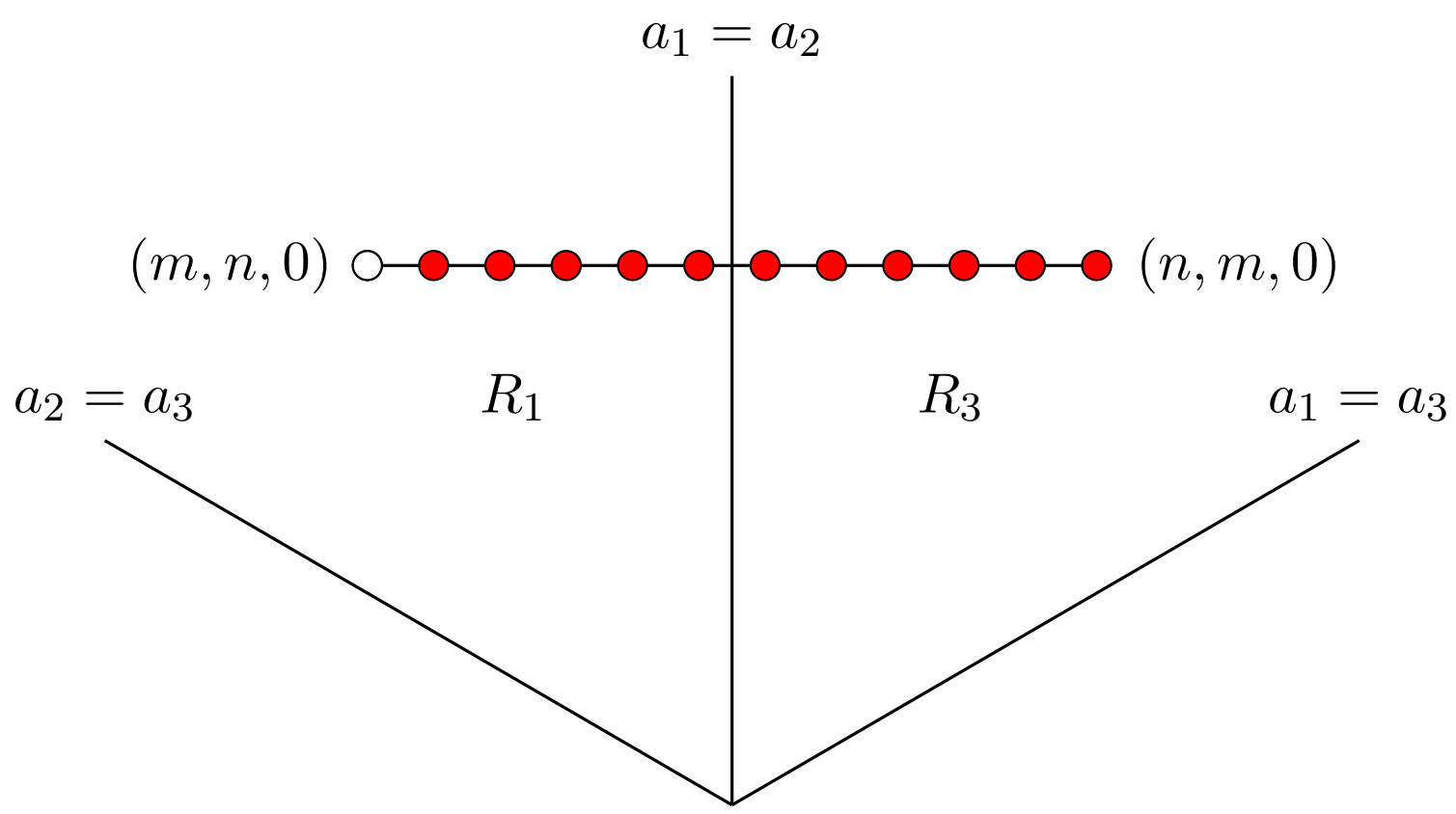

Then reflect the line along $a_2 = a_3$:

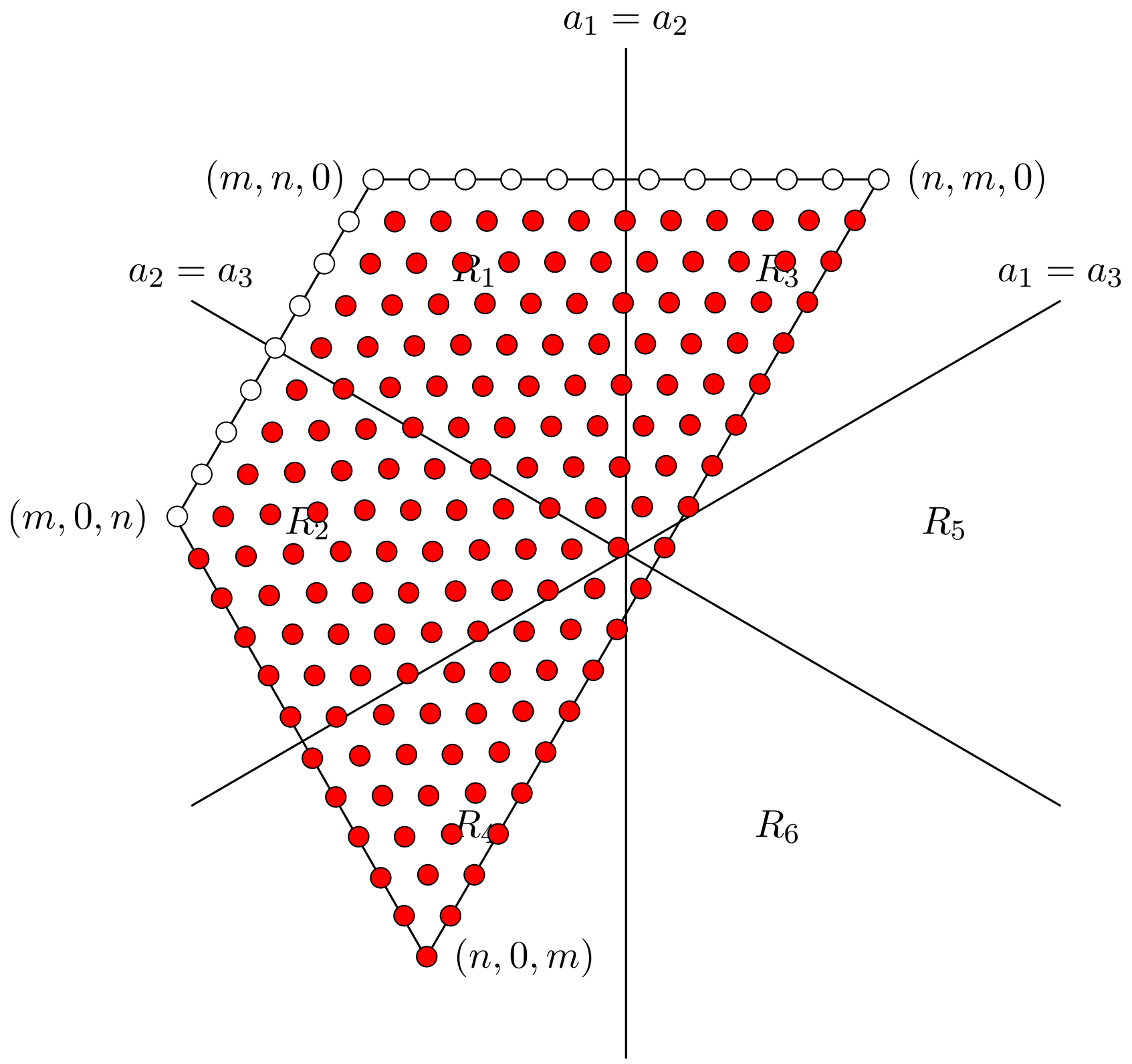

Note that there is another case where the trapezoid does not have any point in $R_5$ or $R_6$:

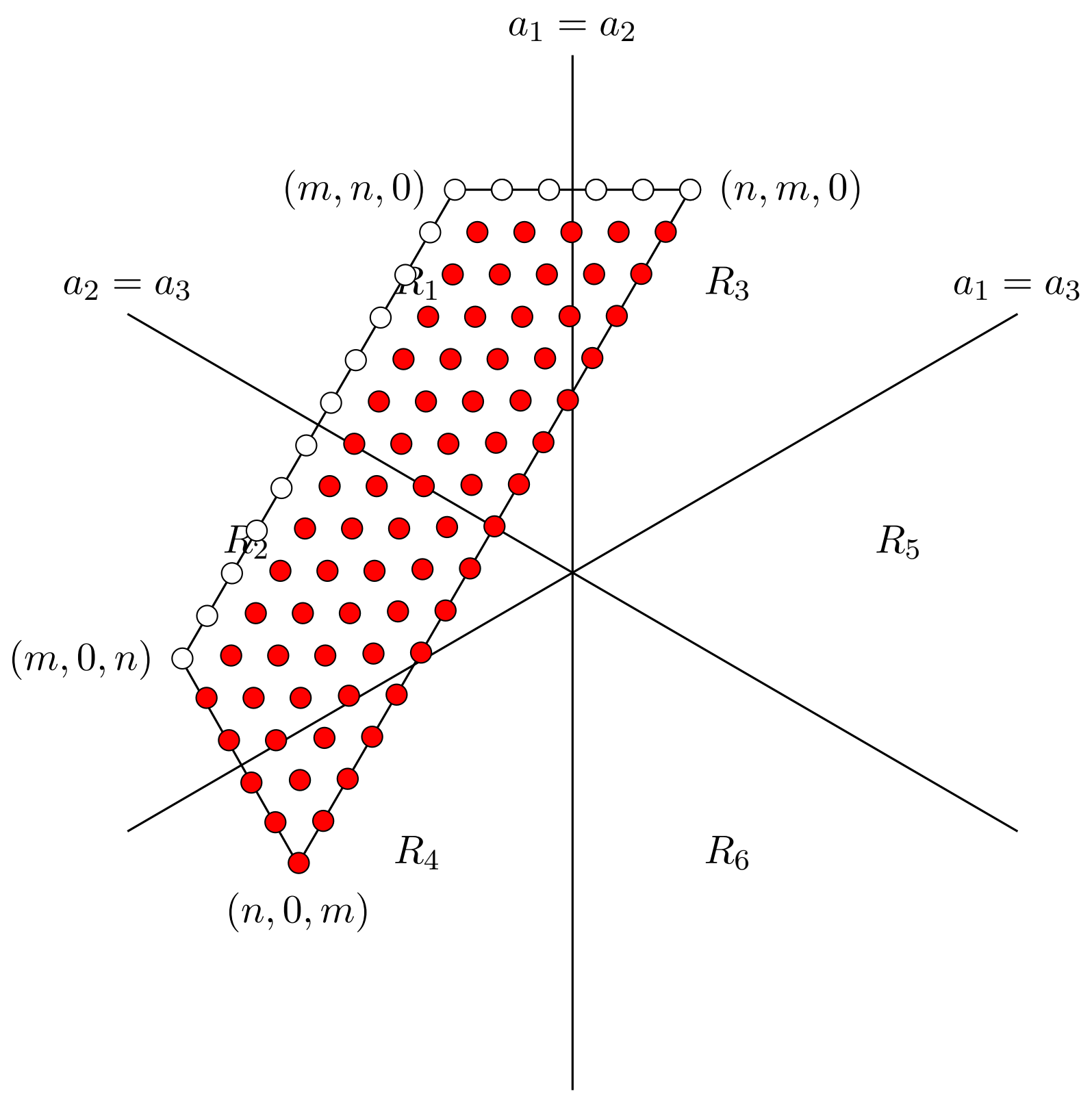

As Case 4 in Section 6.1, the trapezoid has at least a point in $R_5$ or $R_6$ if and only if $m \geq 2n$ (one can prove that by locating the mid-point of the line joining $(n, m, 0)$ and $(n, 0, m)$.)

Again, we can illustrate $\pi_{21} = 1 + \theta_1 + \theta_2 + \theta_{21}$:

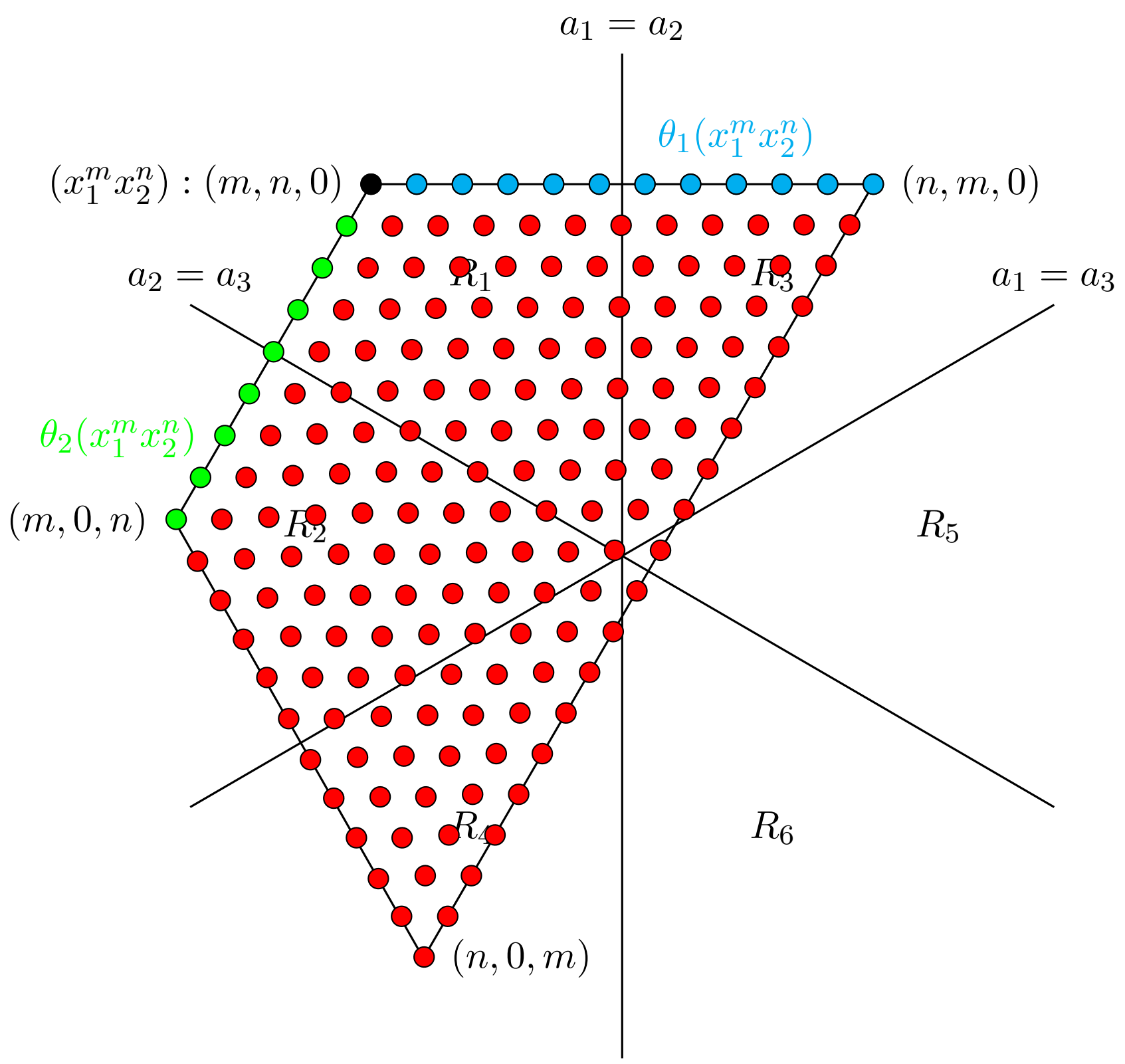

Case 5. $\theta_{12}x^{\lambda_\alpha} = \theta_1\theta_2(x_1^m x_2^n) = \theta_1(x_1^m x_2^{n-1}x_3 + x_1^m x_2^{n-2}x_3^2 + \cdots + x_1^m x_3^n)$

Similar to Case 4., it corresponds to the lattice points in the semi-open trapezoid formed by first reflecting $(m, n, 0)$ along $a_2 = a_3$ followed by reflecting the resulting line along $a_1 = a_2$. The trapezoid has at least a point in $R_4$ or $R_6$ if and only if $2n \geq m$. The trapezoid is as follows (we just shade the region using dotted pattern for convenience, but the actual correspondence should be lattice points in the shaded region including the solid boundaries.):

For $2n < m$:

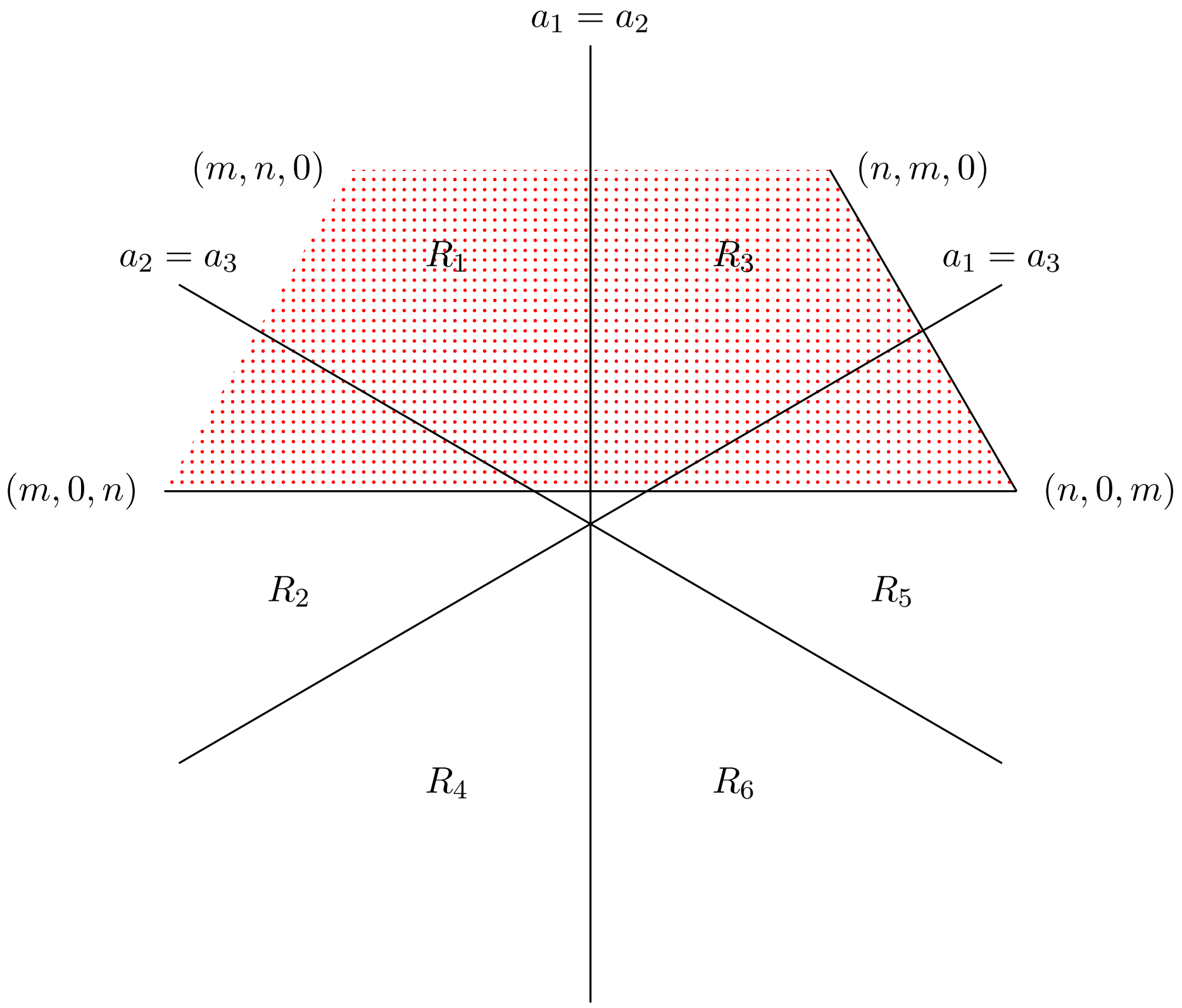

For $2n \geq m$:

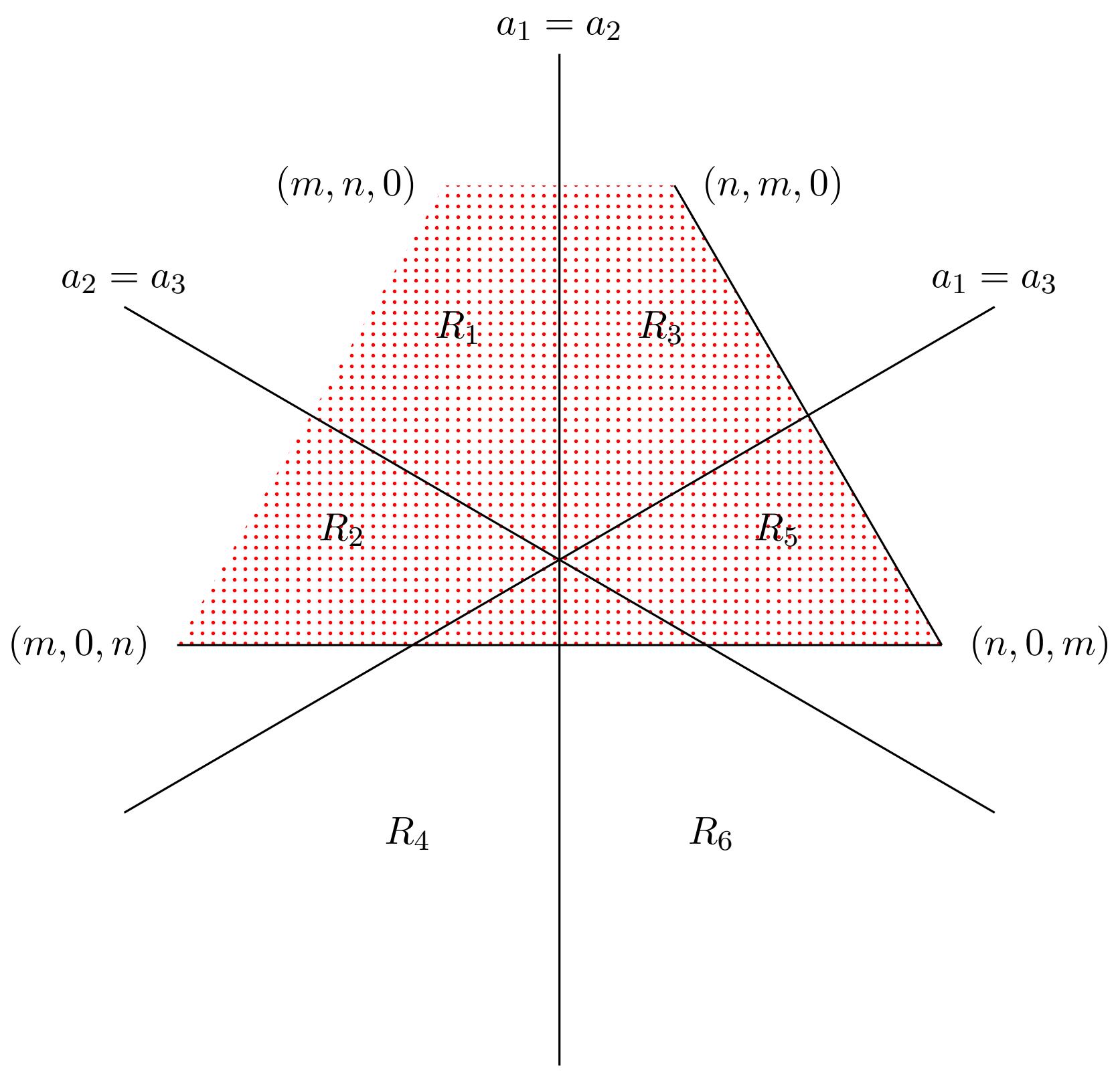

We again illustrate how to decompose $\pi_{12}$ into $1 + \theta_1 + \theta_2 + \theta_{12}$.

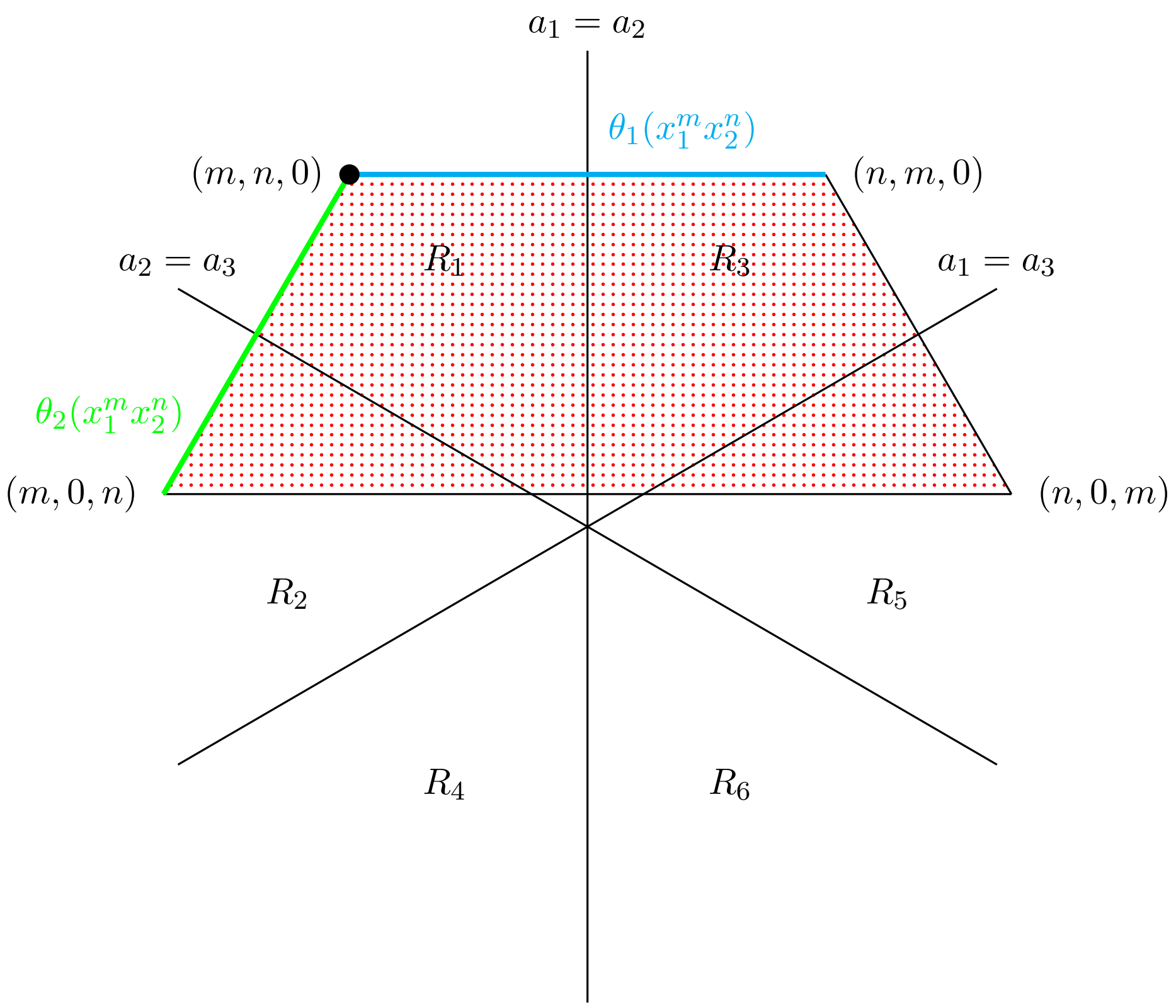

Case 6. $\theta_{121} x^{\lambda_\alpha} = \theta_1\theta_2\theta_1(x_1^m x_2^n) = \theta_1(\theta_2\theta_1(x_1^m x_2^n))$

First recall that by Proposition 3.1 $\theta_i s_i = -\theta_i$, hence $\theta_i s_i + \theta_i = 0$.

For example, when $i = 1$, $\theta_1(x_1^a x_2^b) + \theta_i(x_1^b x_2^a) = 0$ for any and any integers $a \geq b \geq 0$.

We can plot $\theta_1(x_1^b x_2^a)$ as:

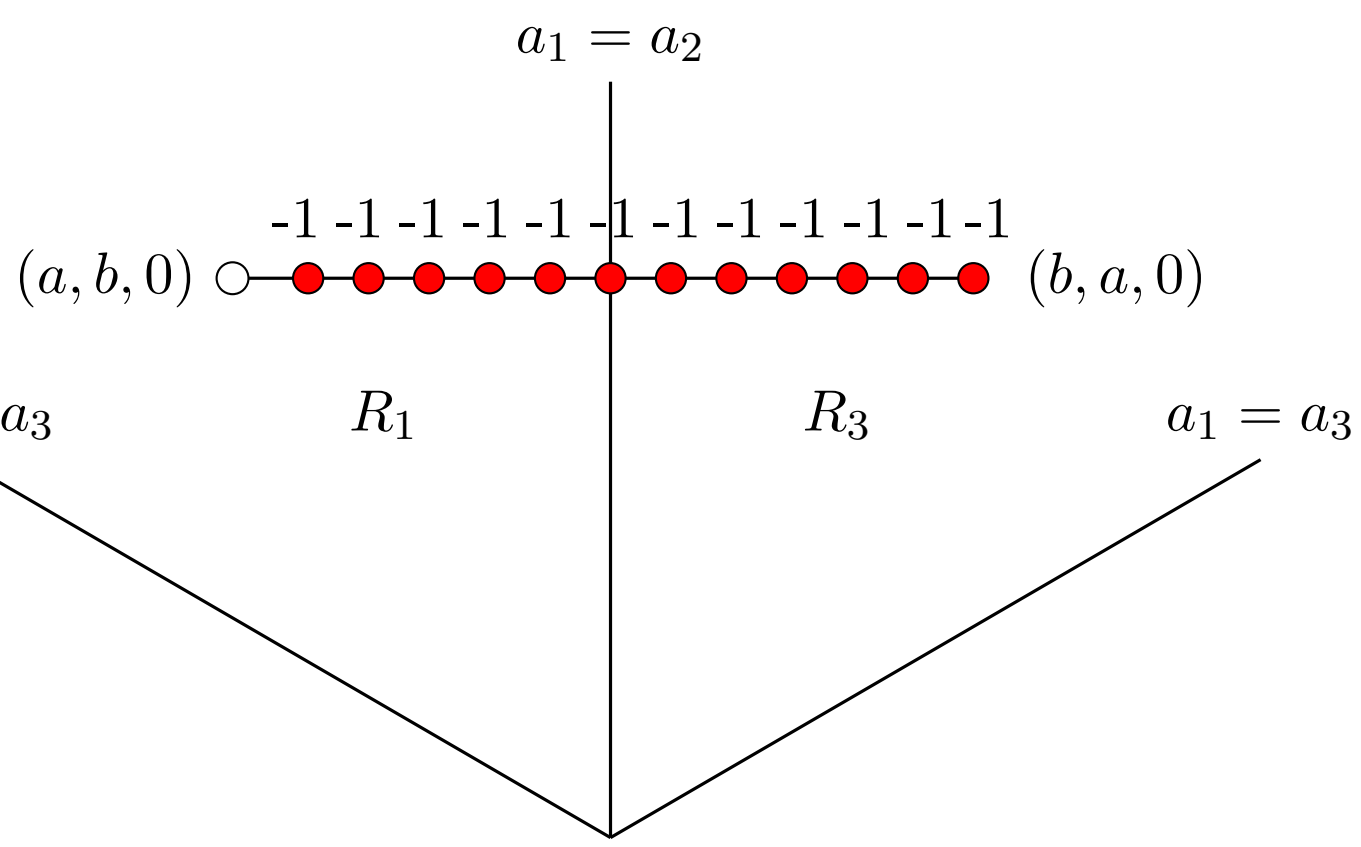

Now consider $\theta_{121}(x_1^m x_2^n) = \theta_1(\theta_{21}(x_1^m x_2^n))$.

By Case 4., we have $\theta_{21}(x_1^m x_2^n)$ as a trapezoid as follows:

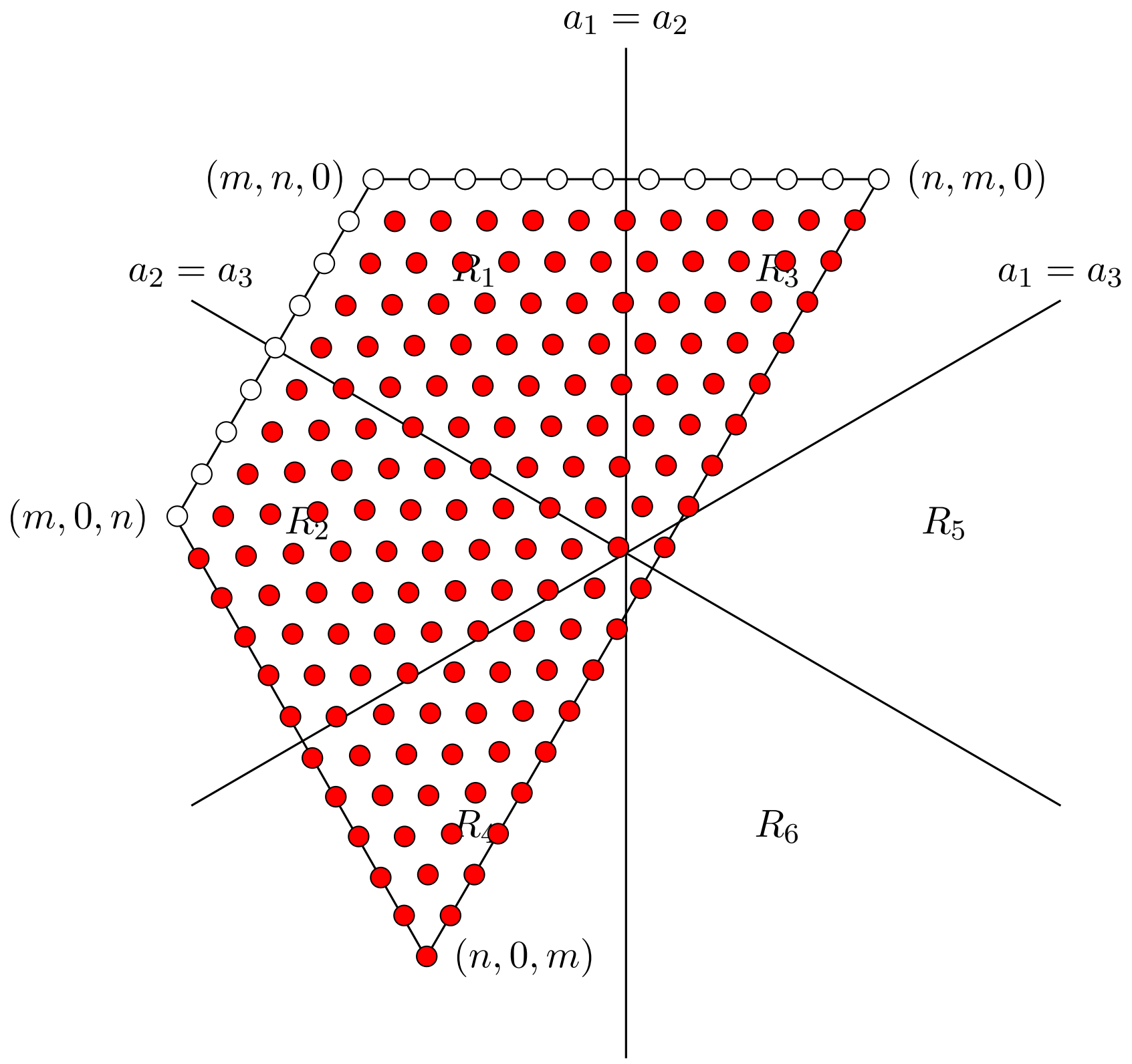

Apply $\theta_1$ to each lattice points in the semi-open trapezoid is equivalent to reflecting the trapezoid along $a_1 = a_2$ and get a hexagon with multiplicities on the points. Recall that each "'$a_1 = a_2'$ - reflection pair in the trapezoid region (i.e. $(a.b, c)$ and $(b, a, c)$) vanishes under $\theta_1$. Hence the multiplicity stays constant along the horizontal line perpendicular to the line $a_1 = a_2$. Here is an example for $m \geq 2n$:

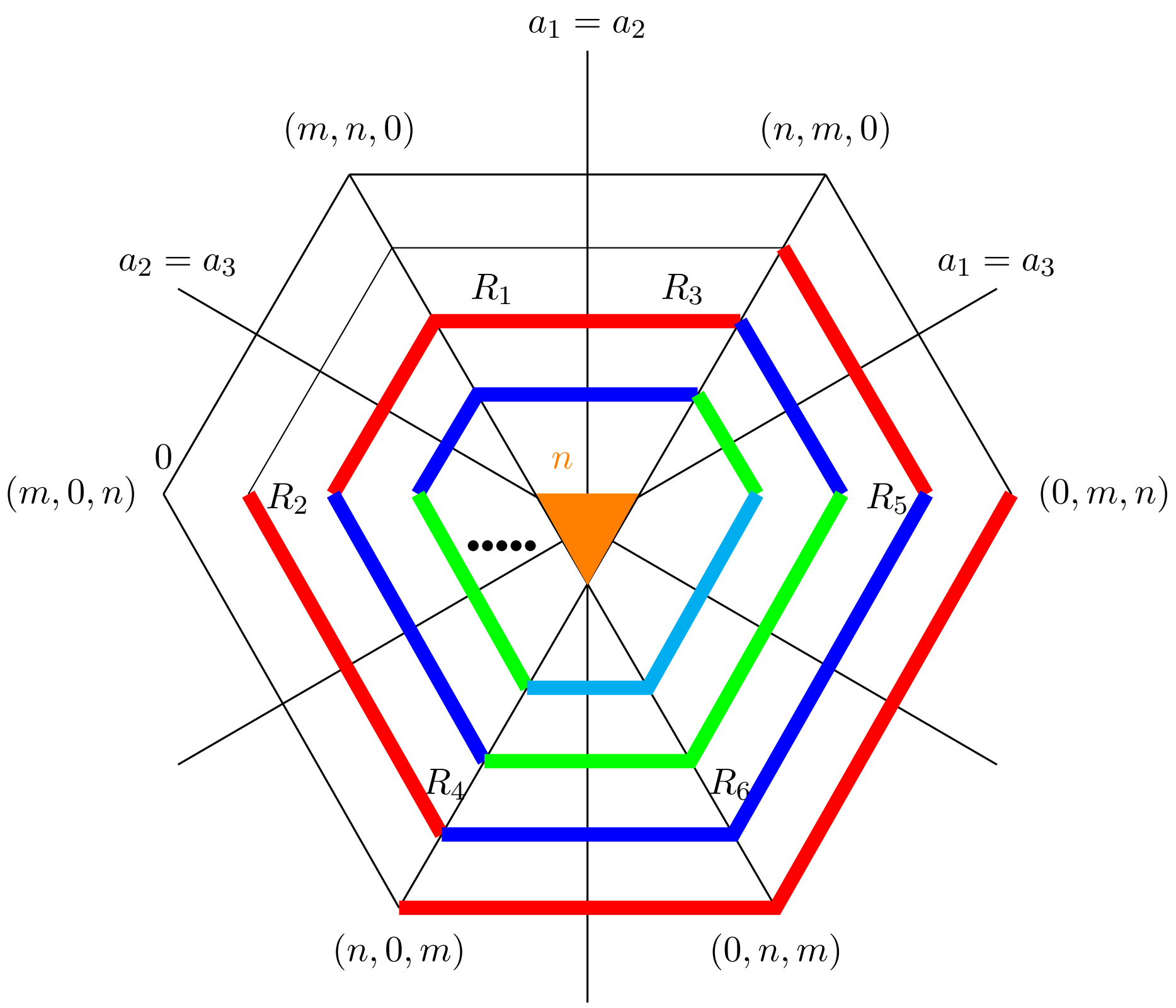

The monomials corresponding to the lattice points on the red boundary have coefficient 1 in the key polynomial, and those monomial corresponding to the lattice points on the blue boundary have coefficient 2 and so on, while those corresponding to the lattice points in the inner most triangle (the orange region) have coefficient $n$ for $m \geq 2n$ (and $m - n$ for $m \leq 2n$).

If we also plot the multiplicity, with $xy$ -plane being the Coxeter arrangement and $z$-coordinates being the multiplicity, we get a polytope similar to the one in Case 6 in Section 6.1 but with different heights (as not all the multiplicity of the monomials in an atom is the same as those in keys).

One can verify that $\theta_{121} = \theta_{212}$ by starting with $\theta_{12}$ instead. Also one can get the decomposition $\pi_{121} = 1 + \theta_1 + \theta_2 + \theta_{12} + \theta_{21} + \theta_{121}$ by putting the figures shown in this section together.

# 7 Products of two Demazure characters

Let $m \geq n \geq 0$ and $k \geq l \geq 0$ be integers.

| | $x_1^k x_2^l$ | $\pi_1(x_1^k x_2^l)$ | $\pi_2(x_1^k x_2^l)$ | $\pi_{21}(x_1^k x_2^l)$ | $\pi_{12}(x_1^k x_2^l)$ | $\pi_{121}(x_1^k x_2^l)$ |
|---|---|---|---|---|---|---|
| $x_1^m x_2^n$ | 4.13 | 4.13 | 4.13 | 4.13 | 4.13 | 4.13/4.18 |
| $\pi_1(x_1^m x_2^n)$ | | (i) | 7.1 | 7.2 | (ii) | 4.18 |
| $\pi_2(x_1^m x_2^n)$ | | | (iii) | (iv) | 7.3 | 4.18 |
| $\pi_{21}(x_1^m x_2^n)$ | | | | (v) | 7.4 | 4.18 |
| $\pi_{12}(x_1^m x_2^n)$ | | | | | (vi) | 4.18 |
| $\pi_{121}(x_1^m x_2^n)$ | | | | | | 4.18 |

Table 2: Decomposition of products of keys into atoms

In the following sections, we will first state the result and one can verify by expanding both sides directly. We will state some other methods (either using operators or polytopes) to verify or interpret the decomposition in the first two sections. These methods are applicable to all cases.

Also, we denote the indicator function as $\mathbb{1}_S = \begin{cases} 1 & \text{if } ((m,n),(k,l)) \in S \\ 0 & \text{otherwise.} \end{cases}$

## 7.1 $\pi_1(x_1^m x_2^n) \times \pi_2(x_1^k x_2^l)$

$$\begin{aligned} \pi_1(x_1^m x_2^n) \times \pi_2(x_1^k x_2^l) = & \sum_{s=0}^{\min\{m-n,k\}} \sum_{t=\max\{0,s-(k-l)\}}^{\min\{l,s+n\}} x_1^{m+k-s} x_2^{n+l+s-t} x_3^t \\ & + \mathbb{1}_{\{m-n>k-l\}} \sum_{t=0}^{\min\{l,(m-n)-(k-l)\}} \theta_1(x_1^{m+l-t} x_2^{k+n} x_3^t) \\ & + \mathbb{1}_{\{l>n\}} \sum_{s=0}^{\min\{l-n,m-n\}} \theta_2(x_1^{m+k-s} x_2^l x_3^{n+s}) \end{aligned}$$

We can also use polytopes to verify the atom positivity of $\pi_1(x_1^m x_2^n) \times \pi_2(x_1^k x_2^l)$. The polytope corresponding to $\pi_1(a_1^m x_2^n)$ is the line:

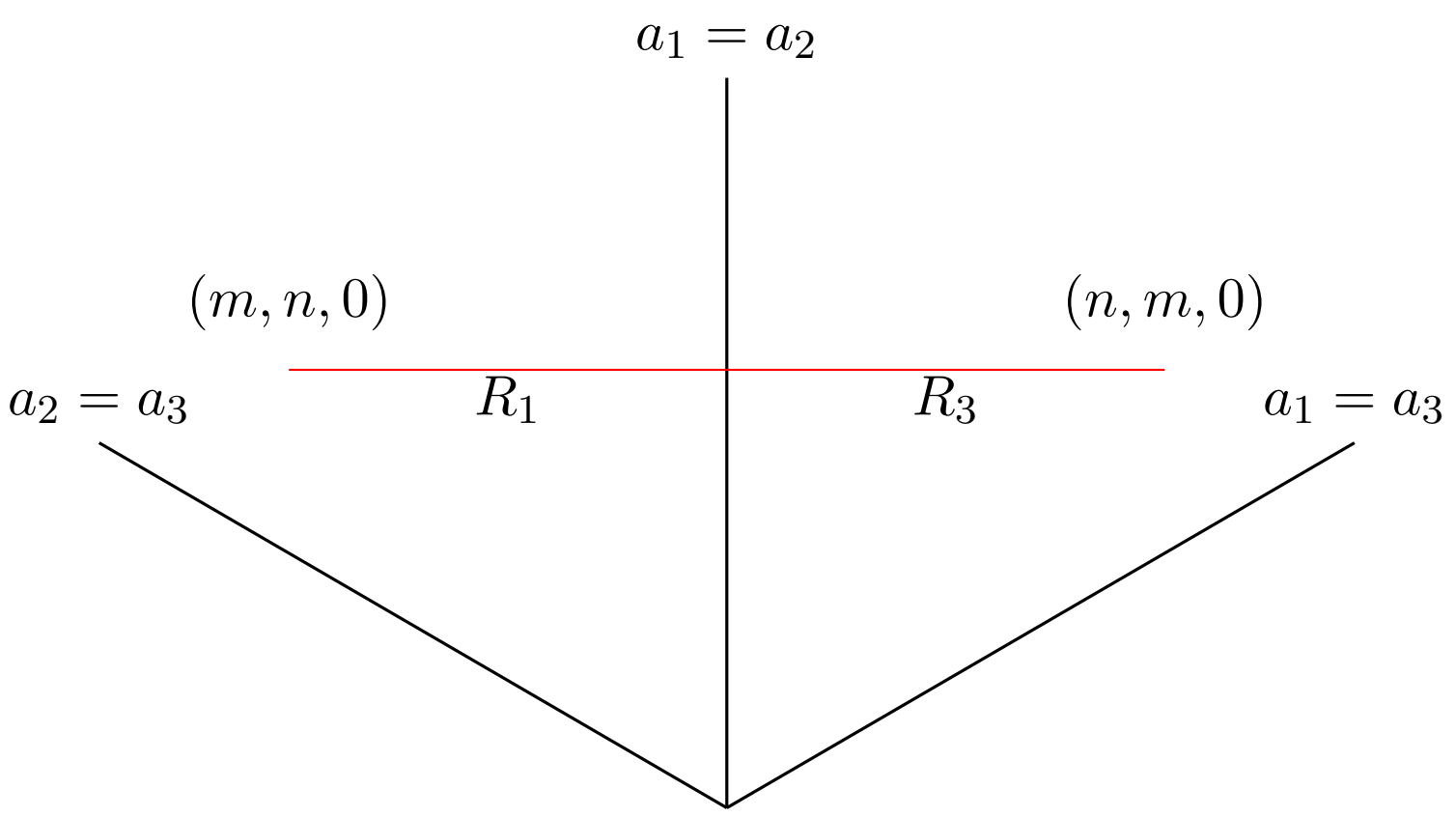

The polytope corresponding to $\pi_2(a_1^k x_2^l)$ is the line:

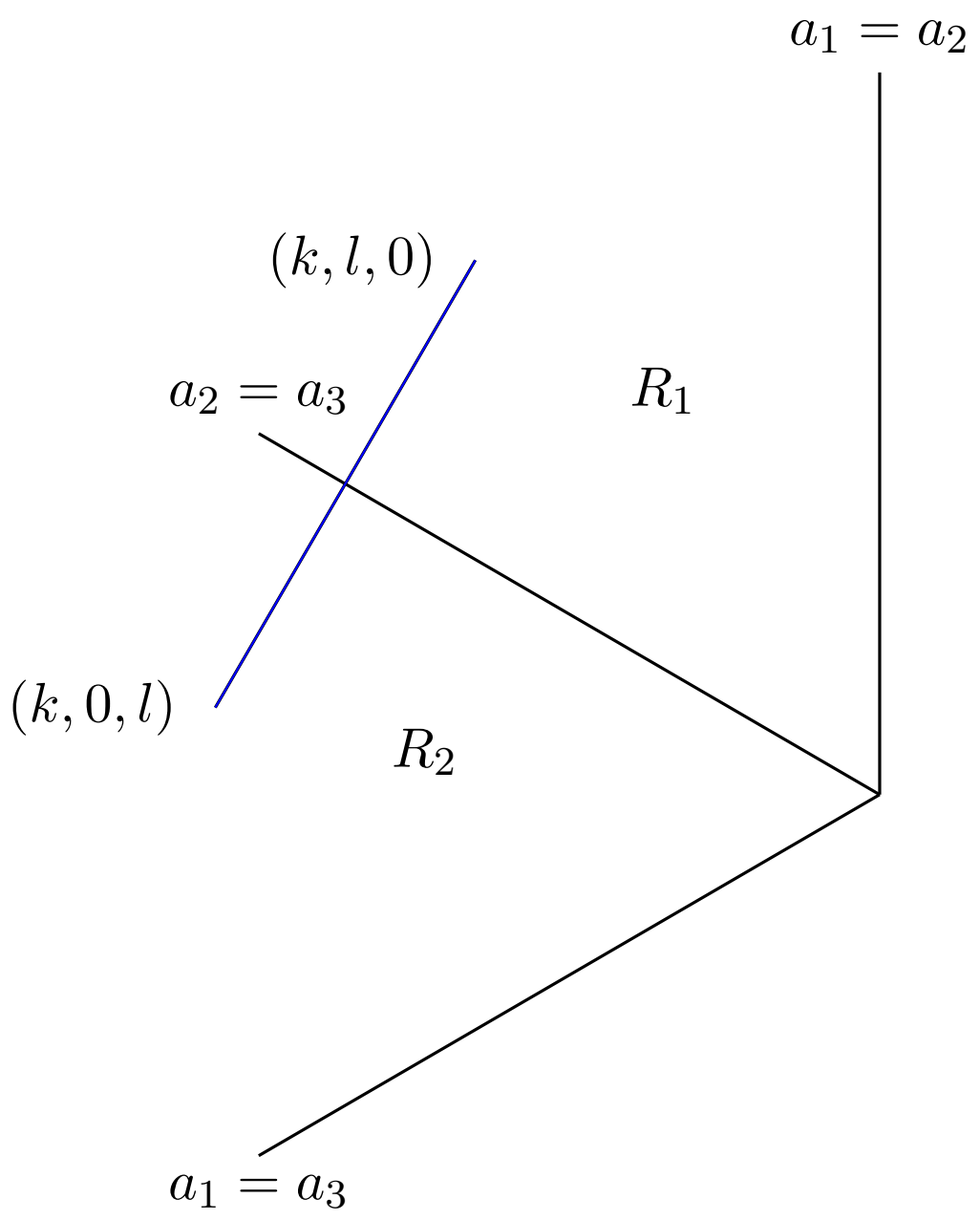

Hence the product is a parallelogram:

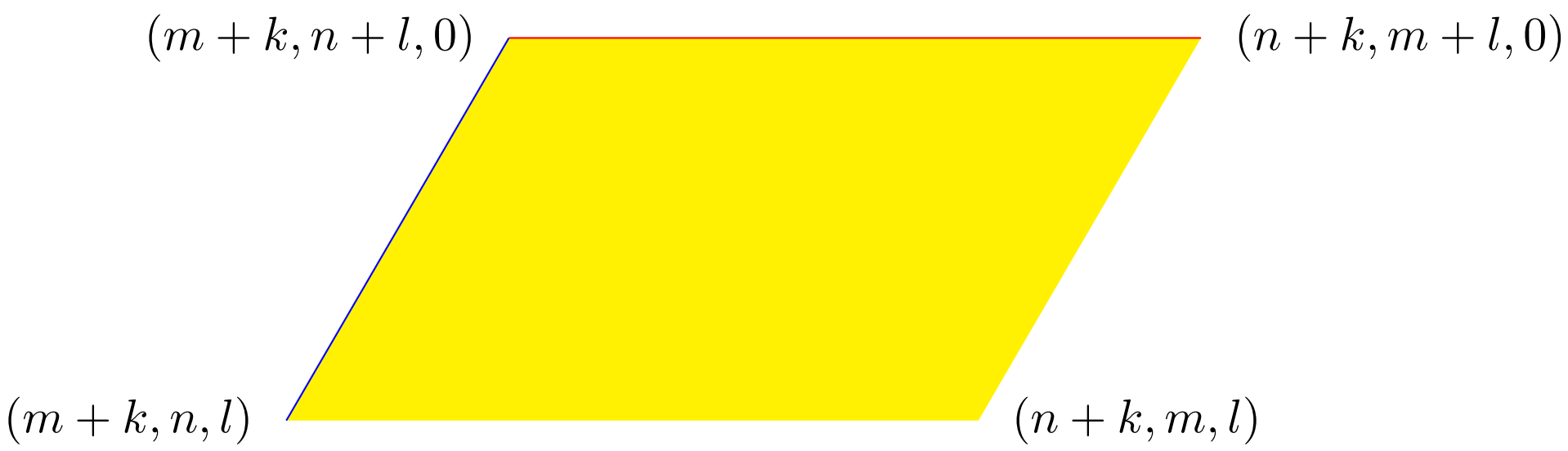

There are different possible positions for the parallelogram:

1. The whole parallelogram lies in $R_1$:

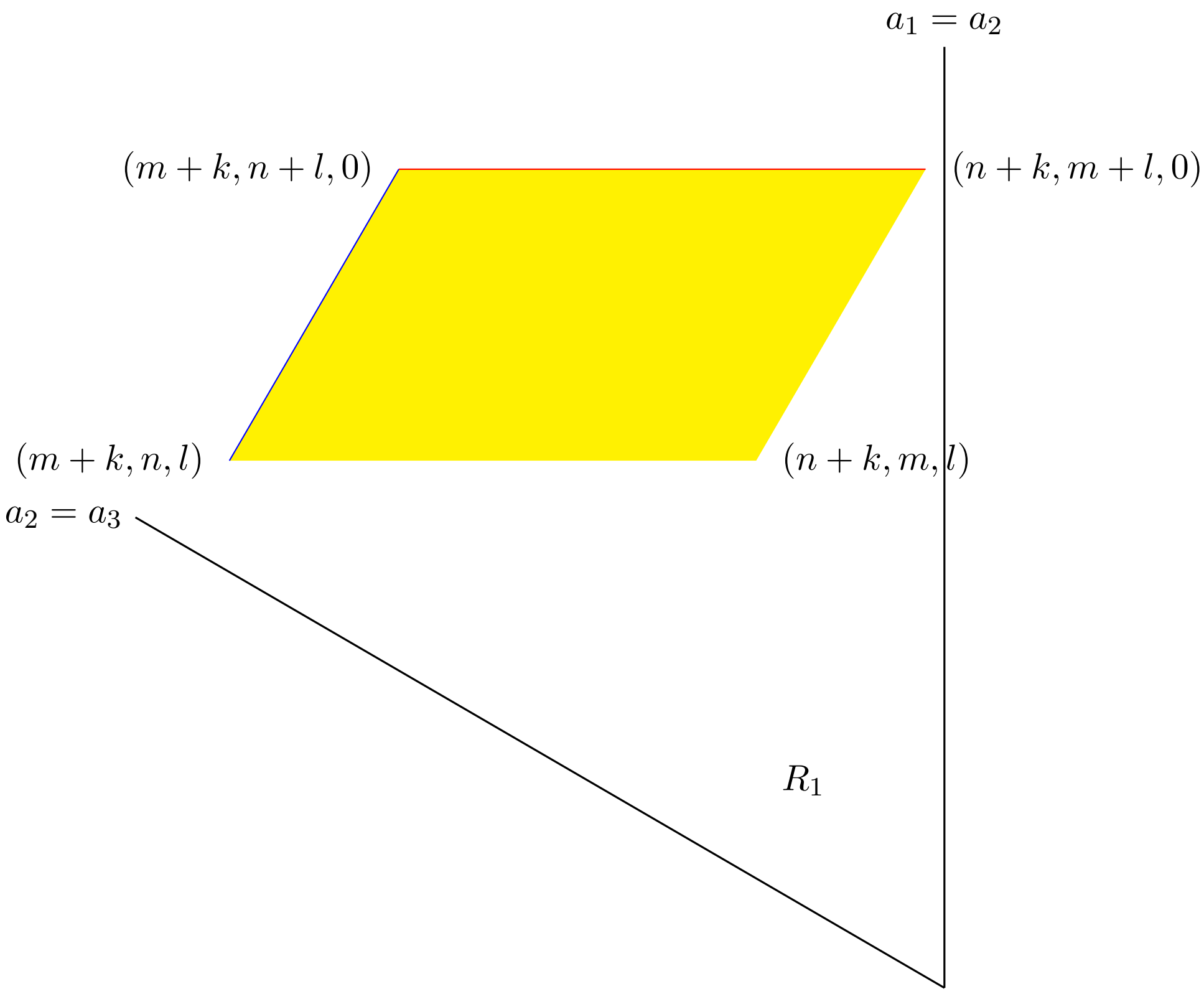

   This case corresponds to $m - n \leq k - l$ and $l \leq n$ in the expansion.

2. The parallelogram lies in $R_1$ and $R_3$:

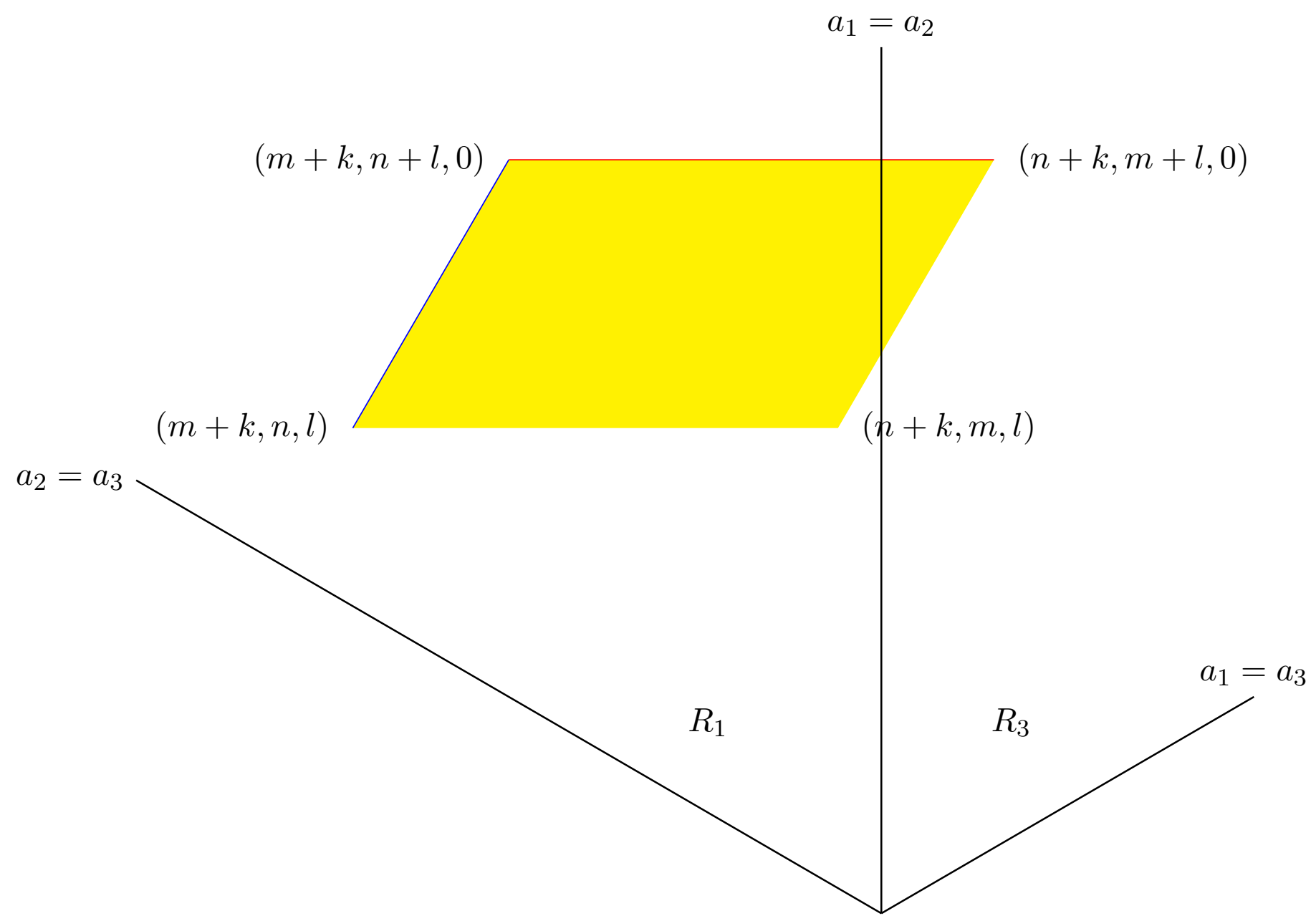

Then we can decompose the parallelogram as:

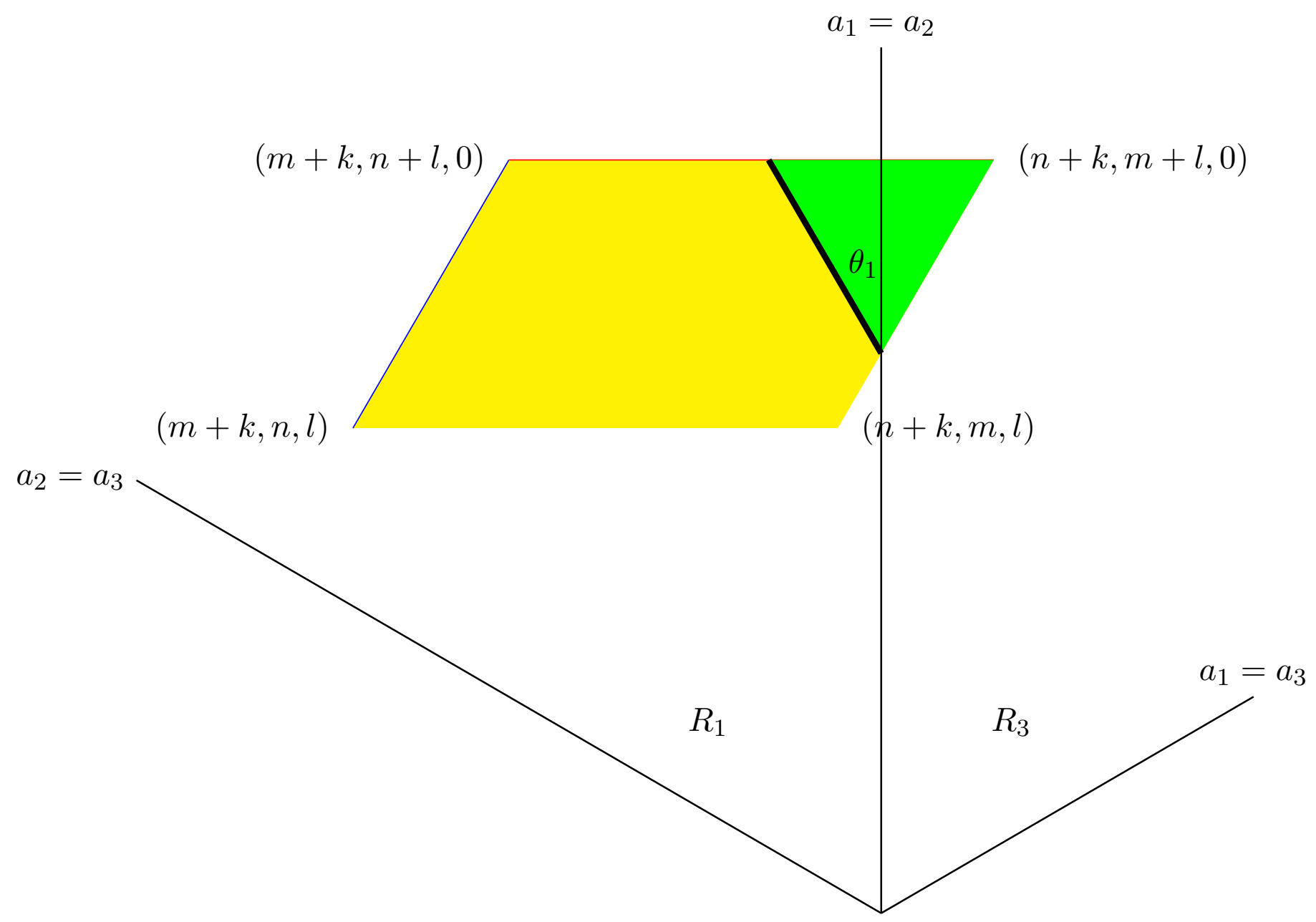

Here the green region is obtained by applying $\theta_1$ on each lattice point on the black line. Notice that the black line is in $R_1$, meaning that every monomial corresponding to the a lattice point is a dominating monomial.

Hence we can decompose the parallelogram into the yellow region in $R_1$ and $\theta_1$ of the black line, all of which are positive sum of atoms.

This corresponds to the case when $m - n > k - l$ and $l \leq n$. Note that there are two different cases for the positions for the line $a_1 = a_2$, either $(n+k, m, l)$ is on the left of the line or on the right of the line. They correspond to the two cases in $\min\{l, (m-n)-(k-l)\}$ in the upper limit of the summation. In fact, one can locate the black line simply by flipping along $a_1 = a_2$ the boundary of the parallelogram which lies in $R_3$.

3. The parallelogram lies in $R_1$ and $R_2$:

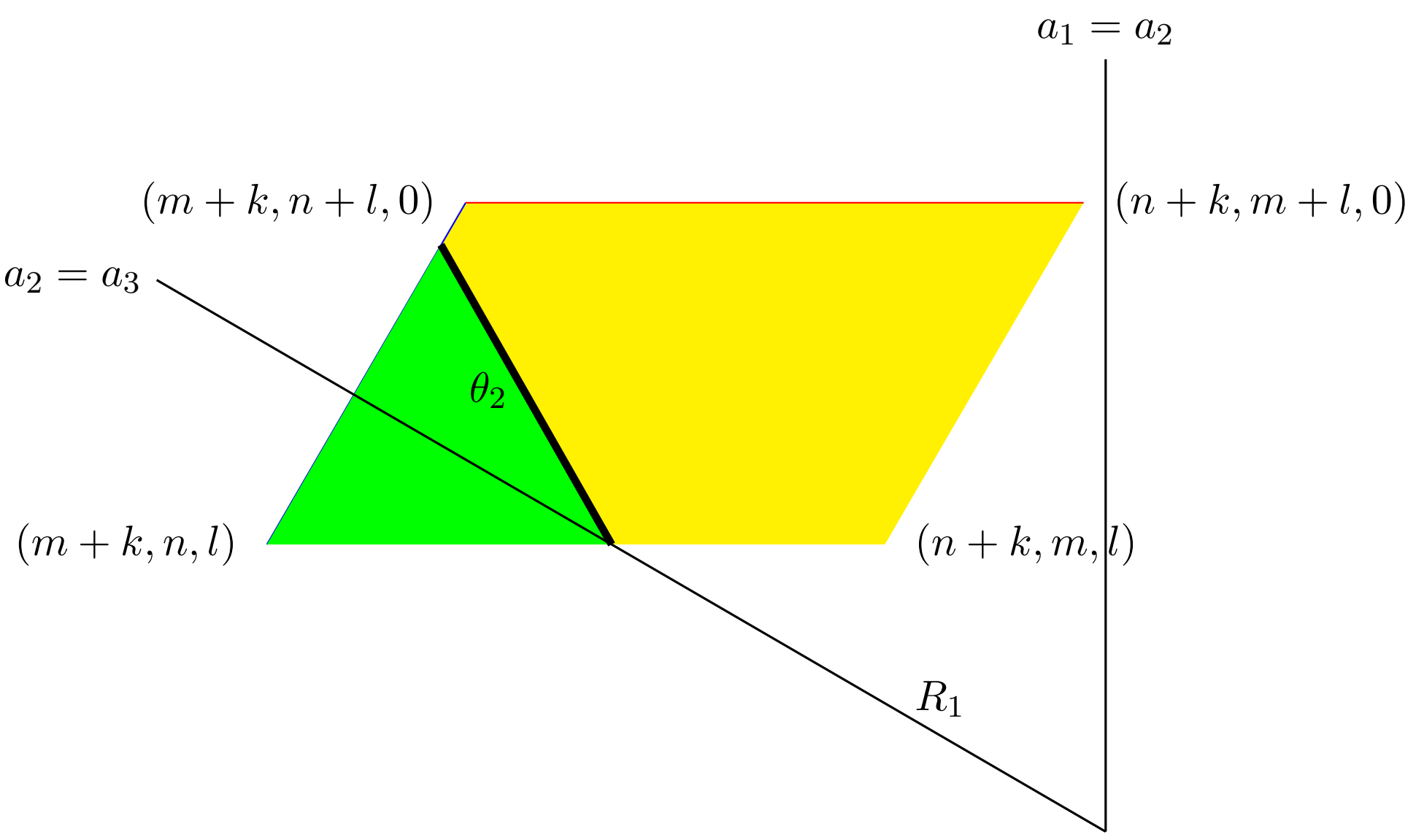

This case corresponds to $m - n \leq k - l$ and $l > n$ in the expansion.

Similar to the previous case when the parallelogram lies in $R_1$ and $R_3$, we can decompose the parallelogram into the yellow region which corresponds to a sum of dominating monomials and the green region which corresponds to $\theta_2$ of a sum of dominating monomials (corresponding to the lattice points on the black line obtained by reflecting along the line $a_2 = a_3$ the 'base' of the parallelogram in $R_2$.

Again there are also two cases: either $(n + k, m, l)$ is above or below the line $a_2 = a_3$, corresponding to the upper limit $\min\{l - n, m - n\}$ of the summation in the third summand in the expansion.

4. The parallelogram lies in $R_1, R_2$ and $R_3$:

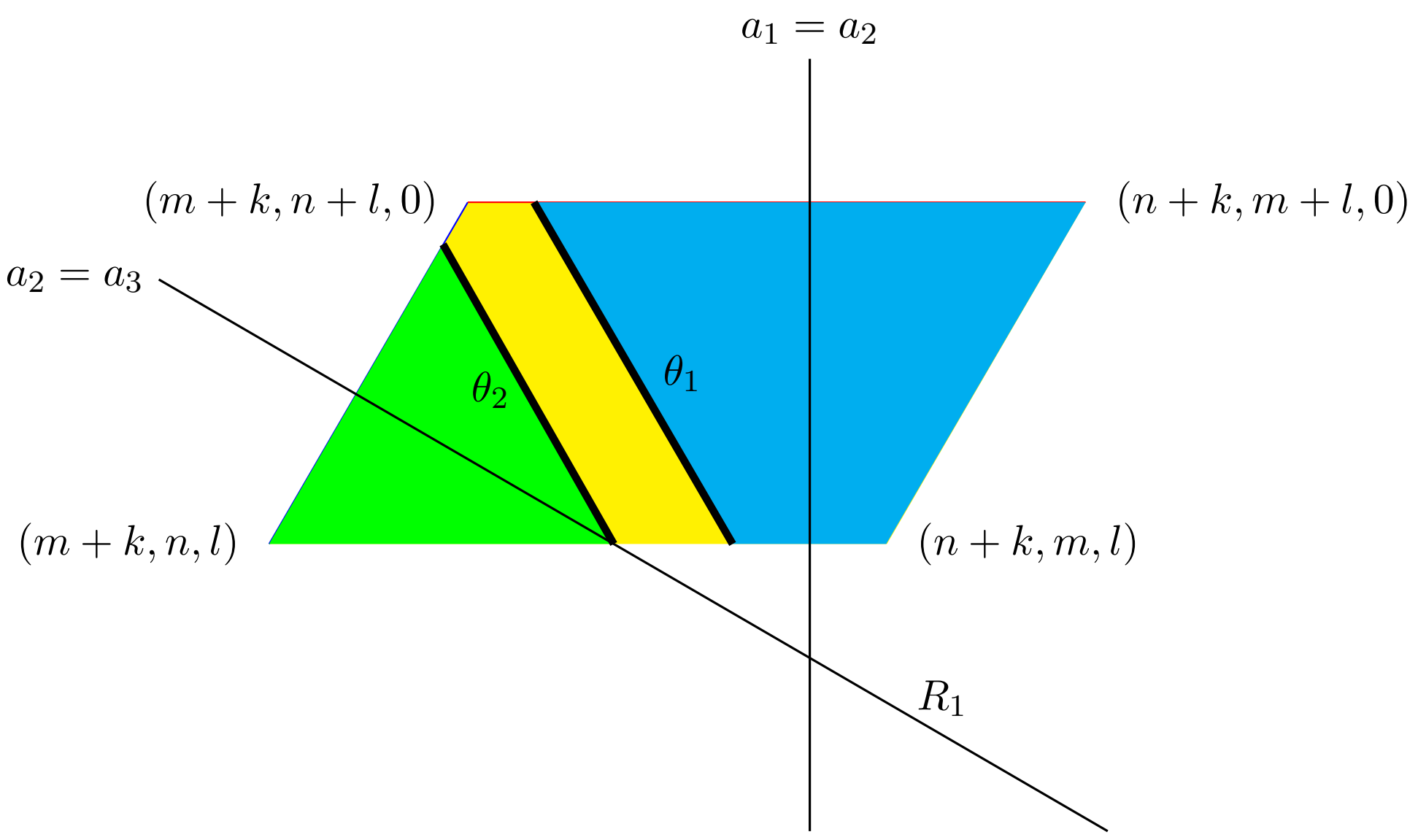

This case corresponds to $m - n > k - l$ and $l > n$ in the expansion and the three regions correspond to the three summands in the expansion.

## 7.2 $\pi_1(x_1^m x_2^n) \times \pi_{21}(x_1^k x_2^l)$

We first write the decomposition of $(x_1^m x_2^n) \times \pi_{121}(x_1^k x_2^l)$ as follows:

$$(x_1^m x_2^n) \times \pi_{121}(x_1^k x_2^l) = A_0(x) + \theta_1 A_1(x) + \theta_2 A_2(x) + \theta_{21} A_{21}(x) + \theta_{12} A_{12}(x) + \theta_{121} A_{121}(x),$$

where $A_I(x) = \sum_{\lambda \in Par} a_\lambda^I x^\lambda$ with $a_\lambda^I \in \mathbb{Z}$ for $I \in \{0, 1, 2, 12, 21, 121\}$.

By Theorem 4.18, $(x_1^m x_2^n) \times \pi_{121}(x_1^k x_2^l)$ is key positive and hence atom positive. i.e. $A_I(x)$ is a sum of dominating monomials with integer coefficients. Then

$$\begin{aligned}\pi_1(x_1^m x_2^n) \times \pi_{21}(x_1^k x_2^l) = \quad & A_0(x) + \sum_{r=0}^{\min\{m,k\}} \sum_{s=\max\{0, r-(n+l)\}}^{\min\{r, m-n, k-l\}} \theta_1(x_1^{m+k-r} x_2^{n+l+s} x_3^{r-s}) \\ & + \theta_2 A_2(x) + \mathbb{1}_{\{\min\{m,k\} \geq n+l\}} \sum_{r=n+l+1}^{\min\{m,k\}} \theta_{12}(x_1^{m+k-r} x_2^r x_3^{n+l}) \\ & + \theta_{21} A_{21}(x) + \mathbb{1}_{\{k > m > n+l\}} \theta_{121}(x_1^k x_2^m x_3^{n+l}).\end{aligned}$$

One can check the coefficients using polytopes as in Section 7.1. We will show a case where $k > m > n + l$ as an example (Note that even $k > m > n + l$ has several subcases). Other cases can be easily deduced similarly.

$\pi_{21}(x_1^k x_2^l)$ corresponds to:

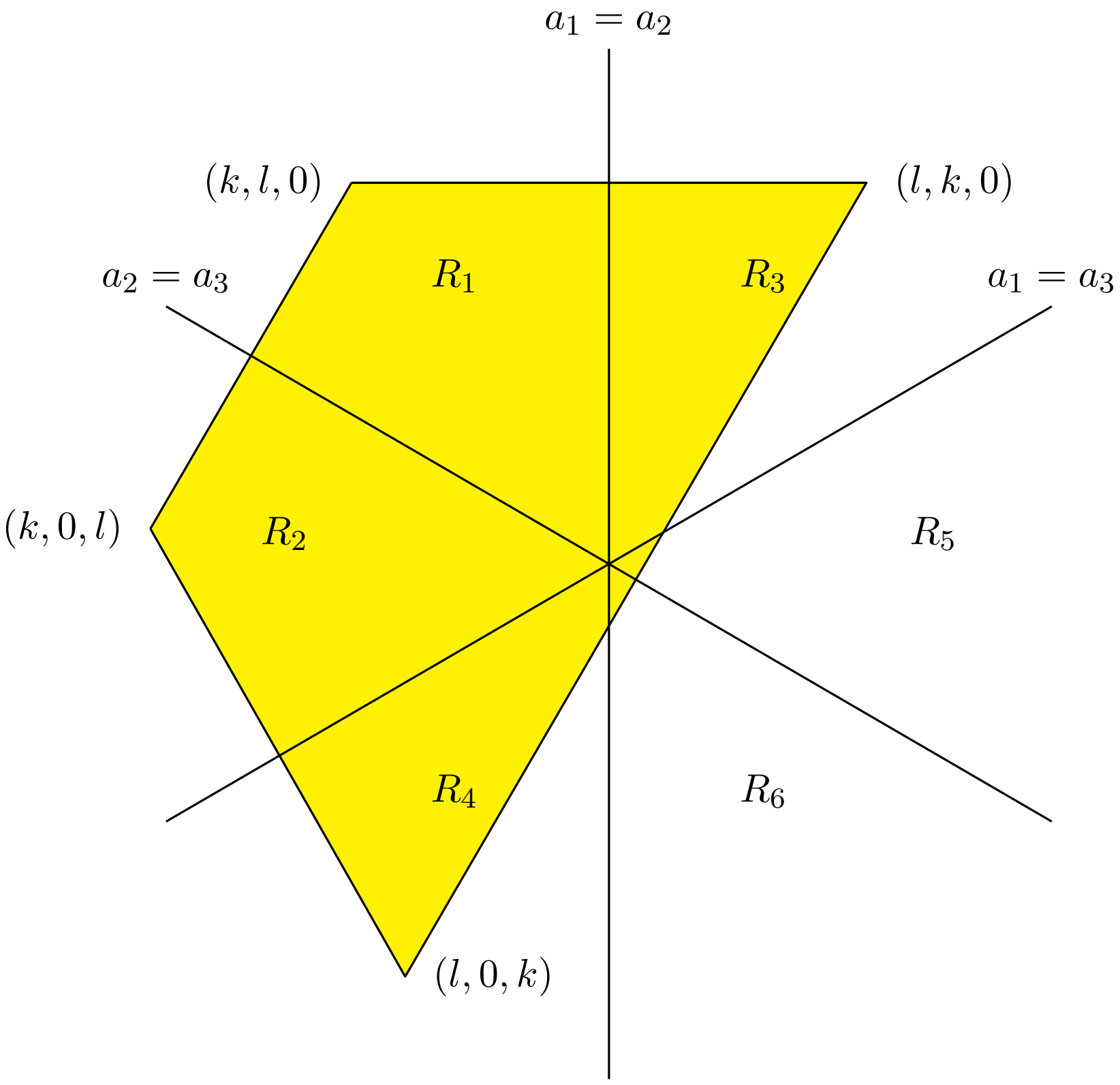

So $\pi_1(x_1^m x_2^n) \times \pi_{21}(x_1^k x_2^l)$ is equivalent to $m-n+1$ trapezoids along the line perpendicular to $a_1 = a_2$ as follows:

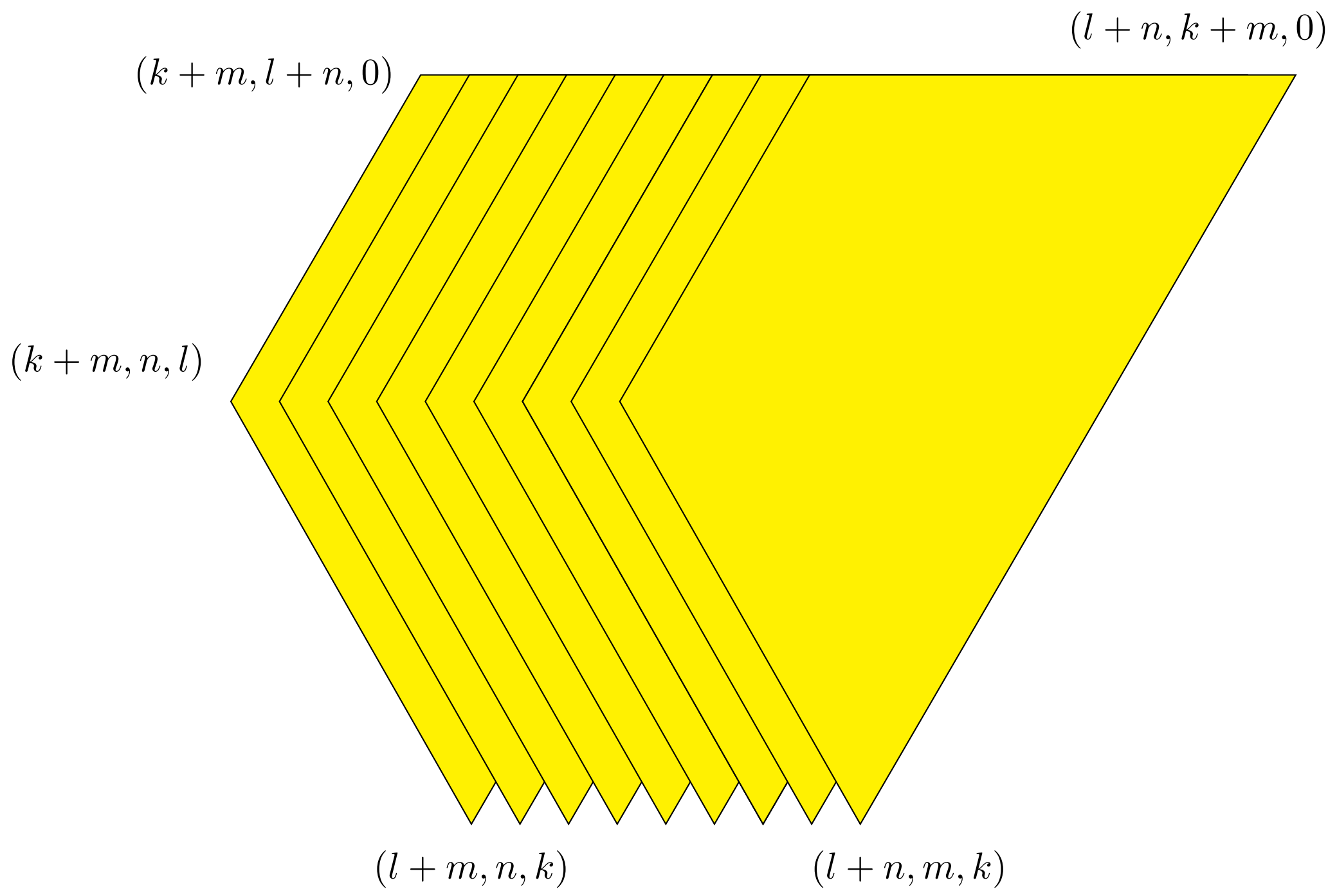

We can draw the product as:

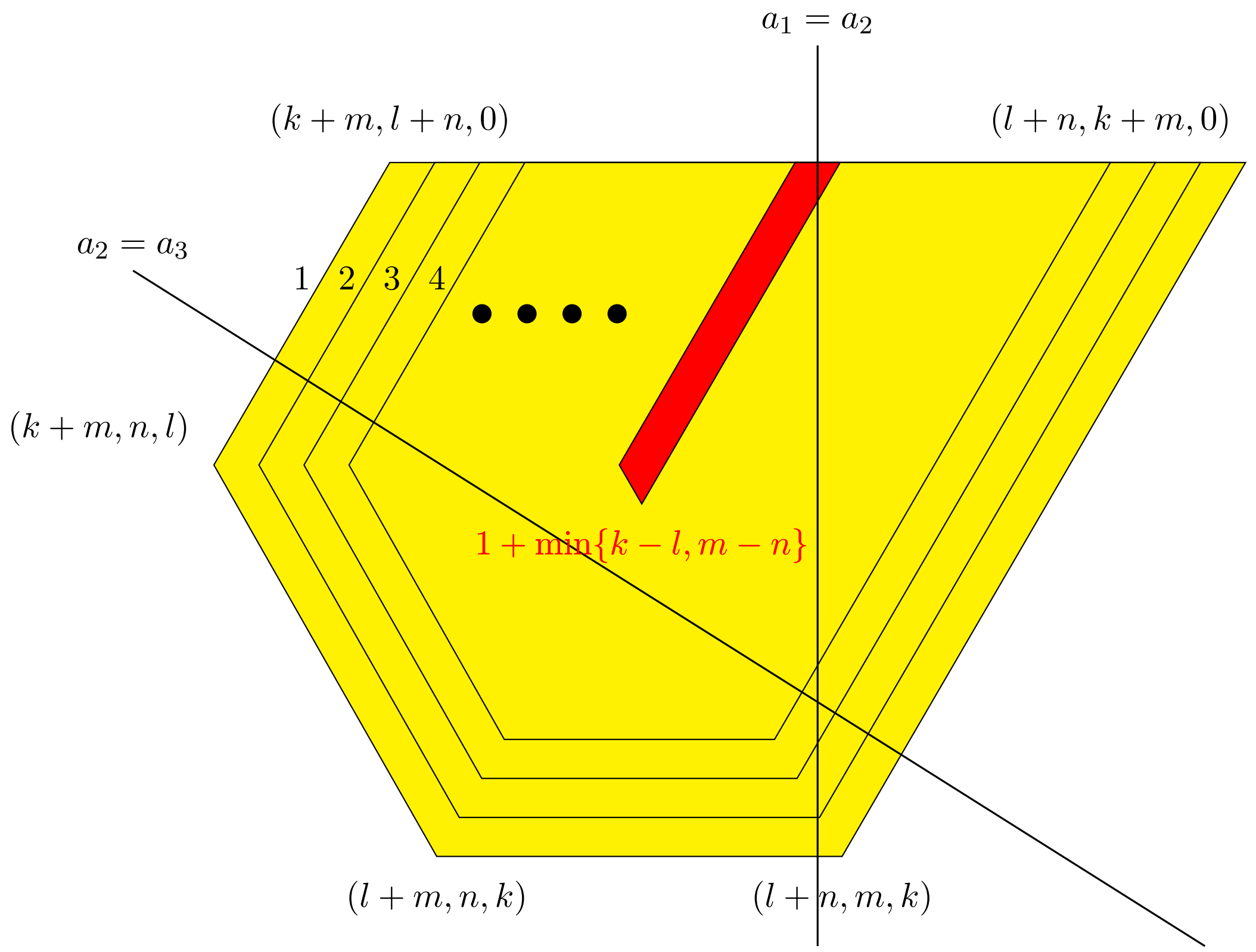

Here the number next to each line represents the multiplicity of the lattice points lying on that line (i.e. the number of trapezoids that the lattice point lies, which is also equal to the coefficient of the monomial corresponding to that lattice point) like what we have shown in Case 6. in Section 6.1. Also all the lattice points in the red region have the maximum multiplicity.

We can now decompose the polytope (with multiplicity) as follows:

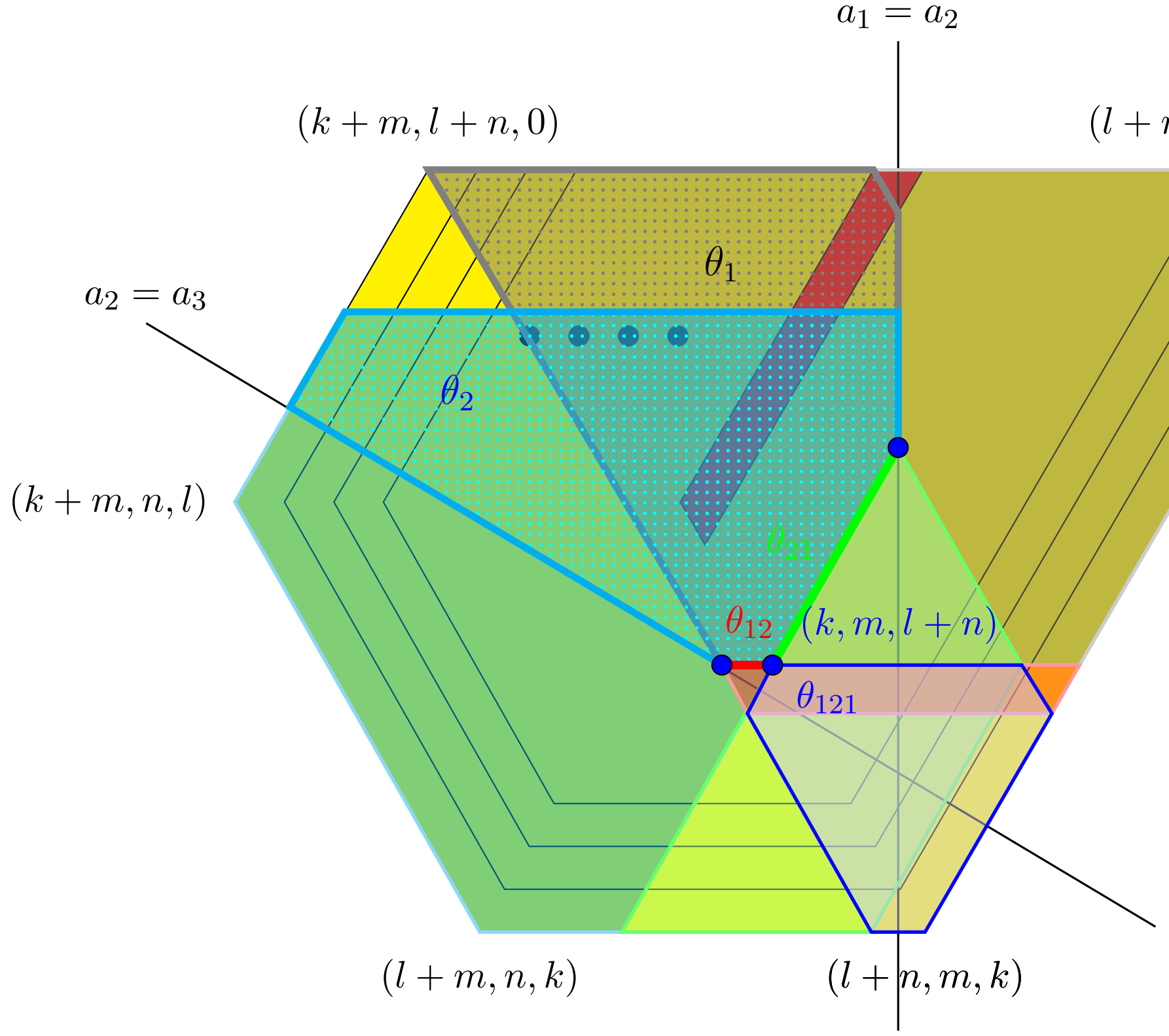

One can check that the multiplicity of each lattice point in the original yellow region is at least the multiplicity of the sum of the multiplicities in all other colored regions. This ensures the remaining points in $R_1$ still corresponds to a positive sum of dominating monomials (i.e. $A_0$ is a positive sum of some dominating monomials).

We can also check some of the coefficients by using operators.

Suppose

$$\begin{aligned} & \pi_1(x_1^m x_2^n) \times \pi_{21}(x_1^k x_2^l) \\ = \quad & B_0(x) + \theta_1 B_1(x) + \theta_2 B_2(x) + \theta_{21} B_{21}(x) + \theta_{12} B_{12}(x) + \theta_{121} B_{121}(x). \end{aligned}$$

By Lemma 3.6 and Proposition 3.1, after applying $\pi_1$ on both sides, we get:

$$\begin{aligned} & \pi_1(\pi_1(x_1^m x_2^n) \times \pi_{21}(x_1^k x_2^l)) \\ = \quad & \pi_1(\pi_{21}(x_1^k x_2^l) \times \pi_1(x_1^m x_2^n)) \\ = \quad & \pi_{121}(x_1^k x_2^l) \times \pi_1(x_1^m x_2^n) + s_1 \pi_{21}(x_1^k x_2^l) \times \theta_1 \pi_1(x_1^m x_2^n) \\ = \quad & \pi_{121}(x_1^k x_2^l) \times \pi_1(x_1^m x_2^n) \end{aligned}$$

and

$$\begin{aligned}&\pi_1(B_0(x)+\theta_1B_1(x)+\theta_2B_2(x)+\theta_{21}B_{21}(x)+\theta_{12}B_{12}(x)+\theta_{121}B_{121}(x))\\=\;&\pi_1B_0(x)+\pi_1\theta_2B_2(x)+\pi_1\theta_{21}B_{21}(x).\end{aligned}$$

Therefore $\pi_1(\pi_1(x_1^mx_2^n)\times\pi_{21}(x_1^kx_2^l))=\pi_1B_0(x)+\pi_1\theta_2B_2(x)+\pi_1\theta_{21}B_{21}(x)$.

Now apply $\pi_1$ on both sides of

$$(x_1^mx_2^n)\times\pi_{121}(x_1^kx_2^l)=A_0(x)+\theta_1A_1(x)+\theta_2A_2(x)+\theta_{21}A_{21}(x)+\theta_{12}A_{12}(x)+\theta_{121}A_{121}(x),$$

we get

$$\begin{aligned}&\pi_1((x_1^mx_2^n)\times\pi_{121}(x_1^kx_2^l))\\=&\pi_1(A_0(x)+\theta_1A_1(x)+\theta_2A_2(x)+\theta_{21}A_{21}(x)+\theta_{12}A_{12}(x)+\theta_{121}A_{121}(x))\\=&\pi_1A_0(x)+\pi_1\theta_2A_2(x)+\pi_1\theta_{21}A_{21}(x).\end{aligned}$$

Since

$$\begin{aligned}&\pi_1((x_1^mx_2^n)\times\pi_{121}(x_1^kx_2^l))\\=&\pi_1(x_1^mx_2^n)\times\pi_{121}(x_1^kx_2^l)+s_1(x_1^mx_2^n)\times\theta_1\pi_{121}(x_1^kx_2^l)\\=&\pi_1(x_1^mx_2^n)\times\pi_{121}(x_1^kx_2^l),\end{aligned}$$

we can conclude that

$$\pi_1A_0(x)+\pi_1\theta_2A_2(x)+\pi_1\theta_{21}A_{21}(x)=\pi_1B_0(x)+\pi_1\theta_2B_2(x)+\pi_1\theta_{21}B_{21}(x).$$

Expand both sides as a sum of atoms:

$$\begin{aligned}&A_0(x)+\theta_1A_0(x)+\theta_2A_2(x)+\theta_1\theta_2A_2(x)+\theta_{21}A_{21}(x)+\theta_1\theta_{21}A_{21}(x)\\=&B_0(x)+\theta_1B_0(x)+\theta_2B_2(x)+\theta_1\theta_2B_2(x)+\theta_{21}B_{21}(x)+\theta_1\theta_{21}B_{21}(x)\end{aligned}$$

and thus

$$\begin{aligned}&A_0(x)+\theta_1A_0(x)+\theta_2A_2(x)+\theta_{12}A_2(x)+\theta_{21}A_{21}(x)+\theta_{121}A_{21}(x)\\=&B_0(x)+\theta_1B_0(x)+\theta_2B_2(x)+\theta_{12}B_2(x)+\theta_{21}B_{21}(x)+\theta_{121}B_{21}(x).\end{aligned}$$

Since the set of all atoms form a basis by item 3 in Theorem 3.8, we have $A_0=B_0, A_2=B_2$ and $A_{21}=B_{21}$.

### 7.3 $\pi_2(x_1^mx_2^n)\times\pi_{12}(x_1^kx_2^l)$

With the same notation in Section 7.2, we have

$$\begin{aligned}&\pi_2(x_1^mx_2^n)\times\pi_{12}(x_1^kx_2^l)\\=&A_0(x)+\theta_1A_1(x)+\sum_{r=0}^{\min\{m,k\}}\sum_{s=\max\{l,n,r,(n+l)-r\}}^{\min\{n+l,m+k-r\}}\theta_2(x_1^{m+k+n+l-s-r}x_2^sx_3^r)+\theta_{12}A_{12}(x)\\&+\mathbb{1}_{\{n+l\geq\max\{m,k\}\}}\sum_{r=m+k-n-l}^{\min\{m,k\}}\theta_{21}(x_1^{n+l}x_2^{m+k-r}x_3^r)+\mathbb{1}_{\{n+l>k>m\}}\theta_{121}(x_1^{n+l}x_2^kx_3^m).\end{aligned}$$

## 7.4 $\pi_{12}(x_1^m x_2^n) \times \pi_{21}(x_1^k x_2^l)$

With the same notation in Section 7.2, we have

$$\begin{aligned}&\pi_{12}(x_1^m x_2^n) \times \pi_{21}(x_1^k x_2^l)\\ =& (1+\theta_1+\theta_2)A_0(x) + \theta_{21}A_1(x) + \theta_{12}A_2(x)\mathbb{1}_{\{m+l>k>n\}} + \sum_{t=0}^{\min\{m-n,l\}} \theta_{121}(x_1^{m+l-t}x_2^k x_3^{n+t}).\end{aligned}$$

## 7.5 All other cases in the table

We will complete the verification for all other unknown cases in the tables by applying Lemma 4.16 and Lemma 3.6 on verified cases.

(i) As

$$\begin{aligned}&\pi_1(x_1^m x_2^n) \times \pi_1(x_1^k x_2^l)\\ =&\ \pi_1(x_1^m x_2^n \times \pi_1(x_1^k x_2^l))\\ =&\ \pi_1(x_1^m x_2^n \times x_1^k x_2^l + x_1^m x_2^n \times \theta_1(x_1^k x_2^l)),\end{aligned}$$

result follows by Theorem 4.13.

(ii)

$$\begin{aligned}&\pi_1(x_1^m x_2^n) \times \pi_{12}(x_1^k x_2^l)\\ =&\ \pi_1(x_1^m x_2^n \times \pi_{12}(x_1^k x_2^l))\\ =&\ \pi_1\Big(x_1^m x_2^n \times (x_1^k x_2^l + \theta_1(x_1^k x_2^l) + \theta_2(x_1^k x_2^l) + \theta_{12}(x_1^k x_2^l))\Big),\end{aligned}$$

result follows by Theorem 4.13.

Alternatively, one can use the fact that $\pi_1(x_1^m x_2^n) \times \pi_{12}(x_1^k x_2^l) = \pi_1(\pi_1(x_1^m x_2^n) \times \pi_2(x_1^k x_2^l)) = \pi_1(\pi_2(x_1^k x_2^l) \times \pi_1(x_1^m x_2^n))$ by putting $f = \pi_2(x_1^k x_2^l)$, $g = \pi_1(x_1^m x_2^n)$ and $i = 1$ in Lemma 3.6 and claim $\pi_1(x_1^m x_2^n) \times \pi_{12}(x_1^k x_2^l)$ is atom positive by applying Lemma 4.16 on $\pi_1(x_1^m x_2^n) \times \pi_2(x_1^k x_2^l)$ which is verified as atom positive in Section 7.1.

(iii)

$$\begin{aligned}&\pi_2(x_1^m x_2^n) \times \pi_2(x_1^k x_2^l)\\ =&\ \pi_2(x_1^m x_2^n \times \pi_2(x_1^k x_2^l))\\ =&\ \pi_2(x_1^m x_2^n \times x_1^k x_2^l + x_1^m x_2^n \times \theta_2(x_1^k x_2^l)),\end{aligned}$$

result follows by Theorem 4.13.

(iv) As

$$\begin{aligned}&\pi_2(x_1^m x_2^n) \times \pi_{21}(x_1^k x_2^l)\\ =&\pi_2(x_1^m x_2^n \times \pi_{21}(x_1^k x_2^l))\\ =&\pi_2(x_1^m x_2^n \times x_1^k x_2^l + x_1^m x_2^n \times \theta_1(x_1^k x_2^l) + x_1^m x_2^n \times \theta_2(x_1^k x_2^l) + x_1^m x_2^n \times \theta_{21}(x_1^k x_2^l))\end{aligned}$$

and result follows by Theorem 4.13.

Similar to Case (ii), one can also use the fact that
$\pi_2(x_1^m x_2^n) \times \pi_{21}(x_1^k x_2^l) = \pi_2(\pi_2(x_1^m x_2^n) \times \pi_1(x_1^k x_2^l)) = \pi_2(\pi_1(x_1^k x_2^l) \times \pi_2(x_1^m x_2^n))$
and claim $\pi_1(x_1^m x_2^n) \times \pi_{12}(x_1^k x_2^l)$ is atom positive.

(v) $\pi_{21}(x_1^m x_2^n) \times \pi_{21}(x_1^k x_2^l) = \pi_2(\pi_1(x_1^m x_2^l) \times \pi_{21}(x_1^k x_2^l))$ and result follows by the case in Section 7.2.

(vi) $\pi_{12}(x_1^m x_2^n) \times \pi_{12}(x_1^k x_2^l) = \pi_1(\pi_2(x_1^m x_2^l) \times \pi_{12}(x_1^k x_2^l))$ and result follows by the case in Section 7.3.